\documentclass[11pt]{article}
\usepackage[a4paper,margin=27mm]{geometry}
\usepackage[T1]{fontenc}
\usepackage{lmodern,amsmath,amssymb,amsthm,mathtools,microtype}
\usepackage[colorlinks=true,linkcolor=blue,citecolor=blue,urlcolor=blue]{hyperref}
\usepackage{enumitem}
\newcommand{\Vzero}{{V^{0}}}
\newcommand{\Rop}{R_{\mathrm{cond}}}
\newcommand{\Bstress}{{B^\Sigma}}
\setlist{itemsep=3pt,topsep=5pt}
\numberwithin{equation}{section}
\newtheorem{theorem}{Theorem}[section]
\newtheorem{proposition}[theorem]{Proposition}
\newtheorem{lemma}[theorem]{Lemma}

\theoremstyle{remark}

\newtheorem*{remark*}{Remark}
\newcommand{\R}{\mathbb R}
\newcommand{\C}{\mathbb C}
\newcommand{\T}{\mathbb T}
\newcommand{\D}{\mathbb D}
\newcommand{\N}{\mathbb N}
\newcommand{\Z}{\mathbb Z}
\newcommand{\E}{\mathbb E}
\newcommand{\cE}{\mathcal E}
\newcommand{\Ccyl}{C_{\mathrm{cyl}}}
\newcommand{\Dom}{\operatorname{Dom}}
\newcommand{\Homeo}{\operatorname{Homeo}}
\newcommand{\Law}{\operatorname{Law}}
\newcommand{\id}{\operatorname{id}}
\newcommand{\Ric}{\operatorname{Ric}}
\newcommand{\supp}{\operatorname{supp}}
\newcommand{\HS}{\mathrm{HS}}
\newcommand{\Tr}{\operatorname{Tr}}
\renewcommand{\Cap}{\operatorname{Cap}}

\DeclareMathOperator{\Rea}{Re}
\DeclareMathOperator{\Ima}{Im}
\allowdisplaybreaks[2]
\hypersetup{
pdftitle={Weak uniqueness and autonomous approximation for the SLE-reversible Weil-Petersson response diffusion},
pdfauthor={Chunhao Cai},
pdfsubject={Self-consistent response diffusion},
pdfkeywords={SLE, weak uniqueness, reversible diffusion, autonomous approximation},
bookmarksnumbered=true,
bookmarksopen=true
}
\title{Weak uniqueness and autonomous approximation\\
for the SLE-reversible Weil--Petersson response diffusion}
\author{Chunhao Cai\\
{\small School of Mathematics (Zhuhai), Sun Yat-Sen University}\\
{\small \texttt{caichh9@mail.sysu.edu.cn}}}
\date{}
\begin{document}
\maketitle
\begin{abstract}
We study the self-consistent common-noise equation associated with the
Weil--Petersson response differential of the rotationally marked SLE
welding probability. For $0<\kappa\leq4$, we prove joint weak uniqueness
among continuous adapted weak solutions whose law at each deterministic
time is the welding probability, and obtain the Markov property in the
original filtration. A conditional Gaussian lift and capacity localization
give the required $L^1$ generator-range density. Identification with the
canonical law gives stationarity and reversibility. We also construct
autonomous Fej\'er and Poisson regularizations with pathwise-unique global
strong solutions. Their joint state-noise laws converge to the canonical
law for every initial probability absolutely continuous with respect to
the welding probability.
\end{abstract}


\setcounter{tocdepth}{2}
\begingroup
\small
\tableofcontents
\endgroup
\clearpage

\section{Introduction}\label{sec:introduction}

The closed response calculus of \cite{Response} associates a
Weil--Petersson gradient and an integration-by-parts formula with the
rotationally marked SLE welding probability. The resulting Dirichlet
form generates the canonical reversible diffusion constructed in
\cite{Diffusion}. That construction gives a common Fourier noise and
a boundary equation whose drift depends on the current welding.

We prove that this self-consistent equation characterizes the
diffusion among weak solutions with the prescribed one-time marginals.
We then recover its joint state-noise law by autonomous regularizations.

\subsection{The self-consistent response equation}

Fix $0<\kappa\leq4$. Write $\D=\{z\in\C:|z|<1\}$,
$S^1=\partial\D$, and $\T=\R/(2\pi\Z)$, identified with $S^1$
by $\theta\mapsto e^{i\theta}$. The state space is
\begin{equation}\label{eq:state}
 \mathcal X=\Homeo_+(S^1),\qquad
 d(h,\widetilde h)=
 \sup_{z\in S^1}|h(z)-\widetilde h(z)|
 +\sup_{z\in S^1}|h^{-1}(z)-\widetilde h^{-1}(z)|.
\end{equation}
Let $\nu$ be the rotationally marked $\mathrm{SLE}_\kappa$ welding
probability of \cite{Response,Diffusion} on the Polish space
\eqref{eq:state}; the rotational mark is retained in the state $h$.

For $n\geq2$ and $\tau\in\{c,s\}$ put
\begin{equation}\label{eq:basis}
 d_n=n^3-n,\qquad \alpha_n=\sqrt{\frac2{d_n}},\qquad
 e_{n,c}(\theta)=\alpha_n\cos(n\theta),\quad
 e_{n,s}(\theta)=\alpha_n\sin(n\theta).
\end{equation}
For real trigonometric polynomials
$v(\theta)=\sum_{|n|\geq2}v_ne^{in\theta}$ and
$w(\theta)=\sum_{|n|\geq2}w_ne^{in\theta}$, set
\[
 (v,w)_H=2\Rea\sum_{n\geq2}d_nv_n\overline{w_n}.
\]
Let $H$ be their completion in this inner product. The fields in
\eqref{eq:basis} form an orthonormal basis of $H$.
We use the index set ${J_H}=\{(n,\tau):n\geq2,\ \tau\in\{c,s\}\}$,
and write $e_j$ when the frequency of $j$ is not needed.

Set
\[
 Q=\frac{\sqrt\kappa}{2}+\frac2{\sqrt\kappa},\qquad
 a_\kappa=\frac{1+6Q^2}{12}.
\]
For a normalized univalent function
$f(z)=z+\sum_{m\geq1}a_mz^{m+1}$ define its Schwarzian derivative
${\operatorname{S}}f=f'''/f'-\frac32(f''/f')^2$ and the coefficients $t_n(f)$ by
\begin{equation}\label{eq:stress}
 a_\kappa{\operatorname{S}}f(z)+\left(\frac{f'(z)}{f(z)}\right)^2
 =z^{-2}+2a_1z^{-1}+\sum_{n\geq2}t_n(f)z^{n-2}.
\end{equation}
Write
\begin{equation}\label{eq:primitive}
 A_f(z)=\sum_{n\geq2}\frac{t_n(f)}{d_n}z^n .
\end{equation}
On the full-measure welding locus, let $f_h^\circ$ be the interior map
in the measurable conformal pair $f_h^\circ=g_h\circ h$ selected in
\cite{Diffusion}, where $g_h$ is the exterior map and
$f_h^\circ(0)=0$, $(f_h^\circ)'(0)=1$.
Write $A_f^*$ for the uniform radial boundary limit of $A_f$.
On the full-measure domain of this selection define
\begin{equation}\label{eq:drift}
 B(h,\theta)=2\Ima A_{f_h^\circ}^*(e^{i\theta})
\end{equation}
and set $B=0$ on the null complement. Then
$B:\mathcal X\to C(\T)$ is measurable.
There is a constant $C_\kappa$ such that, for $0<|\eta|\leq1$,
\begin{align}
 \|B(h,\cdot)\|_\infty&\leq C_\kappa,\notag\\
 |B(h,\theta+\eta)-2B(h,\theta)+B(h,\theta-\eta)|
     &\leq C_\kappa|\eta|.\label{eq:drift-modulus}
\shortintertext{Consequently,}
 |B(h,\theta+\eta)-B(h,\theta)|
     &\leq C_\kappa|\eta|\log(e/|\eta|).\notag
\end{align}
These are the uniform drift bounds of \cite{Diffusion}.

A weak solution consists of a complete filtered probability space
$(\Omega,\mathcal F,(\mathcal F_t),\mathbb P)$ satisfying the usual
conditions, independent standard $(\mathcal F_t)$-Brownian motions
$W=(W^j)_{j\in J_H}$, and an adapted continuous $\mathcal X$-valued
process $X$. Set
\begin{equation}\label{eq:self-consistency}
 \Psi_t=X_t^{-1}\circ X_0,\qquad X_t=X_0\circ\Psi_t^{-1}.
\end{equation}
Let $\Theta_t(u)$ be the jointly continuous angular lift of $\Psi_t$
with $\Theta_0(u)=u$ and
$\Theta_t(u+2\pi)=\Theta_t(u)+2\pi$. Fix the dense label set
$\Lambda=2\pi\mathbb Q$. The equation is
\begin{equation}\label{eq:raw-sde}
 \Theta_t(u)=u+\sqrt2\sum_{j\in {J_H}}\int_0^t
 e_j(\Theta_s(u))\,dW_s^j
 +\int_0^t B(X_s,\Theta_s(u))\,ds,\qquad u\in\Lambda.
\end{equation}
This identity holds simultaneously for $u\in\Lambda$ on each finite
time interval. Its stochastic series converges in
$L^2(\mathbb P;C([0,T],\R))$; the tail estimate is recorded in
Section~\ref{sec:finite-tests}.

We consider weak solutions with the prescribed marginals
\begin{equation}\label{eq:marginal-class}
                         \Law(X_t)=\nu,\qquad t\geq0.
\end{equation}

\subsection{Main results}

Let $D:\Dom(D)\subset L^2(\nu)\to L^2(\nu;H)$ be the closed
response gradient of \cite{Response}, obtained by closing derivatives
along $h\mapsto h\circ\psi_t^v$, where $\psi_t^v$ is the circle flow of
$v(\theta)\partial_\theta$; its cylinder definition is recalled in
Section~\ref{sec:finite-tests}. Its closed form is
\[
 \cE(F,G)=\langle DF,DG\rangle_{L^2(\nu;H)},
 \qquad F,G\in\Dom(D).
\]
The associated generator is $L=-D^*D$ on
$\Dom(L)=\{F\in\Dom(D):DF\in\Dom(D^*)\}$.
Write $(T_t)_{t\geq0}$ for its semigroup.

Denote by $\mathbb Q_\nu$ the joint law of the canonical stationary
diffusion and the Fourier Brownian family reconstructed in
\cite{Diffusion}. Write $C_0([0,\infty),\R)$ for continuous
real paths starting at zero. Joint laws are considered on the Polish space
\begin{equation}\label{eq:path-space}
 {E_{\mathrm{path}}}=C([0,\infty),\mathcal X)\times
       \prod_{j\in {J_H}} C_0([0,\infty),\R),
\end{equation}
where each factor carries uniform convergence on bounded time intervals.

\begin{theorem}\label{thm:weak-uniqueness}
Weak solutions of
\eqref{eq:self-consistency}--\eqref{eq:raw-sde} satisfying
\eqref{eq:marginal-class} exist. Every such solution satisfies
\begin{equation}\label{eq:joint-unique}
 \Law(X_\cdot,W)=\mathbb Q_\nu.
\end{equation}
Moreover, for every bounded measurable $F:\mathcal X\to\R$ and $0\leq s\leq t$,
\begin{equation}\label{eq:conditional-semigroup}
 \E[F(X_t)\mid\mathcal F_s]=T_{t-s}F(X_s)
 \quad\mathbb P\text{-almost surely}.
\end{equation}
\end{theorem}

To state the approximation result, choose the score representatives
as the Fourier coefficients of \eqref{eq:drift}:
\begin{equation}\label{eq:score-versions}
 \begin{aligned}
 \rho_{n,c}(h)&=\frac1{\pi\alpha_n}
       \int_0^{2\pi}B(h,\theta)\cos(n\theta)\,d\theta,\\
 \rho_{n,s}(h)&=\frac1{\pi\alpha_n}
       \int_0^{2\pi}B(h,\theta)\sin(n\theta)\,d\theta.
 \end{aligned}
\end{equation}
They agree $\nu$-almost everywhere with $D^*e_{n,\tau}$.
Write $x_+=\max(x,0)$. The two multiplier families are
\begin{equation}\label{eq:weights}
 \gamma^{N}_n=\left(1-\frac n{N+1}\right)_+,\quad N\geq2,
 \qquad
 \gamma^r_n=r^n,\quad 0<r<1.
\end{equation}
Set $\gamma_{n,c}=\gamma_{n,s}=\gamma_n$.
For either fixed multiplier sequence $\gamma$ define
\begin{equation}\label{eq:weighted-drift}
 B_\gamma(h,\theta)=\sum_{j\in {J_H}}\gamma_j^2\rho_j(h)e_j(\theta).
\end{equation}
If $\operatorname{Fej}_N$ and $\operatorname{Pois}_r$ are respectively the circle Fej\'er
and Poisson convolution operators, then
\begin{equation}\label{eq:matched-convolutions}
 B_{\gamma^N}=\operatorname{Fej}_N^2B,\qquad
 B_{\gamma^r}=\operatorname{Pois}_{r^2}B.
\end{equation}
For a continuous family $\Psi_t^\gamma$ starting at the identity, let
$\Theta_t^\gamma$ be its continuous angular lift, with
$\Theta_0^\gamma(u)=u$ and
$\Theta_t^\gamma(u+2\pi)=\Theta_t^\gamma(u)+2\pi$.
The squared multipliers in $B_\gamma$ match the covariance of the
regularized noise and preserve symmetry of the associated form.
The autonomous equation, imposed simultaneously for $u\in\Lambda$, is
\begin{align}
 X_t^\gamma&=X_0\circ(\Psi_t^\gamma)^{-1},
       \qquad \Theta^\gamma_0(u)=u,\label{eq:autonomous-state}\\
 d\Theta_t^\gamma(u)
   &=\sqrt2\sum_{j\in {J_H}}\gamma_j
       e_j(\Theta_t^\gamma(u))\,dW_t^j
     +B_\gamma(X_t^\gamma,\Theta_t^\gamma(u))\,dt.
       \label{eq:autonomous-sde}
\end{align}
Strong solutions use the completed filtration generated by $(X_0,W)$;
the remaining conventions are those of \eqref{eq:raw-sde}.

For a probability $\mu\ll\nu$, choose a nonnegative density
$g=d\mu/d\nu$ and define the canonical joint law
\begin{equation}\label{eq:initial-weighting}
 d\mathbb Q_\mu=g(X_0)\,d\mathbb Q_\nu.
\end{equation}
The notation $\Longrightarrow$ below means weak convergence of
probability measures on \eqref{eq:path-space}.

\begin{theorem}\label{thm:autonomous-approximation}
For every $\mu\ll\nu$ and each fixed multiplier sequence in
\eqref{eq:weights}, the autonomous equation
\eqref{eq:autonomous-state}--\eqref{eq:autonomous-sde}, with
$\Law(X_0)=\mu$, has a global pathwise-unique strong solution.
For $\mu=\nu$ its state law is stationary and reversible.
For every $\mu\ll\nu$,
\begin{align}
 \Law(X^{\gamma^N}_\cdot,W)&\Longrightarrow\mathbb Q_\mu
       &&(N\to\infty),\notag\\
 \shortintertext{and}
 \Law(X^{\gamma^r}_\cdot,W)&\Longrightarrow\mathbb Q_\mu
       &&(r\uparrow1).\label{eq:main-approximation}
\end{align}
\end{theorem}

The finite-test martingales in Section~\ref{sec:finite-tests} reduce
weak uniqueness to identifying conditional laws.
Section~\ref{sec:geometry} provides the domain and divergence estimates
for conditional compression. Sections~\ref{sec:gaussian-lift}
and~\ref{sec:raw-geometry} represent the welding by two Gaussian fields
and an inverse mark. The marked dynamics in
Section~\ref{sec:gaussian-core} give $L^1$ resolvent-range density;
Section~\ref{sec:conditional-stress} transfers it to the welding law,
identifies the conditional semigroup, and recovers the common noise.
Section~\ref{sec:autonomous} constructs the regularized strong solutions
and identifies their joint limits by Theorem~\ref{thm:weak-uniqueness}.

\section{Response calculus and finite-test martingales}\label{sec:finite-tests}

\subsection{Response identities}

For real smooth periodic functions $a,b$ define
\[
 J_{a,b}(h)=\frac1{2\pi}\int_0^{2\pi}
            a(\theta)b(h(e^{i\theta}))\,d\theta .
\]
Functions on $S^1$ are identified with their periodic angular forms.
The algebra $\Ccyl$ consists of
$F=\phi(J_{a^{(1)},b^{(1)}},\ldots,J_{a^{(m)},b^{(m)}})$ with
$\phi\in C_b^\infty(\R^m)$ and smooth $a^{(\ell)},b^{(\ell)}$.
Let $\psi_t^v$ be the circle flow of $v(\theta)\partial_\theta$.
The response derivative uses the right action:
\[
 D_vF(h)=\left.\frac d{dt}\right|_{t=0}F(h\circ\psi_t^v).
\]
Its closure is $D$. Write $\delta=D^*$ and
$\Gamma(F,G)=(DF,DG)_H$; unlabelled norms are in $L^2(\nu)$ or its
indicated Hilbert-valued version. We use complexifications for
complex-valued coefficients.

For the right response action, the integration-by-parts identity
is proved in Appendix~C of \cite{Response}. For the present law and
$0<\kappa\leq4$, the score and covariance identities are given in
\cite[Sections 4--5]{Response} and \cite{Diffusion}. Thus, for finite
horizontal real Fourier fields $v$ and cylinders $F,G$, writing
$D_j=D_{e_j}$, $\rho_v=D_v^*1$, and $\rho_j=\rho_{e_j}$, we have:
\begin{align}
 \int D_vF\,d\nu&=\int F\rho_v\,d\nu,\qquad
 D_v^*G=\rho_vG-D_vG.\label{eq:ibp}\\
 \shortintertext{The real Fourier scores are}
 \rho_{n,c}&=\alpha_n\Ima t_n,\qquad
 \rho_{n,s}=\alpha_n\Rea t_n,\label{eq:geometric-scores}\\
 \shortintertext{and their covariance is}
 \int t_n\overline{t_m}\,d\nu
     &=\mathbf1_{\{n=m\}}\{a_\kappa d_n+2n\}.\label{eq:ward}\\
 \shortintertext{On cylinders, the generator is}
 LF&=\sum_{j\in{J_H}}(D_j^2F-\rho_jD_jF),\qquad
 \sum_j\|D_j^2F\|_\infty+
 \sum_j\|\rho_jD_jF\|_2<\infty.\label{eq:cylinder-generator}
\end{align}
In \eqref{eq:geometric-scores}, $t_n=t_n(f_h^\circ)$; all identities
hold $\nu$-almost everywhere. Equation \eqref{eq:drift} fixes a bounded
measurable representative of $LF$.

On holomorphic coefficient polynomials the finite responses of
\cite{Response} are
$\ell_n f(z)=z^{n+1}f'(z)$ for $n\geq1$, and, for $n\geq1$, set
\begin{equation}\label{eq:negative-contour}
 \ell_{-n}f(z)=\frac{f(z)^2}{2\pi i}
 \int_{|w|=R}\frac{f'(w)^2w^{1-n}}
 {f(w)^2(f(w)-f(z))}\,dw,\qquad |z|<R<1.
\end{equation}
Each coefficient in this formula is a finite coefficient polynomial.
The actual real responses are
\begin{equation}\label{eq:coefficient-responses}
 D_{\cos(n\cdot)}=\frac i2(\ell_n+\ell_{-n}),\qquad
 D_{\sin(n\cdot)}=\frac12(\ell_n-\ell_{-n}).
\end{equation}
For coefficient polynomials $P_1,P_2$, the same integration by parts gives
\[
 \int D_v(FP_1\overline{P_2})\,d\nu
   =\int FP_1\overline{P_2}\,\rho_v\,d\nu,
 \qquad F\in\Ccyl.
\]
Each fixed polynomial and its responses in fixed finite directions
are bounded.

The form $(\cE,\Dom(D))$ is conservative, strongly local and
quasi-regular in \eqref{eq:state}. Write $\Cap$ for its $1$-capacity.
It has compact sets $K_m\subset \mathcal X$
such that
\begin{equation}\label{eq:compact-nest}
 \Cap(\mathcal X\setminus K_m)\longrightarrow0.
\end{equation}
For $z\in S^1$, let $q_z(h)$ be the angular coordinate of $h^{-1}(z)$.
The inverse-coordinate domain result of \cite{Diffusion} gives, for
smooth periodic $a$,
\begin{align}
 a(q_z)&\in\Dom(L),&
 D_j[a(q_z)]&=-a'(q_z)e_j(q_z).\label{eq:inverse-domain}
 \shortintertext{and}
 L[a(q_z)]&=\frac12a''(q_z)+B(h,q_z)a'(q_z).\label{eq:inverse-generator}
\end{align}

\subsection{Labels and the cylinder martingale problem}

For each label, Doob's inequality bounds the squared
$L^2(\mathbb P;C([0,T],\R))$ norm of the noise tail in
\eqref{eq:raw-sde} by
\begin{equation}\label{eq:series-tail}
 8T\sum_{n>N}\alpha_n^2\longrightarrow0.
\end{equation}

Define the covariance kernel and its evaluation row by
\begin{equation}\label{eq:evaluation-row}
 k_0(\theta)=\sum_{n\geq2}\alpha_n^2\cos(n\theta),\qquad
 k(\theta)=\sum_{j\in{J_H}}e_j(\theta)e_j\in H.
\end{equation}
Then
\begin{equation}\label{eq:kernel-modulus}
 k_0(0)=\frac12,\qquad
 \|k(\theta)-k(\varphi)\|_H^2
 =2\{k_0(0)-k_0(\theta-\varphi)\}
 \leq C\eta^2\log(e/\eta),\quad
 \eta=|\theta-\varphi|\leq1 .
\end{equation}
The estimate follows by splitting the series at
$n=\lfloor\eta^{-1}\rfloor$ and using
$1-\cos(n\eta)\leq\min(2,n^2\eta^2/2)$.

For a deterministic target $e^{iv}$, write $q_t(v)$ for the continuous
angular lift of $X_t^{-1}(e^{iv})$, starting from a measurable choice of
$q_0(v)$.

\begin{lemma}\label{lem:label-extension}
Every weak solution of \eqref{eq:raw-sde} has a common continuous
version of its noise field for all real labels. The equation also
holds at every $\mathcal F_0$-measurable label.
For each deterministic target $e^{iv}$,
\begin{equation}\label{eq:inverse-sde}
 dq_t(v)=\sqrt2\sum_{j\in{J_H}}e_j(q_t(v))\,dW_t^j
                         +B(X_t,q_t(v))\,dt.
\end{equation}
\end{lemma}
\begin{proof}
Set
\[
 I_t(u)=\Theta_t(u)-u-\int_0^t B(X_s,\Theta_s(u))\,ds.
\]
The joint continuity of $\Theta$ and \eqref{eq:drift-modulus} make
$I$ continuous on compact time-label rectangles.
If $u_m\in\Lambda$ tends to $u$, Doob's inequality and
\eqref{eq:kernel-modulus} give
\[
 \E\sup_{t\leq T}\left|
 \sqrt2\sum_j\int_0^t
 [e_j(\Theta_s(u_m))-e_j(\Theta_s(u))]\,dW_s^j
 \right|^2
 \leq8\E\int_0^T
       \|k(\Theta_s(u_m))-k(\Theta_s(u))\|_H^2\,ds
 \longrightarrow0.
\]
The integrands are bounded and converge pointwise. This identifies
$I_t(u)$ with the stochastic series at every fixed $u$.

For an $\mathcal F_0$-measurable label, use finite-valued
$\mathcal F_0$-measurable approximations in $\Lambda$; their indicators
can be inserted into the stochastic integrals. The same estimate
and \eqref{eq:series-tail} pass to the limit.

Choose a measurable lift $q_0(v)$ of $X_0^{-1}(e^{iv})$ and set
$q_t(v)=\Theta_t(q_0(v))$. This proves \eqref{eq:inverse-sde}.
\end{proof}

Set
\[
 C_{\mathrm{inv}}
 =\operatorname{alg}\bigl(\Ccyl,
       \{b(q_{z_1},\ldots,q_{z_m}):
         b\in C^\infty(\T^m),\ z_i\in S^1,\ m\geq1\}\bigr).
\]
For each $F\in C_{\mathrm{inv}}$,
\[
 \|F\|_\infty+\|LF\|_\infty
       +\bigl\|\|DF\|_H\bigr\|_\infty<\infty.
\]

\begin{lemma}\label{lem:raw-martingales}
For a weak solution satisfying \eqref{eq:marginal-class} and $F,G\in{C_{\mathrm{inv}}}$,
\begin{align}
 F(X_t)-F(X_0)&=M_t^F+\int_0^t LF(X_s)\,ds,\label{eq:raw-mp}\\
 M_t^F&=-\sqrt2\sum_j\int_0^t D_jF(X_s)\,dW_s^j,\notag\\
 \shortintertext{Their quadratic covariation is}
 \langle M^F,M^G\rangle_t
     &=2\int_0^t\Gamma(F,G)(X_s)\,ds.\label{eq:raw-bracket}
\end{align}
These are square-integrable martingales on every finite interval.
\end{lemma}
\begin{proof}
For inverse evaluations this is scalar It\^o calculus applied to
\eqref{eq:inverse-sde}, followed by the product rule.

For cylinders put
$A(u)=\int_0^u a(\theta)\,d\theta$ and
$c_a=\int_0^{2\pi}a(\theta)\,d\theta$. Since
$q_t(v+2\pi)=q_t(v)+2\pi$, change of variables and Stieltjes
integration by parts give
\begin{align*}
 \int_0^{2\pi}a(\theta)b(X_t(e^{i\theta}))\,d\theta
 &=\int_0^{2\pi}b(e^{iv})\,d[A(q_t(v))]\\
 &=c_ab(1)-\int_0^{2\pi}A(q_t(v))b'(v)\,dv,
\end{align*}
and hence
\begin{equation}\label{eq:stieltjes}
 J_{a,b}(X_t)=\frac1{2\pi}
 \left[c_ab(1)-\int_0^{2\pi}A(q_t(v))b'(v)\,dv\right].
\end{equation}
The bounds on $a,a',B,b'$ and \eqref{eq:series-tail} justify scalar
It\^o calculus and stochastic Fubini. They give
\begin{align*}
 D_jJ_{a,b}(h)
 &=\frac1{2\pi}\int_0^{2\pi}
       a(q_{e^{iv}}(h))e_j(q_{e^{iv}}(h))b'(v)\,dv,\\
 \shortintertext{and}
 LJ_{a,b}(h)
 &=\frac12J_{a'',b}(h)
 -\frac1{2\pi}\int_0^{2\pi}
       a(\theta)B(h,\theta)\,d(b\circ h)(\theta).
\end{align*}
Finite-dimensional It\^o calculus and the product formula now give
\eqref{eq:raw-mp}--\eqref{eq:raw-bracket} on $C_{\mathrm{inv}}$.

Finally, for every measurable $N$ with $\nu(N)=0$,
\[
 \E\int_0^T\mathbf1_N(X_s)\,ds=T\nu(N)=0 .
\]
\end{proof}

\subsection{Coefficient functions in the closed domain}

Let $dA$ be planar area measure. For $1\leq p<\infty$ write
\[
 A^p=\operatorname{Hol}(\D)\cap L^p(\D,dA/\pi),\qquad
 \|u\|_{A^p}^p=\frac1\pi\int_\D|u(z)|^p\,dA(z).
\]
For a holomorphic function $u$, let $[z^m]u$ denote its coefficient of
$z^m$. A coefficient polynomial is a polynomial in finitely many
$a_m(f)$; its weighted degree assigns weight $m$ to $a_m$.
Antiholomorphic coefficient polynomials are their complex conjugates.
\begin{lemma}\label{lem:coefficient-domain}
Every holomorphic or antiholomorphic coefficient polynomial $P$ belongs
to $\Dom(L)$, and
\begin{equation}\label{eq:polynomial-bounds}
 \|P\|_\infty+\|LP\|_\infty+
       \bigl\|\|DP\|_H\bigr\|_\infty<\infty .
\end{equation}
Its closed response components are the finite coefficient responses
in \eqref{eq:coefficient-responses}.
\end{lemma}
The proof is given in Appendix~\ref{app:coefficient-domain}.

\section{Response domains and divergence}\label{sec:geometry}

Conditional compression requires two domain statements:
$D=D_{\max}$ and $\Dom(\nabla_{H,\max})\subset\Dom(\delta)$,
with divergence bounded by one covariant derivative.
We prove them by adjoining the three low directions, establishing
the full Hodge identity, and then returning to the horizontal metric.
The required Fourier estimates are collected in
Appendix~\ref{app:response-estimates}.

\subsection{The three low directions}

Put $t_1=2a_1$ and
\[
 z_0(\theta)=1,\quad z_c(\theta)=\cos\theta,\quad
 z_s(\theta)=\sin\theta,\qquad Z_a=D_{z_a}\quad(a=0,c,s).
\]
Complex responses are denoted by $\partial_n=D_{e^{in\theta}}$,
$n\in\Z$. The actual response bracket is
\begin{equation}\label{eq:response-bracket}
 [\partial_p,\partial_q]=i(p-q)\partial_{p+q}.
\end{equation}
It holds on the algebra generated by $\Ccyl$ and coefficient
polynomials, by the finite right-action rules.

Define the low scores by
\begin{equation}\label{eq:low-scores}
 \sigma_0=0,\qquad \sigma_c=\Ima t_1,\qquad
 \sigma_s=\Rea t_1.
\end{equation}

\begin{lemma}\label{lem:low-divergences}
For $a\in\{0,c,s\}$ and $F\in\Ccyl$,
\[
 \int Z_aF\,d\nu=\int F\sigma_a\,d\nu.
\]
Each $Z_a$ is closable on $\Ccyl$.
\end{lemma}
\begin{proof}
Write $\sigma_n$ for the complex divergence of $\partial_n$.
For $n\geq2$,
$\sigma_n=i\overline{t_n}$ and $\sigma_{-n}=-it_n$.
Twice applying \eqref{eq:ibp} gives
\[
 \int[\partial_3,\partial_{-2}]F\,d\nu
 =\int F(\partial_3\sigma_{-2}
                       -\partial_{-2}\sigma_3)\,d\nu
 =-5\int F\overline{t_1}\,d\nu.
\]
Set $\tau_f=a_\kappa{\operatorname{S}}f+(f'/f)^2$ and $V(z)=z^{n+1}$.
The positive response is
\[
 \ell_n\tau_f=V\tau_f'+2V'\tau_f+a_\kappa V'''.
\]
Coefficient extraction gives $\ell_3t_2=0$ and
$\ell_2t_3=5t_1$, which prove the last equality.
Equation \eqref{eq:response-bracket} therefore gives
$\int\partial_1F\,d\nu=\int Fi\overline{t_1}\,d\nu$.
Conjugation and the $(2,-2)$ commutator give, respectively,
\[
 \int\partial_{-1}F\,d\nu=-i\int Ft_1\,d\nu,
 \qquad \int\partial_0F\,d\nu=0.
\]
Taking real combinations proves \eqref{eq:low-scores}.

Since $|t_1|\leq4$, for $F\in\Ccyl$ we have
\[
 Z_a^*F=\sigma_aF-Z_aF\in L^2(\nu).
\]
Thus $Z_a^*$ is densely defined, which proves closability.
\end{proof}

Complete the response metric by declaring $z_0,z_c,z_s$ orthonormal
and perpendicular to $H$. Denote this Hilbert space by $\widehat H$.
In the complex Fourier frame $l_n(\theta)=e^{in\theta}$ its weights are
\begin{equation}\label{eq:full-weights}
 \omega_0=1,\qquad \omega_{\pm1}=2,\qquad
 \omega_n=|n|^3-|n|\quad(|n|\geq2).
\end{equation}
Let $(E_A)$ be the full real orthonormal frame, consisting of
$(e_j)_{j\in J_H}$ and $z_0,z_c,z_s$. Define
\[
 D_A=D_{E_A},\qquad
 \widehat\rho_A=
 \begin{cases}
  \rho_j,&E_A=e_j,\\
  \sigma_a,&E_A=z_a.
 \end{cases}
\]
For a finite field $v=\sum_Av_AE_A$, retain the notation
$\rho_v=\sum_Av_A\widehat\rho_A$.

The full gradient $\widehat D$ is the closure of $(D_A)_A$ on
$\Ccyl$. Its form, generator, and divergence are
\[
 \widehat{\cE}(F,G)=\langle\widehat DF,\widehat DG\rangle,
 \qquad \widehat L=-\widehat D^*\widehat D,
 \qquad \widehat\delta=\widehat D^*.
\]

The constant metric connection is defined on finite Fourier fields by
\begin{align}
 \nabla_{l_p}l_q&=ic(p,q)l_{p+q},\label{eq:full-connection}\\
 c(p,q)&=\frac12\left[
 p-q-\frac{(p+2q)\omega_p+(2p+q)\omega_q}{\omega_{p+q}}
 \right].\notag
\end{align}
Its real coefficients are written
$\nabla_{E_A}E_B=\sum_C c_{AB}^{C}E_C$.
The Koszul formula gives
\begin{equation}\label{eq:metric-torsion}
 c_{AB}^{C}=-c_{AC}^{B},\qquad
 \nabla_vw-\nabla_wv=[v,w],\qquad
 \sum_A\nabla_{E_A}E_A=0 .
\end{equation}
The last sum vanishes in every complete sine-cosine cutoff.

\subsection{The corrected score Hessian}\label{sec:score-hessian}

For finite real fields define
\begin{equation}\label{eq:score-hessian}
 {H^\rho}_h(v,w)=D_v\rho_w(h)-\rho_{\nabla_vw}(h),
\end{equation}
using the low scores when a field has low components.
The divergence cocycle
$D_v\rho_w-D_w\rho_v=\rho_{[v,w]}$ and
\eqref{eq:metric-torsion} make this form symmetric.

\begin{lemma}\label{lem:hessian-bound}
For $\nu$-almost every $h$, the form in \eqref{eq:score-hessian}
extends to a symmetric bounded operator on $\widehat H$, and
\begin{equation}\label{eq:uniform-hessian}
 \mathop{\mathrm{ess\,sup}}_{h\sim\nu}
           \|{H^\rho}_h\|_{\widehat H\to\widehat H}<\infty .
\end{equation}
The finite responses of the low scores also satisfy
\[
 \mathop{\mathrm{ess\,sup}}_{h\sim\nu}
       \sum_{a=0,c,s}\sum_A|D_A\sigma_a(h)|^2<\infty.
\]
\end{lemma}
The coefficient proof is given in Appendix~\ref{app:score-hessian}.

\subsection{The curvature contraction}

For finite fields set
\[
 R(v,w)=[\nabla_v,\nabla_w]-\nabla_{[v,w]},\qquad
 R(l_p,l_q)l_k=R(p,q,k)l_{p+q+k}.
\]
The curvature coefficients are
\begin{equation}\label{eq:curvature-coefficients}
 R(p,q,k)=-c(q,k)c(p,q+k)+c(p,k)c(q,p+k)
                                      +(p-q)c(p+q,k).
\end{equation}
The Ricci contraction in the orthonormal frame is denoted by $\Ric$.
Rotation equivariance makes it diagonal on complex Fourier modes.
With the curvature convention above its eigenvalue at frequency $n$ is
\[
 \Ric_n=\frac1{\omega_n}\sum_{k\in\Z}R(k,n,-n),
 \qquad \Ric_{-n}=\Ric_n.
\]
The contraction is absolutely convergent; its Fourier calculation is
given in Appendix~\ref{app:ricci}.

\begin{lemma}\label{lem:ricci-bound}
The curvature contraction converges on finite fields and extends
to a bounded operator on $\widehat H$. Its norm satisfies $\|\Ric\|\leq65/16$.
\end{lemma}

\subsection{Difference quotients for the full form}

For a finite real field $v$ write
$\Phi_t^vh=h\circ\psi_t^v$ and $U_t^vF=F\circ\Phi_t^v$.
For a smooth circle diffeomorphism $g$, write
$\operatorname{Ad}_g v=g_*v=(g'v)\circ g^{-1}$; the star on
$\operatorname{Ad}_g^*$ denotes the Hilbert-space adjoint.
The cylinder algebra is invariant under $U_t^v$, since a change of
variables changes only the smooth source weight in $J_{a,b}$.
For horizontal $v$, the finite-flow density and compression formulas
are proved in \cite[Lemma~6.1 and Proposition~6.2]{Response}.
For a finite real field $v$, including the three low directions,
integration by parts gives, for $F\in\Ccyl$, $F\geq0$,
\[
 \frac d{dt}\int U_t^vF\,d\nu
 =\int\rho_vU_t^vF\,d\nu.
\]
Hence
\[
 \left|\frac d{dt}\int U_t^vF\,d\nu\right|
 \leq\|\rho_v\|_\infty\int U_t^vF\,d\nu.
\]
Gronwall and the fact that cylinders determine the state sigma-algebra yield
\begin{equation}\label{eq:compression}
 e^{-\|\rho_v\|_\infty|t|}\nu
 \leq(\Phi_t^v)_*\nu
 \leq e^{\|\rho_v\|_\infty|t|}\nu.
\end{equation}
In particular $\|U_t^v\|_{2\to2}\leq e^{|t|\|\rho_v\|_\infty/2}$.
Continuity first on bounded continuous cylinders, then density, gives
strong continuity on $L^2(\nu)$.
For $\eta\in C_c^\infty(\R)$ and $f\in L^2(\nu)$ put
\[
 S_\eta^vf=\int_\R\eta(t)U_t^vf\,dt,
 \qquad D_vS_\eta^vf=-\int_\R\eta'(t)U_t^vf\,dt.
\]
If $F_m\in\Ccyl$ and $F_m\to f$ in $L^2$, then
\[
 \|S_\eta^vF_m-S_\eta^vf\|_2
 +\|D_vS_\eta^vF_m-D_vS_\eta^vf\|_2
 \leq C_{v,\eta}\|F_m-f\|_2\longrightarrow0.
\]
Cylinder-valued Riemann sums approximate $S_\eta^vF_m$ in this graph norm.
An approximate identity $\eta_\varepsilon$ gives, for $f$ in the flow
generator domain,
\[
 (S_{\eta_\varepsilon}^vf,D_vS_{\eta_\varepsilon}^vf)
       \longrightarrow(f,D_vf)\quad\text{in }L^2\times L^2.
\]
Thus $\Ccyl$ is a core for the scalar flow generator, also denoted by $D_v$.
Writing $w_t^v=d(\Phi_t^v)_*\nu/d\nu$, the density cocycle gives
\[
 (U_t^v)^*F=w_t^vU_{-t}^vF.
\]
Consequently $F\mapsto(w_{-t}^v)^{1/2}U_t^vF$ is a unitary group
with generator $D_v-\rho_v/2$. Bounded perturbation gives
\[
 \Dom(D_v^*)=\Dom(D_v),\qquad D_v^*F=\rho_vF-D_vF.
\]

\begin{lemma}\label{lem:full-difference-quotient}
If $\lambda>0$, $f\in L^2(\nu)$, and
$u=(\lambda-\widehat L)^{-1}f$, then, for every finite real field $v$,
\begin{equation}\label{eq:directional-regularity}
 D_vu\in\Dom(\widehat D),\qquad
 \lambda\|D_vu\|_2^2+\|\widehat D D_vu\|_2^2
                                  \leq C_{v,\lambda}\|f\|_2^2.
\end{equation}
\end{lemma}
\begin{proof}
Write $\operatorname{Ad}_t=\operatorname{Ad}_{\psi_t^v}$.
For a fixed finite field $v$, the symmetric part of its adjoint
action on $\widehat H$ is bounded. For a shift by a fixed integer $j$,
the high-frequency numerator of the symmetric matrix is
\[
 (j-n)d_{n+j}+(n+2j)d_n=(2n+j)(j^3-j).
\]
After division by $\sqrt{d_nd_{n+j}}$ it is $O_j(n^{-2})$.
Only finitely many entries involve low modes or a change of sign.
Differentiating the squared norm of the smooth pushforward, then
using density, gives
\[
 \|\operatorname{Ad}_t^*\operatorname{Ad}_t-1\|
                              \leq C_v|t|,\qquad |t|\leq1.
\]
The chain rule and \eqref{eq:compression} therefore imply that
$U_t^v$ preserves the minimal full form domain. With
$q_\lambda(a,b)=\lambda\langle a,b\rangle+
\widehat{\cE}(a,b)$,
\begin{align}
 |q_\lambda(U_t^va,U_t^vb)-q_\lambda(a,b)|
       &\leq C_{v,\lambda}|t|\,
                        q_\lambda(a,a)^{1/2}q_\lambda(b,b)^{1/2},
                        \label{eq:form-distortion}\\
 \shortintertext{and the difference quotient satisfies}
 \|(U_t^va-a)/t\|_2&\leq C_v\|\widehat Da\|_2.\notag
\end{align}
Both statements pass from cylinders by their full form density.
Set $v_t=(U_t^vu-u)/t$ and
$w_t=(U_{-t}^v-1)v_t/t$.
Using \eqref{eq:form-distortion} and the resolvent identity,
\[
 \begin{aligned}
 q_\lambda(v_t,v_t)
 &=q_\lambda(u,w_t)
       +O\!\left(q_\lambda(u,u)^{1/2}q_\lambda(v_t,v_t)^{1/2}\right)\\
 &=\langle f,w_t\rangle
       +O\!\left(q_\lambda(u,u)^{1/2}q_\lambda(v_t,v_t)^{1/2}\right).
 \end{aligned}
\]
Also $\|w_t\|_2\leq C_v q_\lambda(v_t,v_t)^{1/2}$, so
\[
 \sup_{0<|t|\leq1}q_\lambda(v_t,v_t)
 \leq C_{v,\lambda}\|f\|_2^2.
\]
The scalar closed flow derivative gives $v_t\to D_vu$ in $L^2$.
Weak compactness in the closed gradient graph proves
\eqref{eq:directional-regularity}.
\end{proof}

\subsection{Smoothing by disk automorphisms}

Let $G$ be the group of orientation-preserving disk automorphisms,
acting on $S^1$, with Haar measure $dg$. It is unimodular and has
Lie algebra spanned by $z_0,z_c,z_s$.
Fix a smooth Riemannian distance $\operatorname{dist}$ on a
neighborhood of its identity.
For $\eta\in C_c^\infty(G)$ define
\[
 S_\eta F=\int_G\eta(g)U_gF\,dg,\qquad U_gF(h)=F(h\circ g).
\]
On compact subsets of $G$ the operators $U_g$ have uniformly bounded
compression by \eqref{eq:compression}.
For $v\in H$ decompose
\begin{equation}\label{eq:low-adjoint-decomposition}
 \operatorname{Ad}_{g^{-1}}v
       =H_gv+\sum_{a=0,c,s}c_a(g,v)z_a .
\end{equation}
The operators $H_g$ are uniformly bounded on compact sets.
The functionals $c_a(g,\cdot)\in H^*$ and their first group
derivatives are bounded there; locally at the identity,
\begin{equation}\label{eq:low-projection-small}
 \|c_a(g,\cdot)\|_{H^*}\leq C\,\operatorname{dist}(g,\id).
\end{equation}
To verify these facts, differentiate the adjoint action. The horizontal
part of each low adjoint is skew for the horizontal metric, and its
low projection on $H$ is supported only at frequency $2$.
In the $(z_c,z_s)$ coordinates the two nonzero matrices are
\[
 \frac{\sqrt3}2
       \begin{pmatrix}0&1\\1&0\end{pmatrix},
 \qquad
 \frac{\sqrt3}2
       \begin{pmatrix}-1&0\\0&1\end{pmatrix}.
\]
Integrating the resulting finite low-component equations proves
\eqref{eq:low-projection-small} and the compact-set bounds.

Let $R_a$ denote differentiation along $g\exp(tz_a)$.
Let $L_a$ denote differentiation along $\exp(tz_a)g$.
For a cylinder function, the response chain rule and Haar
integration by parts give the $H$-valued identity
\begin{equation}\label{eq:group-smoothing-gradient}
 D S_\eta F
 =\int_G\eta(g)H_g^*U_gDF\,dg
 -\sum_a\int_G R_a\bigl(\eta(g)c_a(g,\cdot)\bigr)U_gF\,dg .
\end{equation}
The low derivatives and the full form bound are
\begin{align*}
 Z_aS_\eta F&=-\int_G(L_a\eta)(g)U_gF\,dg.\\
 \shortintertext{Together with \eqref{eq:group-smoothing-gradient}, this yields}
 \|S_\eta F\|_2+\|\widehat DS_\eta F\|_2
       &\leq C_\eta(\|F\|_2+\|DF\|_2).
\end{align*}
Riemann sums for a cylinder input converge in the full form norm.
For $F_m\to F$ in the horizontal form norm, the displayed bound gives
\[
 \|S_\eta(F_m-F)\|_2+\|\widehat DS_\eta(F_m-F)\|_2\longrightarrow0.
\]
Hence $S_\eta:\Dom(D)\to\Dom(\widehat D)$ is bounded.

Apply this to a coefficient polynomial $P$ from
Lemma~\ref{lem:coefficient-domain}. Its finite low responses are bounded.
Integration by parts and the scalar flow core identify them with
$Z_0P,Z_cP,Z_sP$.

For an approximate identity $\eta_\varepsilon$ on $G$, the full
response chain rule gives
\[
 \widehat DS_{\eta_\varepsilon}P
 =\int_G\eta_\varepsilon(g)\operatorname{Ad}_{g^{-1}}^*
             U_g(DP,Z_0P,Z_cP,Z_sP)\,dg.
\]
Strong continuity and the compact-set bounds for the adjoint action yield
\[
 S_{\eta_\varepsilon}P\to P,\qquad
 \widehat D S_{\eta_\varepsilon}P
       \to (D_jP,Z_0P,Z_cP,Z_sP)
       \quad\text{in }L^2 .
\]
Closure therefore gives
\[
 P\in\Dom(\widehat D),\qquad
 \rho_v\in\Dom(\widehat D)
 \quad\text{for every finite real Fourier field }v.
\]

\subsection{Covariant derivatives and two commutators}

Write
\[
 \mathcal V_{\mathrm{fin}}(\widehat H)
 =\left\{\sum_{A\in J_0}F_AE_A:
       J_0\text{ finite},\ F_A\in\Ccyl\right\}.
\]
For $V\in\mathcal V_{\mathrm{fin}}(\widehat H)$, its covariant derivative is
\[
 (\nabla_A V)_B=D_AV_B-\sum_Cc_{AB}^{C}V_C.
\]
The divergence and exterior derivative are, respectively,
\begin{align*}
 \widehat\delta V&=\sum_A(\widehat\rho_AV_A-D_AV_A),\\
 (dV)_{AB}&=(\nabla_AV)_B-(\nabla_BV)_A.
\end{align*}
The exterior norm is $\|dV\|^2=\frac12\sum_{A,B}|(dV)_{AB}|^2$.
Finite integration by parts makes $\nabla$ and $d$ closable. We use
their minimal closures, with
\begin{align*}
 \Dom(\nabla)
 &=\overline{\mathcal V_{\mathrm{fin}}(\widehat H)}
                  ^{\,\|V\|_2+\|\nabla V\|_2},\\
 \shortintertext{and}
 \Dom(d)
 &=\overline{\mathcal V_{\mathrm{fin}}(\widehat H)}
                  ^{\,\|V\|_2+\|dV\|_2}.
\end{align*}
The connection tail of each finite covector is square summable.

For $F\in\Ccyl$ the vector $\widehat DF$ belongs to the minimal
covariant domain, and $d\widehat DF=0$.
In the normalized complex frame $l_p/\sqrt{\omega_p}$, the
covariant Hessian of an inverse angle has coefficient
\[
 \frac{i\{q+c(p,q)\}}{\sqrt{\omega_p\omega_q}}\,
                         e^{i(p+q)\theta}.
\]
For horizontal $p,q$ its numerator has absolute value at most
$C\min(|p|,|q|)$. Thus
\[
 \sum_{|p|,|q|\geq2}
 \frac{|q+c(p,q)|^2}{\omega_p\omega_q}
 \leq C\sum_{2\leq p\leq q}\frac{p^2}{p^3q^3}
 \leq C\sum_{p\geq2}\frac1{p^3}<\infty.
\]
The finitely many low rows and columns have square-summable tails
by \eqref{eq:full-connection}. The Hilbert--Schmidt bound is uniform
in the inverse angle.
The Stieltjes formula \eqref{eq:stieltjes} and cylinder chain rules
give the same assertion for each $F$.
Let $P_M$ project onto the full modes of frequency at most $M$.
The pointwise bound $|(\widehat DF)_n|\leq C_F/\sqrt{\omega_n}$
gives
\begin{equation}\label{eq:cylinder-covariant-closure}
 \|[\nabla,P_M]\widehat DF\|_{L^2(\widehat H\otimes\widehat H)}^2
 \leq\frac{C_F}{M}.
\end{equation}
The equal-sign and opposite-sign sums satisfy
\begin{align*}
 \sum_{2\leq b\leq M,\ a+b>M}\frac1{a^3b}
 &\leq C\left(\frac{\log(e+M)}{M^2}
          +\frac1M\sum_{r\geq1}\frac1{r^2}\right)\leq\frac C M,\\
 \shortintertext{and}
 \sum_{b>M,\ a\geq\max(2,b-M)}\frac1{a^3b}
 &\leq C\sum_{r\geq1}\frac1{(M+r)\max(2,r)^2}
       \leq\frac C M.
\end{align*}
The adjacent low-direction terms satisfy the same bound.
Thus $P_M\widehat DF\to\widehat DF$ in the minimal covariant
graph norm. Torsion freeness gives $d\widehat DF=0$.

Write the connection part as $\Xi_A$:
\[
 \nabla_A=D_A+\Xi_A,\qquad
 (\Xi_AV)_B=-\sum_Cc_{AB}^{C}V_C,\qquad \Xi_A^*=-\Xi_A.
\]
These matrices are constant on the fibre. For a finite field
$v=\sum_Av_AE_A$, write $\Xi_v=\sum_Av_A\Xi_A$.
Choose $\chi\in C^\infty([0,\infty))$ with $0\leq\chi\leq1$,
$\chi=1$ on $[0,1]$, and $\chi=0$ on $[2,\infty)$.
For $R\geq2$ define the finite-rank even Fourier multiplier $Q_R$ by
\begin{equation}\label{eq:slow-cutoff}
 w_R(n)=\chi\!\left(\frac{\log\log(e+|n|)}R\right)
       \quad(|n|\geq2),\qquad w_R(0)=w_R(\pm1)=1 .
\end{equation}
Set
\begin{equation}\label{eq:commutator-definitions}
 \begin{aligned}
 C_{R,A}&=[\Xi_A,Q_R],&
 (C_Rv)_A&=C_{R,A}v,\\
 {T_R^{\nabla}}&=\sum_A
     \bigl([\Xi_A,[\Xi_A,Q_R]]
                         -\widehat\rho_A[\Xi_A,Q_R]\bigr).
 \end{aligned}
\end{equation}
The sums are first evaluated on finite Fourier vectors, with complete
real sine--cosine pairs in the direction index $A$; the resulting
bounded extension is used for $T_R^\nabla$.

\begin{lemma}\label{lem:two-commutator-bounds}
The sums in \eqref{eq:commutator-definitions} have well-defined
limits on finite columns, and
\begin{equation}\label{eq:two-commutator-bounds}
 \|C_R\|_{\widehat H\to\widehat H\otimes\widehat H}
       \leq\frac C R,\qquad
 \mathop{\mathrm{ess\,sup}}_{h\sim\nu}
                  \|{T_R^{\nabla}}(h)\|\leq\frac{C_\kappa}R .
\end{equation}
\end{lemma}

The proof of these two estimates is in
Appendix~\ref{app:commutators}.

\subsection{The full resolvent and first-order domain equality}

Lemmas~\ref{lem:hessian-bound} and \ref{lem:ricci-bound} control
the zero-order term in the Hodge identity. Finite integration by parts,
metricity, and torsion freeness give
\begin{equation}\label{eq:full-hodge}
 \|\widehat\delta V\|_2^2+\|dV\|_2^2
 =\|\nabla V\|_2^2+
       \int_\mathcal X(({H^\rho}_h+\Ric)V,V)_{\widehat H}\,d\nu(h)
\end{equation}
for finite $\Ccyl$-valued covectors.
At every complete real Fourier cutoff, the contractions use
\begin{align*}
 D_A\widehat\rho_B-\sum_Cc_{AB}^C\widehat\rho_C
       &={H^\rho}_h(E_A,E_B),\\
 \shortintertext{and}
 [\nabla_A,\nabla_B]-\nabla_{[E_A,E_B]}&=R(E_A,E_B).\\
 \shortintertext{The two-form norm expands as}
 \|dV\|_2^2
       &=\|\nabla V\|_2^2
          -\Rea\sum_{A,B}\langle(\nabla_AV)_B,(\nabla_BV)_A\rangle.
\end{align*}
Integration by parts, the zero connection trace, and metricity give
\[
 \|\widehat\delta V\|_2^2
       -\Rea\sum_{A,B}\langle(\nabla_AV)_B,(\nabla_BV)_A\rangle
 =\int_\mathcal X(({H^\rho}_h+\Ric)V,V)_{\widehat H}\,d\nu(h).
\]
For each selected component $V_B$, the tails satisfy
\[
 \sum_A\|\Xi_AE_B\|^2<\infty,\qquad
 \sum_{j\in{J_H}}\|\rho_jD_jV_B\|_2
       \leq C_V\sum_{n\geq2}\alpha_n<\infty.
\]
Together with the absolutely convergent Ricci traces in
Lemma~\ref{lem:ricci-bound}, these bounds pass the finite contractions
to the displayed identities and give \eqref{eq:full-hodge}.

Let $\widehat\delta_0$ be the restriction of divergence to
$\mathcal V_{\mathrm{fin}}(\widehat H)$. Equation \eqref{eq:full-hodge} gives
\begin{align*}
 \Dom(\nabla)&\subset\Dom(\overline{\widehat\delta_0}),\\
 \shortintertext{and}
 \|\overline{\widehat\delta_0}V\|_2^2
 &\leq\|\nabla V\|_2^2+C_\kappa\|V\|_2^2,
       \qquad V\in\Dom(\nabla).
\end{align*}
Write $\widehat D_{\max}=\widehat\delta_0^*$.

\begin{lemma}\label{lem:full-resolvent-regularity}
For every $u\in\Dom(\widehat L)$,
$\widehat Du$ belongs to the minimal covariant domain and
\begin{equation}\label{eq:full-resolvent-identity}
 \|\widehat Lu\|_2^2
 =\|\nabla\widehat Du\|_2^2+
     \int_\mathcal X(({H^\rho}_h+\Ric)\widehat Du,\widehat Du)_{\widehat H}
                                                        \,d\nu(h).
\end{equation}
Consequently $\widehat D=\widehat D_{\max}$.
\end{lemma}
\begin{proof}
Fix $\lambda>0$ and $f\in L^2(\nu)$. Put $u=(\lambda-\widehat L)^{-1}f$, $p=\widehat Du$.
We first derive a weak vector resolvent equation for $p$.
Lemma~\ref{lem:full-difference-quotient} puts each $p_A$
in the full minimal scalar form domain. Thus every finite
frequency truncation of $p$ belongs to the minimal covariant domain.
Choose $F_m\in\Ccyl$ converging to $u$ in the full form norm. Then
\[
 \widehat DF_m\longrightarrow p,\qquad
 d\widehat DF_m=0,
\]
so closedness of $d$ gives $p\in\Dom(d)$ and $dp=0$.

For a finite cylinder vector $\phi$, the product rule gives
$\widehat\delta\phi\in\Dom(\widehat D)$: each term is a product
of cylinders, finite derivatives, and finite score coefficients.

Let $J_0$ be the finite support of $\phi$. We use
\[
 \|\widehat D\phi\|_2^2
   =\sum_{B\in J_0}\|\widehat D\phi_B\|_2^2.
\]
All connection sums use complete real Fourier cutoffs.
For $n,m\geq2$, \eqref{eq:full-connection} gives
\[
 \frac{c(n,m)\sqrt{d_{n+m}}}{\sqrt{d_nd_m}}
 =-\frac{(2n+m)\sqrt{d_m}}{\sqrt{d_nd_{n+m}}}
 =O_m(n^{-2}).
\]
The opposite-sign high tail is zero, and the low-column symbols
are $O(n^{-2})$. Consequently, for each full frame index $B$,
\begin{equation}\label{eq:connection-column-bounds}
 m_B:=\sum_A\|\Xi_AE_B\|_{\widehat H}^2<\infty,
 \qquad
 \mathop{\mathrm{ess\,sup}}_h
 \left\|\sum_A\widehat\rho_A(h)\Xi_AE_B\right\|_{\widehat H}
 \leq C_B.
\end{equation}
The second bound uses the $O_B(n^{-3/2})$ score-connection tail
proved in Appendix~\ref{app:commutators}. Rotation equivariance and
skew-adjointness give
\[
 \sum_{A:E_A\in P_M\widehat H}\Xi_A^2E_B
 =-\left(\sum_{A:E_A\in P_M\widehat H}
                 \|\Xi_AE_B\|_{\widehat H}^2\right)E_B
 \longrightarrow-m_BE_B.
\]
Expansion of the covariant Laplacian now gives, componentwise,
\begin{align*}
 \nabla^*\nabla\phi
 &=-\widehat L\phi-2\sum_A\Xi_AD_A\phi
       -\sum_A\Xi_A^2\phi
       +\sum_A\widehat\rho_A\Xi_A\phi.\\
 \shortintertext{The first-order remainder satisfies}
 \left\|\sum_A\Xi_AD_A\phi\right\|_2
 &\leq\left(\sum_{B\in{J_0}}m_B\right)^{1/2}
                                      \|\widehat D\phi\|_2.
\end{align*}
Thus $L_{\mathrm{vec}}\phi=(\nabla^*\nabla+{H^\rho}+\Ric)\phi$ belongs to
$L^2(\nu;\widehat H)$, and
\begin{equation}\label{eq:vector-resolvent-remainder}
 \|L_{\mathrm{vec}}\phi+\widehat L\phi\|_2
 \leq C_{{J_0}}\bigl(\|\widehat D\phi\|_2+\|\phi\|_2\bigr).
\end{equation}
Polarizing \eqref{eq:full-hodge} and using
\eqref{eq:cylinder-covariant-closure} yields
\[
 \langle\widehat DF,L_{\mathrm{vec}}\phi\rangle
 =\langle\widehat DF,\widehat D(\widehat\delta\phi)\rangle
 \quad(F\in\Ccyl).
\]
Approximation of $u$ in the full form norm and the resolvent equation
with test $\widehat\delta\phi$ give
\begin{equation}\label{eq:vector-weak-resolvent}
 \lambda\langle p,\phi\rangle+\langle p,L_{\mathrm{vec}}\phi\rangle
                                      =\langle f,\widehat\delta\phi\rangle.
\end{equation}
For form tests supported on ${J_0}$, define the left pairing by
\[
 \langle p,L_{\mathrm{vec}}\phi\rangle
 =\sum_{B\in{J_0}}\langle\widehat Dp_B,\widehat D\phi_B\rangle
       +\langle p,L_{\mathrm{vec}}\phi+\widehat L\phi\rangle .
\]
Equation \eqref{eq:vector-resolvent-remainder} gives
\begin{align*}
 |\langle p,L_{\mathrm{vec}}\phi\rangle|
 &\leq\left(\sum_{B\in{J_0}}\|\widehat Dp_B\|_2^2\right)^{1/2}
                    \|\widehat D\phi\|_2
   +C_{{J_0}}\|p\|_2(\|\widehat D\phi\|_2+\|\phi\|_2),\\
 \shortintertext{while the divergence obeys}
 \|\widehat\delta\phi\|_2
 &\leq C_{{J_0},\kappa}
                       (\|\widehat D\phi\|_2+\|\phi\|_2).
\end{align*}
Both sides of \eqref{eq:vector-weak-resolvent} therefore extend by
componentwise full form approximation.

To obtain an estimate uniform in $R$, take $\phi=Q_R^2p$ and use
\begin{equation}\label{eq:two-commutator-identity}
 \langle p,\nabla^*\nabla(Q_R^2p)\rangle
 =\|\nabla(Q_Rp)\|_2^2
     -2\langle C_Rp,\nabla(Q_Rp)\rangle
     -\langle p,{T_R^{\nabla}}Q_Rp\rangle.
\end{equation}
For finite vectors it follows by expansion from
$[\nabla^*\nabla,Q_R]=-{T_R^{\nabla}}-2C_R^*\nabla$.
The matrices $C_{R,A}$ are self-adjoint, since
$\Xi_A$ is skew-adjoint.

Let $P_S$ be the full-mode projection of cutoff $S$, with $S$
containing the support of $Q_R$, and first apply the identity to $P_Sp$.
Then $Q_R^2P_Sp=Q_R^2p$ and $Q_RP_Sp=Q_Rp$.
Integration by parts in the components selected by $Q_R$ and
\eqref{eq:vector-resolvent-remainder} justify the principal pairing,
while
\[
 \|C_R(P_Sp-p)\|_2+\|{T_R^{\nabla}}(P_Sp-p)\|_2
       \leq\frac{C_\kappa}R\|P_Sp-p\|_2\longrightarrow0.
\]
Thus $S\to\infty$ gives \eqref{eq:two-commutator-identity}.

Let $a_R=\|\nabla(Q_Rp)\|_2$. Since
\[
 \nabla(Q_R^2p)=Q_R\nabla(Q_Rp)+C_RQ_Rp,
\]
the Hodge identity \eqref{eq:full-hodge} applied to $Q_R^2p$ gives
\[
 \|\widehat\delta(Q_R^2p)\|_2\leq a_R+C_\kappa\|p\|_2.
\]
Insert this and \eqref{eq:two-commutator-identity} into
\eqref{eq:vector-weak-resolvent}.
The boundedness of $({H^\rho}+\Ric)$ and
\eqref{eq:two-commutator-bounds} yield
\[
 a_R^2
 \leq C_\kappa\|p\|_2^2+\frac C R\|p\|_2a_R
                         +\|f\|_2(a_R+C_\kappa\|p\|_2).
\]
Absorbing the terms linear in $a_R$ gives
\[
 \sup_Ra_R^2\leq C_\kappa(\|f\|_2^2+\|p\|_2^2).
\]
Weak closedness puts $p$ in the minimal covariant domain. Moreover,
\begin{align*}
 \|\nabla Q_Rp-\nabla p\|_2
 &\leq\|(Q_R-1)\nabla p\|_2+\frac C R\|p\|_2\longrightarrow0.\\
 \shortintertext{Applying the same argument to $Q_R^2$ gives}
 Q_R^2p&\longrightarrow p
             \quad\text{in the minimal covariant graph norm}.\\
 \shortintertext{Closedness therefore yields}
 \widehat\delta Q_R^2p&\longrightarrow
       \overline{\widehat\delta_0}p
       =\widehat D^*p=-\widehat Lu,\qquad
 dQ_R^2p\longrightarrow dp=0.
\end{align*}
Passing to the limit in \eqref{eq:full-hodge} proves
\eqref{eq:full-resolvent-identity}.

To identify the maximal gradient, observe that
$p=\widehat Du\in\Dom(\overline{\widehat\delta_0})
=\Dom(\widehat D_{\max}^*)$ for every
$u\in\Dom(\widehat L)$.
As $\widehat D\subset\widehat D_{\max}$,
\[
 \widehat L\subset-\widehat D_{\max}^*\widehat D_{\max}.
\]
Both operators are self-adjoint. They coincide, and equality of
their closed forms proves $\widehat D=\widehat D_{\max}$.
\end{proof}

Return to the horizontal gradient.
Let $\delta_0$ be $\delta$ restricted to finite
$\Ccyl$-valued horizontal vectors and put $D_{\max}=\delta_0^*$.

\begin{proposition}\label{prop:horizontal-domain}
The horizontal minimal and maximal derivatives coincide:
\begin{equation}\label{eq:horizontal-domain-equality}
 D=D_{\max}.
\end{equation}
\end{proposition}
\begin{proof}
The single-flow core assertion following \eqref{eq:compression} gives
\[
 \Dom(D_{\max})
 =\left\{u\in\bigcap_{j\in J_H}\Dom(D_j):
             \sum_{j\in J_H}\|D_ju\|_2^2<\infty\right\}.
\]
Let $u\in\Dom(D_{\max})$. Formula
\eqref{eq:group-smoothing-gradient} holds distributionally for $u$,
with $D_{\max}u$ in the first integral.
Put $\Phi_gh=h\circ g$ and $w_g=d(\Phi_g)_*\nu/d\nu$.
The low-flow cocycle and the bounded full gradients of the low scores give
\begin{align*}
 w_{\exp(tz_a)}
 &=\exp\left(\int_0^t U_{-s}^{z_a}\sigma_a\,ds\right),\\
 \|\widehat D w_{\exp(tz_a)}\|_\infty
 &\leq e^{|t|\|\sigma_a\|_\infty}
  \int_{-|t|}^{|t|}\|\operatorname{Ad}_{\exp(sz_a)}\|
                         \|\widehat D\sigma_a\|_\infty\,ds.\\
 \shortintertext{Using finite products of low flows, we obtain}
 \sup_{g\in\supp\eta}
 \bigl(\|w_g\|_\infty+\|w_g^{-1}\|_\infty
                      +\|\widehat Dw_g\|_\infty\bigr)&<\infty.
\end{align*}
For $\phi\in\Ccyl$, the closed product rule gives
\[
 \phi w_g\in\Dom(\widehat D),\qquad
 \widehat D(\phi w_g)=w_g\widehat D\phi+\phi\widehat Dw_g.
\]
Change variables in a finite cylinder-vector test $V$, use the
divergence cocycle and \eqref{eq:low-adjoint-decomposition}, and
integrate the low group derivative by parts. This gives
\begin{align*}
 \langle S_\eta u,\delta_0V\rangle
 ={}&\left\langle
       \int_G\eta(g)H_g^*U_gD_{\max}u\,dg
       -\sum_a\int_G R_a(\eta(g)c_a(g,\cdot))U_gu\,dg,V
                         \right\rangle.\\
 \shortintertext{The low derivatives are}
 Z_aS_\eta u&=-\int_G(L_a\eta)(g)U_gu\,dg.\\
 \shortintertext{Together these formulas imply}
 \|S_\eta u\|_2+\|\widehat D_{\max}S_\eta u\|_2
       &\leq C_\eta(\|u\|_2+\|D_{\max}u\|_2).
\end{align*}
The scalar flow core and the compact-set bounds justify the test
identities and square summation. Thus $S_\eta u\in\Dom(\widehat D_{\max})$.

Choose $\eta_\varepsilon$ supported within distance $C\varepsilon$
of the identity, with integral one.
In local coordinates it has size $O(\varepsilon^{-3})$ and
$\|R_a\eta_\varepsilon\|_1=O(\varepsilon^{-1})$.
By \eqref{eq:low-projection-small},
\[
 \|R_a(\eta_\varepsilon c_a)\|_{L^1(G;H^*)}\leq C,\qquad
 \int_G R_a(\eta_\varepsilon c_a)\,dg=0 .
\]
The second integral in \eqref{eq:group-smoothing-gradient} is
therefore bounded in $L^2(\nu;H)$ by
\[
 C\sup_{\operatorname{dist}(g,\id)\leq C\varepsilon}
                       \|U_gu-u\|_2\longrightarrow0 .
\]
The first integral tends to $D_{\max}u$ by strong continuity
of $U_g$ on vector-valued $L^2$ and of $H_g^*$ on $H$.
Consequently
\[
 S_{\eta_\varepsilon}u\to u,\qquad
 D_{\max}S_{\eta_\varepsilon}u\to D_{\max}u
                                      \quad\text{in }L^2 .
\]
Lemma~\ref{lem:full-resolvent-regularity} puts each mollification
in the full minimal domain, hence in the horizontal minimal
domain. Closedness of $D$ proves $D_{\max}\subset D$.
The reverse inclusion follows from \eqref{eq:ibp}.
\end{proof}

\subsection{Horizontal divergence}\label{sec:horizontal-divergence}

Let $P_H$ project onto the original horizontal
space $H$.
For a horizontal field $V$, put
\begin{align*}
 \|\nabla_H V\|_2^2&=\sum_{j\in{J_H}}\|\nabla_jV\|_2^2,\\
 \shortintertext{and define}
 T_av&=\Xi_vz_a,\qquad K_a=P_H T_a^*P_H.
\end{align*}
The Fourier symbols of $T_a$ have order $-2$, outside finitely many
indices. In particular $T_a$ and $K_a$ are bounded, and for every
fixed finite Fourier direction $p$,
\begin{equation}\label{eq:HD.1}
\|[\Xi_p,K_a]\|_{\widehat H\to\widehat H}<\infty.
\end{equation}
Indeed $\Xi_p$ is a finite shift with order-one symbol, while $K_a$
is a fixed shift of order $-2$. Their commutator has order $-2$:
the leading first-order factor multiplies one finite difference of
the order-$-2$ symbol. The finitely many low rows and columns cause
only a bounded finite-rank correction.

Put $\mathcal V_{\mathrm{fin}}(H)=P_H\mathcal V_{\mathrm{fin}}(\widehat H)$.
For finitely supported $\Phi=(\Phi_j)_{j\in J_H}$ with
$\Phi_j\in\mathcal V_{\mathrm{fin}}(\widehat H)$, set
\[
 \nabla_{H,0}^{*}\Phi
    =P_H\sum_{j\in J_H}(\rho_j\Phi_j-\nabla_j\Phi_j).
\]
The minimal horizontal covariant domain is
\[
 \Dom(\nabla_{H,\min})
 =\overline{\mathcal V_{\mathrm{fin}}(H)}
             ^{\,\|V\|_2+\|\nabla_HV\|_2}.
\]
For $V\in L^2(\nu;H)$ and $G\in L^2(\nu;H\otimes\widehat H)$,
the maximal derivative is defined by
\[
 \nabla_{H,\max}V=G
 \quad\Longleftrightarrow\quad
 \langle G,\Phi\rangle_2
   =\langle V,\nabla_{H,0}^{*}\Phi\rangle_2
 \quad\text{for every finite }\Phi.
\]

\begin{proposition}\label{prop:horizontal-upper-hodge}
The horizontal covariant domains satisfy
\[
 \Dom(\nabla_{H,\min})=\Dom(\nabla_{H,\max})
       \subset\Dom(\delta).
\]
There is $C_\kappa<\infty$ such that, for every
$V\in\Dom(\nabla_{H,\max})$,
\begin{equation}\label{eq:HD.2}
\|\delta V\|_2^2\le
 C_\kappa\big(\|\nabla_H V\|_2^2+\|V\|_2^2\big).
\end{equation}
\end{proposition}

\begin{proof}
Start with a finite horizontal $V$, extended by zero on the low
components. Its full divergence is its horizontal divergence.
In the full Hodge identity, decompose the two-form norm into the
horizontal--horizontal, mixed, and low--low blocks. Metricity gives
\[
 (\nabla_iV)_a=-\langle V,\Xi_i z_a\rangle
                 =-(T_a^*V)_i,
\]
and hence the mixed block is
\[
 \sum_a\|P_H\nabla_aV+K_aV\|_2^2.
\]
Subtracting it from $\sum_a\|\nabla_aV\|_2^2$ cancels the
uncontrolled horizontal parts of the low-direction derivative.
Since $V$ has zero low components, $(1-P_H)\nabla_aV$ is a bounded
zero-order map of $V$. In particular
\begin{align*}
 \sum_a\bigl(\|\nabla_aV\|_2^2-\|P_H\nabla_aV+K_aV\|_2^2\bigr)
 ={}&\sum_a\bigl(\|(1-P_H)\nabla_aV\|_2^2-\|K_aV\|_2^2
                 -2\Rea\langle\nabla_aV,K_aV\rangle\bigr)\\
 \leq{}& C\|V\|_2^2+2\sum_a|\langle\nabla_aV,K_aV\rangle|.
\end{align*}
Drop the other nonnegative two-form terms and use the bounded
Hessian and Ricci operators to obtain
\begin{equation}\label{eq:HD.3}
\|\delta V\|_2^2\le \|\nabla_H V\|_2^2+C_\kappa\|V\|_2^2
       +2\sum_a|\langle \nabla_aV,K_aV\rangle_{L^2}|.
\end{equation}

Each $z_a$ is a finite real linear combination of brackets of fixed
horizontal Fourier fields. In the complex notation this follows
from $[l_2,l_{-2}]=4i l_0$, $[l_3,l_{-2}]=5i l_1$, and their
conjugates; real and imaginary parts give real horizontal brackets.
Choosing a fixed finite collection of real horizontal fields
$p_\ell,q_\ell$ and real constants $b_{a\ell}$, we have
\begin{align*}
 z_a&=\sum_\ell b_{a\ell}[p_\ell,q_\ell],\\
 \shortintertext{and hence}
 \nabla_a&=\sum_\ell b_{a\ell}
       \bigl([\nabla_{p_\ell},\nabla_{q_\ell}]
                         -R(p_\ell,q_\ell)\bigr).
\end{align*}
Covariant integration by parts and the product rule give
\begin{align*}
 \langle[\nabla_p,\nabla_q]V,K_aV\rangle
 &=\langle\nabla_qV,(-\nabla_p+\rho_p)K_aV\rangle
   -\langle\nabla_pV,(-\nabla_q+\rho_q)K_aV\rangle,\\
 \shortintertext{where}
 \nabla_p(K_aV)&=K_a\nabla_pV+[\Xi_p,K_a]V.
\end{align*}
Since $\rho_p,\rho_q$ are bounded, \eqref{eq:HD.1} implies
\[
 |\langle[\nabla_p,\nabla_q]V,K_aV\rangle|
 \leq C_{p,q,\kappa}\bigl(\|\nabla_HV\|_2^2+\|V\|_2^2\bigr).
\]
The fixed curvature operators are bounded by
\eqref{eq:full-connection}--\eqref{eq:curvature-coefficients}.
Combining this bound with the finite bracket representation
therefore yields
\begin{equation}\label{eq:HD.4}
 |\langle\nabla_aV,K_aV\rangle|
 \leq C_\kappa\bigl(\|\nabla_HV\|_2^2+\|V\|_2^2\bigr),
\end{equation}
and proves \eqref{eq:HD.2}. 

To pass to the maximal domain, let $V\in\Dom(\nabla_{H,\max})$.
Apply the finite-rank slow Fourier cutoffs $Q_R$, which preserve $H$,
are the identity on the three low modes, and tend strongly to
$I_{\widehat H}$.
Lemma~\ref{lem:two-commutator-bounds} gives
\begin{align}
 \nabla_j(Q_RV)&=Q_R\nabla_jV+[\Xi_j,Q_R]V,\notag\\
 \shortintertext{and}
 \sum_j\|[\Xi_j,Q_R]V\|_2^2
 &\le C R^{-2}\|V\|_2^2.\label{eq:HD.5}
\end{align}
Skew-adjointness and \eqref{eq:connection-column-bounds} imply
\begin{align*}
 \sum_{j\in{J_H}}\|Q_R\Xi_j\|_{\widehat H\to\widehat H}^2
 &\leq\sum_{B:Q_RE_B\ne0}\sum_{j\in{J_H}}
                            \|\Xi_jE_B\|_{\widehat H}^2
 \leq\sum_{B:Q_RE_B\ne0}m_B<\infty.\\
 \shortintertext{Since}
 D_j(Q_RV)&=Q_R\nabla_jV-Q_R\Xi_jV,\\
 \shortintertext{we conclude that}
 \sum_j\|D_j(Q_RV)\|_2^2
 &\leq2\|\nabla_HV\|_2^2+C(R)\|V\|_2^2.
\end{align*}
Every selected scalar coordinate therefore belongs to
$\Dom(D_{\max})=\Dom(D)$. For fixed $R$, let
$J_R=\{B:Q_RE_B\ne0\}$ and choose
$V_{R,m}\in\mathcal V_{\mathrm{fin}}(H)$ supported on $J_R$ such that
\[
 \|V_{R,m}-Q_RV\|_2^2
   +\sum_{j\in J_H}\|D_j(V_{R,m}-Q_RV)\|_2^2
          \longrightarrow0.
\]
With $M_R=\sum_{B\in J_R}m_B<\infty$, the connection column bound gives
\begin{align*}
 \|\nabla_H(V_{R,m}-Q_RV)\|_2^2
 &\leq2\sum_{j\in J_H}\|D_j(V_{R,m}-Q_RV)\|_2^2
          +2M_R\|V_{R,m}-Q_RV\|_2^2\\
 &\longrightarrow0.
\end{align*}
Finally, \eqref{eq:HD.5} gives
\[
 \|Q_RV-V\|_2+\|\nabla_H(Q_RV-V)\|_2\longrightarrow0,
\]
which proves minimal--maximal equality. For finite cylinder-vector
approximations $V_m\to V$ in this graph norm, \eqref{eq:HD.2} gives
\[
 \|\delta V_m-\delta V_n\|_2^2
 \leq C_\kappa\bigl(\|\nabla_H(V_m-V_n)\|_2^2+\|V_m-V_n\|_2^2\bigr)
 \longrightarrow0.
\]
Closedness of $\delta$ gives $V\in\Dom(\delta)$ and \eqref{eq:HD.2}.

\end{proof}

\section{A conditional Gaussian representation}\label{sec:gaussian-lift}

We couple the welding to two boundary Gaussian fields and use a capacity
tilt to identify their marginal and response calculus.

\subsection{The trace kernel and the capacity}

On the full-measure removable locus, choose the normalized conformal pair
$f=f_h^\circ$, $g$ with
\[
 f(0)=0,\qquad f'(0)=1,\qquad g(\infty)=\infty,
 \qquad f=g\circ h\quad\text{on }S^1.
\]
Set $a=g'(\infty)$ and $c(h)=\log|a|\geq0$, and extend these choices
measurably to the null complement.

Let $\varphi$ be the whole-plane Gaussian free field trace with exterior
harmonic extension zero at infinity \cite[Section 2.1]{BJ}. Denote its
conditional law by $P_h$; its covariance on zero-mass signed measures is
$-\log|x-y|$. Put $U=\varphi(0)$ using the interior harmonic extension.

Let $\omega_0$ and $\omega_\infty$ be harmonic measure from $0$ and
$\infty$. Conformal exhaustion gives, on the present Jordan curves,
\begin{align}
 \iint\log|x-y|\,d\omega_0(x)d\omega_0(y)&=0,\notag\\
 \shortintertext{while the exterior energy is}
 \iint\log|x-y|\,d\omega_\infty(x)d\omega_\infty(y)&=c(h),\notag\\
 \shortintertext{and the mixed energy has the same value:}
 \iint\log|x-y|\,d\omega_0(x)d\omega_\infty(y)&=c(h).\notag\\
 \shortintertext{Hence, since $U$ is paired with $\omega_0-\omega_\infty$,}
 \operatorname{Var}_{P_h}(U)&=-0-c(h)+2c(h)=c(h).
 \label{eq:N.3}
\end{align}
Consequently, the probability kernel
\begin{equation}\label{eq:N.4}
                 dP_h^Q=e^{QU-Q^2c(h)/2}\,dP_h
\end{equation}
is normalized. With $P_0$ denoting subtraction of the circle average, set
\begin{align*}
 2\Rea Y&=P_0(\varphi\circ f+Q\log|f'|),\\
 \shortintertext{and}
 2\Rea Z&=P_0(\varphi\circ g+Q\log|g'|).
\end{align*}
Write $Y(z)=\sum_{n\geq1}Y_nz^n$ and
$Z(w)=\sum_{n\geq1}Z_nw^{-n}$. The conditional trace kernel is
\[
 K_h^Q=\operatorname{Law}_{P_h^Q}(Y,Z).
\]

The covariance identities
\begin{align*}
 \operatorname{Cov}(\varphi(f(z))-U,U)&=-\log|f(z)/z|,\\
 \shortintertext{and}
 \operatorname{Cov}(\varphi(g(w)),U)&=-\log|g(w)/(aw)|
\end{align*}
and the Gaussian shift formula give the conditional means of the
real fields:
\begin{equation}\label{eq:N.6}
 Q\log\left|\frac{zf'(z)}{f(z)}\right|,
 \qquad Q\log\left|\frac{wg'(w)}{g(w)}\right|.
\end{equation}
We use $\mathbb G$ for the centered Neumann probability, normalized by
\[
 \mathbb E_{\mathbb G}Y_n=0,\qquad
 \mathbb E_{\mathbb G}[Y_nY_m]=0,\qquad
 \mathbb E_{\mathbb G}[Y_n\overline{Y_m}]=\frac{\delta_{nm}}n.
\]
On the removable locus, the univalent pre-Schwarzian and
logarithmic-derivative estimates give the uniform bounds
\begin{align*}
 \sup_h\sup_{z\in\D}(1-|z|^2)
 \left|\partial_z\log\frac{zf'(z)}{f(z)}\right|&\leq C,\\
 \shortintertext{and}
 \sup_h\sup_{z\in\D}(1-|z|^2)
 \left|\partial_z\log\frac{z^{-1}g'(z^{-1})}{g(z^{-1})}\right|&\leq C.
\end{align*}
The conditional real covariance is bounded by the Neumann covariance:
for every finite real linear functional $\Lambda$ of either field,
\begin{align*}
 \operatorname{Var}_{K_h^Q}(\Lambda(Y))
       &\leq\operatorname{Var}_{\mathbb G}(\Lambda(Y)),\\
 \shortintertext{and}
 \operatorname{Var}_{K_h^Q}(\Lambda(Z))
       &\leq\operatorname{Var}_{\mathbb G}(\Lambda(Z)).
\end{align*}
Indeed, the inverse covariance form is the sum of the two harmonic
Dirichlet energies.

\subsection{Capacity estimates}

We first relate the capacity to the logarithmic coefficients of $f$.
Put $p(z)=\log(f(z)/z)=\sum_{n\geq1}p_nz^n$ and
$E_p=\sum_{n\ge1}n|p_n|^2$.
The potential of $\omega_0-\omega_\infty$ equals
$-\operatorname{Re}p$ in the interior, up to a constant.
Its total Dirichlet energy, normalized by $1/(2\pi)$, is $c$.
Thus
\begin{equation}\label{eq:N.8}
 \frac{E_p}{2}
 =\frac1{2\pi}\int_\D|p'(z)|^2\,dA(z)
 \leq c.
\end{equation}

The Koebe transform and growth theorem give
\[
 \left|\frac{f'(z)}{f(z)}\right|
 \leq\frac{1+|z|}{|z|(1-|z|)}.
\]
Hence $(1-|z|^2)|p'(z)|\leq6$ for $|z|\geq1/2$.
The maximum principle applied to $p'$ gives the bound $8$ inside
that disk. Therefore
\begin{equation}\label{eq:N.9}
\sup_{z\in\mathbb D}(1-|z|^2)|p'(z)|\le8.
\end{equation}
Write
\[
 \left(\frac{f'(z)}{f(z)}\right)^2-z^{-2}
       =\frac{2p_1}{z}+\tau_{\rm reg}(z),
 \qquad
 \tau_{\rm reg}(z)=2\frac{p'(z)-p_1}{z}+p'(z)^2.
\]
With normalized area $dA/\pi$,
\[
 \left\|\frac{p'-p_1}{z}\right\|_{A^2}^2
 =\sum_{n\ge2}\frac{n^2}{n-1}|p_n|^2\le2E_p.
\]
For the weighted Bergman space
\[
 A^2_2=\operatorname{Hol}(\D)\cap
       L^2\left(\D,(1-|z|^2)^2dA/\pi\right),
\]
these bounds give
\[
 \|\tau_{\rm reg}\|_{A^2_2}^2
 \le 16E_p+128E_p\le288c.
\]
Also
\begin{equation}\label{eq:N.10}
\|\tau_{\rm reg}\|_{A^1}
             \le C(\sqrt c+c).
\end{equation}
For $n\geq2$, write $b_n=[z^{n-2}]\tau_{\rm reg}$; these are the
coefficients also used in Appendix~\ref{app:score-hessian}.
The exterior variation is
\[
 \ell_na=0,\qquad
 \ell_{-n}a=-a\operatorname*{Res}_{u=0}
                 u^{1-n}\left(\frac{f'(u)}{f(u)}\right)^2=-ab_n.
\]
It follows by taking the leading coefficient of the exterior Schiffer
variation, derived in \eqref{eq:KD.17}. The real response combinations give
\[
 D_{n,c}c=\frac{\alpha_n}{2}\operatorname{Im}b_n,\qquad
 D_{n,s}c=\frac{\alpha_n}{2}\operatorname{Re}b_n.
\]
Thus the capacity part of the score is $2Dc$. Parseval for $A^2_2$ gives
\begin{equation}\label{eq:N.11}
 \begin{aligned}
 \|Dc\|_H^2
 &=\frac14\sum_{n\ge2}\alpha_n^2|b_n|^2\\
 &=\frac14\|\tau_{\rm reg}\|_{A^2_2}^2\leq72c.
 \end{aligned}
\end{equation}

To identify the capacity drift, we pair its response with the score.
Write $Q_f(z)=\sum_{n\geq2}t_n(f)z^{n-2}$ for the regular part of
\eqref{eq:stress}. For $0<r<1$, put
\[
 Q_{f,r}(z)=r^2Q_f(rz),\qquad
 \tau_{{\rm reg},r}(z)=r^2\tau_{\rm reg}(rz).
\]
The response normalization and the radial moment $2/d_n$ give
\begin{align*}
 \langle\rho,Dc\rangle_r
 &:=\sum_{n\ge2}r^{2n}
       \sum_{\tau\in\{c,s\}}\rho_{n,\tau}D_{n,\tau}c
   =\Rea\sum_{n\ge2}\frac{r^{2n}t_n\overline{b_n}}{d_n}\\
 &=\frac12\Rea\int_\D
       Q_{f,r}(z)\overline{\tau_{{\rm reg},r}(z)}
                  (1-|z|^2)^2\,\frac{dA(z)}\pi .
\end{align*}
Define the limiting pairing by
\[
 \langle\rho,Dc\rangle
 :=\frac12\Rea\int_\D Q_f(z)\overline{\tau_{\rm reg}(z)}
                  (1-|z|^2)^2\,\frac{dA(z)}\pi .
\]
Bers--$A^1$ duality and \eqref{eq:N.10} give
\begin{equation}\label{eq:N.12}
|\langle\rho,Dc\rangle|
                    \le C_\kappa(\sqrt c+c).
\end{equation}
The Bers bound and radial contraction in $A^1$ imply
\begin{align*}
 \sup_{0<r<1}\sup_{z\in\D}
       (1-|z|^2)^2|Q_{f,r}(z)|&\le C_\kappa,\\
 \shortintertext{while radial contraction gives}
 \|\tau_{{\rm reg},r}\|_{A^1}
       &\le\|\tau_{\rm reg}\|_{A^1},\\
 \|\tau_{{\rm reg},r}-\tau_{\rm reg}\|_{A^1}
       &\longrightarrow0\qquad(r\uparrow1).
\end{align*}
Since $Q_{f,r}(z)\to Q_f(z)$ for each $z\in\D$, dominated
convergence gives, for each fixed welding,
\begin{align*}
 |\langle\rho,Dc\rangle_r-\langle\rho,Dc\rangle|
 &\le C_\kappa\|\tau_{{\rm reg},r}-\tau_{\rm reg}\|_{A^1}\\
 &\quad+\frac12\int_\D
       |Q_{f,r}-Q_f|\,|\tau_{\rm reg}|
                    (1-|z|^2)^2\,\frac{dA(z)}\pi
       \longrightarrow0.\\
 \shortintertext{The same estimates give the uniform bound}
 \sup_{0<r<1}|\langle\rho,Dc\rangle_r|
 &\le C_\kappa(\sqrt c+c).
\end{align*}
Consequently, for every $R<\infty$,
\[
 \left\|\mathbf1_{\{c\le R\}}
       \bigl(\langle\rho,Dc\rangle_r-\langle\rho,Dc\rangle\bigr)
       \right\|_{L^2(\nu)}
       \longrightarrow0 .
\]

For $J=f'/f$, the positive Witt action is
\[
 \ell_nJ=z^{n+1}J'+(n+1)z^nJ.
\]
Extracting the coefficient of $z^{n-2}$ in $\ell_n(J^2)$ gives
$\ell_nb_n=2n$. Hence
\[
 (D_{n,c}^2+D_{n,s}^2)c
 =-\frac{\alpha_n^2}{2}
   (\ell_n\ell_{-n}+\ell_{-n}\ell_n)c
 =n\alpha_n^2=\frac2{n^2-1}.
\]
Summing the complete sine--cosine pairs gives
\begin{equation}\label{eq:N.13}
\sum_jD_j^2c
        =\sum_{n\ge2}\frac{2}{n^2-1}=\frac{3}{2}.
\end{equation}
Write $Lc$ for the following local generator expression; this notation
does not assert that $c\in\Dom(L)$:
\begin{equation}\label{eq:N.14}
L c=\frac{3}{2}-\langle\rho,Dc\rangle,
              \qquad |Lc|\le C_\kappa(1+c).
\end{equation}

For the later capacity localization, we record the generator of a bounded
test $\chi(c)$.
For a smooth bounded $\chi$ whose first and second derivatives have
compact support in $[0,\infty)$, choose $R\ge1$ with
$\supp\chi'\cup\supp\chi''\subset[0,R]$.
Finite-flow differentiation and \eqref{eq:N.11} give
\[
 D_j\chi(c)=\chi'(c)D_jc,\qquad
 \|D\chi(c)\|_2^2
       \le72\sup_{x\ge0}x|\chi'(x)|^2<\infty .
\]
Minimal--maximal equality therefore gives $\chi(c)\in\Dom(D)$.
The local $L^2$ drift limit above and \eqref{eq:N.11}--\eqref{eq:N.13}
yield
\begin{align*}
 \left\|\sum_{n\ge2}\sum_{\tau\in\{c,s\}}
       r^{2n}D_{n,\tau}\chi(c)e_{n,\tau}-D\chi(c)
       \right\|_{L^2(\nu;H)}&\longrightarrow0,\\
 \shortintertext{and}
 \left\|\chi''(c)\left(
       \sum_{n\ge2}\sum_{\tau\in\{c,s\}}r^{2n}|D_{n,\tau}c|^2
                         -\|Dc\|_H^2\right)\right\|_2
       &\longrightarrow0.\\
 \shortintertext{The paired trace also converges:}
 \sum_{n\ge2}\frac{2r^{2n}}{n^2-1}
       &\longrightarrow\frac32 .
\end{align*}
For every real $F\in\Ccyl$, finite-frequency integration by parts
and then absolute convergence at fixed $r<1$ give
\begin{align*}
 -\cE(\chi(c),F)
 &=-\lim_{r\uparrow1}\int
       \sum_{n\ge2}\sum_{\tau\in\{c,s\}}
       r^{2n}D_{n,\tau}\chi(c)D_{n,\tau}F\,d\nu\\
 &=\lim_{r\uparrow1}\int F\left[
       \chi'(c)\left(\sum_{n\ge2}\frac{2r^{2n}}{n^2-1}
                         -\langle\rho,Dc\rangle_r\right)
       +\chi''(c)\sum_{n\ge2}\sum_{\tau\in\{c,s\}}
                         r^{2n}|D_{n,\tau}c|^2\right]d\nu\\
 &=\int F\left[
       \chi'(c)\left(\frac32-\langle\rho,Dc\rangle\right)
                         +\chi''(c)\|Dc\|_H^2\right]d\nu .
\end{align*}
The bracketed expression belongs to $L^\infty(\nu)$ by
\eqref{eq:N.11}--\eqref{eq:N.12}. Both sides are continuous in
the form norm, so the identity extends to every $F\in\Dom(D)$.
The generator definition gives
\begin{equation}\label{eq:N.15}
 \begin{aligned}
 \chi(c)&\in\Dom(L),\\
 L\chi(c)&=\chi'(c)\Big(\frac{3}{2}-\langle\rho,Dc\rangle\Big)
                         +\chi''(c)\|Dc\|_H^2\in L^\infty(\nu).
 \end{aligned}
\end{equation}

We now use the capacity to tilt the welding law. Set
\begin{equation}\label{eq:capacity-tilt}
 \beta=2-\frac{Q^2}{2}\leq0,\qquad
 z_\beta=\int e^{\beta c}\,d\nu\in(0,1],\qquad
 d\nu_\beta=z_\beta^{-1}e^{\beta c}\,d\nu.
\end{equation}
Its finite-direction score is $\rho_v^\beta=\rho_v-\beta D_vc$.
We will identify the field marginal of
\[
                  \pi(dh,dY,dZ)=\nu_\beta(dh)K_h^Q(dY,dZ).
\]
For a function of the fields and the welding, $K$ denotes conditional
expectation given $h$ under $\pi$. Its $H$-valued version is $K_H$.

\subsection{The kernel derivative}
Retain $f=f_h^\circ$ and $p=\log(f/z)$, and set
\[
 u_f(z)=\log(zf'(z)/f(z)),\qquad
 G^{\rm in}(z,w)=\log\frac{f(z)-f(w)}{z-w}-p(z)-p(w).
\]
Write $m(z)=Qu_f(z)/2=\sum m_rz^r$ and
$G^{\rm in}(z,w)=\sum_{r,s\ge1}G^{\rm in}_{rs}z^rw^s$. The interior component of
the tilted trace kernel has Gaussian coordinates $(Y_r)_{r\ge1}$
with
\begin{align}
 KY_r&=m_r,\notag\\
 \shortintertext{with covariances}
 C_{rs}:=K[(Y_r-m_r)(Y_s-m_s)]&=-G^{\rm in}_{rs}/2,\notag\\
 K[(Y_r-m_r)\overline{(Y_s-m_s)}]&=\frac{\delta_{rs}}{2r}.
 \label{eq:KD.1}
\end{align}
Set $Y_{-r}=\bar Y_r$ for $r\geq1$ and $Y_0=0$.
For $n\ge1$ define the centered field derivations $\Vzero_n$ by
\begin{equation}\label{eq:KD.2}
\Vzero_nY_r=nQ\delta_{nr}+(r-n)Y_{r-n},\qquad
 \Vzero_n\bar Y_r=-(r+n)\bar Y_{r+n},
\end{equation}
and set $\Vzero_{-n}P=-\overline{\Vzero_n\bar P}$. Put
\begin{equation}\label{eq:KD.3}
\Sigma_n=Qn^2Y_n-\sum_{k+s=n}ksY_kY_s.
\end{equation}
Its generating Laurent series is
\[
 \Sigma(z)=\sum_{n\geq1}\Sigma_nz^{n-2}
                   =-Y'(z)^2+QY''(z)+\frac{QY'(z)}z.
\]

The exact pointwise identities are
\begin{align}
 \ell_n(KP)&=K(\Vzero_nP)
        -K[(\bar\Sigma_n-K\bar\Sigma_n)P],
 \quad n\ge1,\label{eq:KD.4}\\
\shortintertext{and}
\ell_{-n}(KP)&=K(\Vzero_{-n}P)
        +K[(\Sigma_n-K\Sigma_n)P].\label{eq:KD.5}
\end{align}
The convention \eqref{eq:coefficient-responses} implies
$\ell_n\bar F=-\overline{\ell_{-n}F}$. Thus the real version of
\eqref{eq:KD.4}--\eqref{eq:KD.5} uses
\[
 \Vzero_{\cos n\theta}=\frac i2(\Vzero_n+\Vzero_{-n}),\qquad
 \Vzero_{\sin n\theta}=\frac12(\Vzero_n-\Vzero_{-n}),
\]
and the real linear score family
\[
  \begin{aligned}
 \Sigma_{\cos n\theta}&=\operatorname{Im}\Sigma_n,&
 \Sigma_{\sin n\theta}&=\operatorname{Re}\Sigma_n,\\
 \Sigma_{n,c}&=\alpha_n\operatorname{Im}\Sigma_n,&
 \Sigma_{n,s}&=\alpha_n\operatorname{Re}\Sigma_n\quad(n\geq2).
 \end{aligned}
\]
For the low frame, set
\[
 \Sigma_{z_0}=0,\qquad \Sigma_{z_c}=\operatorname{Im}\Sigma_1,
 \qquad\Sigma_{z_s}=\operatorname{Re}\Sigma_1,
\]
and define the rotation derivation by
\[
 \Vzero_{z_0}Y_r=irY_r,\qquad
 \Vzero_{z_0}\bar Y_r=-ir\bar Y_r.
\]
The derivations and scores extend by real linearity. For every finite
real Fourier field $v$,
\begin{equation}\label{eq:KD.6}
 D_v(KP)=K(\Vzero_vP)-K[(\Sigma_v-K\Sigma_v)P].
\end{equation}
For the constant direction, the normalized maps under rotation are
$f_t(z)=e^{-it}f(e^{it}z)$ and $g_t(w)=e^{-it}g(w)$. The Gaussian
means and covariances therefore give $D_{z_0}(KP)=K(\Vzero_{z_0}P)$,
which is \eqref{eq:KD.6} with $\Sigma_{z_0}=0$.

We verify the kernel identities first on means and covariances.
The positive Kirillov formula $\ell_nf=z^{n+1}f'$ gives
\begin{equation}\label{eq:KD.7}
\ell_nm_r=\frac{Qn}{2}\delta_{nr}
                  +(r-n)m_{r-n}\mathbf1_{r>n},
\end{equation}
\noindent and
\begin{equation}\label{eq:KD.8}
\ell_n C_{rs}
 =(r-n)C_{r-n,s}\mathbf1_{r>n}
 +(s-n)C_{r,s-n}\mathbf1_{s>n}
 -\tfrac12\mathbf1_{r+s=n}.
\end{equation}
Indeed, differentiation of the logarithmic kernel gives
\[
 \ell_n G^{\rm in}(z,w)
 =(z^{n+1}\partial_z+w^{n+1}\partial_w)G^{\rm in}(z,w)
   +\sum_{k=1}^{n-1}z^kw^{n-k}.
\]
Coefficient extraction proves \eqref{eq:KD.8}.

For the negative responses, set
\begin{align}
 J(u)&=f'(u)/f(u),\qquad
 S_z(u)=\frac{1}{u-z}-\frac{1}{u},\notag\\
 \Gamma_z^{\rm in}(u)&=\partial_u G^{\rm in}(z,u)+S_z(u)
       =\frac{f(z)f'(u)}{f(u)(f(u)-f(z))}.
 \label{eq:KD.9}
\end{align}
Take Laurent expansions with $|z|,|w|<|u|$.
The original negative contour formula implies
\begin{align}
 \ell_{-n}p(z)&=\operatorname*{Res}_{u=0}u^{1-n}J(u)\Gamma_z^{\rm in}(u),\notag\\
 \shortintertext{and}
 \ell_{-n}G^{\rm in}(z,w)
 &=\operatorname*{Res}_{u=0}u^{1-n}\Gamma_z^{\rm in}(u)\Gamma_w^{\rm in}(u).
 \label{eq:KD.10}
\end{align}
The second identity follows by substituting the first in
$\delta G^{\rm in}=(f(w)\delta p(z)-f(z)\delta p(w))/(f(z)-f(w))$.
For the second residue, the coefficient expansion is
\[
 [z^r]\Gamma_z^{\rm in}(u)
 =u^{-r-1}-2\sum_{k\ge1}kC_{rk}u^{k-1}.
\]
The product of the two principal terms has zero residue, and the
remaining terms give
\begin{equation}\label{eq:KD.11}
\ell_{-n}C_{rs}
 =(r+n)C_{r+n,s}+(s+n)C_{r,s+n}
             -2\sum_{k+t=n}kt C_{rk}C_{st}.
\end{equation}

Since $u_f(z)=G^{\rm in}(z,z)+p(z)$, \eqref{eq:KD.10} also gives
$\ell_{-n}u_f(z)=\operatorname*{Res}u^{1-n}\Gamma_z^{\rm in}(\Gamma_z^{\rm in}+J)$.
Now
\[
 \partial_u\log\Gamma_z^{\rm in}=\frac{f''}{f'}-J-\frac{f'}{f-f(z)},
 \qquad u_f'(u)=\frac{f''}{f'}-J+1/u.
\]
It follows that
$\Gamma_z^{\rm in}(\Gamma_z^{\rm in}+J-u_f')
=-\partial_u\Gamma_z^{\rm in}-\Gamma_z^{\rm in}/u$.
One residue integration by parts consequently yields
\[
 \ell_{-n}u_f(z)
 =\operatorname*{Res}_{u=0}u^{1-n}u_f'(u)\Gamma_z^{\rm in}(u)
                              -n^2[u^n]G^{\rm in}(z,u).
\]
Coefficient extraction proves
\begin{equation}\label{eq:KD.12}
\ell_{-n}m_r=(r+n)m_{r+n}+Qn^2C_{rn}
                       -2\sum_{k+t=n}kt m_kC_{rt}.
\end{equation}

\subsection{The exterior field and the cross covariances}

Retain the exterior normalization $g(w)=aw+O(1)$, including the
complex phase of $a$ fixed by the rotational mark. In addition to
\eqref{eq:KD.1}, the common tilted
trace has means and covariances
\begin{align}
 m^{\rm out}(w)&=\frac{Q}{2}\log\frac{wg'(w)}{g(w)},\notag\\
 G^{\rm out}(w,v)&=\log\frac{g(w)-g(v)}{a(w-v)},\qquad
 C^{\rm out}(w,v)=-\tfrac12 G^{\rm out}(w,v).
 \label{eq:KD.15}
\end{align}
\noindent and
\begin{equation}\label{eq:KD.16}
 C^{\rm cross}(z,w)=-\frac12\log\big(1-f(z)/g(w)\big),\qquad
 K[(Y-m)\overline{(Z-m^{\rm out})}]=0.
\end{equation}
The exterior mixed covariance is $\delta_{rs}/(2r)$, as for the
interior. To verify \eqref{eq:KD.16}, subtract the interior value at
zero. The covariance of the real harmonic extensions in opposite
components is then $-\log|1-f(z)/g(w)|$; its holomorphic part is
$C^{\rm cross}$.

For \eqref{eq:KD.15}, subtract the exterior Dirichlet Green function
from the whole-plane covariance. The linear Gaussian tilt changes
only the means.

In the complexified right response, $\ell_ng=0$ for $n\ge1$.
For the negative direction write
\begin{equation}\label{eq:KD.17}
F_n(x)=\operatorname*{Res}_{u=0}
       \frac{u^{1-n}x^2f'(u)^2}{f(u)^2(f(u)-x)},\qquad
               \ell_{-n}g(w)=F_n(g(w)).
\end{equation}
The contour in $F_n$ surrounds zero and excludes the preimage of
$x$. The interior contour in \eqref{eq:negative-contour} gives
\[
 \ell_{-n}f(z)=z^{1-n}f'(z)+F_n(f(z)),
\]
where the first term is the residue at $u=z$. The formulas preserve
$f(0)=0,f'(0)=1$ and solve the infinitesimal welding identity. Their
interior residues pass from analytic weldings to the removable locus by
finite analytic quasiconformal deformation.

Set $C_u^{\rm cross}(u,w)=\partial_uC^{\rm cross}(u,w)
=-f'(u)/(2(f(u)-g(w)))$. Substitution of \eqref{eq:KD.17} and
cancellation of the leading exterior coefficient give
\begin{align}
 \ell_{-n}C^{\rm cross}(z,w)
 &=\operatorname*{Res}_{u=0}u^{1-n}
       \Gamma_z^{\rm in}(u)C_u^{\rm cross}(u,w),
 \label{eq:KD.18}\\
 \shortintertext{while the exterior covariance satisfies}
 \ell_{-n}G^{\rm out}(w,v)
 &=4\operatorname*{Res}_{u=0}u^{1-n}
       C_u^{\rm cross}(u,w)C_u^{\rm cross}(u,v),
 \label{eq:KD.19}\\
 \shortintertext{and the exterior mean satisfies}
 \ell_{-n}\log\frac{wg'(w)}{g(w)}
 &=-2\operatorname*{Res}_{u=0}u^{1-n}u_f'(u)C_u^{\rm cross}(u,w)
       +2n^2[u^n]C^{\rm cross}(u,w).
 \label{eq:KD.20}
\end{align}
For \eqref{eq:KD.19}, divide the difference $F_n(x)-F_n(y)$ by $x-y$.
The identity
\[
 \frac{x^2/(s-x)-y^2/(s-y)}{x-y}
 =-1+\frac{s^2}{(s-x)(s-y)}
\]
shows that the constant term cancels the variation of $\log a$.

For \eqref{eq:KD.20}, put $x=g(w)$. Differentiating the Schiffer field gives
\begin{align*}
 \ell_{-n}\log\frac{wg'(w)}{g(w)}
 &=F_n'(x)-\frac{F_n(x)}x\\
 &=2\operatorname*{Res}_{u=0}u^{1-n}
       \bigl(2(C_u^{\rm cross})^2+JC_u^{\rm cross}\bigr).
\end{align*}
Now use
\[
 2(C_u^{\rm cross})^2+J C_u^{\rm cross}+u_f'C_u^{\rm cross}
     =\partial_uC_u^{\rm cross}+C_u^{\rm cross}/u.
\]
Residue integration by parts gives the last term in \eqref{eq:KD.20}.

Define the coefficients by
\[
 C^{\rm cross}_{r,s}=[z^rw^{-s}]C^{\rm cross}(z,w),\qquad
 C^{\rm out}_{r,s}=[w^{-r}v^{-s}]C^{\rm out}(w,v),\qquad
 m^{\rm out}_r=[w^{-r}]m^{\rm out}(w).
\]
In coefficients, these identities become
\begin{align}
 \ell_{-n}C^{\rm cross}_{r,s}
 &=(r+n)C^{\rm cross}_{r+n,s}
       -2\sum_{k+t=n}kt C_{rk}C^{\rm cross}_{t,s},
 \label{eq:KD.21}\\
 \shortintertext{while}
 \ell_{-n}C^{\rm out}_{r,s}
 &=-2\sum_{k+t=n}kt C^{\rm cross}_{k,r}C^{\rm cross}_{t,s},
 \label{eq:KD.22}\\
 \shortintertext{and}
 \ell_{-n}m^{\rm out}_r
 &=Qn^2 C^{\rm cross}_{n,r}
       -2\sum_{k+t=n}kt m_kC^{\rm cross}_{t,r}.
 \label{eq:KD.23}
\end{align}
The positive derivatives of the exterior holomorphic mean and
covariance vanish, while
$\ell_n C^{\rm cross}(z,w)=z^{n+1}\partial_z C^{\rm cross}(z,w)$.
Set $\Vzero_nZ_r=\Vzero_n\bar Z_r=0$.
To pass from these responses to polynomial tests, let $\mathbf y$ be a
finite list of coordinates from both fields and their conjugates. Write
$\mathbf m=K\mathbf y$ and
$\mathbf C_{ij}=K[(\mathbf y_i-\mathbf m_i)(\mathbf y_j-\mathbf m_j)]$.
For finitely many formal variables $\mathbf u$, put
\[
 M(\mathbf u)=Ke^{\mathbf u\cdot\mathbf y}
       =\exp\left(\mathbf u\cdot\mathbf m
                   +\tfrac12\mathbf u^\top\mathbf C\mathbf u\right).
\]
The Gaussian exponential calculation in
Appendix~\ref{app:gaussian-exponential}, with the exterior and cross
responses \eqref{eq:KD.21}--\eqref{eq:KD.23}, gives $\ell_nM$ in the form
prescribed by \eqref{eq:KD.4}. Comparing formal coefficients proves
\eqref{eq:KD.4} for every finite polynomial in both fields; conjugation
and the real response combinations give \eqref{eq:KD.5}--\eqref{eq:KD.6}.

\subsection{The product Gaussian marginal}\label{sec:gaussian-marginal}

Wick's rule identifies the tilted score as $\rho_v^\beta=K\Sigma_v$:
\begin{equation}\label{eq:KD.13}
K\Sigma(z)=\frac{1+6Q^2}{12}{\operatorname{S}}f(z)
       +\frac{Q^2}{4}\big((f'(z)/f(z))^2-z^{-2}\big).
\end{equation}
Here the centered conditional covariance gives
$K[(Y'(z)-m'(z))^2]=-{\operatorname{S}}f(z)/12$. Direct differentiation yields
\[
 u_f''+\frac{u_f'}z-\frac{(u_f')^2}{2}
 ={\operatorname{S}}f+\frac12\left(\frac{f'}f\right)^2-\frac1{2z^2}.
\]
Substitution of $m=Qu_f/2$ in
${\operatorname{S}}f/12-m'^2+Qm''+Qm'/z$ proves \eqref{eq:KD.13}.
Since the response score of $((f'/f)^2-z^{-2})$ is $2Dc$, it gives
$K\Sigma_v=\rho_v-\beta D_vc$ in every fixed direction.

For a fixed field polynomial $P$, Wick's formula makes $KP$ a polynomial
in finitely many coefficients of $f$, $g/a$, $a^{-1}$ and their
conjugates, by the mean and covariance formulas
\eqref{eq:KD.1}, \eqref{eq:KD.15} and \eqref{eq:KD.16}.
Univalent coefficient bounds and $|a^{-1}|=e^{-c}\leq1$
bound this polynomial and its fixed finite responses, so finite-flow
integration by parts applies.
Integrating \eqref{eq:KD.6} against $\nu_\beta$ gives
\begin{equation}\label{eq:KD.14}
\int K(\Vzero_vP)\,d\nu_\beta
                  =\int K(\Sigma_vP)\,d\nu_\beta.
\end{equation}
Direct integration against $\mathbb G$ gives, on finite polynomials,
\begin{align}
 (\partial_{Y_r})^*&=rY_r-\partial_{\bar Y_r},\qquad
 \Vzero_n^*=\Vzero_{-n}+\Sigma_n,\notag\\
 \shortintertext{and therefore}
 \int\Vzero_n P\,d\mathbb G
       &=\int\bar\Sigma_nP\,d\mathbb G.
 \label{eq:KD.24}
\end{align}
The finite-dimensional Gaussian Stein identities give
\begin{align*}
 \mathbb E_{\mathbb G}[\Vzero_nP]
  &=Qn\mathbb E_{\mathbb G}[\partial_{Y_n}P]
       -\sum_{k+t=n}\mathbb E_{\mathbb G}[\partial_{Y_k}\partial_{Y_t}P]\\
  &=\mathbb E_{\mathbb G}[\bar\Sigma_nP].
\end{align*}
The transport terms with $r>n$ cancel those in the antiholomorphic
variables after shifting the index by $n$. All differentiated sums
are finite for each $P$.

To identify the field moments, give
$Y_r$ and $\bar Y_r$ weight $r$, and write $\deg_w$ for the
resulting weighted degree. On polynomials in the antiholomorphic
variables, define
\begin{equation}\label{eq:KD.25}
 {B_n^{\mathrm{raise}}} R
  =\sum_{r\geq1}(r+n)\bar Y_{r+n}\partial_{\bar Y_r}R
    +\left(Qn^2\bar Y_n-\sum_{k+t=n}kt\bar Y_k\bar Y_t\right)R.
\end{equation}
Every sum is finite on each polynomial. Let
$\C[\bar Y]_N$ be the polynomials homogeneous of weighted degree $N$.
The moment recursion will use the spanning identity
\begin{equation}\label{eq:gaussian-raising-span}
 \C[\bar Y]_N
       =\sum_{n=1}^N{B_n^{\mathrm{raise}}}\C[\bar Y]_{N-n},\qquad N\geq1.
\end{equation}
It follows from the Feigin--Fuchs cyclicity statement in
\cite[Section 4.1, equations (4.2)--(4.4)]{BJ}. To check its applicability,
we match the operator and charge in our normalization.
In variables $x_r=-\bar Y_r$, the operator in \eqref{eq:KD.25} is
\[
 \sum_{r\geq1}(r+n)x_{r+n}\partial_{x_r}
      -\sum_{k+t=n}ktx_kx_t-Qn^2x_n.
\]
This is $\mathbf L^{\rm FF}_{-n,\alpha}$ in the cited formulas with
module parameter $\alpha=Q$ and background charge $Q$.
Indeed, the zero-mode term vanishes because $Q-\alpha=0$.
The exceptional parameters in the cited formulas are
$Q-r\sqrt\kappa/2-2s/\sqrt\kappa$ and
$Q+r\sqrt\kappa/2+2s/\sqrt\kappa$, with $r,s\geq1$.
The former are nonpositive and the latter are at least $2Q$.
Thus the module parameter $Q>0$ is nonexceptional, also at
$\kappa=4$.
The cited cyclicity gives a homogeneous basis of descendants of $1$.
Every descendant of weight $N>0$ begins with one raising operator,
which proves \eqref{eq:gaussian-raising-span}.

For a fixed polynomial $S$ in the exterior coordinates and their conjugates,
write $S(Z)=S(Z,\bar Z)$ and define the linear functional
\[
 \Delta_S(P)=\mathbb E_\pi[P(Y)S(Z)]
                 -\mathbb E_{\mathbb G}[P]\mathbb E_\pi[S(Z)].
\]
Since the exterior field is passive in \eqref{eq:KD.14}, the joint
identity and \eqref{eq:KD.24} give
\begin{equation}\label{eq:gaussian-spectator-ward}
 \Delta_S(1)=0,\qquad
 \Delta_S(\Vzero_nP)=\Delta_S(\bar\Sigma_nP),\qquad n\geq1,
\end{equation}
and the conjugate Ward identity also holds. Here $S$ may be complex:
the conjugate identity follows by first testing with $\bar S$.

Let $F$ be holomorphic and homogeneous of weight $M$, and let
$R$ be antiholomorphic and homogeneous of weight $N-n$. The
product rule in \eqref{eq:gaussian-spectator-ward} gives
\begin{equation}\label{eq:gaussian-moment-recursion}
                  \Delta_S(F{B_n^{\mathrm{raise}}}R)
                        =\Delta_S((\Vzero_nF)R).
\end{equation}
For a monomial $F$, every nonzero term on the right has strictly smaller
weight than $M+N$, as the following degrees of nonzero terms show:
\begin{align*}
 \Vzero_nF
  &=Qn\partial_{Y_n}F
    +\sum_{r>n}(r-n)Y_{r-n}\partial_{Y_r}F
    -\sum_{1\leq r<n}(n-r)\bar Y_{n-r}\partial_{Y_r}F.\\
 \shortintertext{After multiplication by $R$, the three types of terms have weights}
 \deg_w\big((\partial_{Y_n}F)R\big)&=M+N-2n,\\
 \deg_w\big(Y_{r-n}(\partial_{Y_r}F)R\big)
                      &=M+N-2n\quad(r>n),\\
 \deg_w\big(\bar Y_{n-r}(\partial_{Y_r}F)R\big)
  &=M+N-2r\quad(1\leq r<n).
\end{align*}
First, taking $F=1$ in \eqref{eq:gaussian-moment-recursion} gives
\[
 \Delta_S({B_n^{\mathrm{raise}}}R)=\Delta_S((\Vzero_n1)R)=0.
\]
Equation \eqref{eq:gaussian-raising-span} identifies every pure
antiholomorphic moment. The pure holomorphic moments follow from
\[
 \Delta_S(F)=\overline{\Delta_{\bar S}(\bar F)}=0
                         \quad(M>0).
\]
Induction on total weight in \eqref{eq:gaussian-moment-recursion}
now gives $\Delta_S(FR)=0$. Thus, for every holomorphic $F$,
antiholomorphic $R$, and exterior polynomial $S$,
\begin{equation}\label{eq:gaussian-spectator-factorization}
 \mathbb E_\pi[F(Y)R(\bar Y)S(Z)]
     =\mathbb E_{\mathbb G}[F(Y)R(\bar Y)]\mathbb E_\pi[S(Z)].
\end{equation}

The moment factorization determines the joint law once finite-coordinate
exponential moments are controlled. For $m\geq1$, let
${\xi^{(m)}}\in\R^{4m}$ list the real and imaginary
parts of $Y_1,Z_1,\ldots,Y_m,Z_m$. For constants $M_m,C_m<\infty$,
the conditional Gaussian bounds give, for every
$\epsilon\in\{-1,1\}^{4m}$,
\[
 |\epsilon\cdot K{\xi^{(m)}}|\leq M_m,\qquad
 K\big[|\epsilon\cdot({\xi^{(m)}}-K{\xi^{(m)}})|^2\big]\leq C_m.
\]
Hence, for $b>0$,
\[
 \mathbb E_\pi e^{b\|{\xi^{(m)}}\|_1}
 \leq\sum_{\epsilon\in\{-1,1\}^{4m}}
                \mathbb E_\pi e^{b\epsilon\cdot{\xi^{(m)}}}
 \leq2^{4m}e^{bM_m+b^2C_m/2}<\infty.
\]
The corresponding bound for the product of $\mathbb G$ and the $Z$
marginal follows from the same finite-coordinate bounds.
Absolute convergence and
\eqref{eq:gaussian-spectator-factorization} therefore give, for
$t\in\R^{4m}$,
\begin{align*}
 \mathbb E_\pi e^{it\cdot{\xi^{(m)}}}
  &=\sum_{\ell\in\{0,1,2,\ldots\}^{4m}}\frac{(it)^\ell}{\ell!}
                         \mathbb E_\pi[{\xi^{(m)}}^\ell]\\
  &=\mathbb E_{\mathbb G\otimes\operatorname{Law}_\pi(Z)}
                              e^{it\cdot{\xi^{(m)}}}.
\end{align*}
Thus
\[
       \operatorname{Law}_\pi(Y,Z)
                        =\mathbb G\otimes\operatorname{Law}_\pi(Z).
\]

Inversion identifies the remaining exterior marginal. Let $\iota(z)=1/z$
and define
\[
 f^\#(z)=\frac{a}{g(1/z)},\qquad
 g^\#(w)=\frac{a}{f(1/w)},\qquad
 h^\#=\iota\circ h^{-1}\circ\iota.
\]
These maps satisfy $(f^\#)'(0)=1$, $(g^\#)'(\infty)=a$, and
$c(h^\#)=c(h)$. With the logarithmic conformal-radius disintegration in
\cite[Appendix~B.2.3]{Response}, inversion acts by
\[
 (s,h)\longmapsto(-s-c(h),h^\#).
\]
Invariance of the SLE loop measure \cite{ZhanLoop} and uniformity of the
rotational mark therefore give, for every bounded measurable $H$
and every interval $J\subset\R$ of length one,
\begin{align*}
 \int H(h^\#)\,d\nu(h)
 &=\int_\R\int\mathbf1_J(-s-c(h))H(h^\#)\,d\nu(h)\,ds\\
 &=\int_\R\int\mathbf1_J(s)H(h)\,d\nu(h)\,ds
   =\int H(h)\,d\nu(h).
\end{align*}
Since $c(h^\#)=c(h)$, the same identity holds for $\nu_\beta$.

The conditional means are interchanged because
\[
 \log\frac{z(f^\#)'(z)}{f^\#(z)}
    =\log\frac{wg'(w)}{g(w)}\Big|_{w=1/z}.
\]
For the covariance kernels,
\[
 G^\#(z,u)=G^{\rm out}(1/z,1/u),
\]
where $G^\#$ is the interior logarithmic kernel for $f^\#$ in
\eqref{eq:KD.1}; the cross kernel is interchanged in the same way.
Thus, for every bounded cylinder function $\Phi$ of the two fields,
\[
 K_{h^\#}^Q\Phi(Y,Z)=K_h^Q\Phi(Z,Y),\qquad
 \mathbb E_\pi\Phi(Y,Z)=\mathbb E_\pi\Phi(Z,Y).
\]
Combining this symmetry with the $Y$ marginal and the factorization gives
\begin{equation}\label{eq:KD.26}
 \int K_h^{Q}(dY,dZ)\,d\nu_\beta(h)
                       =\mathbb G(dY)\mathbb G(dZ).
\end{equation}

Let $q(h)$ be the angular coordinate of $h^{-1}(1)$, considered in
$\mathbb T$. Since $D_vq=-v(q)$, the kernel identities extend by the
product rule to $P(Y,Z)b(q)$ with $b$ smooth and periodic. The field
and mark derivation is $V_v=\Vzero_v-v(q)\partial_q$.

\subsection{The joint finite-flow density}\label{sec:joint-finite-flow}

Let
\[
 P_0u=u-\frac1{2\pi}\int_0^{2\pi}u(\theta)\,d\theta
\]
be the projection onto zero-average distributions.
For an orientation-preserving smooth circle diffeomorphism $\psi$,
write $\psi'$ for its angular derivative and define $T_\psi Y$ by
the real field $P_0(2\operatorname{Re}Y\circ\psi+Q\log\psi')$.
For a real trigonometric vector field $v$ with flow $\psi_t$, define
\begin{equation}\label{eq:KD.27}
{R_t^{\pi}}(h,Y,Z)=(h\circ\psi_t,T_{\psi_t}Y,Z).
\end{equation}
Fix $s>0$ and let $\dot H^{-s}(\T)$ be the real Sobolev space
of zero-average distributions. Identifying each field with its real
boundary distribution, \eqref{eq:KD.27} is a measurable group on
$\mathcal X\times\dot H^{-s}(\T)\times\dot H^{-s}(\T)$.
On a test $F_0(h)P(Y,Z)$ its generator is
$X_v=D_v+\Vzero_v$. The base finite-flow integration by parts and
\eqref{eq:KD.6} give
\begin{equation}\label{eq:KD.28}
\int X_v(F_0P)\,d\pi
                  =\int\Sigma_v(Y)F_0P\,d\pi.
\end{equation}
Here $D_v$ acts on $F_0$, $\Vzero_v$ acts on $P$, and
\eqref{eq:KD.6} differentiates the conditional expectation.

Set
\begin{align*}
 J_t(Y)&=\int_0^t\Sigma_v(T_{\psi_s}Y)\,ds,\\
 \shortintertext{and define the two flow weights by}
 w_t(Y)&=e^{-J_t(Y)},\\
 r_t(Y)&=w_t(T_{\psi_{-t}}Y)^{-1}.
\end{align*}
The joint finite-flow identity is
\begin{equation}\label{eq:KD.29}
 ({R_t^{\pi}})_*\pi=r_t(Y)\pi.
\end{equation}
To exponentiate \eqref{eq:KD.28}, we use the estimates proved in
Appendix~\ref{app:finite-flow-estimates}: for each $1\leq p<\infty$,
there is $t_p>0$ such that
\begin{equation}\label{eq:finite-flow-exponential}
 \sup_{|t|\leq t_p}
 \left(\|e^{J_t}\|_{L^p(\mathbb G)}
      +\|e^{-J_t}\|_{L^p(\mathbb G)}
      +\|\Sigma_v(T_{\psi_t}Y)e^{-J_t}\|_{L^p(\mathbb G)}\right)
 <\infty.
\end{equation}
Take $F_0\in\Ccyl$ and a finite field polynomial $P$, and set
$F=F_0(h)P(Y,Z)$ and $G_t=(F\circ {R_t^{\pi}})w_t$.
The cocycle identity gives
\[
                 \partial_tG_t=X_vG_t-\Sigma_v(Y)G_t.
\]
For sufficiently small $|t|$, the appendix also gives tests $G_{t,k}$
in the algebra of \eqref{eq:KD.28} such that
\begin{equation}\label{eq:finite-flow-graph}
 \begin{aligned}
 &\|G_{t,k}-G_t\|_{L^2(\pi)}
   +\|X_vG_{t,k}-X_vG_t\|_{L^1(\pi)}\\
 &\qquad+\|\Sigma_v(G_{t,k}-G_t)\|_{L^1(\pi)}\longrightarrow0.
 \end{aligned}
\end{equation}
This approximation and the time differentiability proved there allow
\eqref{eq:KD.28} to pass to $G_t$:
\[
 \frac d{dt}\int G_t\,d\pi
       =\int(X_vG_t-\Sigma_vG_t)\,d\pi=0.
\]
Integration in time yields
\[
 \int(F\circ{R_t^{\pi}})w_t\,d\pi=\int F\,d\pi.
\]
\begin{samepage}
The tests determine the joint measure, proving \eqref{eq:KD.29} for
small $t$. The identity
\[
 r_{t+s}(Y)=r_t(Y)r_s(T_{\psi_{-t}}Y)
\]
extends it by composition to every finite $t$.
\end{samepage}

For the base flow define its density by
\[
 b_t^\beta(h)=\exp\left(\int_{-t}^0
                  \rho_v^\beta(h\circ\psi_s)\,ds\right),
 \qquad (h\mapsto h\circ\psi_t)_*\nu_\beta
                    =b_t^\beta\nu_\beta.
\]
Disintegrating \eqref{eq:KD.29} now gives, for each fixed $t$ and
$\nu_\beta$-almost every $h$,
\begin{align*}
 (T_{\psi_t},\id)_*K_h^Q(dY,dZ)
   &=\frac{r_t(Y)}{b_t^\beta(h\circ\psi_t)}
                      K_{h\circ\psi_t}^Q(dY,dZ),\\
 \int r_t(Y)\,K_{h\circ\psi_t}^Q(dY,dZ)
   &=b_t^\beta(h\circ\psi_t).
\end{align*}

\section{Geometric tests along every weak solution}\label{sec:raw-geometry}

We extend the finite-test martingale identities to the conditional field moments.

\subsection{The geometric test algebra}

Retain the marked normalization of Section~\ref{sec:gaussian-lift} and set
\begin{equation}\label{eq:JL.1}
g_*=g/a=w+\gamma_0+\sum_{r\ge1}\gamma_rw^{-r},
 \qquad f_*=f/a=a^{-1}f,\qquad c=\log|a|.
\end{equation}
With $a_m=[z^{m+1}]f$ and $|a^{-1}|=e^{-c}\leq1$, define
\[
 C_{\mathrm{coef}}
 =\C[a_m,\bar a_m,\gamma_r,\bar\gamma_r,a^{-1},\overline{a^{-1}}:
        m\ge1,\ r\ge0].
\]
Each polynomial uses only finitely many generators and is bounded.

For the exterior assertion, apply the Koebe growth theorem to
$z\mapsto1/g_*(1/z)$, which is normalized univalent in the disk because
$0$ lies inside the curve.  It gives $|g_*(w)|\le4$ on the boundary.
Consequently
\begin{equation}\label{eq:JL.2}
\|f_*\|_\infty\le4,\qquad |\gamma_r|\le4.
\end{equation}
Bounds for finitely many interior coefficients are the usual univalent
coefficient bounds.

Let $\mathcal G_{\mathrm{geo}}$ consist of the real generators
\begin{align*}
 \mathcal G_{\mathrm{geo}}={}&\Ccyl
 \cup\{\Rea R,\Ima R:R\in C_{\mathrm{coef}}\}\\
 &\cup\{b(q_{z_1},\ldots,q_{z_m}):
               b\in C^\infty(\T^m),\ z_i\in S^1,\ m<\infty\}\\
 &\cup\{\chi(c):\chi\in C_b^\infty([0,\infty)),
                  \ \supp\chi'\cup\supp\chi''\text{ compact}\}.
\end{align*}
Define
\[
 C_{\mathrm{geo}}
 =\{\Phi(R_1,\ldots,R_m):
       R_i\in\mathcal G_{\mathrm{geo}},\ \Phi\in C_b^\infty(\R^m),\ m<\infty\}.
\]
Since the generators are bounded, this class is an algebra and contains
each generator. We use the same symbol for its complexification when
testing complex polynomials.

\begin{lemma}\label{lem:raw-geometric-tests}
For every $F\in C_{\mathrm{geo}}$,
\[
 F\in\Dom(L),\qquad
 \|F\|_\infty+\|LF\|_\infty
       +\|DF\|_{L^\infty(\nu;H)}<\infty.
\]
Along every weak solution satisfying
\eqref{eq:marginal-class}, $F(X_\cdot)$ admits a continuous modification
satisfying \eqref{eq:raw-mp}--\eqref{eq:raw-bracket}.
Moreover, for every field polynomial $P$ and smooth periodic $b$,
\begin{equation}\label{eq:JL.3}
                   K[P(Y,Z)b(q)]\in C_{\mathrm{geo}}.
\end{equation}
\end{lemma}

\begin{proof}
The coefficient-polynomial description in
Section~\ref{sec:gaussian-marginal} proves \eqref{eq:JL.3}.
Fix a weak solution in the stated class. Lemma~\ref{lem:raw-martingales}
already treats cylinders and inverse evaluations. It remains to treat
conformal coefficients and capacity cutoffs on the same filtered
probability space, with the same Brownian family.

We first obtain response bounds for the added conformal coefficients.
Write $J=f'/f$. The exterior Schiffer formula \eqref{eq:KD.17} gives
\begin{equation}\label{eq:JL.4}
\ell_n(a^{-1})=\ell_ng_*=0\quad(n\ge1),\qquad
 \ell_{-n}(a^{-1})=a^{-1}\operatorname*{Res}_{u=0}u^{1-n}J(u)^2.
\end{equation}
Indeed $\ell_{-n}a=-a\operatorname*{Res}u^{1-n}J^2$.
Using $\tau_{\rm reg}$ from the capacity estimates, for $n\ge2$ we obtain
\begin{equation}\label{eq:JL.5}
\sum_{n\ge2}(\ell_{-n}(a^{-1}))z^{n-2}=a^{-1}\tau_{\rm reg}(z).
\end{equation}
The capacity estimate \eqref{eq:N.10} implies
\begin{equation}\label{eq:JL.6}
\|a^{-1}\tau_{\rm reg}\|_{A^1}
       \le C e^{-c}(\sqrt c+c)\le C.
\end{equation}

For $g_*$, cancellation of its leading coefficient is essential.
Expand the Schiffer field at infinity before extracting a fixed
exterior coefficient.  The term proportional to $g$ cancels
$g_*\ell_{-n}a/a$ and leaves
\begin{equation}\label{eq:JL.7}
\ell_{-n}g_*(w)=-\operatorname*{Res}_{u=0}u^{1-n}
 \left\{\frac{f_*'(u)^2}{f_*(u)}
       +\sum_{j\ge0}\frac{f_*(u)^jf_*'(u)^2}{g_*(w)^{j+1}}\right\}.
\end{equation}
Each fixed coefficient of $w^{-r}$ uses only finitely many summands.
For $\gamma_0$ the analytic kernel, after discarding its Laurent
principal part, is
\begin{equation}\label{eq:JL.8}
\Gamma_0^{\rm out}(z)
=-\left(\frac{f_*'(z)^2}{f_*(z)}-\frac{a^{-1}}z\right).
\end{equation}
For $r\ge1$ it is
\begin{equation}\label{eq:JL.9}
\Gamma_r^{\rm out}(z)=-\sum_{j=0}^{r-1}[w^{-r}]g_*(w)^{-j-1}
                         f_*(z)^jf_*'(z)^2.
\end{equation}
In both formulas
$\Gamma_r^{\rm out}=\sum_{n\ge2}(\ell_{-n}\gamma_r)z^{n-2}$.
The area-change formula and \eqref{eq:JL.2} give
\begin{align}
 \int_{\mathbb D}|f_*'|^2dA
       &=|f_*(\mathbb D)|\le16\pi,\notag\\
 \shortintertext{and, with the weight $|w|^{-1}$,}
 \int_{\mathbb D}\frac{|f_*'|^2}{|f_*|}dA
       &=\int_{f_*(\mathbb D)}\frac{dA(w)}{|w|}\le8\pi.
 \label{eq:JL.10}
\end{align}
Thus $\|\Gamma_r^{\rm out}\|_{A^1}\le C_r$ uniformly over the normalized pairs.

Combine these exterior kernels with the interior kernels in
Lemma~\ref{lem:coefficient-domain}. For a holomorphic polynomial in
finitely many coefficient generators, the derivation rule reads
\[
 R=P(\xi_1,\ldots,\xi_m),\qquad
 \Gamma_R^{\mathrm{resp}}
   =\sum_{i=1}^m(\partial_iP)(\xi)\Gamma_{\xi_i}^{\mathrm{resp}}.
\]
Here the kernel for $a^{-1}$ is \eqref{eq:JL.5}, and those for
$\gamma_r$ are \eqref{eq:JL.8}--\eqref{eq:JL.9}. In particular,
\[
 \sup_h\|\Gamma_R^{\mathrm{resp}}(h)\|_{A^1}
 \leq\sum_{i=1}^m\sup_h|(\partial_iP)(\xi(h))|
                 \sup_h\|\Gamma_{\xi_i}^{\mathrm{resp}}(h)\|_{A^1}
 <\infty.
\]
The embedding \eqref{eq:bergman-embedding} then bounds $\|DR\|_H$.
Conjugation and the product rule give the same bound for the real algebra.

To pass these formulas to the given weak solution, use the prescribed-environment
approximations $h_t^r=X_0\circ(\Psi_t^r)^{-1}$ and their normalized
pairs $(f^r,g^r)$ constructed in Appendix~\ref{app:coefficient-domain}.
By \eqref{eq:prescribed-pair-convergence}, each $R\in C_{\mathrm{coef}}$ has a
continuous modification satisfying
\[
 \mathbb P(R_t=R(X_t))=1,\qquad
 \mathbb E\int_0^T|R_t-R(X_t)|^2\,dt=0.
\]
The same construction gives zero limiting boundary area for
$dt\,d\mathbb P$-almost every $(t,\omega)$.

The required kernel convergence follows from boundary convergence and
zero boundary area. Suppose normalized Jordan pairs $(f_k,g_k)$ converge boundary-uniformly
to $(f,g)$ and the limiting boundary has area zero. Put
$a_k=g_k'(\infty)$, $g_{*,k}=g_k/a_k$ and $f_{*,k}=f_k/a_k$.
Boundary convergence gives
\[
 a_k\to a,\qquad a_k^{-1}\to a^{-1},\qquad
 f_{*,k}\to f_*\text{ uniformly on }\overline\D.
\]
Winding numbers then give
\[
 \mathbf1_{f_{*,k}(\D)}\to\mathbf1_{f_*(\D)}
                         \quad dA\text{-almost everywhere}.
\]
Local analytic convergence and the area identity imply
\[
 \|f_{*,k}'\|_{A^2}^2=\frac{|f_{*,k}(\D)|}{\pi}
       \longrightarrow\frac{|f_*(\D)|}{\pi}=\|f_*'\|_{A^2}^2.
\]
As in \eqref{eq:strong-kernel-convergence}, weak convergence and
convergence of norms yield
\begin{equation}\label{eq:JL.11}
f_{*,k}'\longrightarrow f_*'\quad\text{in }A^2.
\end{equation}
Moreover, $f_{*,k}'^2/f_{*,k}\to f_*'^2/f_*$ pointwise off zero.
By change of area and domination in $4\mathbb D$,
\[
 \int_\D\left|\frac{f_{*,k}'^2}{f_{*,k}}\right|dA
 =\int_{f_{*,k}(\D)}\frac{dA(w)}{|w|}
 \longrightarrow\int_{f_*(\D)}\frac{dA(w)}{|w|}
 =\int_\D\left|\frac{f_*'^2}{f_*}\right|dA.
\]
The complex-valued form of Scheff\'e's lemma therefore gives
\begin{equation}\label{eq:JL.12}
f_{*,k}'^2/f_{*,k}\longrightarrow f_*'^2/f_*\quad\text{in }L^1.
\end{equation}
For every integer $j\geq0$, \eqref{eq:JL.2} gives
\begin{align*}
 \|f_{*,k}^jf_{*,k}'^2-f_*^jf_*'^2\|_{A^1}
 &\leq4^j(\|f_{*,k}'\|_{A^2}+\|f_*'\|_{A^2})
                       \|f_{*,k}'-f_*'\|_{A^2}\\
 &\quad+\|f_{*,k}^j-f_*^j\|_\infty\|f_*'\|_{A^2}^2
             \longrightarrow0.
\end{align*}
Equations \eqref{eq:JL.8}--\eqref{eq:JL.12} and convergence of the
fixed exterior coefficients therefore give
\[
 \|\Gamma_r^{\rm out}(f_k,g_k)-\Gamma_r^{\rm out}(f,g)\|_{A^1}\longrightarrow0
                         \qquad\text{for every fixed }r.
\]

For \eqref{eq:JL.5}, it remains to prove strong Dirichlet convergence of
$p_k=\log(f_k/z)$.
Let $\Omega=f(\mathbb D)$ and choose a fixed $0<\delta<1/4$.
All normalized interior domains contain $\delta\mathbb D$ by Koebe.
The precise renormalized area identity is
\begin{equation}\label{eq:JL.13}
\int_{\mathbb D}|p'(z)|^2dA(z)
   =\int_{\Omega\setminus\delta\mathbb D}\frac{dA(w)}{|w|^2}
                                  +2\pi\log\delta.
\end{equation}
Indeed, the cross term in $|1/z+p'|^2$ has zero angular integral.
Changing variables $w=f(z)$ gives
\begin{align*}
 \int_{\Omega\setminus f(\varepsilon\D)}\frac{dA(w)}{|w|^2}
 &=\int_{\varepsilon<|z|<1}\left|\frac{f'(z)}{f(z)}\right|^2dA(z)\\
 &=2\pi\log(1/\varepsilon)
       +\int_{\varepsilon<|z|<1}|p'(z)|^2dA(z).
\end{align*}
Since $f(z)=z+O(z^2)$, replacing $f(\varepsilon\D)$ by
$\varepsilon\D$ changes the integral by $o(1)$. Thus
\[
 \int_{\Omega\setminus f(\varepsilon\D)}\frac{dA(w)}{|w|^2}
 =\int_{\Omega\setminus\delta\D}\frac{dA(w)}{|w|^2}
      +2\pi\log(\delta/\varepsilon)+o(1).
\]
Subtract $2\pi\log(1/\varepsilon)$ and let
$\varepsilon\downarrow0$ to obtain \eqref{eq:JL.13}.

The domains $\Omega_k$ lie in a common bounded set along the convergent
sequence.  Their indicators converge almost everywhere, so \eqref{eq:JL.13}
implies convergence of the Dirichlet norms.  Local analytic convergence
and this norm convergence imply
\begin{equation}\label{eq:JL.14}
p_k'\longrightarrow p'\quad\text{in }A^2.
\end{equation}
The downsampling operator $h\mapsto(h-h(0))/z$ has $A^2$ norm at
most $\sqrt2$. Therefore
\begin{align*}
 \|\tau_{{\rm reg},k}-\tau_{\rm reg}\|_{A^1}
 &\leq\left(2\sqrt2+\|p_k'\|_{A^2}+\|p'\|_{A^2}\right)
                            \|p_k'-p'\|_{A^2}\longrightarrow0.\\
 \shortintertext{Since $a_k^{-1}\to a^{-1}$, we also have}
 \|a_k^{-1}\tau_{{\rm reg},k}-a^{-1}\tau_{\rm reg}\|_{A^1}
 &\leq |a_k^{-1}|\|\tau_{{\rm reg},k}-\tau_{\rm reg}\|_{A^1}
        +|a_k^{-1}-a^{-1}|\|\tau_{\rm reg}\|_{A^1}
                    \longrightarrow0.
\end{align*}
For a fixed holomorphic $R=P(\xi_1,\ldots,\xi_m)$, write
$\xi^{(k)}=\xi(f_k,g_k)$ and $\xi=\xi(f,g)$. The preceding product formula gives
\begin{align*}
 &\|\Gamma_R^{\mathrm{resp}}(f_k,g_k)
             -\Gamma_R^{\mathrm{resp}}(f,g)\|_{A^1}\\
 &\quad\leq\sum_{i=1}^m
       |(\partial_iP)(\xi^{(k)})-(\partial_iP)(\xi)|
          \|\Gamma_{\xi_i}^{\mathrm{resp}}(f_k,g_k)\|_{A^1}\\
 &\qquad+\sum_{i=1}^m |(\partial_iP)(\xi)|
       \|\Gamma_{\xi_i}^{\mathrm{resp}}(f_k,g_k)
               -\Gamma_{\xi_i}^{\mathrm{resp}}(f,g)\|_{A^1}
 \longrightarrow0.
\end{align*}
Together with \eqref{eq:bergman-embedding}, this proves
\[
 \|DR(f_k,g_k)-DR(f,g)\|_H\longrightarrow0.
\]
The kernel and gradient bounds are uniform over all states by
\eqref{eq:polynomial-kernel}, \eqref{eq:JL.6} and \eqref{eq:JL.10}.
Conjugation and the product rule give the gradient conclusion for
every real polynomial in $C_{\mathrm{coef}}$.

To identify the limiting It\^o drift, we also control the paired second
responses. The extended holomorphic coefficient algebra is preserved by each
fixed Witt response.  This follows directly by fixed coefficient
extraction in \eqref{eq:JL.4} and \eqref{eq:JL.7}; only nonnegative powers of $a^{-1}$ occur.
All fixed first and second response coefficients are consequently
bounded.

Under the zero mode,
\begin{equation}\label{eq:JL.15}
\ell_0(a^{-1})=a^{-1},\qquad \ell_0g_*=0,
\end{equation}
as is also checked by $f_t(z)=e^{-it}f(e^{it}z)$,
$g_t=e^{-it}g$, under right rotation.  Positive modes annihilate
$a^{-1}$ and every $\gamma_r$, and annihilate a fixed interior polynomial once
their index exceeds its weighted degree.  The response bracket gives,
for all sufficiently large $n$ and each fixed holomorphic $R$,
\begin{equation}\label{eq:JL.16}
\ell_nR=0,\qquad \ell_n\ell_{-n}R=2n\ell_0R,
 \qquad
 (D_{n,c}^2+D_{n,s}^2)R=-\frac{2}{n^2-1}\ell_0R.
\end{equation}
Thus its paired trace is absolutely uniformly summable. For fixed
holomorphic coefficient polynomials $R_1,R_2$, the product rule gives
\begin{align*}
 \sum_{\tau=c,s}D_{n,\tau}^2(R_1\bar R_2)
 ={}&\bar R_2\sum_{\tau=c,s}D_{n,\tau}^2R_1
       +R_1\sum_{\tau=c,s}D_{n,\tau}^2\bar R_2\\
 &+2\sum_{\tau=c,s}D_{n,\tau}R_1D_{n,\tau}\bar R_2.
\end{align*}
The last tail satisfies
\[
 \sum_{n>N,\tau}|D_{n,\tau}R_1D_{n,\tau}\bar R_2|
 \leq\left(\sum_{n>N,\tau}|D_{n,\tau}R_1|^2\right)^{1/2}
      \left(\sum_{n>N,\tau}|D_{n,\tau}R_2|^2\right)^{1/2}.
\]
By strong $H$ convergence and the uniform gradient bound,
\[
 \mathbb E\int_0^T
      \|DR(h_s^r)-DR(X_s)\|_H^2\,ds\longrightarrow0.
\]
For holomorphic $R$, let $b_R$ be the bounded expression in
\eqref{eq:polynomial-drift}, with the extended kernel $\Gamma_R^{\mathrm{resp}}$.
For real products use the stochastic product rule. Write
\[
 b_R^r(h,x)=\sum_{n\geq2,\tau}r^{2n}D_{n,\tau}^2R(h)
       -\sum_{n\geq2,\tau}r^n\rho_{n,\tau}(x)D_{n,\tau}R(h).
\]
The prescribed-environment approximation has drift
$b_R^r(h_s^r,X_s)$, so the scores are evaluated at the given candidate.

For $x$ in the removable locus, put
$Q_x=Q_{f_x^\circ}$ and $Q_{x,r}(z)=r^2Q_x(rz)$.
The Bers bound gives
\begin{align*}
 &\left|\int_\D\left[
       \overline{Q_{x,r}}\Gamma_R^{\mathrm{resp}}(h)
             -\overline{Q_x}\Gamma_R^{\mathrm{resp}}(x)\right]
                     (1-|z|^2)^2\,\frac{dA(z)}\pi\right|\\
 &\quad\leq C_\kappa\|\Gamma_R^{\mathrm{resp}}(h)-\Gamma_R^{\mathrm{resp}}(x)\|_{A^1}\\
 &\qquad+\left|\int_\D\overline{(Q_{x,r}-Q_x)}\Gamma_R^{\mathrm{resp}}(x)
                     (1-|z|^2)^2\,\frac{dA(z)}\pi\right|.
\end{align*}
The second term tends to zero by dominated convergence for each
$x$; the first tends to zero along the geometric approximation.
All terms are uniformly bounded. Consequently
\begin{align*}
 \mathbb E\int_0^T|b_R^r(h_s^r,X_s)-b_R(X_s)|^2\,ds
                                  &\longrightarrow0,\\
 \shortintertext{and}
 \mathbb E\int_0^T\sum_{n,\tau}
       |r^nD_{n,\tau}R(h_s^r)-D_{n,\tau}R(X_s)|^2\,ds
                                  &\longrightarrow0.
\end{align*}
Doob's inequality gives
\begin{align*}
 &\mathbb E\sup_{t\leq T}\left|\sqrt2\sum_{n,\tau}\int_0^t
       \big[r^nD_{n,\tau}R(h_s^r)-D_{n,\tau}R(X_s)\big]
                                         \,dW_s^{n,\tau}\right|^2\\
 &\qquad\leq8\mathbb E\int_0^T\sum_{n,\tau}
       |r^nD_{n,\tau}R(h_s^r)-D_{n,\tau}R(X_s)|^2\,ds
                            \longrightarrow0.
\end{align*}
The regularized It\^o identity therefore yields, on the original
filtered space,
\[
 R(X_t)-R(X_0)
   =-\sqrt2\sum_j\int_0^tD_jR(X_s)\,dW_s^j
                         +\int_0^tb_R(X_s)\,ds .
\]
Since $R$ and $b_R$ are bounded, the generator criterion
\eqref{eq:canonical-generator-limit}, applied to this identity for
the canonical candidate, gives
$R\in\Dom(L)$ and $LR=b_R$.

Write $M^R$ for the martingale part of $R(X_\cdot)$.
For a smooth periodic $b$, define $M^{b(q)}$ in the same way.
The scalar inverse-mark formula and the
stochastic product rule give
\begin{align*}
 d\langle M^R,M^{b(q)}\rangle_t
    &=2\langle DR,D[b(q)]\rangle_H(X_t)\,dt.\\
 \shortintertext{The corresponding generator identity is}
 L(Rb(q))&=b(q)LR+R L[b(q)]
                           +2\langle DR,D[b(q)]\rangle_H.
\end{align*}
The same calculation applies to finitely many inverse evaluations.

Put $a_t=g_t'(\infty)$ and $c_t=\log|a_t|$.
Choose $R_0<\infty$ with
$\supp\chi'\cup\supp\chi''\subset[0,R_0]$. Then
$\chi(-\log|a^{-1}|)$ is constant for $|a^{-1}|<e^{-R_0}$ and
extends smoothly across zero. Its raw It\^o identity is
\[
 d\chi(c_t)=-\sqrt2\chi'(c_t)\sum_jD_jc(X_t)\,dW_t^j
       +\left[\chi'(c_t)Lc(X_t)
                     +\chi''(c_t)\|Dc(X_t)\|_H^2\right]dt.
\]
For the local formula for $c$, use
$\sigma_m^c=\inf\{t:c_t\geq m\}$.
Continuity and nonvanishing of the leading coefficient imply
$\sigma_m^c\uparrow\infty$ almost surely. The gradient and drift
are those of \eqref{eq:N.11}--\eqref{eq:N.14}.

For a bounded smooth $\Phi$ and finitely many real generators
$R_1,\ldots,R_m$, the chain rule gives
\begin{align*}
 D\Phi(R_1,\ldots,R_m)&=\sum_i\partial_i\Phi\,DR_i,\\
 \shortintertext{and}
 L\Phi(R_1,\ldots,R_m)&=\sum_i\partial_i\Phi\,LR_i
                    +\sum_{i,k}\partial_{ik}\Phi\,
                                      \langle DR_i,DR_k\rangle_H.
\end{align*}
These quantities are bounded on the ranges of the generators.
\end{proof}

\subsection{Reconstruction from the fields and an inverse mark}\label{sec:reconstruction}

A single inverse mark fixes the rotational ambiguity in the welding.
We use the two-sided quantum-length identity for chordal
$\mathrm{SLE}_\kappa$ from \cite[Theorem~1.2]{PS}. By the
conformal Markov property of
the rooted SLE loop measure, conditional on the past up to a nontrivial
stopping time before the loop returns to its root, the remaining arc
has chordal $\mathrm{SLE}_\kappa$ law \cite[Theorem~4.1(ii)]{ZhanLoop}.
Thus the identity applies to compact subarcs away from the root.
In a conformal chart, the covariance kernels of the independent
whole-plane field and the free-boundary GFF used in \cite{PS} differ
by a function harmonic in each variable. After fixing additive
constants, the fields can therefore be coupled to differ by a
harmonic function. The GMC coordinate-change formula
\cite[Lemma~2.4]{PS} and covariance under harmonic shifts show that
this changes both side measures by the same boundary weight, so their
equality is preserved. At $\kappa=4$ we use the critical version of
\cite[Theorem~1.2]{PS}. Exhausting the loop away from the root and
using nonatomicity gives the identity on the whole loop.

Let $\mu_Y,\mu_Z$ be the boundary Gaussian multiplicative chaos measures
of $2\Rea Y$ and $2\Rea Z$, divided by their respective total masses. For $\kappa<4$
the exponential parameter is $\sqrt\kappa/2$; at $\kappa=4$ use the
critical derivative construction. The total masses are finite and
positive. Adding constants to the fields leaves these normalized
measures unchanged. The scale disintegration in
\cite[Appendix~B.2.3, equations (B.47)--(B.49)]{Response}
identifies $\nu$ with the unit-interior-conformal-radius shape law
together with an independent uniform exterior rotation. Since scaling
does not change the welding and additive constants do not change the
normalized chaos measures, the matching identity passes to $\nu$.
Thus, under $\nu(dh)P_h(d\varphi)$, matching gives
\begin{equation}\label{eq:JL.17}
 h_*\mu_Y=\mu_Z.
\end{equation}
The actual lift is equivalent to the trace probability:
\begin{gather*}
 \frac{d[\nu_\beta(dh)P_h^Q(d\varphi)]}
      {d[\nu(dh)P_h(d\varphi)]}
   =z_\beta^{-1}\exp\left(\beta c(h)+QU-Q^2c(h)/2\right)
       \in(0,\infty)\\
       \nu(dh)P_h(d\varphi)\text{-almost surely}.
\end{gather*}
Thus \eqref{eq:JL.17} also holds under the actual lift.

For $0\leq x\leq2\pi$, define
\[
 F_Y(x)=2\pi\mu_Y([0,x]),\qquad F_Z(x)=2\pi\mu_Z([0,x]).
\]
Extend both functions by $F(x+2\pi)=F(x)+2\pi$.
The normalized chaos measures are nonatomic and have full support,
also at criticality. Thus $F_Y,F_Z$ are lifts of circle homeomorphisms.

For the inverse mark $q=h^{-1}(1)$ in angular coordinates,
\eqref{eq:JL.17} now gives
\begin{equation}\label{eq:JL.18}
 h(x)=F_Z^{-1}\big(F_Y(x)-F_Y(q)\big)\pmod{2\pi}.
\end{equation}
Denote this measurable reconstruction by $h_{\mathrm{rec}}(Y,Z,q)$,
with any measurable extension across the $\pi$-null exceptional set.
Normalized GMC covariance gives
\[
 \mu_{T_\psi Y}=(\psi^{-1})_*\mu_Y.
\]
Substituting this in \eqref{eq:JL.18} yields
\[
 h_{\mathrm{rec}}(T_\psi Y,Z,\psi^{-1}q)
                      =h_{\mathrm{rec}}(Y,Z,q)\circ\psi.
\]

The projection to the marked fields and the reconstruction satisfy
\begin{align*}
 \Pi(h,Y,Z)&=(Y,Z,q(h)),\\
 h_{\mathrm{rec}}\circ\Pi(h,Y,Z)&=h
 \quad\pi\text{-almost surely}.
\end{align*}
From now on, retain the symbol $\pi$ for the pushforward $\Pi_*\pi$.
Define
\begin{align*}
 I:L^2(\nu_\beta)&\longrightarrow L^2(\pi),&
 IF&=F\circ h_{\mathrm{rec}},\\
 I_H:L^2(\nu_\beta;H)&\longrightarrow L^2(\pi;H),&
 I_HU&=U\circ h_{\mathrm{rec}}.
\end{align*}
These are isometries, and their adjoints satisfy
\begin{align*}
 K&=I^*,& KI&=\mathrm{Id},\\
 \shortintertext{and}
 K_H&=I_H^*,& K_HI_H&=\mathrm{Id}.
\end{align*}
The conditional expectation $K$ agrees with the kernel already defined, with
$q=q(h)$ fixed when conditioning on $h$. The field marginal is
$\mathbb G\otimes\mathbb G$; the conditional distribution of $q$
is kept as part of $\pi$.

\section{Gaussian dynamics and the marked generator}\label{sec:gaussian-core}

We first determine the Gaussian field path law, recover its driving noise,
and control the inverse mark driven by that noise. These three steps
give the marked generator range used in the conditional compression.

\subsection{The Gaussian operator core}

On $L^2(\mathbb G)$, let $\mathrm{Pol}_0$ be the algebra of finite
polynomials in $(Y_n,\bar Y_n)_{n\geq1}$, equivalently the finite
complex Wick-polynomial algebra. The horizontal gradient has components
$(\Vzero_j)_{j\in{J_H}}$ defined in Section~\ref{sec:gaussian-lift}.
In terms of the real field $\phi=2\operatorname{Re}Y$, its response in
direction $v$ is $P_0(v\phi'+Qv')$.
Denote its closure by $\Vzero$ and set $A_0=\Vzero^*\Vzero$.

Closability follows from the Gaussian polynomial adjoints in
\eqref{eq:KD.24}. Indeed, if $f_m\to0$ and $\Vzero f_m\to U$ in
their respective $L^2$ spaces, then for every finite polynomial vector $G$,
\[
 \langle U,G\rangle
  =\lim_m\langle f_m,\Vzero^*G\rangle=0,
 \qquad \Vzero^*G\in L^2(\mathbb G),
\]
so $U=0$ by density.

\begin{proposition}\label{prop:gaussian-core}
The algebra $\mathrm{Pol}_0$ is an operator core for $A_0$.
\end{proposition}

See Appendix~\ref{app:gaussian-core} for the proof.

\subsection{The Gaussian Markov semigroup}\label{sec:smooth-approximation}

For a finite list $x\in\R^d$ of real field coordinates, let $\mathbb G_d$
be its Gaussian law. Write $\mathbb N_0=\{0,1,2,\ldots\}$.
For $k,l\in\mathbb N_0$, $1\leq p<\infty$, and
multiindices $\eta\in\mathbb N_0^d$, set
\[
 \|F\|_{k,l,p}^p
 =\sum_{|\eta|\leq k}\sum_{b=0}^l
   \int\sup_{q\in\T}|\partial_x^\eta\partial_q^bF(x,q)|^p\,d\mathbb G_d(x).
\]
For functions independent of $q$ this is the Gaussian Sobolev norm.

\begin{lemma}\label{lem:smooth-polynomial-approximation}
Every smooth finite-coordinate function of polynomial growth,
with all derivatives of polynomial growth, can be approximated by
field polynomials and trigonometric polynomials in the preceding
Sobolev norms, for each finite $p$ and each fixed finite derivative order.
\end{lemma}
The proof is given in Appendix~\ref{app:smooth-polynomial-approximation}.

For $P\in\mathrm{Pol}_0$, let $d$ be its ordinary polynomial degree.
The Fock calculation in Appendix~\ref{app:gaussian-core} puts $A_0P$
in the sum of Wiener chaoses of orders at most $d+2$.

For a symmetric kernel $u$, write $I_r(u)$ for its homogeneous Wiener
integral of order $r$. Wick's formula and Cauchy--Schwarz in each
contraction give, for integers $m\geq1$,
\[
 \mathbb E|I_r(u)|^{2m}\leq C_{m,r}\|u\|^{2m},\qquad
 \mathbb E|I_r(u)|^2=r!\|u\|^2.
\]
Summing the finitely many chaos orders and using H\"older's inequality
therefore gives, for every finite $p\geq2$,
\[
 \|A_0P\|_{L^p(\mathbb G)}
       \leq C_{p,d}\|A_0P\|_{L^2(\mathbb G)}<\infty.
\]
The same estimate applies to the drift and covariance coefficients
of each fixed finite coordinate list.

Let $F$ satisfy Lemma~\ref{lem:smooth-polynomial-approximation} and
choose polynomial approximants $F_m$. H\"older's inequality gives
\begin{align*}
 \|A_0(F_m-F_n)\|_2&\leq C_F\|F_m-F_n\|_{2,0,4}\longrightarrow0,\\
 \shortintertext{and}
 \|\Vzero(F_m-F_n)\|_2&\leq C_F\|F_m-F_n\|_{1,0,4}\longrightarrow0.
\end{align*}
Since $F_m\to F$ in $L^2$, closedness gives
\[
 F\in\Dom(A_0),\qquad
 A_0F_m\longrightarrow A_0F\quad\text{in }L^2(\mathbb G).
\]
Let $P\in\mathrm{Pol}_0$ be real and let $\Phi$ be a smooth normal
contraction. Choose $\chi\in C_c^\infty(\R)$, $0\leq\chi\leq1$,
equal to one near zero, and set
\[
 \Phi_R(t)=\int_0^t\chi(u/R)\Phi'(u)\,du.
\]
Each $\Phi_R(P)$ satisfies the preceding growth conditions, so
\[
 \Vzero\Phi_R(P)=\chi(P/R)\Phi'(P)\Vzero P.
\]
Since $|\Phi_R(P)|+|\Phi(P)|\leq2|P|$ and $|\Phi'|\leq1$,
dominated convergence gives
\[
 \|\Phi_R(P)-\Phi(P)\|_2
  +\|\chi(P/R)\Phi'(P)\Vzero P-\Phi'(P)\Vzero P\|_2
       \longrightarrow0.
\]
Closedness of $\Vzero$ yields
\[
 \Vzero\Phi(P)=\Phi'(P)\Vzero P,\qquad
 \|\Vzero\Phi(P)\|_2^2\leq\|\Vzero P\|_2^2.
\]
Consequently the closed form is Markov. Its symmetric Markov
semigroup is $S_t=e^{-tA_0}$, and
\[
 A_0 1=0\quad\Longrightarrow\quad S_t1=1.
\]
\subsection{The conditional field law}\label{sec:conditional-field-law}

Consider an adapted process $Y_t$ with continuous Fourier coordinates
and $\Law(Y_t)=\mathbb G$ for every deterministic $t$.
For $F\in\mathrm{Pol}_0$, assume that
\[
 M_t^F=F(Y_t)-F(Y_0)+\int_0^t A_0F(Y_r)\,dr
\]
is a martingale in the original filtration, with
\[
 d\langle M^F,M^G\rangle_t
   =2\langle\Vzero F,\Vzero G\rangle_H(Y_t)\,dt.
\]

Fix a bounded real $\mathcal F_s$-measurable random variable $\zeta$ and define the signed
density $g_t$ by
\[
 g_t\,d\mathbb G=\mathbb E[\zeta\,\mathbf1_{Y_{s+t}\in\cdot}].
\]
Then $\|g_t\|_\infty\le\|\zeta\|_\infty$, and the martingale identities imply
\[
 \langle g_t,F\rangle-\langle g_0,F\rangle
       =-\int_0^t\langle g_r,A_0F\rangle\,dr,
                      \qquad F\in\mathrm{Pol}_0.
\]
If $F_m\to F$ in the $A_0$ graph norm, the error in this identity is
bounded by
\[
 \|\zeta\|_\infty
       \bigl(2\|F_m-F\|_2+t\|A_0(F_m-F)\|_2\bigr)\longrightarrow0.
\]
Proposition~\ref{prop:gaussian-core} therefore extends it to
$F\in\Dom(A_0)$. For $F\in\Dom(A_0^2)$,
graph-norm Riemann sums with the backward test $S_{t-r}F$ give
\begin{align*}
 \langle g_t,F\rangle-\langle g_0,S_tF\rangle
 &=\int_0^t\left(-\langle g_r,A_0S_{t-r}F\rangle
              +\langle g_r,\partial_rS_{t-r}F\rangle\right)\,dr\\
 &=0,\qquad \partial_rS_{t-r}F=A_0S_{t-r}F.
\end{align*}
For $F\in L^2(\mathbb G)$ use
$F_\epsilon=(\operatorname{Id}+\epsilon A_0)^{-2}F\in\Dom(A_0^2)$ and
$\|F_\epsilon-F\|_2\to0$. Hence
\[
                  \langle g_t,F\rangle=\langle g_0,S_tF\rangle.
\]
Thus $g_t=S_tg_0$.

To obtain the identity at every time, let $P$ be a polynomial.
The Gaussian marginals give uniform integrability through
\[
 \sup_t\mathbb E|P(Y_t)|^p=\|P\|_{L^p(\mathbb G)}^p<\infty
 \qquad(1\leq p<\infty).
\]
Coordinate continuity therefore implies
\[
 \mathbb E|P(Y_t)-P(Y_r)|^2\longrightarrow0\qquad(t\to r).
\]
Polynomial density and $\sup_t\|g_t\|_2\leq\|\zeta\|_\infty$
extend this continuity to all $L^2$ pairings.
Since $\zeta$ was arbitrary,
\begin{equation}\label{eq:GW.1}
\mathbb E[f(Y_{s+t})\mid\mathcal F_s]=S_tf(Y_s)
                         \qquad (f\in L^2(\mathbb G)).
\end{equation}
This determines the state path law.

\subsection{Injectivity of the field response}

For $\mathbb G$-almost every field $\phi$, the real linear response map
\begin{equation}\label{eq:GW.2}
H\longrightarrow C^\infty(S^1)',\qquad
 v\longmapsto P_0(v\phi'+Qv')
\end{equation}
has zero kernel. Here $H\subset H^{3/2}(S^1)$ is the horizontal
WP space and $P_0$ removes the constant Fourier coefficient.

Fix small $\epsilon,\delta>0$ with $\epsilon+\delta<1/2$.
Almost surely $\phi\in H^{-\epsilon}$.
Define $v\phi'=(v\phi)'-v'\phi$.
The one-dimensional Sobolev product estimates give
\[
 \|v\phi\|_{H^{-\epsilon}}
  +\|v'\phi\|_{H^{-\epsilon-\delta}}
 \leq C_{\epsilon,\delta}\|v\|_{H^{3/2}}\|\phi\|_{H^{-\epsilon}}.
\]
Thus \eqref{eq:GW.2} is continuous into
$H^{-1-\epsilon-\delta}(S^1)$ for each such field.

Suppose its value is zero for a nonzero $v$. Since $v$ is continuous,
there is an open interval $I$ on which it is bounded away from zero.
Let $a_v\in\R$ denote the constant Fourier mode of $v\phi'+Qv'$.
On $I$ the kernel equation gives
\[
 P_0(v\phi'+Qv')=0
 \quad\Longrightarrow\quad v\phi'+Qv'=a_v
 \quad\Longrightarrow\quad (v\phi)'=v'\phi-Qv'+a_v.
\]
The Sobolev product and local reciprocal estimates yield
\[
 (v\phi)'\in H^{-\epsilon-\delta}_{\rm loc}(I)
 \quad\Longrightarrow\quad
 v\phi\in H^{1-\epsilon-\delta}_{\rm loc}(I).
\]
Division by $v$, which is bounded away from zero on smaller intervals,
then gives
\[
 \phi\in H^{1-\epsilon-\delta}_{\rm loc}(I)\subset L^2_{\rm loc}(I).
\]

For a nonzero real smooth cutoff $b$ in a fixed interval, let $C_N$
be the covariance matrix of the real Fourier components of $b\phi$
with frequencies $|k|\leq N$. With normalized circle measure,
\[
 \|C_N\|_{\rm op}\leq\|b\|_\infty^2,\qquad
 \operatorname{Tr}C_N
 =\sum_{|k|\leq N}\sum_{n\ne0}\frac{|\widehat b(k-n)|^2}{|n|}.
\]
Choose an integer $m$ with $\widehat b(m)\ne0$. Then
\[
 \operatorname{Tr}C_N
 \geq |\widehat b(m)|^2\sum_{1\leq|n|\leq N-|m|}\frac1{|n|}
 \longrightarrow\infty.
\]
The finite Gaussian determinant identity gives
\[
 \mathbb E\exp\left(-\sum_{|k|\leq N}|\widehat{b\phi}(k)|^2\right)
   =\det(\operatorname{Id}+2C_N)^{-1/2}
   \le\exp(-c\operatorname{Tr}C_N)\longrightarrow0.
\]
Monotone convergence yields
\[
 \mathbb E e^{-\|b\phi\|_2^2}=0,
 \qquad \mathbb G(\|b\phi\|_2=\infty)=1.
\]
A countable collection of cutoffs in rational intervals excludes
local $L^2$ on every interval, proving \eqref{eq:GW.2} simultaneously
for all $v\in H$.

\subsection{Recovery of the Brownian motion}\label{sec:field-noise-recovery}

Choose a countable list of real Fourier-coordinate tests $F_a$ and set
\[
 v_a(Y)=\Vzero F_a(Y),\qquad
 M_t^a=F_a(Y_t)-F_a(Y_0)+\int_0^t A_0F_a(Y_r)\,dr.
\]
Their brackets are
$d\langle M^a,M^b\rangle_t=2\langle v_a(Y_t),v_b(Y_t)\rangle_H\,dt$.
Each $v_a$ belongs to $H$ almost surely, since
$\int\|v_a\|_H^2\,d\mathbb G<\infty$.
If $u\in H$ is orthogonal to every $v_a(Y)$, all real Fourier
coefficients of $P_0(u\phi'+Qu')$ vanish. The response distribution is
therefore zero, and injectivity in \eqref{eq:GW.2} gives $u=0$.
Hence
\[
 \overline{\operatorname{span}\{v_a(Y):a\geq1\}}=H
                   \quad\mathbb G\text{-almost surely}.
\]
For $a,b\leq N$, let
\[
 G_N(Y)_{ab}=\langle v_a(Y),v_b(Y)\rangle_H.
\]
Write $G_N^\dagger$ for its measurable Moore--Penrose inverse and
$P_N(Y)$ for the orthogonal projection onto the span of
$v_1(Y),\ldots,v_N(Y)$. For $j\in J_H$, define
\begin{align}
 c_a^{j,N}(Y)&=\frac1{\sqrt2}\sum_{b=1}^N
       (G_N(Y)^\dagger)_{ab}\langle v_b(Y),e_j\rangle_H,\notag\\
 \shortintertext{and set}
 W_t^{j,N}&=\sum_{a=1}^N\int_0^t c_a^{j,N}(Y_r)\,dM_r^a.
 \label{eq:field-noise-approximation}
\end{align}
Use the common cutoff $\max_a|c_a^{j,N}|\leq m$ to define each integral.
Its martingale-norm error tends to zero because
\[
 \mathbb E\int_0^T\mathbf1_{\{\max_a|c_a^{j,N}(Y_r)|>m\}}
                      \|P_N(Y_r)e_j\|_H^2\,dr\longrightarrow0.
\]
Finite Gram inversion gives
\begin{align*}
 d\langle W^{j,N},W^{k,M}\rangle_t
     &=\langle P_N(Y_t)e_j,P_M(Y_t)e_k\rangle_H\,dt.\\
 \shortintertext{Doob's inequality consequently gives}
 \mathbb E\sup_{t\leq T}|W_t^{j,N}-W_t^{j,M}|^2
     &\leq4T\int\|(P_N-P_M)e_j\|_H^2\,d\mathbb G
          \longrightarrow0.
\end{align*}
The limit $W^j$ is a continuous martingale in the original filtration,
and dominated convergence gives
\begin{equation}\label{eq:field-noise-brackets}
 \langle W^j,W^k\rangle_t=\delta_{jk}t,
 \qquad d\langle M^a,W^j\rangle_t
              =\sqrt2\langle v_a(Y_t),e_j\rangle_H\,dt.
\end{equation}
L\'evy's characterization makes every finite subfamily a standard
Brownian motion. Moreover,
\[
 d\left\langle M^a-\sqrt2\sum_j\int_0^\cdot
                        \Vzero_jF_a(Y_r)\,dW_r^j\right\rangle_t=0.
\]
Thus this noise realizes the field martingales. If a noise was already
given, the same Gram calculation yields
\[
 W_t^{j,N}=\int_0^t\langle P_N(Y_r)e_j,dW_r\rangle_H.
\]
Thus
\[
 \mathbb E\sup_{t\leq T}|W_t^{j,N}-W_t^j|^2
       \leq4T\int\|(P_N-1)e_j\|_H^2\,d\mathbb G\longrightarrow0.
\]
Fix $j,N,T$ and set
\[
 \mathbf c_t=(c_a^{j,N}(Y_t))_{a=1}^N,\qquad
 \mathbf M_t=(M_t^a)_{a=1}^N,\qquad
 \Gamma_t=2G_N(Y_t).
\]
For integers $K,m\geq1$, use the bounded adapted coefficients
\begin{align*}
 \mathbf c_t^{(K)}
   &=\mathbf c_t\mathbf1_{\{\max_a|c_a^{j,N}(Y_t)|\leq K\}},\\
 \shortintertext{and their continuous past averages}
 \mathbf c_t^{(K,m)}
   &=m\int_{(t-1/m)_+}^t\mathbf c_r^{(K)}\,dr.
\end{align*}
Since $\mathbb E\int_0^T\operatorname{Tr}\Gamma_t\,dt<\infty$,
Lebesgue differentiation and the preceding cutoff estimate give
\[
 \lim_{K\to\infty}\lim_{m\to\infty}
 \mathbb E\int_0^T
  \left\langle\mathbf c_t^{(K,m)}-\mathbf c_t,
       \Gamma_t(\mathbf c_t^{(K,m)}-\mathbf c_t)\right\rangle_{\R^N}\,dt=0.
\]
For $t_\ell=\ell T/2^L$, set
\[
 J_{K,m,L}(t)=\sum_{\ell=0}^{2^L-1}
  \left\langle\mathbf c_{t_\ell}^{(K,m)},
     \mathbf M_{t_{\ell+1}\wedge t}-\mathbf M_{t_\ell\wedge t}
  \right\rangle_{\R^N}.
\]
Continuity of the past averages, It\^o's isometry and Doob's inequality yield
\[
 \lim_{K\to\infty}\lim_{m\to\infty}\lim_{L\to\infty}
 \mathbb E\sup_{t\leq T}|J_{K,m,L}(t)-W_t^{j,N}|^2=0.
\]
Each sum is a measurable field-path functional. Passing to these limits
determines the joint law of $(Y,W)$.

For an $\mathcal F_0$-measurable spectator $Z$ and bounded $f,g$,
\[
 \mathbb E[f(Y_t)g(Z)]
   =\mathbb E[S_tf(Y_0)g(Z)].
\]
Iteration of \eqref{eq:GW.1} likewise determines the joint law of
$(Y,Z,W)$ from that of $(Y_0,Z)$.
\subsection{An Osgood modulus for the Gaussian stress}\label{sec:stress-osgood}

Use the stress coefficients $\Sigma_n$ of \eqref{eq:KD.3} under
$\mathbb G$, and the standard complex coordinates $\xi_n=\sqrt nY_n$.
For $N\geq2$ set
\begin{equation}\label{eq:GO.1}
 \Bstress_N(Y,\theta)=2\operatorname{Im}\sum_{n=2}^N
                      \frac{\Sigma_n}{d_n}e^{in\theta}.
\end{equation}
By the real-score definitions following \eqref{eq:KD.3} and the basis
normalization \eqref{eq:basis}, this is equivalently
\[
 \Bstress_N(Y,\theta)=\sum_{n=2}^N\bigl(
   \Sigma_{n,c}(Y)e_{n,c}(\theta)
  +\Sigma_{n,s}(Y)e_{n,s}(\theta)\bigr),
 \qquad \alpha_n^2=\frac2{d_n}.
\]
Thus its limit is the real-frame stress $\sum_j\Sigma_j(Y)e_j(\theta)$
used in \eqref{eq:MS.4}.
There is a measurable random variable $\Lambda_Q(Y)$ with finite moments
of every finite order such that $\Bstress_N$ along dyadic $N$ converges
uniformly to a continuous periodic function $\Bstress$, and
\begin{align}
 \|\Bstress(Y,\cdot)\|_\infty&\le \Lambda_Q(Y),\notag\\
 \shortintertext{and, for $\delta=\operatorname{dist}_{S^1}(x,y)\le1$,}
 |\Bstress(Y,x)-\Bstress(Y,y)|
 &\le \Lambda_Q(Y)\delta\log(e/\delta).\label{eq:GO.2}
\end{align}
In fact $\Lambda_Q$ can be chosen with a small exponential moment.

The estimate \eqref{eq:GO.2} is proved in
Appendix~\ref{app:stress-modulus}.

\subsubsection{The inverse-mark equation}

Let $Y_t$ be any progressively measurable field process such that
\begin{equation}\label{eq:GO.7}
\int_0^T \Lambda_Q(Y_t)\,dt<\infty
                    \quad\text{almost surely for every }T.
\end{equation}
If each $Y_t$ has the Neumann Gaussian marginal, then
\[
 \mathbb E\int_0^T \Lambda_Q(Y_t)\,dt
       =T\int \Lambda_Q(Y)\,d\mathbb G(Y)<\infty,
\]
so \eqref{eq:GO.7} holds. Consider the circle equation driven by the same
cylindrical Brownian motion as the field,
\begin{equation}\label{eq:GO.8}
dq_t=-\sqrt2\sum_j e_j(q_t)\,dW_t^j+\Bstress(Y_t,q_t)\,dt.
\end{equation}
Equation \eqref{eq:GO.8} has pathwise uniqueness, conditional on the given
adapted pair $(Y,W)$ and its initial mark.

For $\delta=\operatorname{dist}_{S^1}(x,y)\leq1$,
\eqref{eq:kernel-modulus} gives
\begin{equation}\label{eq:GO.9}
\sum_j|e_j(x)-e_j(y)|^2
 =\|k(x)-k(y)\|_H^2\leq C\delta^2\log(e/\delta).
\end{equation}
For two solutions with equal initial marks, choose angular lifts
$q_t^1,q_t^2$ with the same initial value. Fix a sufficiently small
$\varepsilon>0$ and set
\[
 \Delta_t=|q_t^1-q_t^2|^2,\qquad
 \tau_\varepsilon=\inf\{t:\Delta_t\geq\varepsilon\},
\]
It\^o's formula gives, before
$\tau_\varepsilon$,
\begin{align*}
 d\Delta_t&=dM_t+b_t\,dt,\\
 \shortintertext{where the drift and bracket satisfy}
 b_t&\leq C(1+\Lambda_Q(Y_t))\Delta_t\log(C/\Delta_t),\\
 d\langle M\rangle_t&\leq C\Delta_t^2\log(C/\Delta_t)\,dt.
\end{align*}
Define the clock and its inverse by
\[
 \eta_t=\int_0^t(1+\Lambda_Q(Y_u))\,du,\qquad
 t_s=\inf\{t:\eta_t>s\},\qquad
 \widetilde\Delta_s=\Delta_{t_s\wedge\tau_\varepsilon}.
\]
The clock is finite on compact intervals and $\eta_t\geq t$.
The stopped time-changed martingale has bounded bracket density on
each finite $s$ interval. Taking expectations and using concavity yields
\[
 \mathbb E\widetilde\Delta_s\leq C\int_0^s
  \mathbb E\widetilde\Delta_r\log(C/\mathbb E\widetilde\Delta_r)\,dr.
\]
Since $\int_0^\varepsilon dr/[r\log(C/r)]=\infty$, Osgood's
inequality gives
\[
 \mathbb E\widetilde\Delta_s=0\quad(s\geq0),\qquad
 \mathbb P(\Delta_{t\wedge\tau_\varepsilon}=0\text{ for all }t)=1,
 \qquad \mathbb P(\tau_\varepsilon<\infty)=0.
\]
Hence $q^1=q^2$ almost surely on every compact time interval.

For the marked Witt derivation
$V_j=\Vzero_j-e_j(q)\partial_q$ with divergence $\Sigma_j(Y)$,
the formal generator applied to a smooth mark test is exactly
\begin{equation}\label{eq:GO.10}
\sum_j(V_j^2-\Sigma_jV_j)b(q)
            =\frac{1}{2}b''(q)+\Bstress(Y,q)b'(q).
\end{equation}
Indeed complete sine/cosine pairs give
\[
 \sum_je_j(q)e_j'(q)=0,\qquad \sum_je_j(q)^2=\frac12.
\]
For each $p<\infty$, the stress tail satisfies
\[
 \|\Bstress-\Bstress_N\|_{L^p(\mathbb G;C(\T))}
 \leq C\sum_{n>N}\alpha_n\|\alpha_n\Sigma_n\|_{L^p(\mathbb G)}
 \leq C(p,Q)\sum_{n>N}\sqrt{\frac2{d_n}}\longrightarrow0.
\]
This justifies \eqref{eq:GO.10} under every coupling with Gaussian
field marginals.

\subsection{The marked generator range}\label{sec:marked-range}

On the probability $\pi$ of Section~\ref{sec:reconstruction},
let $\mathrm{Pol}$ be the algebra of finite field polynomials times
trigonometric polynomials in $q$. The finite-flow identity gives
\begin{equation}\label{eq:MS.1}
 V_j=\Vzero_j-e_j(q)\partial_q,\qquad
 V_j^*=-V_j+\Sigma_j(Y).
\end{equation}
All Fourier sums below are taken in complete sine/cosine pairs. Set
\[
 J_{H,N}=\{(n,\tau):2\leq n\leq N,\ \tau\in\{c,s\}\}.
\]
On $\mathrm{Pol}$ the generator is defined by
\begin{align*}
 L_{\pi,N}p&=\sum_{j\in J_{H,N}}(V_j^2-\Sigma_jV_j)p,\\
 \shortintertext{and}
 L_\pi p&=\lim_{N\to\infty}L_{\pi,N}p
                  \quad\text{in }L^2(\pi),\qquad p\in\mathrm{Pol}.
\end{align*}
Write $j\leq N$ for $j\in J_{H,N}$. For $p=P(Y,Z)b(q)$, with $P$ a
field polynomial and $b$ a trigonometric polynomial, the paired product rule gives
\begin{align*}
 L_{\pi,N}(Pb)
  &=-b\sum_{j\le N}(\Vzero_j)^*\Vzero_jP
     +P\left(\sum_{n=2}^N\alpha_n^2\,b''
                          +\Bstress_N b'\right)
     -2b'\sum_{j\le N}e_j\Vzero_jP.
\end{align*}
Here $A_0$ acts only on $Y$, and the Gaussian field calculations give
\[
 \sum_{j\le N}(\Vzero_j)^*\Vzero_jP\longrightarrow A_0P
                 \quad\text{in }L^2(\mathbb G\otimes\mathbb G),
 \qquad \sum_{n=2}^\infty\alpha_n^2=\tfrac12.
\]
The two remaining tails satisfy
\begin{align*}
 \left\|\sum_{j>N}e_j\Vzero_jP\right\|_{L^2(\pi)}
  &\le\left(\sum_{n>N}\alpha_n^2\right)^{1/2}
            \|\Vzero P\|_{L^2(\mathbb G\otimes\mathbb G;H)}
       \longrightarrow0,\\
 \shortintertext{and}
 \|P(\Bstress_N-\Bstress)\|_{L^2(\pi)}
  &\le\|P\|_{L^4(\mathbb G\otimes\mathbb G)}
          \left\|\sup_q|\Bstress_N-\Bstress|\right\|_{L^4(\mathbb G)}
       \longrightarrow0.
\end{align*}
Thus $L_{\pi,N}p\to L_\pi p$ in $L^2(\pi)$, with
\[
 L_\pi(Pb)=-bA_0P
       +P(\tfrac12b''+\Bstress b')
       -2b'\sum_j e_j\Vzero_jP.
\]

\begin{proposition}\label{prop:marked-range}
For every $\lambda>0$,
\begin{equation}\label{eq:MS.2}
\overline{(\lambda-L_{\pi})\mathrm{Pol}}^{\,L^1(\pi)}=L^1(\pi).
\end{equation}
\end{proposition}

\begin{proof}
Fix $\lambda>0$. By $L^1$--$L^\infty$ duality, it suffices to exclude a
bounded function $u$ with
\[
 \int u(\lambda-L_\pi)F\,d\pi=0
                 \qquad(F\in\mathrm{Pol}).
\]
Superposition will turn the annihilator into a probability curve;
we first verify the coordinate martingale problem and its uniqueness.
Encode $(Y,Z,q)$ by the real and imaginary Fourier coordinates of
the two fields and the two real coordinates $\cos q,\sin q$.
Denote these coordinate functions by $x_i$, $i\geq1$.
Equip this countable-coordinate space with its product sigma-algebra and topology.
Its drift and covariance coefficients are
\[
 \mu_i=L_\pi x_i,\qquad
 \sigma_{ik}=L_\pi(x_ix_k)-x_iL_\pi x_k-x_kL_\pi x_i.
\]
For field coordinates, the Gaussian Fock calculations give
\begin{equation}\label{eq:MS.3}
\mu_i\in\bigcap_{p<\infty}L^p(\mathbb G),\qquad
 \sigma_{ik}=2\langle \Vzero x_i,\Vzero x_k\rangle_H
       \in\bigcap_{p<\infty}L^p(\mathbb G).
\end{equation}
The passive coordinates have zero drift and covariance. For a
smooth mark function $f$,
\begin{equation}\label{eq:MS.4}
L_{\pi}f(q)=\tfrac12 f''(q)+\Bstress(Y,q)f'(q),
 \quad
 \Bstress(Y,q)=\sum_j\Sigma_j(Y)e_j(q).
\end{equation}
The Gaussian stress result gives
$\sup_q|\Bstress(Y,q)|\in L^p(\mathbb G)$ for every finite $p$.
Mixed covariances obey
\begin{equation}\label{eq:MS.5}
 |\sigma_{i,f(q)}|
 \leq \sqrt2\,\|f'\|_\infty\|\Vzero x_i\|_H,
\end{equation}
because $\sum_je_j(q)^2=1/2$. Thus every coordinate drift and
covariance is integrable under every curve of measures whose
densities relative to $\pi$ are uniformly bounded.

For $F\in C_b^\infty(\R^d\times\T)$ on a finite field-coordinate
list, take $F_m\in\mathrm{Pol}$
as in Lemma~\ref{lem:smooth-polynomial-approximation}. Equations
\eqref{eq:MS.3}--\eqref{eq:MS.5} and H\"older's inequality give
\begin{align*}
 \|F_m-F\|_{L^1(\pi)}
  +\|L_\pi F_m-L_\pi F\|_{L^1(\pi)}
 &\leq C_F\|F_m-F\|_{2,2,2}\longrightarrow0,\\
 \shortintertext{and}
 \|V(F_m-F)\|_{L^2(\pi;H)}
 &\leq C_F\|F_m-F\|_{1,1,4}\longrightarrow0.
\end{align*}
These estimates use only the Gaussian field marginal, since for
every integrable $G$,
\[
 \int|G(Y,Z,q)|\,d\pi
 \leq\int\sup_{q\in\T}|G(Y,Z,q)|\,
                  d(\mathbb G\otimes\mathbb G)(Y,Z).
\]
They also extend a bounded annihilator identity to these tests.

Fixed-direction integration by parts and paired coordinate convergence give
\[
 \int L_\pi F\,d\pi
 =-\lim_{N\to\infty}\sum_{j\in J_{H,N}}
      \langle V_jF,V_j1\rangle_{L^2(\pi)}=0.
\]
Hence the constant curve $\pi_t=\pi$ solves its Fokker--Planck equation.

Let $x=(x_1,\ldots,x_d)$ be a finite coordinate list,
$f\in C_b^2(\R^d)$, and let $f_\epsilon$ be standard smooth
mollifications. For $\Phi=f(x)$ and $\Phi_\epsilon=f_\epsilon(x)$,
\begin{align*}
 L_\pi\Phi
 &=\sum_{i=1}^d\mu_i\partial_i f(x)
       +\frac12\sum_{i,k=1}^d\sigma_{ik}\partial_{ik}f(x),\\
 \shortintertext{and \eqref{eq:MS.3}--\eqref{eq:MS.5} give}
 \|\Phi_\epsilon-\Phi\|_{L^1(\pi)}
   +\|L_\pi\Phi_\epsilon-L_\pi\Phi\|_{L^1(\pi)}&\longrightarrow0.\\
 \shortintertext{If $\sup_{t\leq T}\|d\pi_t/d\pi\|_\infty\leq M$, then}
 &\sup_{t\leq T}\left|\int(\Phi_\epsilon-\Phi)\,d\pi_t\right|
   +\int_0^T\!\int|L_\pi\Phi_\epsilon-L_\pi\Phi|\,d\pi_t\,dt\\
 &\qquad\leq M\|\Phi_\epsilon-\Phi\|_{L^1(\pi)}
             \!+MT\|L_\pi\Phi_\epsilon-L_\pi\Phi\|_{L^1(\pi)}
                  \longrightarrow0.
\end{align*}
Thus the Fokker--Planck identities extend to all $C_b^2$
finite-coordinate tests, including $(\cos q,\sin q)$.

We apply \cite[Theorem 3.3, with Definitions 3.1--3.2]{Qiao} on
$\R^\N$ with its product topology. Choose measurable versions of all
coefficients on their common full-measure set and extend them by zero
elsewhere. Every finite covariance matrix is symmetric and nonnegative:
\[
 \sum_{i,k=1}^d r_i\sigma_{ik}r_k
   =2\left\|\sum_{i=1}^d r_iVx_i\right\|_H^2\geq0
   \qquad(r\in\R^d).
\]

For a curve with $\sup_{t\leq T}\|d\pi_t/d\pi\|_\infty\leq M$,
\eqref{eq:MS.3}--\eqref{eq:MS.5} give, for each $i,k$,
\[
 \int_0^T\int(|\mu_i|+|\sigma_{ik}|)\,d\pi_t\,dt
 \leq MT\int(|\mu_i|+|\sigma_{ik}|)\,d\pi<\infty.
\]
The curves used below are weakly continuous probability curves.
The preceding approximation verifies their Fokker--Planck identities
on all required $C_b^2$ finite-coordinate tests.
The superposition theorem therefore supplies a probability on
$C([0,T];\R^\N)$ with these marginals and the coordinate martingales.

Let $i_c,i_s$ be the indices of the two mark coordinates. Marginal
support and path continuity give
\[
 \mathbb P\bigl(x_{i_c}(t)^2+x_{i_s}(t)^2=1
               \text{ for all }t\in[0,T]\bigr)
 =\mathbb P\left(\bigcap_{t\in\mathbb Q\cap[0,T]}
          \{x_{i_c}(t)^2+x_{i_s}(t)^2=1\}\right)=1.
\]
For every passive field coordinate $x_i$,
\[
 \mu_i=\sigma_{ii}=0
 \quad\Longrightarrow\quad
 \langle x_i-x_i(0)\rangle_t=0
 \quad\Longrightarrow\quad x_i(t)=x_i(0)\quad\text{almost surely}.
\]

To identify the marked law, consider two such martingale solutions with the same initial
law and Gaussian $Y$ marginals. For a polynomial field test $F$,
let $x\in\R^d$ be its finite real-coordinate list and choose
$\chi\in C_c^\infty(\R^d)$ equal to one near zero. Put
$F_R=\chi(x/R)F$. Gaussian moments
and \eqref{eq:MS.3} give
\[
 \sup_{t\leq T}\mathbb E|F_R(Y_t)-F(Y_t)|
 +\mathbb E\int_0^T\bigl(|A_0(F_R-F)(Y_t)|
                 +\|\Vzero(F_R-F)(Y_t)\|_H^2\bigr)\,dt
 \longrightarrow0.
\]
Thus the bounded cylinder martingale identities pass to $F$.
The backward-test calculation of \eqref{eq:GW.1} applies on $[0,T]$.
Apply \eqref{eq:field-noise-approximation}--\eqref{eq:field-noise-brackets}
to the field coordinate martingales $M_i$. This constructs $W$ in the
full filtration from the field path, with
\[
 dM_i=\sqrt2\sum_j\Vzero_jx_i(Y_t)\,dW_t^j,
 \qquad \langle W^j,W^k\rangle_t=\delta_{jk}t.
\]
Let $q_t$ be a continuous angular lift and set
\[
 M_t^q=q_t-q_0-\int_0^t\Bstress(Y_r,q_r)\,dr.
\]
The encoded mark generator and its mixed covariances give
\begin{align*}
 d\langle M^q\rangle_t&=dt,\\
 d\langle M^q,M_i\rangle_t
 &=-2\langle k(q_t),\Vzero x_i(Y_t)\rangle_H\,dt.\\
 \shortintertext{Applying the Gram coefficients gives}
 d\langle M^q,W^{j,N}\rangle_t
 &=-\sqrt2\langle k(q_t),P_N(Y_t)e_j\rangle_H\,dt.\\
 \shortintertext{Passing to the limit yields}
 d\langle M^q,W^j\rangle_t&=-\sqrt2 e_j(q_t)\,dt.
\end{align*}
Since $\|k(q)\|_H^2=1/2$,
\[
 d\left\langle M^q+\sqrt2\sum_j\int_0^\cdot
                              e_j(q_r)\,dW_r^j\right\rangle_t
          =(1-2\|k(q_t)\|_H^2)\,dt=0.
\]
Thus
\begin{equation}\label{eq:MS.6}
dq_t=-\sqrt2\sum_j e_j(q_t)\,dW_t^j+\Bstress(Y_t,q_t)\,dt.
\end{equation}

Let $q^1,q^2$ denote two candidate mark paths with the same initial
joint law of $(Y_0,Z,q_0)$. For every bounded cylinder function
$\Phi$ of the future field segment $Y_{[s,T]}$, iterating
\eqref{eq:GW.1} gives
\[
 \mathbb E[\Phi(Y_{[s,T]})\mid\mathcal F_s]
     =\mathbb E[\Phi(Y_{[s,T]})\mid Y_s].
\]
Thus the joint law of $(Y_{[0,T]},Z,q_0)$ is common to both
solutions. Write $\mathbb E_r$ for expectation under the $r$th
candidate law.

For $r=1,2$, let ${K^{\mathrm{path}}}_r$ be the conditional law
of $q^r_{[0,T]}$ given $(Y_{[0,T]},Z,q_0)$. If $A$ depends only on
$q^r_{[0,s]}$ and $\Phi$ is a bounded function of $Y_{[s,T]}$,
the preceding identity and a monotone-class argument give
\begin{align*}
 &\mathbb E_r[\mathbf1_A\Phi\mid Y_{[0,s]},Z,q_0]\\
 &\qquad=\mathbb E_r[\mathbf1_A\mid Y_{[0,s]},Z,q_0]
         \mathbb E_r[\Phi\mid Y_{[0,s]},Z,q_0].\\
 \shortintertext{Thus conditioning on the future does not change a past mark event:}
 &{K^{\mathrm{path}}}_r(A\mid Y_{[0,T]},Z,q_0)
       ={K^{\mathrm{path}}}_r(A\mid Y_{[0,s]},Z,q_0).
\end{align*}
This is the compatibility condition for conditional-product coupling
in \cite[Lemma~2.11]{Kurtz}. We therefore couple the two marks
conditionally independently given the common base path, with their
respective conditional laws ${K^{\mathrm{path}}}_1$ and
${K^{\mathrm{path}}}_2$. The coupled probability is
\[
 d\widehat{\mathbb P}
  =\operatorname{Law}(Y_{[0,T]},Z,q_0)(dY,dZ,dq_0)
       {K^{\mathrm{path}}}_1(dq^1\mid Y,Z,q_0)
       {K^{\mathrm{path}}}_2(dq^2\mid Y,Z,q_0).
\]
Let $\widehat{\mathcal F}_s$ be generated by
$(Y_{[0,s]},Z,q^1_{[0,s]},q^2_{[0,s]})$, with its usual augmentation.
Compatibility says that for every bounded $G$ depending on
$q^2_{[0,s]}$,
\[
 {K^{\mathrm{path}}}_2(G\mid Y_{[0,T]},Z,q_0)
   ={K^{\mathrm{path}}}_2(G\mid Y_{[0,s]},Z,q_0).
\]
Consequently, for bounded $F$ measurable with respect to
$(Y_{[0,s]},Z,q^1_{[0,s]})$,
\[
 \widehat{\mathbb E}[(W_t^j-W_s^j)FG]
 =\mathbb E_1\!\left[(W_t^j-W_s^j)F
   {K^{\mathrm{path}}}_2(G\mid Y_{[0,s]},Z,q_0)\right]=0,
\]
where the last equality uses the Brownian martingale property of $W$
in the first candidate's filtration. The symmetric argument is not
needed for $W$, which is a functional of the common base path.
Its quadratic covariations are unchanged by the coupling, so
$\langle W^j,W^k\rangle_t=\delta_{jk}t$. L\'evy's characterization
therefore makes $W$ Brownian in $\widehat{\mathcal F}_t$.

For $r=1,2$, let $M_i^r$ be the $i$th coordinate martingale and set
\[
 N_t^{r,ik}=M_i^r(t)M_k^r(t)
              -\int_0^t\sigma_{ik}(x_u^r)\,du.
\]
For bounded $F$ measurable with respect to
$(Y_{[0,s]},Z,q^1_{[0,s]})$ and bounded $G$ measurable with respect to
$q^2_{[0,s]}$, the same kernel identity gives
\begin{align*}
 \widehat{\mathbb E}[(M_i^1(t)-M_i^1(s))FG]
 &=\mathbb E_1\!\left[(M_i^1(t)-M_i^1(s))F
   {K^{\mathrm{path}}}_2(G\mid Y_{[0,s]},Z,q_0)\right]=0,\\
 \widehat{\mathbb E}[(N_t^{1,ik}-N_s^{1,ik})FG]
 &=\mathbb E_1\!\left[(N_t^{1,ik}-N_s^{1,ik})F
   {K^{\mathrm{path}}}_2(G\mid Y_{[0,s]},Z,q_0)\right]=0.
\end{align*}
The symmetric identities hold for the second solution. Hence both
coordinate martingale problems, including their covariance identities,
remain valid in the coupled filtration. Since each marginal satisfies
\eqref{eq:MS.6} and the recovered $W$ is the same base-path functional,
both marks solve \eqref{eq:MS.6} with this common Brownian family.
Pathwise uniqueness gives $q^1=q^2$ almost surely, so the two joint
laws coincide.

To eliminate the annihilator, assume first that $u$ is real.
Testing with constants gives
$\int u\,d\pi=0$. Testing with field polynomials independent of
the mark, and using the Gaussian scalar operator core tensored
with the passive field, gives
\begin{equation}\label{eq:MS.7}
\mathbb E_\pi[u\mid Y,Z]=0.
\end{equation}
Indeed this conditional expectation annihilates the dense range
of $\lambda+A_0\otimes I_Z$ in $L^2(\mathbb G\otimes\mathbb G)$.

Fix $T>0$ and choose $\varepsilon>0$ so small that
$\varepsilon e^{\lambda T}\|u\|_\infty<1/2$. Then
\begin{equation}\label{eq:MS.8}
\pi_t=(1+\varepsilon e^{\lambda t}u)\pi,
                          \qquad 0\le t\le T,
\end{equation}
is a curve of positive probability measures.
It is continuous in total variation, since
\[
 \|\pi_t-\pi_s\|_{\rm TV}
 =\varepsilon|e^{\lambda t}-e^{\lambda s}|\|u\|_{L^1(\pi)}
 \longrightarrow0\qquad(t\to s).
\]
For $F\in\mathrm{Pol}$,
the annihilator identity and $\int L_\pi F\,d\pi=0$ yield
\begin{align*}
 \frac{d}{dt}\int F\,d\pi_t
 &=\varepsilon\lambda e^{\lambda t}\int uF\,d\pi
  =\varepsilon e^{\lambda t}\int uL_\pi F\,d\pi
  =\int L_\pi F\,d\pi_t.
\end{align*}
For bounded $G(Y,Z)$, \eqref{eq:MS.7} gives
\begin{align*}
 \int G(Y,Z)\,d\pi_t
 &=\int G(Y,Z)\,d\pi
  +\varepsilon e^{\lambda t}\int G(Y,Z)\mathbb E_\pi[u\mid Y,Z]\,d\pi\\
 &=\int G\,d(\mathbb G\otimes\mathbb G).
\end{align*}
The extension to $C_b^2$ tests permits superposition for
$(\pi_t)_{0\le t\le T}$. Its density in \eqref{eq:MS.8} is
uniformly bounded relative to $\pi$, so the coordinate
integrability follows from \eqref{eq:MS.3}--\eqref{eq:MS.5}.
Denote the resulting path law by $\mathbb P^{\mathrm{curve}}$.
It solves the marked martingale problem, has initial law $\pi_0$,
and has time marginals $\pi_t$.

Superpose the stationary curve $\pi$ and denote its path law by
$\mathbb P^0$. Define
\[
 \frac{d\mathbb P^\varepsilon}{d\mathbb P^0}
      =1+\varepsilon u(Y_0,Z_0,q_0).
\]
Write $\mathbb E^0,\mathbb E^\varepsilon$ for the corresponding expectations.
For a coordinate martingale $M$ and bounded
$\mathcal F_s$-measurable $\zeta$,
\[
 \mathbb E^\varepsilon[(M_t-M_s)\zeta]
 =\mathbb E^0[(M_t-M_s)\zeta(1+\varepsilon u(Y_0,Z_0,q_0))]=0.
\]
For bounded $F,G$, \eqref{eq:GW.1} and \eqref{eq:MS.7} give
\begin{align*}
 \mathbb E^\varepsilon[F(Y_t)G(Z)]
 &=\mathbb E^0[F(Y_t)G(Z)]
  +\varepsilon\mathbb E^0[u(Y_0,Z,q_0)S_tF(Y_0)G(Z)]\\
 &=\int F\,d\mathbb G\int G\,d\mathbb G.
\end{align*}
Thus $\mathbb P^\varepsilon$ has initial law $\pi_0$ and
field marginals $\mathbb G\otimes\mathbb G$ at every time, as does
$\mathbb P^{\mathrm{curve}}$. The marked-path uniqueness proved above
therefore gives
$\mathbb P^\varepsilon=\mathbb P^{\mathrm{curve}}$. In particular,
its time marginals are \eqref{eq:MS.8}.

Taking time marginals contracts the total variation of signed
measures. Thus
\begin{equation}\label{eq:MS.9}
 \varepsilon e^{\lambda t}\|u\,\pi\|_{\rm TV}
     =\|\pi_t-\pi\|_{\rm TV}
     \le\varepsilon\|u\,\pi\|_{\rm TV}.
\end{equation}
For any $t>0$ these inequalities imply $u=0$. Apply the argument
to real and imaginary parts for a complex annihilator. This proves
\eqref{eq:MS.2}.
\end{proof}

\section{Conditional compression and weak uniqueness}\label{sec:conditional-stress}

Conditional compression transfers the marked range to the base.
We control its correction, remove the capacity tilt, and determine
the conditional state and noise laws.

\subsection{Conditional score bounds}

Retain $I$, $K$, $\pi$ and $\nu_\beta$ from
Section~\ref{sec:reconstruction}. For $i,j\in{J_H}$ define
\[
 \eta_j=\Sigma_j-I\rho_j^\beta,\qquad
 H^{\rm field}_{ij}=\Vzero_i\Sigma_j-\Sigma_{\nabla_{e_i}e_j},
 \qquad \Rop g=(K(\eta_jg))_j.
\]
The full connection, including its low components, is used in
$H^{\rm field}$.
In the next lemma $a=(a_j)$ and $b=(b_{ij})$ are finitely supported
complex arrays, and $K$ acts on each displayed scalar expression.

\begin{lemma}\label{lem:conditional-stress-bounds}
There is $C_\kappa<\infty$ such that, $\nu_\beta$-almost everywhere,
\begin{align}
 K\left|\sum_j a_j\eta_j\right|^2
      &\leq C_\kappa\sum_j|a_j|^2,\label{eq:conditional-score-bound}\\
 \shortintertext{and}
 K\left|\sum_{i,j}b_{ij}
       \bigl(\eta_i\eta_j-IK(\eta_i\eta_j)\bigr)\right|^2
 &+K\left|\sum_{i,j}b_{ij}
       \bigl(H^{\rm field}_{ij}-IKH^{\rm field}_{ij}\bigr)\right|^2
       \leq C_\kappa\sum_{i,j}|b_{ij}|^2.
       \label{eq:conditional-matrix-bound}
\end{align}
\end{lemma}
The proof is given in Appendix~\ref{app:stress-covariance}.
Its covariance and Bloch-mean hypotheses hold uniformly for $K_h^Q$
by Section~\ref{sec:gaussian-lift}.

Conditional Cauchy--Schwarz gives
\begin{equation}\label{eq:B.3}
 \|\Rop g\|_{L^2(\nu_\beta;H)}\le C_*\|g\|_{L^2(\pi)}.
\end{equation}
Indeed, with the supremum over finite sequences $a$,
\begin{align*}
 \|\Rop g(h)\|_H^2
 &=\sup_{\|a\|_{\ell^2}\leq1}
   \left|K_h^Q\left[(g-K_h^Qg)\sum_j\overline{a_j}\eta_j\right]\right|^2\\
 &\leq C_*^2\operatorname{Var}_{K_h^Q}(g)
 \leq C_*^2K_h^Q|g|^2.\\
 \shortintertext{Integrating over the welding gives}
 \|\Rop g\|_{L^2(\nu_\beta;H)}^2
 &\leq C_*^2\int K_h^Q|g|^2\,d\nu_\beta(h)
 =C_*^2\|g\|_{L^2(\pi)}^2.
\end{align*}

\subsection{The tilted gradient}

The bounded finite-direction scores $\rho_j^\beta$ give a densely
defined adjoint on cylinder-valued vectors. Thus the cylinder gradient
is closable in $L^2(\nu_\beta)$. Denote its closure by $D_\beta$ and set
\[
 A_\beta=D_\beta^*D_\beta,\qquad L_\beta=-A_\beta,\qquad
 \|F\|_{\Dom(D_\beta)}^2=\|F\|_{L^2(\nu_\beta)}^2
                              +\|D_\beta F\|_{L^2(\nu_\beta;H)}^2.
\]

For the tilted Hodge estimate and finite fields $v,w$, set
$\operatorname{Hess}c(v,w)=D_vD_wc-D_{\nabla_vw}c$.
The capacity-score calculation in Appendix~\ref{app:score-hessian},
summarized in \eqref{eq:capacity-hessian-bound}, gives
\begin{equation}\label{eq:DA.2}
\|\operatorname{Hess}c\|_{\widehat H\to\widehat H}
                                      \le C.
\end{equation}
Consequently
\[
 \|H^\rho-\beta\operatorname{Hess}c+\Ric\|_{\widehat H\to\widehat H}
 \leq \|H^\rho\|+|\beta|\|\operatorname{Hess}c\|+\|\Ric\|
 \leq C_\beta,
\]
and every fixed tilted score component is bounded.

For a finite horizontal covector $U$, extended by zero on the low
components, the tilted Hodge identity is
\begin{equation}\label{eq:tilted-hodge}
 \|D_\beta^*U\|_{L^2(\nu_\beta)}^2+\|dU\|_{L^2(\nu_\beta)}^2
 =\|\nabla U\|_{L^2(\nu_\beta)}^2
  +\int\bigl((H^\rho-\beta\operatorname{Hess}c+\Ric)U,U\bigr)_{\widehat H}
                                                        \,d\nu_\beta.
\end{equation}
The horizontal--low decomposition used in
Proposition~\ref{prop:horizontal-upper-hodge} therefore gives
\begin{align*}
 \|D_\beta^*U\|_2^2
 &\leq \|\nabla_HU\|_2^2+C_\beta\|U\|_2^2
       +2\sum_{a=0,c,s}|\langle\nabla_aU,K_aU\rangle_{L^2(\nu_\beta)}|.
\end{align*}
For fixed horizontal modes $p,q$, tilted integration by parts gives
\begin{align*}
 \langle[\nabla_p,\nabla_q]U,K_aU\rangle_{\nu_\beta}
   &=\langle\nabla_qU,(-\nabla_p+\rho_p)K_aU\rangle_{\nu_\beta}\\
   &\quad-\langle\nabla_pU,(-\nabla_q+\rho_q)K_aU\rangle_{\nu_\beta}
          +\Delta_{p,q,a}^{\beta},\\
 \shortintertext{where the change of score contributes}
 \Delta_{p,q,a}^{\beta}
   &=\langle\nabla_qU,(\rho_p^\beta-\rho_p)K_aU\rangle_{\nu_\beta}
     -\langle\nabla_pU,(\rho_q^\beta-\rho_q)K_aU\rangle_{\nu_\beta},\\
 \shortintertext{and}
 |\Delta_{p,q,a}^{\beta}|
   &\leq\|K_a\|
      \bigl(\|\rho_p^\beta-\rho_p\|_\infty
                +\|\rho_q^\beta-\rho_q\|_\infty\bigr)
      \|\nabla_HU\|_2\|U\|_2\\
   &\leq C_{p,q,\beta}\bigl(\|\nabla_HU\|_2^2+\|U\|_2^2\bigr).
\end{align*}
Here $\rho_p^\beta-\rho_p=-\beta D_pc$, and both fixed score components
are bounded. Expanding $\nabla_p(K_aU)$ as in
Proposition~\ref{prop:horizontal-upper-hodge} and using \eqref{eq:HD.1}
therefore gives
\[
 |\langle[\nabla_p,\nabla_q]U,K_aU\rangle_{L^2(\nu_\beta)}|
 \leq C_{p,q,\beta}
       \bigl(\|\nabla_HU\|_2^2+\|U\|_2^2\bigr),
\]
The finite bracket representation of the low directions and the bounded
curvature corrections yield, after closure on the minimal horizontal covariant domain,
\begin{equation}\label{eq:DA.3}
\|D_\beta^*U\|_{L^2(\nu_\beta)}
 \le C_\beta(\|\nabla_HU\|_2+\|U\|_2).
\end{equation}

\subsection{Minimal and maximal base gradients after the tilt}
\label{sec:tilted-domains}

Proposition~\ref{prop:horizontal-domain} identifies the minimal and
maximal gradients under $\nu$. Let $D_{\beta,\max}$ be the adjoint
of the tilted divergence on finite cylinder-valued fields:
\[
 \delta_{\beta,0}U=\sum_j(\rho_j^\beta U_j-D_jU_j),\qquad
 D_{\beta,\max}=\delta_{\beta,0}^*.
\]
Write $w_\beta=z_\beta^{-1}e^{\beta c}$ and choose a
smooth cutoff $\chi_R(c)$, equal to one on $c\le R$ and zero on
$c\ge2R$, with derivative $O(R^{-1})$.  Equation \eqref{eq:N.11} gives
\begin{equation}\label{eq:DA.4}
\|D\chi_R(c)\|_H\le C R^{-1/2}.
\end{equation}
Set $\varphi_R=\chi_R/w_\beta$. Lemma~\ref{lem:raw-geometric-tests}
gives $\varphi_R\in C_{\mathrm{geo}}\subset\Dom(D)$, with bounded value and
gradient. Original cylinder approximations also converge in the
tilted form norm, since
\[
 \|F\|_{L^2(\nu_\beta)}^2+\|DF\|_{L^2(\nu_\beta;H)}^2
 \leq\|w_\beta\|_\infty
       \bigl(\|F\|_{L^2(\nu)}^2+\|DF\|_{L^2(\nu;H)}^2\bigr).
\]
Thus $\varphi_R\in\Dom(D_\beta)$. For a finite cylinder-valued $U$,
approximation of $\varphi_R$ gives
$\varphi_RU\in\Dom(\overline{\delta_{\beta,0}})$ and
\[
 \overline{\delta_{\beta,0}}(\varphi_RU)
   =w_\beta^{-1}\bigl(\chi_R\delta_0U
                             -\langle D\chi_R,U\rangle_H\bigr).
\]
For $u\in\Dom(D_{\beta,\max})$, this gives
\begin{align*}
 \langle\chi_Ru,\delta_0U\rangle_{L^2(\nu)}
 &=\langle u,\overline{\delta_{\beta,0}}(\varphi_RU)\rangle_{L^2(\nu_\beta)}
                  +\langle uD\chi_R,U\rangle_{L^2(\nu;H)}\\
 &=\langle\chi_RD_{\beta,\max}u+uD\chi_R,U\rangle_{L^2(\nu;H)}.\\
 \shortintertext{Local equivalence of the two measures gives}
 \|\chi_Ru\|_{L^2(\nu)}^2
   +\|\chi_RD_{\beta,\max}u+uD\chi_R\|_{L^2(\nu;H)}^2
 &\leq C_{\beta,R}\bigl(\|u\|_{L^2(\nu_\beta)}^2
                      +\|D_{\beta,\max}u\|_{L^2(\nu_\beta;H)}^2\bigr).
\end{align*}
Proposition~\ref{prop:horizontal-domain} and the preceding cylinder
approximation yield
\[
 \chi_Ru\in\Dom(D)\cap\Dom(D_\beta),\qquad
 D_\beta(\chi_Ru)=\chi_RD_{\beta,\max}u+uD\chi_R.
\]
Finally,
\begin{align*}
 \|\chi_Ru-u\|_{L^2(\nu_\beta)}&\longrightarrow0,\\
 \shortintertext{and}
 \|D_\beta(\chi_Ru)-D_{\beta,\max}u\|_{L^2(\nu_\beta;H)}
 &\leq\|(\chi_R-1)D_{\beta,\max}u\|_2
          +CR^{-1/2}\|u\|_2\longrightarrow0.
\end{align*}
Closedness proves $D_\beta=D_{\beta,\max}$.
Define $\Dom(\nabla_{H,\beta,\min})$ and
$\Dom(\nabla_{H,\beta,\max})$ as in Section~\ref{sec:horizontal-divergence},
using $\nu_\beta$ and the tilted weak adjoints.
For $U\in\Dom(\nabla_{H,\beta,\max})$, the slow Fourier cutoffs $Q_L$
in \eqref{eq:slow-cutoff} satisfy
\begin{align*}
 \nabla_H(Q_LU)&=Q_L\nabla_HU+([\Xi_j,Q_L]U)_{j\in J_H},\\
 \shortintertext{with}
 \|([\Xi_j,Q_L]U)_j\|_{L^2(\nu_\beta;\HS(H,\widehat H))}
   &\leq CL^{-1}\|U\|_{L^2(\nu_\beta;H)},\\
 \shortintertext{so that}
 \|Q_LU-U\|_2+\|\nabla_H(Q_LU)-\nabla_HU\|_2&\longrightarrow0.
\end{align*}
Moreover, the connection column bound gives
\begin{align*}
 D_j(Q_LU)&=Q_L\nabla_jU-Q_L\Xi_jU,\\
 \shortintertext{and hence}
 \sum_{j\in J_H}\|D_j(Q_LU)\|_2^2
   &\leq2\|\nabla_HU\|_2^2+C(L)\|U\|_2^2.
\end{align*}
Every selected component of $Q_LU$ therefore belongs to
$\Dom(D_{\beta,\max})=\Dom(D_\beta)$.
Choose $U_{L,k}\in\mathcal V_{\mathrm{fin}}(H)$ on the same finite
Fourier support, approximating these components in the $D_\beta$ form norm.
For fixed $L$, the connection column bound implies
\begin{align*}
 &\|U_{L,k}-Q_LU\|_2
  +\|\nabla_H(U_{L,k}-Q_LU)\|_2
 \\
 &\qquad\leq C(L)\bigl(\|U_{L,k}-Q_LU\|_2
                 +\|D_\beta(U_{L,k}-Q_LU)\|_2\bigr)
 \longrightarrow0.
\end{align*}
Here $D_\beta$ acts componentwise. A diagonal choice in $L,k$ and
closedness of $D_\beta^*$ in \eqref{eq:DA.3} give
\[
 \Dom(\nabla_{H,\beta,\min})
   =\Dom(\nabla_{H,\beta,\max})
   =:\Dom(\nabla_{H,\beta})\subset\Dom(D_\beta^*).
\]
Moreover,
\[
 1\in\Dom(D_\beta),\qquad D_\beta1=0.
\]
The conservative Markov form therefore has resolvent bound
\[
 \|(\lambda+A_\beta)^{-1}\|_{L^p\to L^p}\leq\lambda^{-1}
 \qquad(\lambda>0,\ 1\leq p\leq\infty).
\]

\subsection{The marked form and base pullbacks}\label{sec:marked-form}

Let $C_{\mathrm{joint}}$ be the smallest algebra containing
$\mathrm{Pol}$ and $I\Ccyl$ that is stable under bounded smooth scalar
compositions with bounded derivatives. On this algebra set
\begin{equation}\label{eq:DA.7}
 V_j=\Vzero_j-e_j(q)\partial_q,
 \qquad V_j(IF)=I(D_jF).
\end{equation}
Bounded base gradients and Gaussian moments give
\[
 \|g\|_{L^p(\pi)}+\|Vg\|_{L^p(\pi;H)}<\infty,
                  \qquad g\in C_{\mathrm{joint}},\quad1\leq p<\infty.
\]
The pullback identity in \eqref{eq:DA.7} follows from finite-flow
equivariance.

The exact joint Radon--Nikodym cocycle gives finite-direction
integration by parts between any two tests in $C_{\mathrm{joint}}$,
with divergence $\Sigma_j(Y)$.  In particular, for $p\in\mathrm{Pol}$,
\begin{equation}\label{eq:DA.8}
\langle V_jg,p\rangle
        =\langle g,(-V_j+\Sigma_j)p\rangle,
          \qquad g\in C_{\mathrm{joint}},
\end{equation}
and $(-V_j+\Sigma_j)p\in L^2(\pi)$. Smooth finite-coordinate density
and Lemma~\ref{lem:smooth-polynomial-approximation} give
$\overline{\mathrm{Pol}}^{\,L^2(\pi)}=L^2(\pi)$. If
$g_k\to0$ in $L^2(\pi)$ and $Vg_k\to U$ in $L^2(\pi;H)$, then
\[
 \langle U_j,p\rangle
   =\lim_k\langle g_k,(-V_j+\Sigma_j)p\rangle=0,
       \qquad j\in{J_H},\quad p\in\mathrm{Pol},
 \qquad U=0.
\]
The closure $V$ therefore defines the Markov form with generator
$A_\pi=V^*V$. Closing \eqref{eq:DA.7} gives
\begin{equation}\label{eq:DA.9}
I\Dom(D_\beta)\subset\Dom(V),
                           VI=I_HD_\beta.
\end{equation}

The paired calculation in Section~\ref{sec:marked-range} gives
$L_{\pi,N}p\to L_\pi p$ in $L^2(\pi)$.
For $p=P(Y,Z)b(q)$, let $V_N$ retain the directions $j\in J_{H,N}$. Then
\begin{align*}
 \|V_Np-Vp\|_{L^2(\pi;H)}
 &\le\|b\|_\infty
       \|\mathbf1_{\{j>N\}}\Vzero P\|_{L^2(\mathbb G\otimes\mathbb G;H)}\\
 &\quad+\|b'\|_\infty\|P\|_{L^2(\mathbb G\otimes\mathbb G)}
             \left(\sum_{n>N}\alpha_n^2\right)^{1/2}
 \longrightarrow0.
\end{align*}
For $g\in C_{\mathrm{joint}}$, \eqref{eq:DA.8} now passes to the limit as
\[
 \langle Vp,Vg\rangle
    =\lim_{N\to\infty}\langle V_Np,V_Ng\rangle
    =-\lim_{N\to\infty}\langle L_{\pi,N}p,g\rangle
    =-\langle L_\pi p,g\rangle.
\]
The form core extends this identity to every $g\in\Dom(V)$.
Consequently $p\in\Dom(A_\pi)$ and $A_\pi p=-L_\pi p$.
Linearity proves the claim for every $p\in\mathrm{Pol}$.

For polynomial form density, we verify the contraction hypothesis of
Lemma~\ref{lem:mixed-core} and use the marked range.
Let $g\in\mathrm{Pol}$ be real, and let $\chi$ be a bounded smooth real
function with bounded derivatives, $\chi(0)=0$ and $0\leq\chi'\leq1$.
Choose a finite list $x=(x_1,\ldots,x_d)$ of real field coordinates
used by $g$, and put $f=\chi(g)$. Equation \eqref{eq:MS.3} and the
polynomial growth of the derivatives of $g$ give
\begin{align*}
 \|Vg\|_{L^4(\pi;H)}
 &\leq\sum_{a=1}^d
    \|\Vzero x_a\|_{L^8(\mathbb G\otimes\mathbb G;H)}
    \left\|\sup_q|\partial_{x_a}g|\right\|_{L^8(\mathbb G\otimes\mathbb G)}
    +2^{-1/2}\left\|\sup_q|\partial_qg|\right\|_{L^4(\mathbb G\otimes\mathbb G)}
    <\infty.
\end{align*}
The chain rule in Appendix~\ref{sec:mixed-core} gives
\[
 f\in\Dom(A_\pi),\qquad
 A_\pi f=\chi'(g)A_\pi g-\chi''(g)\|Vg\|_H^2.
\]
Every derivative of $f$ has polynomial growth uniformly in $q$.
Lemma~\ref{lem:smooth-polynomial-approximation} gives real
$p_k\in\mathrm{Pol}$ such that, for multiindices
$\zeta\in\{0,1,\ldots\}^d$ and integers $s\geq0$,
\[
 \varepsilon_k=
 \max_{|\zeta|+s\leq2}
 \left\|\sup_{q\in\mathbb T}
  |\partial_x^\zeta\partial_q^s(p_k-f)|
                  \right\|_{L^4(\mathbb G\otimes\mathbb G)}
 \longrightarrow0.
\]
Set $r_k=p_k-f$ and $L_\pi r_k=-A_\pi r_k$, since
$p_k,f\in\Dom(A_\pi)$. The product and chain rules give
\begin{align*}
 Vr_k
 &=\sum_{a=1}^d(\partial_{x_a}r_k)\Vzero x_a
                 -(\partial_qr_k)(e_j(q))_j.\\
 \shortintertext{For the generator,}
 L_\pi r_k
 &=\sum_{a=1}^d(L_\pi x_a)\partial_{x_a}r_k
   +\sum_{a,b=1}^d
       \langle\Vzero x_a,\Vzero x_b\rangle_H\partial_{x_ax_b}r_k\\
 &\quad+\tfrac12\partial_{qq}r_k+\Bstress\,\partial_qr_k
       -2\sum_{a=1}^d\partial_{x_aq}r_k
                      \sum_j e_j(q)\Vzero_jx_a.
\end{align*}
Moreover,
\[
 \sup_q\left|\sum_j e_j(q)\Vzero_jx_a\right|
 \leq2^{-1/2}\|\Vzero x_a\|_H,\qquad
 \sup_q|\Bstress|\in L^2(\mathbb G).
\]
The field marginal of $\pi$ is $\mathbb G\otimes\mathbb G$.
H\"older's inequality and the coefficient moment bounds
\eqref{eq:MS.3}--\eqref{eq:MS.5} therefore give
\[
 \|r_k\|_{L^2(\pi)}
 +\|Vr_k\|_{L^2(\pi;H)}
 +\|L_\pi r_k\|_{L^1(\pi)}
 \leq C_x\varepsilon_k\longrightarrow0,
\]
where $C_x<\infty$ depends only on the fixed coordinate list $x$ and $Q$.

Thus
\begin{equation}\label{eq:marked-contraction-approximation}
 \|p_k-f\|_{L^2(\pi)}
 +\|V(p_k-f)\|_{L^2(\pi;H)}
 +\|A_\pi(p_k-f)\|_{L^1(\pi)}
                         \longrightarrow0.
\end{equation}

Proposition~\ref{prop:marked-range} and
\eqref{eq:marked-contraction-approximation} verify the hypotheses of
Lemma~\ref{lem:mixed-core} for $\mathrm{Pol}$.

To obtain form density, let $u\in\Dom(V)$ be real. Choose smooth
bounded truncations $\chi_m$, with
$|\chi_m'|\leq1$, equal to the identity on $[-m,m]$. Put
\[
 u_{m,n}=n(n+A_\pi)^{-1}\chi_m(u).
\]
The Markov resolvent gives
$u_{m,n}\in\Dom(A_\pi)\cap L^\infty(\pi)$. Spectral calculus
and the gradient chain rule give, with
$\|v\|_{\Dom(V)}^2=\|v\|_2^2+\|Vv\|_2^2$,
\[
 \lim_{n\to\infty}\|u_{m,n}-\chi_m(u)\|_{\Dom(V)}=0,
 \qquad
 \lim_{m\to\infty}\|\chi_m(u)-u\|_{\Dom(V)}=0.
\]
Choose $n(m)$ and, by Lemma~\ref{lem:mixed-core}, $p_m\in\mathrm{Pol}$
with
\[
 \|u_{m,n(m)}-\chi_m(u)\|_{\Dom(V)}<m^{-1},\qquad
 \|p_m-u_{m,n(m)}\|_{\Dom(V)}<m^{-1}.
\]
Then $p_m\to u$ in $\Dom(V)$.

\subsection{Compression identities and the lower-order adjoint}\label{sec:closed-compression}

For $g\in C_{\mathrm{joint}}$ and a finite cylinder-valued horizontal vector $U$,
the joint integration-by-parts identity gives
\begin{align*}
 \langle K_HVg,U\rangle_{L^2(\nu_\beta;H)}
 &=\langle Vg,I_HU\rangle_{L^2(\pi;H)}\\
 &=\sum_j\langle g,\Sigma_jIU_j-ID_jU_j\rangle_{L^2(\pi)}\\
 &=\langle Kg,\delta_{\beta,0}U\rangle_{L^2(\nu_\beta)}
       +\langle\Rop g,U\rangle_{L^2(\nu_\beta;H)}.
\end{align*}
Consequently
\begin{align*}
 Kg&\in\Dom(D_{\beta,\max}),\qquad
 D_{\beta,\max}Kg=K_HVg-\Rop g,\\
 \shortintertext{with}
 \|D_{\beta,\max}Kg\|_2&\leq\|Vg\|_2+C_*\|g\|_2.
\end{align*}
The tilted minimal--maximal equality yields
\begin{equation}\label{eq:DA.10}
D_\beta Kg=K_HVg-\Rop g,\qquad \Rop I=0.
\end{equation}
For $g_n\in C_{\mathrm{joint}}$ converging to $g\in\Dom(V)$ in the field form norm,
\begin{align*}
 \|Kg_n-Kg\|_2&\leq\|g_n-g\|_2\longrightarrow0,\\
 \shortintertext{and}
 \|D_\beta Kg_n-D_\beta Kg_m\|_2
 &\leq\|V(g_n-g_m)\|_2+C_*\|g_n-g_m\|_2\longrightarrow0.
\end{align*}
Closedness of $D_\beta$ extends \eqref{eq:DA.10} to $\Dom(V)$.

The bounded density $w_\beta$ transfers original scalar form
approximations to the tilted form norm. Thus
Lemma~\ref{lem:raw-geometric-tests} and \eqref{eq:DA.9} give, for $i,j\in J_H$,
\begin{align*}
 \Sigma_j&\in\mathrm{Pol},\\
 \rho_j^\beta=K\Sigma_j&\in C_{\mathrm{geo}}\cap\Dom(D_\beta),\\
 \shortintertext{and hence}
 I\rho_j^\beta&\in\Dom(V).\\
 \shortintertext{The pullback and centered-score derivatives are}
 VI\rho_j^\beta&=I_HD_\beta\rho_j^\beta,\\
 \eta_j&\in\Dom(V),\\
 \shortintertext{and}
 V_i\eta_j&=\Vzero_i\Sigma_j-ID_i\rho_j^\beta.
\end{align*}
The bounded base gradients and Gaussian moments give, for
$1\leq p<\infty$,
\[
 \|\eta_j\|_{L^p(\pi)}+\|V\eta_j\|_{L^p(\pi;H)}<\infty.
\]
For $g\in C_{\mathrm{joint}}$, the moment bounds and the closed product rule give
\[
 \eta_jg\in\Dom(V),\qquad
 V_i(\eta_jg)
       =g(\Vzero_i\Sigma_j-ID_i\rho_j^\beta)+\eta_jV_i g.
\]
In addition,
\[
 \|\eta_i\eta_jg\|_2
       \leq\|\eta_i\|_6\|\eta_j\|_6\|g\|_6<\infty.
\]
For $F\in\Ccyl$, finite joint integration by parts and closure give
\[
 IF\in\Dom(V_i^*),\qquad V_i^*IF=\Sigma_iIF-ID_iF.
\]
Consequently, for $g\in C_{\mathrm{joint}}$,
\[
 \langle K(\eta_jg),\rho_i^\beta F-D_iF\rangle_{L^2(\nu_\beta)}
 =\langle V_i(\eta_jg)-\eta_i\eta_jg,IF\rangle_{L^2(\pi)}.
\]
For $i,j\in{J_H}$ and $a\in\{0,c,s\}$ set
\[
 H^{\rm base}_{ij}=D_i\rho_j^\beta-\rho^\beta_{\nabla_{e_i}e_j},
 \qquad \eta_a=\Sigma_{z_a}-I\rho_{z_a}^\beta.
\]
The low scores and \eqref{eq:KD.1} give
\begin{equation}\label{eq:GH.low-scores}
\begin{aligned}
 \eta_0&=0,\qquad
 \eta_c=Q\operatorname{Im}(Y_1-I(m_1)),\qquad
 \eta_s=Q\operatorname{Re}(Y_1-I(m_1)),\\
 K(\eta_c^2+\eta_s^2)&=\frac{Q^2}{2}.
\end{aligned}
\end{equation}
For $g\in C_{\mathrm{joint}}$, the preceding identity and the product rule give
\begin{align*}
 D_i(\Rop g)_j
 &=K\big[(\Vzero_i\Sigma_j-ID_i\rho_j^\beta
                                    -\eta_i\eta_j)g\big]
                                      +K(\eta_jV_i g).\\
 \shortintertext{The covariant derivative is therefore obtained from}
 (\nabla_i\Rop g)_j
 &=D_i(\Rop g)_j-\sum_{k\in{J_H}}c_{ij}^kK(\eta_kg),\\
 \shortintertext{using the connection decomposition}
 \Sigma_{\nabla_{e_i}e_j}-I\rho^\beta_{\nabla_{e_i}e_j}
 &=\sum_{k\in{J_H}}c_{ij}^k\eta_k
                           +\sum_{a=0,c,s}c_{ij}^a\eta_a.
\end{align*}
Consequently
\begin{equation}\label{eq:GH.6}
\begin{aligned}
 (\nabla_i\Rop g)_j
 &=K\big[(H^{\rm field}_{ij}-IH^{\rm base}_{ij}
                                     -\eta_i\eta_j)g\big]\\
 &\quad+K(\eta_jV_i g)+\sum_{a=0,c,s}c_{ij}^aK(\eta_ag).
\end{aligned}
\end{equation}
Differentiating $K\Sigma_j=\rho_j^\beta$ yields
\begin{align*}
 D_i\rho_j^\beta
 &=K(\Vzero_i\Sigma_j)-K(\eta_i\eta_j),\\
 \shortintertext{so that}
 KH^{\rm field}_{ij}-H^{\rm base}_{ij}&=K(\eta_i\eta_j).
\end{align*}
The supremum below is over finite complex matrices $b=(b_{ij})$.
Conditional Cauchy--Schwarz and \eqref{eq:conditional-matrix-bound} imply
\begin{equation}\label{eq:GH.centered-bound}
\begin{aligned}
 &\sum_{i,j\in{J_H}}
   \left|K\big[(H^{\rm field}_{ij}-IH^{\rm base}_{ij}
                                      -\eta_i\eta_j)g\big]\right|^2\\
 &=\sup_{\sum_{i,j}|b_{ij}|^2\leq1}
   \left|K\left[g\sum_{i,j}\overline{b_{ij}}
      \big(H^{\rm field}_{ij}-IKH^{\rm field}_{ij}
                   -\eta_i\eta_j+IK(\eta_i\eta_j)\big)\right]\right|^2\\
 &\leq C K|g|^2.
\end{aligned}
\end{equation}
The conditional covariance bound \eqref{eq:conditional-score-bound} gives
\begin{equation}\label{eq:GH.response-bound}
 \int\sum_{i,j\in{J_H}}|K(\eta_jV_i g)|^2\,d\nu_\beta
       \leq C_*^2\int\sum_{i\in{J_H}}K|V_i g|^2\,d\nu_\beta
       =C_*^2\|Vg\|_{L^2(\pi;H)}^2.
\end{equation}
By metricity and \eqref{eq:connection-column-bounds}, for $a=c,s$,
\[
 \sum_{i,j\in{J_H}}|c_{ij}^a|^2
 =\sum_{i\in{J_H}}\|P_H\Xi_i z_a\|_{\widehat H}^2
 \leq\sum_{i\in{J_H}}\|\Xi_i z_a\|_{\widehat H}^2
 <\infty.
\]
Thus \eqref{eq:GH.low-scores} and conditional Cauchy--Schwarz give
\begin{align}
 |K(\eta_cg)|^2+|K(\eta_sg)|^2
 &\leq K(\eta_c^2+\eta_s^2)K|g|^2
       =\frac{Q^2}{2}K|g|^2.\notag\\
 \shortintertext{Summing the connection coefficients gives}
 \sum_{i,j\in J_H}
       \left|\sum_{a=0,c,s}c_{ij}^aK(\eta_ag)\right|^2
 &\leq\left(\sum_{i,j\in J_H}\sum_{a=c,s}|c_{ij}^a|^2\right)
                     \sum_{a=c,s}|K(\eta_ag)|^2\notag\\
 &\leq C K|g|^2.\notag\\
 \shortintertext{After integration,}
 \int\sum_{i,j\in J_H}
       \left|\sum_{a=0,c,s}c_{ij}^aK(\eta_ag)\right|^2\,d\nu_\beta
 &\leq C\|g\|_{L^2(\pi)}^2.\label{eq:GH.low-bound}
\end{align}
The low output components satisfy
\begin{align*}
 (\nabla_i\Rop g)_a
 &=-\langle\Rop g,\Xi_i z_a\rangle_{\widehat H}
                                      =-(T_a^*\Rop g)_i.\\
 \shortintertext{The boundedness of $T_a$ therefore gives}
 \sum_{a=0,c,s}\sum_{i\in{J_H}}
              \|(\nabla_i\Rop g)_a\|_{L^2(\nu_\beta)}^2
 &\leq\sum_{a=0,c,s}\|T_a\|^2\|\Rop g\|_{L^2(\nu_\beta;H)}^2
       \leq C\|g\|_{L^2(\pi)}^2.
\end{align*}
Combining these bounds with \eqref{eq:GH.6} gives
\begin{equation}\label{eq:GH.7}
 \|\nabla_H\Rop g\|_{L^2(\nu_\beta;\HS(H,\widehat H))}
       \leq C\big(\|Vg\|_{L^2(\pi;H)}+\|g\|_{L^2(\pi)}\big),
       \qquad g\in C_{\mathrm{joint}}.
\end{equation}
By \eqref{eq:B.3}, \eqref{eq:GH.7} and the domain equality in
Section~\ref{sec:tilted-domains}, \eqref{eq:DA.3} gives, for $g\in C_{\mathrm{joint}}$,
\[
\begin{aligned}
 \Rop g&\in\Dom(\nabla_{H,\beta})\cap\Dom(D_\beta^*),\\
 \|D_\beta^*\Rop g\|_2
 &\leq C_\beta\big(\|\nabla_H\Rop g\|_2+\|\Rop g\|_2\big)\\
 &\leq C_\beta\big(\|Vg\|_2+\|g\|_2\big).
\end{aligned}
\]
For $g\in\Dom(V)$, choose $g_n\in C_{\mathrm{joint}}$ with
$\|g_n-g\|_2+\|Vg_n-Vg\|_2\to0$. Then
\begin{align*}
 \|\Rop g_n-\Rop g\|_2&\leq C_*\|g_n-g\|_2\longrightarrow0,\\
 \shortintertext{and the covariant derivatives satisfy}
 \|\nabla_{H,\beta}(\Rop g_n-\Rop g_m)\|_2
  &\leq C\bigl(\|V(g_n-g_m)\|_2+\|g_n-g_m\|_2\bigr)
                                                   \longrightarrow0,\\
 \shortintertext{while}
 \|D_\beta^*\Rop g_n-D_\beta^*\Rop g_m\|_2
 &\leq C_\beta\big(\|V(g_n-g_m)\|_2+\|g_n-g_m\|_2\big)
                                                   \longrightarrow0.
\end{align*}
Closedness of $\nabla_{H,\beta}$ and $D_\beta^*$ therefore gives
\begin{align}
 \Rop&:\Dom(V)\longrightarrow
          \Dom(\nabla_{H,\beta})\cap\Dom(D_\beta^*),\notag\\
 \shortintertext{with}
 \|\Rop g\|_2+\|\nabla_{H,\beta}\Rop g\|_2
             +\|D_\beta^*\Rop g\|_2
   &\leq C_\beta\bigl(\|Vg\|_2+\|g\|_2\bigr).
                                                    \label{eq:DA.11}
\end{align}

For $p\in\mathrm{Pol}$ and $F\in\Dom(D_\beta)$,
\begin{align*}
 \langle K_HVp,D_\beta F\rangle_{L^2(\nu_\beta;H)}
 &=\langle Vp,I_HD_\beta F\rangle_{L^2(\pi;H)}\\
 &=\langle Vp,VIF\rangle_{L^2(\pi;H)}
   =\langle A_\pi p,IF\rangle_{L^2(\pi)}\\
 &=\langle KA_\pi p,F\rangle_{L^2(\nu_\beta)}.
\end{align*}
Hence
\[
 K_HVp\in\Dom(D_\beta^*),\qquad
 D_\beta^*K_HVp=KA_\pi p.
\]
Combining this with \eqref{eq:DA.10} and \eqref{eq:DA.11} gives
\begin{equation}\label{eq:DA.12}
KA_{\pi}p-D_\beta^*\Rop p=A_\beta Kp.
\end{equation}

\subsection{Geometric tests in the tilted generator domain}

Write $A=-L=D^*D$. For the capacity cutoffs $\chi_R(c)$ above,
Lemma~\ref{lem:raw-geometric-tests} gives
$\chi_Rw_\beta\in C_{\mathrm{geo}}\subset\Dom(D)$ and
\[
 D(\chi_Rw_\beta)=\beta\chi_Rw_\beta Dc+w_\beta D\chi_R.
\]
For $\beta<0$, the bound $\|Dc\|^2\leq72c$ gives
\[
 \|\beta w_\beta Dc\|_{L^\infty(\nu;H)}^2
 \leq72\beta^2z_\beta^{-2}\sup_{r\geq0}r e^{2\beta r}<\infty.
\]
Equation \eqref{eq:DA.4} gives
\begin{align*}
 \|\chi_Rw_\beta-w_\beta\|_{L^2(\nu)}^2
 &\leq\|w_\beta\|_\infty^2\nu(c\geq R)\longrightarrow0,\\
 \shortintertext{and}
 \|D(\chi_Rw_\beta)-\beta w_\beta Dc\|_{L^2(\nu;H)}
 &\leq\|\mathbf1_{\{c\geq R\}}\beta w_\beta Dc\|_{L^2(\nu;H)}
       +C\|w_\beta\|_\infty R^{-1/2}\longrightarrow0.
\end{align*}
Thus $w_\beta\in\Dom(D)$ with $Dw_\beta=\beta w_\beta Dc$.
For $\beta=0$, this follows directly from $w_0=1$.

Let $F\in C_{\mathrm{geo}}$ be real. Lemma~\ref{lem:raw-geometric-tests} gives
$F\in\Dom(A)$ and bounded $F$, $AF$, and $\|DF\|_H$.
Choose real $F_k\in\Ccyl$ converging to $F$ in the original form norm.
Then
\begin{align*}
 &\|F_k-F\|_{L^2(\nu_\beta)}^2
       +\|DF_k-DF\|_{L^2(\nu_\beta;H)}^2\\
 &\qquad\leq\|w_\beta\|_\infty
       \left(\|F_k-F\|_{L^2(\nu)}^2
             +\|DF_k-DF\|_{L^2(\nu;H)}^2\right)\longrightarrow0.
\end{align*}
Hence $F\in\Dom(D_\beta)$ and $D_\beta F=DF$. For real $G\in\Ccyl$,
$w_\beta G\in\Dom(D)$, and hence
\begin{align*}
 \int\langle DF,DG\rangle_H\,d\nu_\beta
 &=\langle DF,D(w_\beta G)\rangle_{L^2(\nu;H)}
       -\int G\langle DF,Dw_\beta\rangle_H\,d\nu\\
 &=\int\bigl(AF-\beta\langle Dc,DF\rangle_H\bigr)G\,d\nu_\beta.
\end{align*}
For $\beta<0$,
\[
 \|\beta\langle Dc,DF\rangle_H\|_{L^2(\nu_\beta)}^2
 \leq72\beta^2\|DF\|_\infty^2z_\beta^{-1}
                   \int c e^{\beta c}\,d\nu<\infty;
\]
for $\beta=0$ the term vanishes. The cylinder form core therefore gives
\begin{equation}\label{eq:geometric-tilted-domain}
 F\in\Dom(A_\beta),\qquad
 A_\beta F=AF-\beta\langle Dc,DF\rangle_H,\qquad
 D_\beta F\in L^4(\nu_\beta;H).
\end{equation}
Complexification gives the same domain identity for complex tests.

\subsection{Final proof of Theorem~\ref{thm:weak-uniqueness}}
\label{sec:final-weak-proof}

For $0<\kappa\leq4$, the canonical stationary construction in
\cite[Sections 3--4]{Diffusion} provides a continuous adapted state,
the response normalization used here, and the Brownian family in
\eqref{eq:self-consistency}--\eqref{eq:raw-sde}.
It remains to identify every weak solution with marginals \eqref{eq:marginal-class}.
We first establish the range density that will determine its conditional laws.

For the generator $A_\beta$ and the test class $C_{\mathrm{geo}}$ of
Lemma~\ref{lem:raw-geometric-tests}, we have
\begin{equation}\label{eq:tilted-range}
 \overline{(\lambda+A_\beta)C_{\mathrm{geo}}}^{\,L^1(\nu_\beta)}
                              =L^1(\nu_\beta),\qquad\lambda>0.
\end{equation}
Indeed, let a bounded real function $a$ annihilate this range and put
$u=Ia$. Since $Kp\in C_{\mathrm{geo}}$ for $p\in\mathrm{Pol}$,
\eqref{eq:DA.12} gives
\[
 \langle u,(\lambda+A_\pi)p\rangle_{\pi}
                     =\langle a,D_\beta^*\Rop p\rangle_{\nu_\beta}.
\]
Define the real functional
\[
 \ell(g)=\langle a,D_\beta^*\Rop g\rangle_{\nu_\beta},
 \qquad g\in\Dom(V).
\]
Equation \eqref{eq:DA.11} gives its form-norm bound:
\[
 |\ell(g)|\leq C_\beta\|a\|_2(\|Vg\|_2+\|g\|_2).
\]
Proposition~\ref{prop:marked-range} and Lemma~\ref{lem:mixed-core}
give $u\in\Dom(V)$. Polynomial form density then gives
\[
 \lambda\langle u,g\rangle+\langle Vu,Vg\rangle=\ell(g),
                              \qquad g\in\Dom(V)\text{ real}.
\]
Taking $g=u$ yields
\[
 \lambda\|u\|_2^2+\|Vu\|_2^2
   =\langle a,D_\beta^*\Rop u\rangle=0,
\]
because $\Rop I=0$ as a bounded operator on the whole $L^2$ space.
Hence $a=0$. The real and imaginary parts and Hahn--Banach prove
\eqref{eq:tilted-range}.

We transfer this range identity to the original welding probability.
Write $A=-L=D^*D$ and retain $A_\beta=D_\beta^*D_\beta$.
The density in the direction needed below is
\[
 w=\frac{d\nu}{d\nu_\beta}=z_\beta e^{-\beta c},\qquad
 \vartheta=\log w,\qquad D\vartheta=-\beta Dc.
\]
The score change and the generator product rule give
\begin{equation}\label{eq:generator-weight-change}
 AF=A_\beta F-\langle D\vartheta,DF\rangle_H,
                         \qquad F\in C_{\mathrm{geo}}.
\end{equation}
For real $g\in C_{\mathrm{geo}}$ and a contraction $\chi$ as in
\eqref{eq:mixed-contraction-approximation}, we have
$\chi(g)\in C_{\mathrm{geo}}$.
Equation \eqref{eq:geometric-tilted-domain} puts both functions in
$\Dom(A_\beta)$ with $L^4$ gradients. Thus the constant sequence
$p_k=\chi(g)$ verifies the contraction hypothesis of
Lemma~\ref{lem:mixed-core}.

Together with \eqref{eq:tilted-range}, this verifies the hypotheses of
Lemma~\ref{lem:mixed-core} for $(A_\beta,C_{\mathrm{geo}})$.

The transfer will also use the capacity estimates
\begin{equation}\label{eq:capacity-transfer-bounds}
 c<\infty\quad\nu\text{-almost everywhere},\qquad
 \|Dc\|_H^2\leq72c,
 \qquad |Lc|\leq C_\kappa(1+c).
\end{equation}
Their raw It\^o and canonical-domain interpretations were proved in
Lemma~\ref{lem:raw-geometric-tests}.

Fix $\lambda>0$ and suppose that $a\in L^\infty(\nu)$ is real and
\begin{equation}\label{eq:original-annihilator}
                 \int a(\lambda+A)F\,d\nu=0,
                                      \qquad F\in C_{\mathrm{geo}}.
\end{equation}
Choose a smooth function $\chi:[0,\infty)\to[0,1]$ equal to one on
$[0,1]$ and zero on $[2,\infty)$, and set
$\chi_N=\chi(c/N)$, $b_N=\chi_Nwa$. Inserting $\chi_NF$ into
\eqref{eq:original-annihilator} gives
\begin{align}
 \langle b_N,(\lambda+A_\beta)F\rangle_{\nu_\beta}
   &=\langle f_N,F\rangle_{\nu_\beta}
       +\langle U_N,DF\rangle_{\nu_\beta},
                                      \label{eq:localized-annihilator}\\
 \shortintertext{where}
 f_N&=-wa A\chi_N,
 &U_N&=wa(2D\chi_N+\chi_ND\vartheta).\notag
\end{align}
Indeed, the two identities used in this expansion are
\begin{align*}
 A_\beta(\chi_NF)&=\chi_NA_\beta F+FA_\beta\chi_N
                         -2\langle D\chi_N,DF\rangle,\\
 \shortintertext{and}
 A_\beta\chi_N-\langle D\vartheta,D\chi_N\rangle&=A\chi_N.
\end{align*}

Equation \eqref{eq:capacity-transfer-bounds} gives
\[
 \|b_N\|_\infty+\|f_N\|_\infty+\|U_N\|_2<\infty.
\]
Hence the right-hand side is continuous in the form norm:
\[
 |\langle f_N,F\rangle+\langle U_N,DF\rangle|
 \leq(\|f_N\|_2+\|U_N\|_2)\|F\|_{\Dom(D_\beta)}.
\]
Lemma~\ref{lem:mixed-core} yields
\begin{equation}\label{eq:localized-energy}
                              b_N\in\Dom(D_\beta).
\end{equation}
Choose a smooth capacity cutoff $\psi_N(c)$ equal to one on $c\leq N$
and zero on $c\geq2N$. Then
\begin{align*}
 \psi_Na&=(\psi_N/w)b_{4N},\\
 D_\beta(\psi_Na)
   &=\frac{\psi_N}{w}D_\beta b_{4N}
      +\frac{b_{4N}}w\bigl(D\psi_N-\psi_ND\vartheta\bigr).
\end{align*}
The multiplier $\psi_N/w$ and its gradient are bounded. Local
equivalence of $\nu$ and $\nu_\beta$, followed by
Proposition~\ref{prop:horizontal-domain}, gives
\[
                              \psi_Na\in\Dom(D).
\]
Strong locality gives the compatibility relation
\begin{align*}
 \mathbf1_{\{c\leq N\}}D(\psi_{2M}a)
   &=\mathbf1_{\{c\leq N\}}D(\psi_{2N}a),\qquad M\geq N,\\
 \shortintertext{which defines the local gradient by}
 \mathbf1_{\{c\leq N\}}Da
   &:=\mathbf1_{\{c\leq N\}}D(\psi_{2N}a),\qquad N\geq1.
\end{align*}
The definition is independent of the chosen cutoffs.

The local gradient gives a capacity cutoff estimate.
Fix $N$ and $p\in\Ccyl$, and put $F=\chi_Np$. The product rules give
\begin{align*}
 DF&=\chi_NDp+pD\chi_N,\\
 AF&=\chi_NAp+pA\chi_N-2\langle D\chi_N,Dp\rangle_H.\\
 \shortintertext{In particular, all three quantities are supported in the same capacity set:}
 \mathbf1_{\{c>2N\}}(F,DF,AF)&=0
                                      \quad\nu\text{-almost everywhere}.
\end{align*}
Choose the preceding capacity cutoff $\psi_{2N}$ equal to one on
$c\leq3N$ and zero on $c\geq4N$. The local gradient identity gives
\[
 \psi_{2N}a\in\Dom(D),\qquad
 \mathbf1_{\{c\leq2N\}}D(\psi_{2N}a)
       =\mathbf1_{\{c\leq2N\}}Da.
\]
Consequently,
\begin{align*}
 \int aAF\,d\nu
 &=\int(\psi_{2N}a)AF\,d\nu\\
 &=\int\langle D(\psi_{2N}a),DF\rangle\,d\nu
 =\int\langle Da,DF\rangle\,d\nu.
\end{align*}
Hence \eqref{eq:original-annihilator} becomes
\begin{equation}\label{eq:local-weak-energy}
               \lambda\int aF\,d\nu+
                     \int\langle Da,DF\rangle\,d\nu=0.
\end{equation}
Choose real $p_k\in\Ccyl$ with $p_k\to\chi_Na$ in $\Dom(D)$, and set
$F_k=\chi_Np_k$. Then
\[
 F_k\in C_{\mathrm{geo}},\qquad
 \mathbf1_{\{c>2N\}}(F_k,DF_k,AF_k)=0\quad\nu\text{-almost everywhere},
\]
\noindent and
\[
 \|F_k-\chi_N^2a\|_{\Dom(D)}
    \leq C_N\|p_k-\chi_Na\|_{\Dom(D)}\longrightarrow0.
\]
Since $\int_{\{c\leq2N\}}\|Da\|_H^2d\nu<\infty$, this passes
\eqref{eq:local-weak-energy} to $F=\chi_N^2a$ and gives
\begin{align}
 \lambda\|\chi_Na\|_2^2+
       \int\chi_N^2\|Da\|_H^2\,d\nu
  &=-2\int\chi_Na\langle Da,D\chi_N\rangle_H\,d\nu.\notag\\
 \shortintertext{Young's inequality and the cutoff bound give}
 \lambda\|\chi_Na\|_2^2+
       \frac12\int\chi_N^2\|Da\|_H^2\,d\nu
  &\leq2\|a\|_\infty^2\int\|D\chi_N\|_H^2\,d\nu
   \leq\frac{288\|\chi'\|_\infty^2\|a\|_\infty^2}{N}.
                                     \label{eq:capacity-caccioppoli}
\end{align}
The last inequality uses
$\|D\chi_N\|^2\leq72N^{-2}\|\chi'\|_\infty^2
c\mathbf1_{\{N\leq c\leq2N\}}$. Since $\chi_N\to1$ almost
everywhere,
\[
 \lambda\|a\|_2^2
 =\lim_{N\to\infty}\lambda\|\chi_Na\|_2^2
 \leq\lim_{N\to\infty}
       \frac{288\|\chi'\|_\infty^2\|a\|_\infty^2}{N}=0.
\]
Real and imaginary parts and Hahn--Banach now prove
\begin{equation}\label{eq:original-range}
       \overline{(\lambda-L)C_{\mathrm{geo}}}^{\,L^1(\nu)}=L^1(\nu),
                                           \qquad\lambda>0.
\end{equation}

Now let $(X,W)$ be any weak solution satisfying
\eqref{eq:marginal-class}, on its given filtered probability space.

To identify the conditional state law, fix $s\geq0$ and a bounded real
$\mathcal F_s$-measurable random variable $\zeta$. Define signed
densities $g_t$ by
\[
      g_t\,d\nu=\mathbb E[\zeta\mathbf1_{\{X_{s+t}\in\cdot\}}],
                     \qquad\|g_t\|_\infty\leq\|\zeta\|_\infty.
\]
Lemma~\ref{lem:raw-geometric-tests} gives, for $F\in C_{\mathrm{geo}}$,
\[
 \langle g_t,F\rangle-\langle g_0,F\rangle
                     =\int_0^t\langle g_r,LF\rangle\,dr.
\]
The canonical evolution $T_tg_0$ satisfies the same identity because
$C_{\mathrm{geo}}\subset\Dom(L)$. Put $u_t=g_t-T_tg_0$. For $\lambda>0$,
define $a_\lambda\in L^\infty(\nu)$ through its pairings by
\begin{align*}
 \langle a_\lambda,\phi\rangle
    &=\int_0^\infty e^{-\lambda t}\langle u_t,\phi\rangle\,dt,
                                    \qquad \phi\in L^1(\nu),\\
 \shortintertext{with}
 \|a_\lambda\|_\infty&\leq\frac{2\|\zeta\|_\infty}{\lambda}.
\end{align*}
For $F\in C_{\mathrm{geo}}$, the weak evolution and integration by parts give
\begin{align*}
 \langle a_\lambda,(\lambda-L)F\rangle
   &=\int_0^\infty e^{-\lambda t}
       \left(\lambda\langle u_t,F\rangle
                    -\frac d{dt}\langle u_t,F\rangle\right)\,dt\\
   &=\langle u_0,F\rangle=0.
\end{align*}
Equation~\eqref{eq:original-range} implies $a_\lambda=0$.
Scalar Laplace uniqueness gives $\langle u_t,F\rangle=0$ for almost
every $t$. For $F\in\Ccyl$, path continuity and semigroup continuity give
\[
 \lim_{r\to t}\mathbb E[\zeta F(X_{s+r})]
       =\mathbb E[\zeta F(X_{s+t})],\qquad
 \lim_{r\to t}\|T_rg_0-T_tg_0\|_2=0.
\]
Hence, for every $t\geq0$,
\begin{align*}
 \mathbb E[\zeta F(X_{s+t})]
 &=\langle T_tg_0,F\rangle
  =\langle g_0,T_tF\rangle
  =\mathbb E[\zeta T_tF(X_s)],\qquad F\in\Ccyl.
\end{align*}
A countable measure-determining subalgebra and the monotone class
theorem extend the identity to bounded measurable $F$. Varying $\zeta$ gives
\[
                    \mathbb E[F(X_{s+t})\mid\mathcal F_s]
                                      =T_tF(X_s)
                    \quad\mathbb P\text{-almost surely}.
\]
This proves \eqref{eq:conditional-semigroup} and identifies the state
path law with the canonical law.

It remains to identify the noise jointly with the state. Choose
countably many distinct dense targets $z_i\in S^1$, and let $q_i(t)$
be the continuous angular lift of $X_t^{-1}(z_i)$, with initial value
in $[0,2\pi)$. By Lemma~\ref{lem:label-extension}, the state path
determines the martingales
\[
 M_i(t)=q_i(t)-q_i(0)-\int_0^tB(X_r,q_i(r))\,dr
       =\sqrt2\sum_j\int_0^te_j(q_i(r))\,dW_r^j.
\]
Let $P_N^{\mathrm{rec}}(t)$ be the orthogonal projection in $H$ onto
$\operatorname{span}\{k(q_i(t)):1\leq i\leq N\}$, with $k$ from
\eqref{eq:evaluation-row}. As in the proof of \cite[Lemma~4.4]{Diffusion}, Fourier uniqueness
at distinct points $q_i$ modulo $2\pi$ gives
\[
 \sum_{i=1}^Na_i k(q_i)=0
 \quad\Longrightarrow\quad
 \sum_{i=1}^Na_ie^{-inq_i}=0\qquad(|n|\geq2).
\]
The signed atomic measure therefore has only three Fourier modes:
\[
 \sum_{i=1}^Na_i\delta_{q_i}
 =\sum_{n=-1}^1\left(\sum_{i=1}^Na_ie^{-inq_i}\right)
                     e^{in\theta}\frac{d\theta}{2\pi}.
\]
It is also absolutely continuous, so all $a_i$ vanish.
Thus the Gram matrices are invertible.
For each $j\in{J_H}$, define
\begin{equation}\label{eq:noise-reconstruction}
 W^{j,N,\mathrm{rec}}_t
 =\int_0^t\langle P_N^{\mathrm{rec}}(r)e_j,dW_r\rangle_H.
\end{equation}
To express this integral through the state path, set
\[
 G_N^{\mathrm{rec}}(r)_{ik}=\langle k(q_i(r)),k(q_k(r))\rangle_H,\qquad
 v_{j,N}(r)=(e_j(q_i(r)))_{i=1}^N.
\]
With $a_{j,N}^{\mathrm{rec}}(r)=(G_N^{\mathrm{rec}}(r))^{-1}v_{j,N}(r)$,
Gram inversion gives
\[
 W^{j,N,\mathrm{rec}}_t=\frac1{\sqrt2}\sum_{i=1}^N
              \int_0^t[a_{j,N}^{\mathrm{rec}}(r)]_i\,dM_i(r).
\]
For fixed $j,N,T$, the Gram coefficient $a_{j,N}^{\mathrm{rec}}(t)$ is
continuous in the inverse-coordinate paths. Set
\[
 \sigma_K=\inf\{t\geq0:
     \|a_{j,N}^{\mathrm{rec}}(t)\|_{\R^N}\geq K\}\wedge T.
\]
Before $\sigma_K$ the coefficient is bounded, and
\[
 d\langle M_i,M_k\rangle_t=2G_N^{\mathrm{rec}}(t)_{ik}\,dt,
 \qquad \operatorname{Tr}(2G_N^{\mathrm{rec}}(t))=N.
\]
The measurable predictable approximation in
Section~\ref{sec:field-noise-recovery}, applied to the state path and
these inverse-coordinate martingales, therefore determines the joint
law of $(X,W^{j,N,\mathrm{rec}}_{\cdot\wedge\sigma_K})$ from the state path law.
Since $\sigma_K\uparrow T$ almost surely, the same holds without stopping.

For every $t$, the points $X_t^{-1}(z_i)$ are dense in $S^1$.
Continuity of the evaluation map and the embedding
$H\subset C(S^1)$ imply $P_N^{\mathrm{rec}}(t)e_j\to e_j$.
Doob's inequality and
dominated convergence give
\[
 \mathbb E\sup_{t\leq T}|W_t^{j,N,\mathrm{rec}}-W_t^j|^2
 \leq4\mathbb E\int_0^T
       \|(\operatorname{Id}_H-P_N^{\mathrm{rec}}(r))e_j\|_H^2\,dr
 \longrightarrow0.
\]
Thus each Brownian coordinate, jointly with the state path and every
finite collection of other coordinates, is determined by the state
law. This proves \eqref{eq:joint-unique} and completes the proof.

\section{Autonomous approximation}\label{sec:autonomous}

Theorem~\ref{thm:weak-uniqueness} fixes the target joint law.
We construct autonomous strong approximations and prove convergence to that law.

\subsection{The regularized symmetric forms}

Throughout this section $\gamma$ is one fixed sequence from
\eqref{eq:weights}; constants may depend on that sequence.
The global score versions \eqref{eq:score-versions} satisfy
$|\rho_{n,\tau}(h)|\leq C_\kappa/\alpha_n$.
Finite support or exponential decay therefore gives
\begin{equation}\label{eq:fixed-weight-bounds}
 \sup_{h\in \mathcal X}\sum_j\gamma_j^2|\rho_j(h)|^2<\infty,\qquad
 \sum_j\gamma_j^2\|e_j\|_{C^m(\T)}^2<\infty\quad(m\geq0).
\end{equation}
In particular $B_\gamma(h,\cdot)$ has bounded derivatives of
every order, uniformly in $h$. We also use the unweighted choice
$\gamma_j=1$, $B_\gamma=B$ in the limiting equation, and write
\[
 k^\gamma(\theta)=\sum_{j\in J_H}\gamma_je_j(\theta)e_j\in H,
 \qquad k^1=k
\]
for the evaluation rows, with $k$ from \eqref{eq:evaluation-row}.
Since $0\leq\gamma_j\leq1$, \eqref{eq:kernel-modulus} bounds their
common modulus of continuity.

On $\Ccyl$ define
\begin{equation}\label{eq:weighted-form}
 \cE_\gamma(F,G)=\int_\mathcal X\sum_j\gamma_j^2D_jF D_jG\,d\nu .
\end{equation}
For its closure, use the notation $\Dom(\cE_\gamma)$ for the
form domain, $L_\gamma$ for the nonpositive generator,
$T_t^\gamma$ for the semigroup and $\Cap_\gamma$ for the $1$-capacity.
The following lemma proves that this closure exists.

\begin{lemma}\label{lem:inverse-target-extension}
For either a fixed multiplier sequence $\gamma$ from
\eqref{eq:weights}, or the unweighted sequence
$\gamma_j\equiv1$, suppose the inverse-coordinate equation
\begin{equation}\label{eq:inverse-target-dense}
 dq_t(v)=\sqrt2\sum_j\gamma_je_j(q_t(v))\,dW_t^j
                    +B_\gamma(X_t,q_t(v))\,dt
\end{equation}
holds for a countable dense set of deterministic targets.
Then it holds at every $\mathcal F_0$-measurable target.
Consequently $\Psi_t=X_t^{-1}\circ X_0$ satisfies the
corresponding relative-flow equation.
\end{lemma}
\begin{proof}
The inverse-coordinate lifts are continuous jointly in the target
and time. For targets $v_m\to v$ in the dense set, Doob's
inequality and \eqref{eq:kernel-modulus} give
\[
\begin{aligned}
 &\mathbb E\sup_{t\le T}\left|\sqrt2\sum_j\int_0^t\gamma_j
       [e_j(q_s(v_m))-e_j(q_s(v))]\,dW_s^j\right|^2\\
 &\qquad\le8\mathbb E\int_0^T
       \|k^\gamma(q_s(v_m))-k^\gamma(q_s(v))\|_H^2\,ds
       \longrightarrow0.
\end{aligned}
\]
The drift term passes to the limit by \eqref{eq:drift-modulus}
when $\gamma_j\equiv1$, and by \eqref{eq:fixed-weight-bounds}
for the fixed multiplier sequences in \eqref{eq:weights}.
Thus \eqref{eq:inverse-target-dense} extends to each deterministic
target. For an $\mathcal F_0$-measurable circle-valued target $V$,
approximate $V$ by finite-valued $\mathcal F_0$-measurable targets
in the dense set. Choose angular lifts compatibly in local charts;
adding an integer multiple of $2\pi$ does not change the coefficients
or drift. The event indicators are $\mathcal F_0$-measurable and hence
predictable, so the same estimate passes to the limit in the It\^o
integrals. Taking $V=X_0(e^{iu})$, with the lift chosen so that
$q_0(V)=u$, proves the relative-flow equation. A countable
intersection gives a common event for the source labels and
integer time horizons.
\end{proof}

\begin{lemma}\label{lem:weighted-weak-law}
The form $\cE_\gamma$ is closable. Its closure is conservative,
strongly local and quasi-regular. For $u\in\Dom(D)$ and every measurable
$E\subset\mathcal X$,
\begin{equation}\label{eq:form-capacity-comparison}
 \Dom(D)\subset\Dom(\cE_\gamma),\qquad
 \cE_\gamma(u,u)\leq\cE(u,u),\qquad
 \Cap_\gamma(E)\leq\Cap(E).
\end{equation}
Its stationary diffusion can be equipped with a Brownian family
so as to solve \eqref{eq:autonomous-state}--\eqref{eq:autonomous-sde}.
\end{lemma}
\begin{proof}
Suppose $F_m\to0$ in $L^2$ and $(\gamma_jD_jF_m)_j\to U$ in
$L^2(\nu;\ell^2)$. Equation~\eqref{eq:ibp} gives, for $G\in\Ccyl$,
\[
 \langle U_j,G\rangle
 =\lim_m\gamma_j\langle F_m,\rho_jG-D_jG\rangle=0.
\]
Density of $\Ccyl$ yields $U=0$, proving closability.
For a smooth contraction $\chi$ and a cylinder $F$,
\[
 \cE_\gamma(\chi(F),\chi(F))
  =\int|\chi'(F)|^2\sum_j\gamma_j^2|D_jF|^2\,d\nu
                                      \leq\cE_\gamma(F,F).
\]
Closure gives the Dirichlet property; $D_jF=0$ wherever $F$
is locally constant gives strong locality.
Since $0\leq\gamma_j\leq1$, the two form assertions in
\eqref{eq:form-capacity-comparison} follow. For open $E$, the
variational definition of capacity gives the last inequality;
outer approximation extends it to all sets. Thus the compact nest
$K_m$ of \eqref{eq:compact-nest} satisfies
\[
 \Cap_\gamma(\mathcal X\setminus K_m)
 \leq\Cap(\mathcal X\setminus K_m)\longrightarrow0.
\]
Together with the countable continuous separating cylinder core,
this proves quasi-regularity. Also
\[
 1\in\Dom(\cE_\gamma),\qquad \cE_\gamma(1,1)=0,
 \qquad T_t^\gamma1=1.
\]
The local form--process correspondence \cite[Section~2]{AMR}
gives a continuous conservative diffusion whose stationary law
is symmetric and has marginal $\nu$.

To identify the equation, set
\[
 k_\gamma(\theta)=\sum_{n\geq2}\gamma_n^2\alpha_n^2\cos(n\theta).
\]
The original minimal-form approximations of $a(q_z)$ also
converge in the weighted form norm by
\eqref{eq:form-capacity-comparison}. Scalar integration by parts
against cylinders gives
\begin{equation}\label{eq:weighted-inverse-generator}
 a(q_z)\in\Dom(L_\gamma),\qquad
 L_\gamma a(q_z)
       =k_\gamma(0)a''(q_z)+B_\gamma(h,q_z)a'(q_z).
\end{equation}
Indeed, for $F\in\Ccyl$,
\begin{align*}
 \sum_j\gamma_j^2e_je_j'&=0,\qquad
 \sum_j\gamma_j^2e_j^2=k_\gamma(0),\\
 \shortintertext{so integration by parts gives}
 \cE_\gamma(a(q_z),F)
 &=-\int\{k_\gamma(0)a''(q_z)
                    +B_\gamma(h,q_z)a'(q_z)\}F\,d\nu.
\end{align*}
The series converge uniformly by \eqref{eq:fixed-weight-bounds}.
Form-core approximation extends this identity to all form tests.
For $u_\ell=a_\ell(q_{z_\ell})$, $1\leq\ell\leq m$, the same argument gives
\[
 L_\gamma\prod_{\ell=1}^m u_\ell
 =\sum_{\ell=1}^m(L_\gamma u_\ell)\prod_{k\ne\ell}u_k
   +2\sum_{\ell<k}\sum_j\gamma_j^2D_ju_\ell D_ju_k
                                  \prod_{a\ne\ell,k}u_a\in L^\infty(\nu).
\]

The generator identity and the Markov property give, for
$u=a(q_z)$, the continuous square-integrable martingale
\[
 M_t^u=u(X_t)-u(X_0)-\int_0^tL_\gamma u(X_r)\,dr.
\]
The semigroup identity gives, for $0\leq s\leq t$,
\[
 \E[M_t^u-M_s^u\mid\mathcal F_s]
  =\left(T_{t-s}^\gamma u-u
           -\int_0^{t-s}T_r^\gamma L_\gamma u\,dr\right)(X_s)=0.
\]
The product identity above yields
\[
 d\langle M^u,M^v\rangle_t
       =2\sum_j\gamma_j^2D_ju(X_t)D_jv(X_t)\,dt.
\]
Localization in angular charts gives inverse-lift martingales
$M^{q_z}$.
For Fej\'er weights choose $m>2N$ distinct targets $z_1,\ldots,z_m$
and define the matrix
\[
 ({U^{\mathrm{ev}}}_t)_{a,j}=\sqrt2\gamma_je_j(q_{z_a}(X_t)),
 \qquad 1\leq a\leq m,\quad j=(n,\tau),\quad2\leq n\leq N.
\]
A nonzero trigonometric polynomial of degree at most $N$ has at
most $2N$ zeros. Thus $U_t^{\mathrm{ev}}v=0$ implies $v=0$.
The left inverse defines
\[
 dW_t^j=\sum_{a=1}^m
          [ ((U_t^{\mathrm{ev}})^{\top}U_t^{\mathrm{ev}})^{-1}
                                (U_t^{\mathrm{ev}})^{\top}]_{j,a}
                                                \,dM_t^{q_{z_a}}.
\]
Its quadratic covariation is
\[
 d\langle W^i,W^j\rangle_t=\mathbf1_{\{i=j\}}\,dt.
\]
Localize where the left-inverse norm is bounded; append unused
Brownian coordinates independently.
For Poisson weights and $v\in\ell^2$, square summability gives
uniform convergence, and the dense inverse targets give
\begin{align*}
 \sum_jv_j\gamma_je_j(q_z)=0\quad\text{for every dense target }z
 &\ \Longrightarrow\ \sum_jv_j\gamma_je_j(\theta)=0\quad(\theta\in\T)\\
 &\ \Longrightarrow\ v=0.
\end{align*}
Thus the rows span $\ell^2$. To construct the Brownian family,
write $q_i(t)=q_{z_i}(X_t)$ and, as in
\eqref{eq:noise-reconstruction}, put
\[
 G_N^\gamma(t)_{ik}
   =\langle k^\gamma(q_i(t)),k^\gamma(q_k(t))\rangle_H,
 \qquad
 v_{j,N}^\gamma(t)=(\gamma_j e_j(q_i(t)))_{i=1}^N.
\]
The finite Gram matrices are invertible by the same Fourier
uniqueness argument, since every Poisson weight is positive.
Let $P_N^\gamma(t)$ project onto the first $N$ weighted rows,
and define
\[
 W_t^{j,N}
 =\frac1{\sqrt2}\sum_{i=1}^N\int_0^t
       [(G_N^\gamma)^{-1}v_{j,N}^\gamma]_i(s)
                                      \,dM_s^{q_{z_i}}.
\]
The vector stochastic integral is well-defined: its bracket
density is
\[
 v_{j,N}^{\gamma\,\top}(G_N^\gamma)^{-1}v_{j,N}^\gamma
      =\|P_N^\gamma(t)e_j\|_H^2\leq1.
\]
Thus no uniform bound on the inverse Gram matrices is needed.
The product martingales give
\[
 d\langle W^{j,N},W^{k,M}\rangle_t
       =\langle P_N^\gamma(t)e_j,P_M^\gamma(t)e_k\rangle_H\,dt.
\]
For $M\geq N$, the projections are nested, and Doob's
inequality gives
\[
 \E\sup_{t\leq T}|W_t^{j,M}-W_t^{j,N}|^2
  \leq4\E\int_0^T
      \|(P_M^\gamma(t)-P_N^\gamma(t))e_j\|_H^2\,dt
      \longrightarrow0.
\]
Here $P_N^\gamma(t)\to I$ strongly for every $t$, since the
inverse images of the dense targets remain dense; dominated
convergence applies because the integrand is at most $1$.
Hence $W^{j,N}$ converges in $L^2(\mathbb P;C([0,T]))$ to
a continuous martingale $W^j$. Its limits satisfy
\[
 \langle W^j,W^k\rangle_t=\mathbf1_{\{j=k\}}t;
\]
hence $(W^j)_j$ is a Brownian family. For each inverse-coordinate martingale, the same projection
identity and $P_N^\gamma(t)\to I$ give
\[
 d\langle M^{q_z},W^j\rangle_t
       =\sqrt2\,\gamma_je_j(q_z(X_t))\,dt.
\]
Moreover, the stochastic series
\[
 S_t^z=\sqrt2\sum_j\int_0^t
                \gamma_je_j(q_z(X_s))\,dW_s^j
\]
converges in $L^2$ on every bounded time interval, since
$\sum_j\gamma_j^2e_j(q)^2=k_\gamma(0)<\infty$ uniformly in $q$.
The displayed cross-bracket identity then yields
\[
 d\langle S^z\rangle_t
 =d\langle M^{q_z},S^z\rangle_t
 =d\langle M^{q_z}\rangle_t
 =2k_\gamma(0)\,dt.
\]
Consequently $\langle M^{q_z}-S^z\rangle_t=0$.
Since $M^{q_z}-S^z$ starts at zero, it vanishes identically.
Lemma~\ref{lem:inverse-target-extension} now gives the
self-consistent equation. This is the required
stationary symmetric weak solution.
\end{proof}

\subsection{Dependence on a smooth source change}

Fix a conformally removable initial welding $h_0$.
Such initial weldings form a set of full $\nu$-measure for $0<\kappa\leq4$,
with the normalization used in \cite{Response,Diffusion}.
For an integer $s\geq6$, let $H^s_{\mathrm{per}}$ be the real
Sobolev space of $2\pi$-periodic functions, and use the open set
of degree-one lifts
\begin{equation}\label{eq:sobolev-chart}
 {U_s}=\{\Psi(u)=u+p(u):
              p\in H^s_{\mathrm{per}},\ \inf_u(1+p'(u))>0\}.
\end{equation}
On the removable orbit of $h_0$, write $B_\gamma^{\mathrm{geom}}$
for the normalized conformal-map formula
\eqref{eq:weighted-drift}, regardless of the chosen extension
of that formula outside the full-measure state locus.

Define the drift and noise coefficients on ${U_s}$ by
\begin{align}
 b_\gamma(\Psi)&=
       B_\gamma^{\mathrm{geom}}(h_0\circ\Psi^{-1},\Psi(\cdot)),\notag\\
 \sigma_\gamma(\Psi)&=
       \bigl(\sqrt2\gamma_je_j(\Psi(\cdot))\bigr)_{j\in{J_H}}.
                                      \label{eq:feedback-coefficients}
\end{align}

\begin{lemma}\label{lem:feedback-chart}
The map $\Psi\mapsto B_\gamma^{\mathrm{geom}}(h_0\circ\Psi^{-1},\cdot)$
is locally continuously differentiable from ${U_s}$
to $C^{s+1}(\T)$. Every $\Psi_*\in U_s$ has a chart ball
$\mathcal B\subset U_s$ such that, for $\Psi,\widetilde\Psi\in\mathcal B$,
\begin{align*}
 \|b_\gamma(\Psi)-b_\gamma(\widetilde\Psi)\|_{H^s}
   &\leq C_{\gamma,h_0,\mathcal B}\|\Psi-\widetilde\Psi\|_{H^s},\\
 \shortintertext{and}
 \|\sigma_\gamma(\Psi)-\sigma_\gamma(\widetilde\Psi)\|_{\HS(\ell^2,H^s_{\mathrm{per}})}
   &\leq C_{\gamma,h_0,\mathcal B}\|\Psi-\widetilde\Psi\|_{H^s}.
\end{align*}
\end{lemma}
The proof is given in Appendix~\ref{app:feedback-chart}.

\subsection{Regularity of every relative flow}

\begin{lemma}\label{lem:regularized-flow-smoothness}
Every weak solution of the regularized equation has a relative
flow with continuous paths in ${U_s}$ for every integer
$s\geq6$.
\end{lemma}
\begin{proof}
For this one solution fix its progressive field
$b_t(\theta)=B_\gamma(X_t,\theta)$, extended periodically to $\R$.
For every integer $m\geq0$, \eqref{eq:fixed-weight-bounds} gives
\[
 \sup_t\|b_t^{(m)}\|_\infty<\infty,\qquad
 \sum_j\gamma_j^2\|e_j^{(m)}\|_\infty^2<\infty.
\]
Consider the periodic semimartingale field
\[
 F_t(x)=\sqrt2\sum_j\gamma_je_j(x)W_t^j
                         +\int_0^t b_a(x)\,da.
\]
On every finite time interval its martingale series converges in
$L^2(\Omega;C([0,T];H^{m+1}_{\mathrm{per}}))$ for each $m$.
Sobolev embedding and the bounded drift derivatives therefore make
$F$ a continuous $C^m$-semimartingale field for every $m$.
Its local characteristics are
\[
 a(x,y)=2\sum_j\gamma_j^2e_j(x)e_j(y),\qquad b_t(x),
\]
with bounded spatial derivatives of every order. Complete
sine/cosine pairs give the vanishing It\^o--Stratonovich correction
\[
 \frac12\left.\partial_xa(x,y)\right|_{y=x}
             =\sum_j\gamma_j^2e_j'(x)e_j(x)=0.
\]
Thus \cite[Theorem~2.7.7]{KunitaLectures}, for general
semimartingale fields with progressive local characteristics,
applies to the prescribed-field It\^o equation. It gives a global
flow $\widehat\Theta_t$ of $C^m$ diffeomorphisms for every $m$;
boundedness of the coefficients prevents explosion. Periodicity
and uniqueness give, on one common event over integer $m$ and $T$,
\[
 \widehat\Theta_t(u+2\pi)=\widehat\Theta_t(u)+2\pi.
\]
The flow starts at the identity, and its spatial derivative is
jointly continuous and never vanishes. Hence, for every $T<\infty$,
\[
 \inf_{t\leq T,\,u\in[0,2\pi]}
                       \partial_u\widehat\Theta_t(u)>0.
\]
Continuity in $C^s$ on the circle implies continuity in $H^s$,
so the induced circle flow has paths in $C([0,T],U_s)$.

The original raw relative flow and $\widehat\Theta$ solve
the same prescribed-field scalar equation. Scalar uniqueness
identifies them at the countable dense labels, and continuity
identifies the full flow.
\end{proof}

\subsection{Uniform compact containment and joint tightness}

For each $\gamma$, let $\mathbb Q_\nu^\gamma$ be the stationary joint weak law
of Lemma~\ref{lem:weighted-weak-law}, and let $\E_\nu^\gamma$ be its expectation.
The bounds below are uniform over both multiplier families.
The positive convolution formulas \eqref{eq:matched-convolutions}
give
\begin{equation}\label{eq:uniform-weighted-spatial}
 \|B_\gamma(h,\cdot)\|_\infty\leq C_\kappa,\qquad
 k_\gamma(0)\leq k_0(0)=\frac12 .
\end{equation}

\begin{lemma}\label{lem:weighted-compact-containment}
For every $T>0$ and $\varepsilon>0$ there is a compact
$K\subset \mathcal X$ such that
\[
 \sup_\gamma\mathbb Q_\nu^\gamma
       (X_t^\gamma\notin K\text{ for some }0\leq t\leq T)
                                                        <\varepsilon .
\]
\end{lemma}
\begin{proof}
By Lemma~\ref{lem:weighted-weak-law}, $X^\gamma$ is the stationary
conservative diffusion associated with the quasi-regular symmetric form
$\cE_\gamma$ on the Polish space $\mathcal X$. Its co-process is the
same conservative diffusion. Thus the Lyons--Zheng decomposition
\cite[Theorem~3.1 and Section~4(i)(a)]{Trutnau} applies directly to every
real $u\in\Dom(\cE_\gamma)$, evaluated in its quasi-continuous version.
On $[0,T]$, let $r_T\omega(s)=\omega(T-s)$ and write
$\widehat M_t^u=M_t^u\circ r_T$, where $M^u$ is the forward Fukushima
martingale. The processes $M^u$ and $\widehat M^u$ are continuous
square-integrable martingales in the forward and reversed natural
filtrations, respectively, and
\[
 u(X_t^\gamma)-u(X_0^\gamma)
    =\frac12\{M_t^u-\widehat M_T^u+\widehat M_{T-t}^u\},
 \qquad 0\leq t\leq T.
\]
The martingale energy identity \cite[Section~5.2]{FOT} and stationarity
give, with brackets in the corresponding forward and reversed filtrations,
\[
 \E_\nu^\gamma\langle M^u\rangle_T
       =\E_\nu^\gamma\langle\widehat M^u\rangle_T
       =2T\cE_\gamma(u,u).
\]
Doob's inequality in the two filtrations therefore gives
\begin{equation}\label{eq:weighted-maximal-form}
 \E_\nu^\gamma\sup_{t\leq T}|u(X_t^\gamma)|^2
       \leq2\|u\|_2^2+40T\cE_\gamma(u,u).
\end{equation}
Take the compact nest $K_m$ of \eqref{eq:compact-nest}.
By the variational definition of capacity there are
$u_m\in\Dom(D)$ with $u_m\geq1$ quasi-everywhere on
$\mathcal X\setminus K_m$ and
\[
 \|u_m\|_2^2+\cE(u_m,u_m)
             \leq2\Cap(\mathcal X\setminus K_m)+m^{-1}.
\]
The original quasi-continuous version is also a weighted
quasi-continuous version: any original capacity nest remains
one for the weighted form by
\eqref{eq:form-capacity-comparison}.
Original exceptional sets have zero weighted capacity.
Thus \eqref{eq:weighted-maximal-form} and
$\cE_\gamma(u_m,u_m)\leq\cE(u_m,u_m)$ imply
\[
 \sup_\gamma\mathbb Q_\nu^\gamma
   (\exists t\leq T:X_t^\gamma\notin K_m)
 \leq(2+40T)\{2\Cap(\mathcal X\setminus K_m)+m^{-1}\}
                                                   \longrightarrow0 .
\]
\end{proof}

\begin{lemma}\label{lem:joint-tightness}
The stationary joint laws $\mathbb Q_\nu^\gamma$ are tight on ${E_{\mathrm{path}}}$.
\end{lemma}
\begin{proof}
For a fixed inverse target and a smooth periodic test $a$,
\eqref{eq:weighted-inverse-generator} and
\eqref{eq:uniform-weighted-spatial} give bounded drift and
bounded martingale quadratic-variation density, uniformly
in $\gamma$. The fourth-moment martingale inequality yields
\begin{equation}\label{eq:uniform-inverse-increments}
 \E|a(q_z(X_t^\gamma))-a(q_z(X_b^\gamma))|^4
                  \leq C_{a,T}|t-b|^2,\qquad 0\leq b,t\leq T .
\end{equation}
Thus every finite family of inverse sine and cosine coordinates
is tight in its continuous path space.

These coordinates for a dense set of targets separate points
of $\mathcal X$. On a compact $K\subset \mathcal X$, for each $\eta>0$ there
are finitely many such coordinates $f_1,\ldots,f_m$ and
$\delta>0$ such that
\[
 h,\widetilde h\in K,\quad
 \max_{\ell\leq m}|f_\ell(h)-f_\ell(\widetilde h)|<\delta
                          \quad\Longrightarrow\quad
 d(h,\widetilde h)<\eta .
\]
Otherwise compactness would give two limiting points with
all coordinates equal and distance at least $\eta$.
Combine this finite-coordinate criterion with
Lemma~\ref{lem:weighted-compact-containment} and
\eqref{eq:uniform-inverse-increments}. For every $\eta>0$,
\[
 \lim_{\delta\downarrow0}\sup_\gamma
 \mathbb Q_\nu^\gamma\left(
   \sup_{\substack{b,t\leq T\\|t-b|\leq\delta}}
       d(X_t^\gamma,X_b^\gamma)>\eta\right)=0 .
\]
Compact containment and the metric Arzel\`a--Ascoli criterion
give tightness in $C([0,T],\mathcal X)$. A diagonal choice of compact
sets over integer $T$ gives tightness on $C([0,\infty),\mathcal X)$.

The countable Brownian family has the same product Wiener
law for every $\gamma$, which is a probability on a Polish
space and hence is tight. Compact products then give joint
tightness on ${E_{\mathrm{path}}}$.
\end{proof}

Along either approximation limit, let $K_\gamma$ be the positive
convolution kernel in \eqref{eq:matched-convolutions}, normalized by
$\int_\T K_\gamma(u)\,du/(2\pi)=1$.
The common spatial modulus in \eqref{eq:drift-modulus} gives a bounded
function $\omega$ with $\omega(\delta)\to0$ as $\delta\downarrow0$ and
\[
 |B(h,\theta-u)-B(h,\theta)|
   \leq\omega(\operatorname{dist}_{S^1}(u,0)).
\]
Consequently
\begin{align}
 \sup_{h,\theta}|B_\gamma(h,\theta)-B(h,\theta)|
   &\leq\int_\T K_\gamma(u)
         \omega(\operatorname{dist}_{S^1}(u,0))\,\frac{du}{2\pi}
                                                 \longrightarrow0,\notag\\
 \shortintertext{and}
 \sup_\theta\sum_j|(\gamma_j-1)e_j(\theta)|^2
   &=\sum_{n\geq2}(\gamma_n-1)^2\alpha_n^2\longrightarrow0 .
                                        \label{eq:uniform-coefficient-limits}
\end{align}
The first limit uses the approximate-identity property of $K_\gamma$;
the second uses domination by $\sum_n\alpha_n^2<\infty$.

\subsection{Proof of strong well-posedness in Theorem~\ref{thm:autonomous-approximation}}\label{sec:strong-wellposedness}

Fix one multiplier $\gamma$ and an initial probability $\mu\ll\nu$.
Put $c_j(h)=\gamma_j\rho_j(h)/\sqrt2$ and
$C_\gamma=\sup_h\sum_j|c_j(h)|^2<\infty$.

For any candidate and finite $T$, Novikov's condition follows from
\[
 \E\exp\left(\frac12\int_0^T\sum_j|c_j(X_s)|^2\,ds\right)
                                      \leq e^{C_\gamma T/2}.
\]
For $m\geq2$ define the positive truncated density by
\begin{align*}
 \log R_T^{\gamma,m}
   ={}&-\sum_{n=2}^m\sum_{\tau\in\{c,s\}}
           \int_0^Tc_{n,\tau}(X_s)\,dW_s^{n,\tau}\\
     &-\frac12\int_0^T\sum_{n=2}^m\sum_{\tau\in\{c,s\}}
                    |c_{n,\tau}(X_s)|^2\,ds.
\end{align*}
The exponential martingale estimate gives
\[
 \sup_{m\geq2}\E[(R_T^{\gamma,m})^p]
       \leq \exp\left(\frac{p(p-1)}2C_\gamma T\right),
                                              \qquad p>1.
\]
The It\^o isometry gives
\[
 \E\left|\sum_{n>m,\tau}\int_0^Tc_{n,\tau}(X_s)
                                      \,dW_s^{n,\tau}\right|^2
 \leq \frac{C_\kappa^2T}{2}
          \sum_{n>m,\tau}\frac{\gamma_n^2}{\alpha_n^2}
                                                    \longrightarrow0.
\]
Thus $R_T^{\gamma,m}\to R_T^\gamma$ in $L^1$, where
\begin{equation}\label{eq:remove-drift-density}
 R_T^\gamma=\exp\left(-\sum_j\int_0^Tc_j(X_s)\,dW_s^j
                       -\frac12\int_0^T\sum_j|c_j(X_s)|^2\,ds\right).
\end{equation}
Finite-coordinate Girsanov and the $L^1$ limit make every finite family
of $\widetilde W$ Brownian under $d\widetilde{\mathbb P}_T=R_T^\gamma\,d\mathbb P$, where
\begin{align*}
 \widetilde W_t^j&=W_t^j+\int_0^tc_j(X_s)\,ds,\\
 \shortintertext{and}
 \E[R_T^\gamma\mid\mathcal F_0]&=1.
\end{align*}
Since
\[
 \sqrt2\sum_j\gamma_je_j(\theta)c_j(h)
   =\sum_j\gamma_j^2\rho_j(h)e_j(\theta)=B_\gamma(h,\theta),
\]
the drift cancels and
\begin{equation}\label{eq:zero-drift-flow}
 d\Theta_t(u)=\sqrt2\sum_j\gamma_je_j(\Theta_t(u))\,d\widetilde W_t^j.
\end{equation}
For $x,y\in\R$ and every $m\geq0$,
\begin{align*}
 \sum_j\gamma_j^2|e_j(x)-e_j(y)|^2
        &\leq |x-y|^2\sum_j\gamma_j^2\|e_j'\|_\infty^2,\\
 \shortintertext{and}
 \sum_j\gamma_j^2\|e_j^{(m)}\|_\infty^2&<\infty.
\end{align*}
These bounds and the prescribed-field argument of
Lemma~\ref{lem:regularized-flow-smoothness} give a unique strong smooth
flow for \eqref{eq:zero-drift-flow}. With $X=X_0\circ\Psi^{-1}$, this
determines the joint law of $(X,\widetilde W)$ from $\mu$.

Conversely, on this zero-drift law use
\[
 (R_T^\gamma)^{-1}
 =\exp\left(\sum_j\int_0^Tc_j(X_s)\,d\widetilde W_s^j
                       -\frac12\int_0^T\sum_j|c_j(X_s)|^2\,ds\right).
\]
It recovers
\begin{align*}
 W_t^j&=\widetilde W_t^j-\int_0^tc_j(X_s)\,ds,\\
 \shortintertext{and, for $0\leq t\leq T$,}
 \E_{\mathbb P}[R_T^\gamma\mid\mathcal F_t]&=R_t^\gamma,\\
 \shortintertext{while}
 \E_{\widetilde{\mathbb P}_T}
       [(R_T^\gamma)^{-1}\mid\mathcal F_t]&=(R_t^\gamma)^{-1}.
\end{align*}
Thus the zero-drift law determines a unique joint law on ${E_{\mathrm{path}}}$.
For $\mu=g\nu$, time-zero weighting preserves the label martingales and
Brownian independence, so this law is
\[
 d\mathbb Q_\mu^\gamma=g(X_0)\,d\mathbb Q_\nu^\gamma,
 \qquad \E_\mu^\gamma[\,\cdot\,]=\int(\,\cdot\,)\,d\mathbb Q_\mu^\gamma.
\]
For $\nu(N)=0$, stationarity under $\mathbb Q_\nu^\gamma$ gives
\begin{align}
 \mathbb Q_\mu^\gamma(X_t^\gamma\in N)
   &=\E_\nu^\gamma[g(X_0)\mathbf1_N(X_t^\gamma)]=0,\notag\\
 \shortintertext{and hence}
 \E_\mu^\gamma\int_0^T\mathbf1_N(X_s^\gamma)\,ds
   &=\int_0^T\mathbb Q_\mu^\gamma(X_s^\gamma\in N)\,ds=0.
                                      \label{eq:weighted-null-occupation}
\end{align}

For two solutions sharing $X_0$ and $W$, weak uniqueness,
\eqref{eq:weighted-null-occupation} and
Lemmas~\ref{lem:regularized-flow-smoothness} and
\ref{lem:feedback-chart} give, for $i=1,2$,
\[
 \Psi^i\in C([0,T],U_s),\qquad
 X_t^i=X_0\circ(\Psi_t^i)^{-1}.
\]
The same occupation identity gives
\[
 \E\int_0^T\|B_\gamma(X_t^i,\cdot)
               -B_\gamma^{\mathrm{geom}}(X_t^i,\cdot)\|_\infty^2\,dt=0.
\]
Removability of $X_0$ persists under smooth source composition.

Disintegrate the joint law given $X_0=h_0$ as $\mathbb P_{h_0}$.
The countable label equations, their stochastic-integral limits and the
zero-occupation identity persist for $\mu$-almost every $h_0$.
For finite $J_1\subset J_H$, $a\in\R^{J_1}$, bounded measurable $r$,
and bounded joint-past $G$ up to $b\leq t$, Brownian independence gives
\begin{align*}
 &\E\left[r(X_0)G
       \exp\left(i\sum_{j\in J_1}a_j(W_t^j-W_b^j)\right)\right]
   =e^{-(t-b)|a|^2/2}\E[r(X_0)G],\\
 \shortintertext{and hence}
 &\E_{h_0}\left[G
       \exp\left(i\sum_{j\in J_1}a_j(W_t^j-W_b^j)\right)\right]
   =e^{-(t-b)|a|^2/2}\E_{h_0}G.
\end{align*}
Countable determining tests, continuity and a monotone class give the
second identity simultaneously for $\mu$-almost every $h_0$.
Thus $W$ is Brownian in the usual augmentation $(\mathcal F_t^{h_0})$
of the joint natural filtration under $\mathbb P_{h_0}$.

In the lift coordinate \eqref{eq:sobolev-chart}, difference-quotient
derivatives are adapted. For either lift, Lemma~\ref{lem:feedback-chart}
and compactness of the path image give
\begin{align*}
 &A_t^\gamma
   :=\int_0^t\left(
       \|\sigma_\gamma(\Psi_a)\|_{\HS(\ell^2,H^s_{\mathrm{per}})}^2
        +\|b_\gamma(\Psi_a)\|_{H^s}^2\right)\,da<\infty,\\
 \shortintertext{and, with $\tau_m=\inf\{t:A_t^\gamma\geq m\}\wedge m$,}
 &\tau_m\uparrow\infty,\\
 \shortintertext{and}
 &\mathbb E_{h_0}\int_0^{T\wedge\tau_m}\left(
       \|\sigma_\gamma(\Psi_a)\|_{\HS(\ell^2,H^s_{\mathrm{per}})}^2
        +\|b_\gamma(\Psi_a)\|_{H^s}^2\right)\,da\leq m.
\end{align*}
The stopped $H^s_{\mathrm{per}}$-valued residual is
\[
 \mathfrak r_t^{(m)}
   =\Psi_{t\wedge\tau_m}-\id
      -\int_0^{t\wedge\tau_m}\sigma_\gamma(\Psi_a)\,dW_a
      -\int_0^{t\wedge\tau_m}b_\gamma(\Psi_a)\,da.
\]
Bounded evaluation and the common label identities give
\[
 \operatorname{ev}_u\mathfrak r_t^{(m)}=0
        \quad(u\in\Lambda,\ t\geq0)
 \quad\Longrightarrow\quad
 \mathfrak r_t^{(m)}=0\quad\text{in }H^s_{\mathrm{per}}.
\]
Thus, as $m\to\infty$,
\begin{equation}\label{eq:hilbert-feedback-equation}
 \Psi_t-\id=\int_0^t\sigma_\gamma(\Psi_a)\,dW_a
                         +\int_0^tb_\gamma(\Psi_a)\,da
                         \quad\text{in }H^s_{\mathrm{per}}.
\end{equation}

For lifts agreeing at a finite $(\mathcal F_t^{h_0})$-stopping time $\sigma$,
let $\tau\geq\sigma$ stop either lift at exit from a common smaller chart ball.
Set
\[
 F_{\sigma,\tau}(t)
 =\E_{h_0}\left[\sup_{0\leq v\leq t}
  \|\Psi^1_{(\sigma+v)\wedge\tau}
       -\Psi^2_{(\sigma+v)\wedge\tau}\|_{H^s}^2
                 \,\middle|\,\mathcal F_\sigma^{h_0}\right].
\]
The local Lipschitz bounds, Hilbert BDG and Cauchy--Schwarz give
\begin{align*}
 F_{\sigma,\tau}(t)&\leq C\int_0^tF_{\sigma,\tau}(b)\,db,\\
 \shortintertext{and conditional Gronwall yields}
 F_{\sigma,\tau}(t)&=0.
\end{align*}
A countable chart cover and compactness of the continuous path images
extend these comparisons to every finite horizon. Hence
\[
 \mathbb P(\Psi^1=\Psi^2)
   =\int\mathbb P_{h_0}(\Psi^1=\Psi^2)\,\mu(dh_0)=1.
\]

The self-consistency relation gives $X^1=X^2$.

The input and output spaces are Polish:
\begin{gather*}
 \mathcal X\times\prod_jC_0([0,\infty),\R),\\
 \shortintertext{and}
 C([0,\infty),\mathcal X).
\end{gather*}
Use \eqref{eq:autonomous-state}--\eqref{eq:autonomous-sde} at the countable labels
as the joint-law constraints.
Let $\mathbb W$ denote the product Wiener probability and let
${K^{\mathrm{path}}}_\gamma(dX\mid h,w)$ be the conditional output law of a
weak solution given $(X_0,W)=(h,w)$.

Write $X_{\leq t},W_{\leq t}$ for the paths restricted to $[0,t]$.
For bounded measurable $f$, Brownian independence gives
the compatibility identity of \cite[Lemma~2.12]{Kurtz}:
\begin{equation}\label{eq:regularized-compatibility}
 \mathbb E[f(X_{\leq t})\mid X_0,W]
       =\mathbb E[f(X_{\leq t})\mid X_0,W_{\leq t}]
                 \quad\text{almost surely}.
\end{equation}
Joint compatibility keeps $W$ Brownian in the coupled filtration;
the preceding comparison therefore applies. Thus
\cite[Lemma~2.10 and Theorem~1.5]{Kurtz} gives a Borel map
$\Phi_\gamma$ with
\begin{equation}\label{eq:strong-regularized-map}
 {K^{\mathrm{path}}}_\gamma(dX\mid h,w)=\delta_{\Phi_\gamma(h,w)}(dX),
 \qquad X^\gamma_\cdot=\Phi_\gamma(X_0,W).
\end{equation}
By \cite[Proposition~2.13]{Kurtz}, the output up to $t$ is measurable
with respect to the completed $\sigma(X_0,W_{\leq t})$.
Thus $\Phi_\gamma$ realizes an adapted global strong solution on every
input space with law $\mu\otimes\mathbb W$.
Choose $\Phi_\gamma$ under $\nu$ and use the same version under $\mu=g\nu$.
Joint weak uniqueness gives
\[
 \Law(\Phi_\gamma(X_0,W),W)=\mathbb Q_\mu^\gamma.
\]
For $\mu=\nu$, the state law is therefore stationary and reversible.

\subsection{Proof of convergence in Theorem~\ref{thm:autonomous-approximation}}
\label{sec:final-approximation-proof}

Take $\mu=\nu$ and let
$\gamma_m=\gamma^{N_m}$ with $N_m\to\infty$, or
$\gamma_m=\gamma^{r_m}$ with $r_m\uparrow1$.
Lemma~\ref{lem:joint-tightness} gives a subsequence, still indexed by $m$,
with $\mathbb Q_\nu^{\gamma_m}\Longrightarrow\mathbb Q_*$ on ${E_{\mathrm{path}}}$.
Write $(X,W)$ for the coordinate process under $\mathbb Q_*$.
Continuity of fixed-time evaluation gives \eqref{eq:marginal-class}:
\[
 \mathbb Q_*\circ X_t^{-1}=\nu,\qquad t\geq0.
\]

Separability of $C(\T)$ and scalar $L^2(\nu)$ density give continuous
maps $\widetilde B_\ell:\mathcal X\to C(\T)$ with
$\|\widetilde B_\ell-B\|_{L^2(\nu;C(\T))}\to0$.
Set
\[
 B_\ell(h)=\frac{C_\kappa\widetilde B_\ell(h)}
                       {C_\kappa\vee\|\widetilde B_\ell(h)\|_\infty}.
\]
Then
\begin{equation}\label{eq:continuous-drift-approximation}
 \|B_\ell\|_\infty\leq C_\kappa,\qquad
 \|B_\ell-B\|_{L^2(\nu;C(\T))}
       \leq2\|\widetilde B_\ell-B\|_{L^2(\nu;C(\T))}\longrightarrow0.
\end{equation}
Stationarity gives, uniformly in $m$,
\begin{align}
 \E\int_0^T
   \|B_\ell(X_s^{\gamma_m},\cdot)
                         -B(X_s^{\gamma_m},\cdot)\|_\infty^2\,ds
 &=T\|B_\ell-B\|_{L^2(\nu;C(\T))}^2 .\label{eq:occupation-error}
\end{align}
The same equality holds for $X$. For a smooth periodic $a$,
Cauchy--Schwarz and \eqref{eq:uniform-coefficient-limits} give
\begin{align*}
 &\E\sup_{t\leq T}\left|\int_0^t
  [B_{\gamma_m}-B_\ell](X_s^{\gamma_m},q_z(X_s^{\gamma_m}))
                  a'(q_z(X_s^{\gamma_m}))\,ds\right|^2\\
 &\qquad\leq2T^2\|a'\|_\infty^2
  \left(\|B_{\gamma_m}-B\|_\infty^2
                    +\|B-B_\ell\|_{L^2(\nu;C(\T))}^2\right).\\
 \shortintertext{For the limiting process,}
 &\E\sup_{t\leq T}\left|\int_0^t
 [B-B_\ell](X_s,q_z(X_s))a'(q_z(X_s))\,ds\right|^2
 \leq T^2\|a'\|_\infty^2\|B-B_\ell\|_{L^2(\nu;C(\T))}^2.
\end{align*}
The integral with fixed $B_\ell$ is continuous on the state path space;
take $m\to\infty$ and then $\ell\to\infty$ in these bounds.

In the joint natural filtration, for smooth periodic $a$ and a target $z$
in the fixed countable dense set, define
\[
 M_t^a=a(q_z(X_t))-a(q_z(X_0))
       -\int_0^t\left\{\frac12a''(q_z(X_s))
                     +B(X_s,q_z(X_s))a'(q_z(X_s))\right\}\,ds .
\]
For $0\leq b\leq t$ and a bounded continuous function $H$ of finitely
many past state and noise coordinates, the prelimit identities,
the drift passage above, and \eqref{eq:uniform-coefficient-limits} yield
\begin{align*}
 \E[H(M_t^a-M_b^a)]&=0.\\
 \shortintertext{For the squared martingale,}
 \E\left[H\left\{(M_t^a)^2-(M_b^a)^2
              -2\int_b^t\sum_j|D_j[a(q_z)](X_s)|^2\,ds\right\}\right]&=0.\\
 \shortintertext{The mixed product with the noise satisfies}
 \E\left[H\left\{M_t^aW_t^j-M_b^aW_b^j
               +\sqrt2\int_b^tD_j[a(q_z)](X_s)\,ds\right\}\right]&=0.
\end{align*}
Uniform integrability follows from the bounded inverse-coordinate
drifts and the Brownian moments on $[0,T]$.
A monotone class proves that $M^a$ is a continuous martingale with
\begin{align}
 d\langle M^a\rangle_t
   &=2\sum_j|D_j[a(q_z)](X_t)|^2\,dt,\notag\\
 \shortintertext{and}
 d\langle M^a,W^j\rangle_t
   &=-\sqrt2D_j[a(q_z)](X_t)\,dt .
                                      \label{eq:limit-cross-brackets}
\end{align}
The same passage for $W^j$ and $W^iW^j-\mathbf1_{\{i=j\}}t$
makes $W$ Brownian in the joint filtration by L\'evy's characterization.

The uniform Fourier tail gives
\[
 \E\sup_{t\leq T}\left|\sqrt2\sum_{n>N,\tau}\int_0^t
               D_{n,\tau}[a(q_z)](X_s)\,dW_s^{n,\tau}\right|^2
       \leq8T\|a'\|_\infty^2\sum_{n>N}\alpha_n^2\longrightarrow0.
\]
By \eqref{eq:limit-cross-brackets},
\begin{align*}
 d\left\langle M^a+\sqrt2\sum_j\int_0^\cdot
                       D_j[a(q_z)](X_s)\,dW_s^j\right\rangle_t
 &=2\sum_j|D_j[a(q_z)](X_t)|^2\,dt
   +2\sum_j|D_j[a(q_z)](X_t)|^2\,dt\\
 &\quad-4\sum_j|D_j[a(q_z)](X_t)|^2\,dt=0.
\end{align*}
The continuous martingale starts at zero and therefore vanishes.

Taking $a=\sin,\cos$ in local angular charts gives the inverse-lift SDE
on one event for all targets in the dense set.
Lemma~\ref{lem:inverse-target-extension}, with $\gamma_j=1$,
extends it to time-zero measurable targets. At
$X_0(e^{iu})$ this is \eqref{eq:raw-sde}, with
\[
 \Psi_t=X_t^{-1}\circ X_0,\qquad X_t=X_0\circ\Psi_t^{-1}.
\]
Continuity of composition and inversion makes $\Psi$ adapted and continuous.
After the usual augmentation, $(X,W)$ is an admissible weak solution.
Theorem~\ref{thm:weak-uniqueness} gives $\mathbb Q_*=\mathbb Q_\nu$;
tightness therefore yields
\begin{align*}
 \mathbb Q_\nu^{\gamma^N}&\Longrightarrow\mathbb Q_\nu
                  &&(N\to\infty),\\
 \shortintertext{and}
 \mathbb Q_\nu^{\gamma^r}&\Longrightarrow\mathbb Q_\nu
                  &&(r\uparrow1).
\end{align*}

For $\mu=g\nu$, write $\gamma\to1$ for either approximation limit and choose
bounded continuous $g_\ell:\mathcal X\to\R$ with
$\|g_\ell-g\|_{L^1(\nu)}\to0$.
For bounded continuous $\Phi$ on ${E_{\mathrm{path}}}$ and
$\mathbb Q=\mathbb Q_\nu^\gamma$ or $\mathbb Q_\nu$, the initial marginal $\nu$ gives
\[
 \left|\int \Phi\{g(X_0)-g_\ell(X_0)\}\,d\mathbb Q\right|
                 \leq\|\Phi\|_\infty\|g-g_\ell\|_{L^1(\nu)}.
\]
For fixed $\ell$, stationary joint convergence applies to the bounded
continuous function $\Phi g_\ell(X_0)$. The two errors above give
\[
 \limsup_{\gamma\to1}\left|
   \int\Phi g(X_0)\,d\mathbb Q_\nu^\gamma
       -\int\Phi g(X_0)\,d\mathbb Q_\nu\right|
 \leq2\|\Phi\|_\infty\|g-g_\ell\|_{L^1(\nu)}
 \xrightarrow[\ell\to\infty]{}0.
\]
This proves \eqref{eq:main-approximation} and completes the proof of Theorem~\ref{thm:autonomous-approximation}.

\clearpage
\appendix
\section{Response geometry and regularized source charts}\label{app:response-estimates}

\subsection{Proof of Lemma~\ref{lem:coefficient-domain}}\label{app:coefficient-domain}

\begin{proof}
It suffices to consider holomorphic polynomials.

For a normalized univalent $f$ put $s_f(z)=\sqrt{z/f(z)}$, with
$s_f(0)=1$. Its odd exterior transform
$F_f^{\mathrm{ext}}(\zeta)=\zeta s_f(\zeta^{-2})$ is univalent.
The area and growth theorems for univalent functions \cite{Pommerenke} give
\begin{equation}\label{eq:area-bounds}
 \sum_{m\geq1}(2m-1)|[z^m]s_f|^2\leq1,\qquad
 \|s_f\|_\infty\leq2,\qquad
 \int_\D|s_f'(z)|^2\,dA(z)\leq\pi .
\end{equation}
Set
\[
 g_\ell(f;z)=s_f(z)^{2\ell}
                  (s_f(z)-2zs_f'(z))^2
            =f'(z)^2\left(\frac z{f(z)}\right)^{\ell+3}.
\]
Equation \eqref{eq:area-bounds} implies
$\int_\D|g_\ell(f;z)|\,dA(z)\leq16\pi4^\ell$.
Coefficient extraction in \eqref{eq:negative-contour} gives
\begin{align}
 \ell_{-n}a_m
 &=\sum_{\ell=0}^{m-1}c_{m\ell}(f)
              [z^{n+\ell+1}]g_\ell(f;z),\label{eq:kernel-coefficients}\\
 \shortintertext{where}
 c_{m\ell}(f)&=[z^{m+1}]f(z)^{\ell+2}.\notag
\end{align}
All $c_{m\ell}$ in a fixed expression are uniformly bounded.
For $g\in A^1$ put
\[
 ({R_\ell^{\mathrm{tail}}} g)(z)
 =z^{-\ell-3}
       \left(g-\sum_{j=0}^{\ell+2}[z^j]g\,z^j\right)
\]
and use submean and Cauchy estimates to obtain
\[
 \sum_{j=0}^{\ell+2}|[z^j]g|
 +\sup_{|z|\leq1/2}|{R_\ell^{\mathrm{tail}}} g(z)|
 \leq C_\ell\|g\|_{A^1}.
\]
On the remaining annulus,
\[
 \int_{1/2<|z|<1}|{R_\ell^{\mathrm{tail}}} g(z)|\,dA(z)
 \leq C_\ell\|g\|_{A^1}.
\]
Thus $\|{R_\ell^{\mathrm{tail}}} g\|_{A^1}\leq C_\ell\|g\|_{A^1}$.
Consequently
\begin{equation}\label{eq:polynomial-kernel}
 {K_P^{\mathrm{resp}}}(f;z)
 =\sum_m(\partial_mP)(f)
       \sum_{\ell=0}^{m-1}c_{m\ell}(f){R_\ell^{\mathrm{tail}}} g_\ell(f;z)
 =\sum_{n\geq2}(\ell_{-n}P)(f)z^{n-2}
\end{equation}
satisfies $\sup_f\int_\D|{K_P^{\mathrm{resp}}}(f;z)|\,dA(z)<\infty$.
For any integrable holomorphic $H$, submean estimates and integration
give
\begin{align}
 \sup_{z\in\D}(1-|z|^2)^2|H(z)|
      &\leq C\int_\D|H|\,dA.\notag\\
 \shortintertext{Multiplication by $|H|$ and integration then give}
 \int_\D|H|^2(1-|z|^2)^2\,dA
      &\leq C\left(\int_\D|H|\,dA\right)^2.\label{eq:bergman-embedding}
\end{align}
Using
$\int_\D |z|^{2n-4}(1-|z|^2)^2\,dA=2\pi/d_n$,
we obtain
\begin{equation}\label{eq:negative-response-energy}
 \sup_f\sum_{n\geq2}\frac{|\ell_{-n}P(f)|^2}{d_n}<\infty .
\end{equation}
By \eqref{eq:coefficient-responses} and
\eqref{eq:negative-response-energy},
\begin{align*}
 \ell_nP&=0,\qquad n>\deg_wP,\\
 \shortintertext{and hence}
 \sup_f\sum_{n\geq2,\tau}|D_{n,\tau}P(f)|^2
 &=\sup_f\sum_{n\geq2}
       \frac{|\ell_nP(f)|^2+|\ell_{-n}P(f)|^2}{d_n}
       <\infty.
\end{align*}

To identify the generator, choose an integer $M$ at least the weighted degree of
$P$, and put $\ell_0=\sum_m ma_m\partial_m$. For $n>M$ the contour
identities imply
\begin{equation}\label{eq:trace-tail-polynomials}
 \begin{aligned}
 \ell_nP&=0,\qquad \ell_n\ell_{-n}P=2n\ell_0P,\\
 &\text{and}\qquad
 \sum_{\tau\in\{c,s\}}D_{n,\tau}^2P
     =-\frac{2}{n^2-1}\ell_0P.
 \end{aligned}
\end{equation}
Indeed, the variation $\delta f=z^{n+1}f'$ gives
\[
 \delta g_\ell=z^{n+1}g_\ell'+(2n-\ell-1)z^ng_\ell,
 \qquad
 [z^{n+\ell+1}]\delta g_\ell=2n[z^{\ell+1}]g_\ell.
\]
The zero-mode residue in \eqref{eq:kernel-coefficients} is $ma_m$.
This proves \eqref{eq:trace-tail-polynomials} first on $a_m$
and then by the derivation rule.

Define the bounded candidate drift
\begin{align}
 b_P(f)
 ={}&-\sum_{n\geq2}
       \frac{(\ell_n\ell_{-n}+\ell_{-n}\ell_n)P(f)}{d_n}
       -\sum_{n=2}^{M}\frac{t_n(f)\ell_nP(f)}{d_n}\notag\\
 &+\frac1{2\pi}\int_\D
       \overline{Q_f(z)}\,{K_P^{\mathrm{resp}}}(f;z)
                          (1-|z|^2)^2\,dA(z).
       \label{eq:polynomial-drift}
\end{align}
The first sum converges uniformly absolutely by
\eqref{eq:trace-tail-polynomials}. The uniform estimate
$|Q_f(z)|(1-|z|^2)^2\leq C_\kappa$ from \cite{Diffusion}
and \eqref{eq:polynomial-kernel} bound the last integral.
For $0<r<1$, the coefficient pairing is
\begin{align*}
 \sum_{n\geq2}\frac{r^n\overline{t_n(f)}\ell_{-n}P(f)}{d_n}
 &=\frac1{2\pi}\int_\D
       \overline{r^2Q_f(rz)}\,{K_P^{\mathrm{resp}}}(f;z)
                         (1-|z|^2)^2\,dA(z).\\
 \shortintertext{The uniform stress bound and dominated convergence give}
 \lim_{r\uparrow1}\sum_{n\geq2}
       \frac{r^n\overline{t_n(f)}\ell_{-n}P(f)}{d_n}
 &=\frac1{2\pi}\int_\D
       \overline{Q_f(z)}\,{K_P^{\mathrm{resp}}}(f;z)
                         (1-|z|^2)^2\,dA(z).
\end{align*}

For a weak solution satisfying \eqref{eq:marginal-class}, put
$b_t=B(X_t,\cdot)$ on its original filtered space and use its Brownian
family $W$. The Osgood estimate in \cite[Section~4]{Diffusion} is
formulated for an adapted field on this filtration, so it applies to
this $b_t$ without conditioning on the environment. The angular lift
of the regularized circle flow solves
\begin{align*}
 d\Theta_t^r(u)&=\sqrt2\sum_{n\geq2,\tau}r^ne_{n,\tau}(\Theta_t^r(u))
                                      \,dW_t^{n,\tau}\\
 &\quad+(\operatorname{Pois}_rb_t)(\Theta_t^r(u))\,dt,
  \qquad\Theta_0^r(u)=u.
\end{align*}
Put $h_t^r=X_0\circ(\Psi_t^r)^{-1}$. By the fixed-environment Osgood
stability estimate of \cite[Section 4]{Diffusion}, for every
$\epsilon>0$ and a deterministic short block $[0,T]$,
\[
 \mathbb P\left(\sup_{t\leq T}d(h_t^r,X_t)>\epsilon\right)
                         \longrightarrow0\qquad(r\uparrow1).
\]
Equation \eqref{eq:drift-modulus} gives the uniform Zygmund bound
for $b_t$ and the associated bounded Beltrami extension; $\Law(X_0)=\nu$
places the initial welding in the removable locus almost surely. These
are the hypotheses of \cite[Lemmas 5.2--5.8 and Section 6]{Diffusion}.
Hence the boundary paths of the conformal pairs
$(\widetilde f^r,\widetilde g^r)$, normalized by
$\widetilde g^r(w)=w+O(w^{-1})$, are tight in
$C([0,T]\times S^1)^2$, and subsequential limits remain injective.

Take any two sequences $r_k,s_k\uparrow1$. Joint tightness gives a
subsequential limit of the state and both conformal-pair sequences.
Write its state coordinate as $X$ and its pair coordinates as
$(\widetilde f^i,\widetilde g^i)$, $i=1,2$. The convergence
$h^{r_k},h^{s_k}\to X$ in probability and the boundary-matching
identities pass to this limit, so both limiting pairs weld the same
$X_t$. For rational $t$, removability and uniqueness of the normalized
conformal pair therefore give
\[
 (\widetilde f_t^1,\widetilde g_t^1)
       =(\widetilde f_t^2,\widetilde g_t^2)\quad\text{almost surely},
 \qquad t\in\mathbb Q\cap[0,T],
\]
where $\Law(X_t)=\nu$ was used.
Continuity extends this identity to the entire time interval:
\[
 \mathbb P\left(
   (\widetilde f_t^1,\widetilde g_t^1)
      =(\widetilde f_t^2,\widetilde g_t^2)
          \text{ for every }t\in[0,T]\right)=1.
\]
Tightness and uniqueness of joint subsequential limits therefore give,
for every $\epsilon>0$,
\[
 \lim_{r,s\uparrow1}\mathbb P\left(
   \|\widetilde f^r-\widetilde f^s\|_{C([0,T]\times S^1)}
  +\|\widetilde g^r-\widetilde g^s\|_{C([0,T]\times S^1)}
                                      >\epsilon\right)=0.
\]
For the resulting limit $(\widetilde f,\widetilde g)$, Cauchy's formula
and nonvanishing of $\widetilde f_t'(0)$ give
\[
 \mathbb P\left(\inf_{t\leq T}|\widetilde f_t'(0)|<m^{-1}\right)
                              \longrightarrow0\qquad(m\to\infty).
\]
Consequently the affine normalization
\[
 f_t^r=\frac{\widetilde f_t^r-\widetilde f_t^r(0)}
                         {(\widetilde f_t^r)'(0)},\qquad
 g_t^r=\frac{\widetilde g_t^r-\widetilde f_t^r(0)}
                         {(\widetilde f_t^r)'(0)}
\]
passes to the limit boundary-uniformly in probability. Deterministic
restarts cover every finite interval, since
\[
 \Law(X_{kT})=\nu,\qquad k\in\mathbb N_0.
\]
At each fixed time the selected maps of $X_t$ equal the resulting
continuous pair. Thus, on every finite interval,
\begin{equation}\label{eq:prescribed-pair-convergence}
 (f^r,g^r)\longrightarrow(f,g)
       \quad\text{in probability in }C([0,T]\times S^1)^2.
\end{equation}
Removability and Fubini also give zero boundary area for
$dt\,d\mathbb P$-almost every $(t,\omega)$.

For the remaining argument take the canonical stationary solution.
In particular, \eqref{eq:prescribed-pair-convergence} gives
$f_t^r\to f_t$ in probability in $C(\overline\D)$ at each time.

The removable limiting curve has zero area. For normalized maps
$f_m\to f$ in $C(\overline\D)$ with this limiting boundary,
$s_{f_m}\to s_f$ uniformly. If $V_m,V$ are the omitted sets of the
odd exterior transforms, winding numbers and boundedness give
\begin{align*}
 \mathbf1_{V_m}&\longrightarrow\mathbf1_V
       &&dA\text{-almost everywhere},\\
 \shortintertext{and hence, by dominated convergence,}
 |V_m|&\longrightarrow|V|.
\end{align*}
The area theorem now yields
\[
 \sum_{j\geq1}(2j-1)|[z^j]s_{f_m}|^2
 =1-|V_m|/\pi
 \longrightarrow1-|V|/\pi
 =\sum_{j\geq1}(2j-1)|[z^j]s_f|^2.
\]
Together with coefficientwise convergence, this gives strong convergence
of the weighted coefficient sequences. Consequently,
\begin{align*}
 \|s_{f_m}'-s_f'\|_{A^2}^2
       &=\sum_{j\geq1}j|[z^j](s_{f_m}-s_f)|^2\\
       &\leq\sum_{j\geq1}(2j-1)|[z^j](s_{f_m}-s_f)|^2
       \longrightarrow0.
\end{align*}
Equation \eqref{eq:area-bounds} now gives
\[
 \|g_\ell(f_m)-g_\ell(f)\|_{A^1}
 \leq C_\ell
       (\|s_{f_m}-s_f\|_\infty+\|s_{f_m}'-s_f'\|_{A^2})
 \longrightarrow0.
\]
The bounded coefficient factors and ${R_\ell^{\mathrm{tail}}}$ in
\eqref{eq:polynomial-kernel} yield
\begin{equation}\label{eq:strong-kernel-convergence}
 g_\ell(f_m;\cdot)\to g_\ell(f;\cdot)\text{ in }L^1(\D),
 \qquad\text{and}\qquad
 {K_P^{\mathrm{resp}}}(f_m;\cdot)\to{K_P^{\mathrm{resp}}}(f;\cdot)
                       \text{ in }L^1(\D).
\end{equation}
Moreover, \eqref{eq:bergman-embedding} implies
\[
 \sum_{n\geq2}\frac{|\ell_{-n}P(f_m)-\ell_{-n}P(f)|^2}{d_n}
 =\frac1{2\pi}\int_\D
       |{K_P^{\mathrm{resp}}}(f_m;z)-{K_P^{\mathrm{resp}}}(f;z)|^2
                          (1-|z|^2)^2\,dA(z)
 \longrightarrow0.
\]
For a fixed regularization, the Beltrami source charts used in
\cite{Diffusion} make each interior coefficient twice continuously
differentiable. It\^o's formula gives
\[
 dP(f_t^r)=-\sqrt2\sum_{n,\tau}r^n
                D_{n,\tau}P(f_t^r)\,dW_t^{n,\tau}
 +\left[\sum_{n,\tau}r^{2n}D_{n,\tau}^2P(f_t^r)
       -\sum_{n,\tau}r^n\rho_{n,\tau}(X_t)
                         D_{n,\tau}P(f_t^r)\right]dt .
\]
The prescribed score pairing is
\begin{align*}
 \sum_{n,\tau}r^n\rho_{n,\tau}(X_t)D_{n,\tau}P(f_t^r)
 ={}&\sum_{n=2}^{M}\frac{r^nt_n(f_t)\ell_nP(f_t^r)}{d_n}\\
 &-\frac1{2\pi}\int_\D
       \overline{r^2Q_{f_t}(rz)}\,{K_P^{\mathrm{resp}}}(f_t^r;z)
                              (1-|z|^2)^2\,dA(z).\\
 \shortintertext{The stress factor satisfies the uniform bound}
 \sup_{r<1,z\in\D}
       |r^2Q_{f_t}(rz)|(1-|z|^2)^2&\leq C_\kappa.
\end{align*}
Denote the bracketed It\^o drift above by $b_t^{P,r}$.
The trace tail \eqref{eq:trace-tail-polynomials}, strong kernel
convergence, and dominated convergence give
\begin{align*}
 \E\int_0^T|b_s^{P,r}-b_P(f_s)|^2\,ds&\longrightarrow0,\\
 \shortintertext{and}
 \E\int_0^T\sum_{n,\tau}
       |r^nD_{n,\tau}P(f_s^r)-D_{n,\tau}P(f_s)|^2\,ds
       &\longrightarrow0.\\
 \shortintertext{Doob's inequality then gives}
 \E\sup_{t\leq T}\left|
       \sqrt2\sum_{n,\tau}\int_0^t
       (r^nD_{n,\tau}P(f_s^r)-D_{n,\tau}P(f_s))\,dW_s^{n,\tau}
                      \right|^2&\longrightarrow0.
\end{align*}
Consequently
\[
 P(f_t)-P(f_0)
 =-\sqrt2\sum_j\int_0^tD_jP(f_s)\,dW_s^j
                      +\int_0^t b_P(f_s)\,ds .
\]
It remains to identify this drift with the closed generator.
For bounded state functions $F$ and $b$, a martingale identity with
drift $b$ along the canonical solution gives, by conditioning at time
zero, the $L^2(\nu)$ identity
\[
 T_tF-F=\int_0^tT_sb\,ds.
\]
Strong continuity gives
\begin{equation}\label{eq:canonical-generator-limit}
 \left\|\frac{T_tF-F}{t}-b\right\|_2
 \leq\frac1t\int_0^t\|T_sb-b\|_2\,ds\longrightarrow0.
\end{equation}
Thus $F\in\Dom(L)$ and $LF=b$. Applying this to $F=P$ and $b=b_P$
gives $P\in\Dom(L)$ and $LP=b_P$.
Let $D_j^{\mathrm{fin}}P$ denote the finite coefficient response
in \eqref{eq:coefficient-responses}. For each fixed direction and
$\phi\in\Ccyl$,
\[
 \langle (DP)_j,\phi\rangle
       =\langle P,\rho_j\phi-D_j\phi\rangle
       =\langle D_j^{\mathrm{fin}}P,\phi\rangle .
\]
The second equality is the finite polynomial integration-by-parts
input stated above.
Density identifies the components. This proves the lemma.
\end{proof}

\subsection{Proof of Lemma~\ref{lem:hessian-bound}}\label{app:score-hessian}
\begin{proof}
Fix a normalized univalent $f$.
Complexify real bilinear forms bilinearly and set
$e_n^+=\alpha_n l_n$, $n\geq2$.
Their Hermitian norms are $\sqrt2$.
Write
\[
 S_n=[z^{n-2}]{\operatorname{S}}f,\qquad
 b_n=[z^{n-2}]\left((f'/f)^2-z^{-2}\right),\qquad
 \beta_{nm}=\frac{d_n(n+2m)}{d_{n+m}},
\]
and define
$R^S_{nm}=\ell_{-m}S_n-\beta_{nm}S_{n+m}$ and
$R^b_{nm}=\ell_{-m}b_n-\beta_{nm}b_{n+m}$.
Substitution in \eqref{eq:full-connection}, together with
\eqref{eq:coefficient-responses}, gives
\begin{align}
 {H^\rho}(e_n^+,e_m^+)
    &=\alpha_n\alpha_m
          \overline{a_\kappa R^S_{nm}+R^b_{nm}}.\label{eq:hessian-block}
\intertext{The mixed block is tridiagonal:}
 {H^\rho}(e_n^+,\overline{e_m^+})
    &=0\qquad (|n-m|\geq2),\notag\\
 {H^\rho}(e_n^+,\overline{e_n^+})
    &=2a_\kappa+\frac4{n^2-1},\notag\\
 {H^\rho}(e_n^+,\overline{e_{n+1}^+})
    &=\frac{\alpha_n\alpha_{n+1}}2(2n+1)t_1.\label{eq:mixed-block}
\end{align}
For the mixed entries, the positive stress response is
\[
 \ell_nt_m=
 \begin{cases}
  (m+n)t_{m-n},&m>n,\\
  a_\kappa d_n+2n,&m=n,\\
  0,&m<n.
 \end{cases}
\]
For $m\geq n+2$, the connection coefficient satisfies
$c(n,-m)=n+m$, so the two terms cancel explicitly:
\begin{align*}
 {H^\rho}(e_n^+,\overline{e_m^+})
 &=\alpha_n\alpha_m
       \bigl(\ell_nt_m-c(n,-m)t_{m-n}\bigr)\\
 &=\alpha_n\alpha_m\bigl((m+n)-(m+n)\bigr)t_{m-n}=0.
\end{align*}
The entries with $n\geq m+2$ vanish by conjugate symmetry.
Since $|t_1|\leq4$, the matrix in \eqref{eq:mixed-block} is bounded.

For the Schwarzian part, consider the Grunsky kernel
\[
 G_f(z,w)=\frac{f'(z)f'(w)}{(f(z)-f(w))^2}
                  -\frac1{(z-w)^2}
       =\sum_{p,q\geq0}g_{pq}z^pw^q.
\]
Its normalized coefficient matrix is
\[
 (G_f^{\mathrm{Gr}})_{jk}=\frac{g_{j-1,k-1}}{\sqrt{jk}},
 \qquad j,k\geq1.
\]
The Grunsky inequality gives $\|{G_f^{\mathrm{Gr}}}\|_{\ell^2\to\ell^2}\leq1$
and $G_f(z,z)={\operatorname{S}}f(z)/6$ \cite{Pommerenke}.
It applies to arbitrary normalized univalent $f$: apply the smooth
case to $f(rz)/r$ and pass to every finite matrix as $r\uparrow1$.
Differentiating \eqref{eq:negative-contour} gives
\[
 \ell_{-m}{\operatorname{S}}f(z)
 =\frac6{2\pi i}\int_{|w|=R}
       \frac{f'(z)^2f'(w)^2w^{1-m}}{(f(w)-f(z))^4}\,dw .
\]
Indeed the third derivative of $y^2/[A^2(A-y)]$ is
$6/(A-y)^4$. Squaring $(z-w)^{-2}+G_f(z,w)$ and taking
coefficients yields
\begin{align}
 \ell_{-m}S_n
 ={}&12\sum_{\ell=0}^{n-2}(\ell+1)g_{n-2-\ell,m+\ell}
       +6\sum_{p=0}^{n-2}\sum_{q=0}^{m-2}
              g_{pq}g_{n-2-p,m-2-q}.\label{eq:schwarzian-residues}
\end{align}
Let $(\varepsilon_j)_{j\geq1}$ be the standard basis of $\ell^2(\N)$.
Define $J_{\mathrm{Gr}}:\ell^2(\{2,3,\ldots\})
\to\ell^2(\N)\otimes\ell^2(\N)$ by
\[
 J_{\mathrm{Gr}}\varepsilon_n=
       \sum_{\substack{j,k\geq1\\j+k=n}}\sqrt{\frac{6jk}{d_n}}\,
                           \varepsilon_j\otimes\varepsilon_k,
 \qquad n\geq2.
\]
The normalized quadratic part of \eqref{eq:schwarzian-residues}
is $2J_{\mathrm{Gr}}^*({G_f^{\mathrm{Gr}}}\otimes{G_f^{\mathrm{Gr}}})
J_{\mathrm{Gr}}$.
The disjoint frequency levels give
\begin{align*}
 \langle J_{\mathrm{Gr}}\varepsilon_n,
                    J_{\mathrm{Gr}}\varepsilon_m\rangle
 &=\delta_{nm}\frac6{d_n}\sum_{j=1}^{n-1}j(n-j)
   =\delta_{nm}.\\
 \shortintertext{Hence}
 \|2J_{\mathrm{Gr}}^*({G_f^{\mathrm{Gr}}}\otimes
             {G_f^{\mathrm{Gr}}})J_{\mathrm{Gr}}\|
 &\leq2\|{G_f^{\mathrm{Gr}}}\|^2\leq2.
\end{align*}

To bound the linear part in trace norm, for a finite matrix
$M=(M_{jk})_{j,k\geq1}$ define
\[
 ({P_{\mathrm{coef}}} M)_{jk}
   =\frac{6\sqrt{jk}}{d_{j+k}}
             \sum_{r+s=j+k}\sqrt{rs}\,M_{rs}.
\]
With $\xi(z)_j=\sqrt j\,z^{j-1}$ and
$F_M(z)=\sum_{j,k}\sqrt{jk}M_{jk}z^{j+k-2}$,
\[
 {P_{\mathrm{coef}}} M=\frac3\pi\int_\D F_M(z)(1-|z|^2)^2
                         \xi(\bar z)\otimes\xi(\bar z)\,dA(z).
\]
The equality follows from the radial moment $2/d_{j+k}$.
Since $\|\xi(z)\|^2=(1-|z|^2)^{-2}$, Cauchy--Schwarz for two
Bergman functions proves the rank-one estimate
$\|{P_{\mathrm{coef}}}(x\otimes y)\|_1\leq3\|x\|\|y\|$.
Singular-value decomposition gives
\begin{equation}\label{eq:hankel-trace-bound}
 \|{P_{\mathrm{coef}}} M\|_1\leq3\|M\|_1 .
\end{equation}
The subscript $1$ here denotes trace norm.

For finite sequences $u,v$ supported on $n\geq2$, extend them by
zero elsewhere and set
\[
 M(u,v)_{jk}=\sqrt{jk}\sum_{\ell\geq1}\ell u_{j+\ell}v_{k-\ell},
 \qquad
 R_{\mathrm{lin}}(u,v)=(1-{P_{\mathrm{coef}}})(M(u,v)+M(u,v)^{\top}).
\]
The finite sum identity
\[
 \sum_{j=1}^{n+m-1}j(n+m-j)
              \{(n-j)_++(j-m)_+\}=\frac{d_n(n+2m)}6
\]
identifies $R_{\mathrm{lin}}$ with the linear correction in $R^S$.
Also $R_{\mathrm{lin}}(u,v)=R_{\mathrm{lin}}(v,u)$, because
$\beta_{nm}-\beta_{mn}=n-m$ and the corresponding bracketed
positive-part expressions differ by $n-m$.

Use dyadic blocks $u^{(k)},v^{(\ell)}$ on
$[2^k,2^{k+1})$ and $[2^\ell,2^{\ell+1})$.
For $k\leq\ell$, let $N=2^{k+1}$ and factor $M=L_uR_v$, with
\[
 (L_u)_{jr}=\sqrt{jr}\,u_{j+r},\qquad
 (R_v)_{rp}=\mathbf1_{\{r<N\}}\sqrt{rp}\,v_{p-r}.
\]
Writing $\|u\|_d^2=\sum_{n\geq2}d_n|u_n|^2$, one obtains
\begin{align*}
 \|L_u\|_{\HS}^2&=\frac16\|u^{(k)}\|_d^2,
\shortintertext{whereas}
 \|R_v\|_{\HS}^2
 &\leq\sum_m|v_m|^2
                 \left(\frac{mN^2}2+\frac{N^3}3\right)
 \leq\frac{56}9\,2^{2(k-\ell)}\|v^{(\ell)}\|_d^2 .
\end{align*}
Schatten H\"older, \eqref{eq:hankel-trace-bound}, and symmetry give
\[
 \|R_{\mathrm{lin}}(u^{(k)},v^{(\ell)})\|_1
 \leq8\sqrt{28/27}\,2^{-|k-\ell|}
                       \|u^{(k)}\|_d\|v^{(\ell)}\|_d.
\]
Since
\begin{align*}
 \sum_{j\in\mathbb Z}2^{-|j|}
 &=1+2\sum_{j\geq1}2^{-j}=3,\\
 \shortintertext{convolution of the dyadic norms gives}
 \|R_{\mathrm{lin}}(u,v)\|_1
 &\leq8\sqrt{28/27}\left(\sum_{j\in\mathbb Z}2^{-|j|}\right)
                \|u\|_d\|v\|_d\\
 &=24\sqrt{28/27}\,\|u\|_d\|v\|_d.
\end{align*}
Apply this with $u_n=\alpha_n x_n$, $v_n=\alpha_n y_n$.
The linear Schwarzian pairing is
$6\Tr({G_f^{\mathrm{Gr}}}R_{\mathrm{lin}}(u,v))$, and
$\|u\|_d=\sqrt2\|x\|$, $\|v\|_d=\sqrt2\|y\|$.
Together with the quadratic part this proves
\begin{equation}\label{eq:schwarzian-operator-bound}
 \|(\alpha_n\alpha_mR^S_{nm})_{n,m\geq2}\|
 \leq2+288\sqrt{28/27}<300 .
\end{equation}

For the remaining stress term, let
$f_{\mathrm{odd}}(z)=\sqrt{f(z^2)}$, normalized by
$f_{\mathrm{odd}}'(0)=1$.
It is odd and univalent. The Schwarzian chain rule and
\eqref{eq:negative-contour}, symmetrized under $w\mapsto-w$, give
\[
 S_{2n}(f_{\mathrm{odd}})=4S_n(f)+\frac32b_n(f),\qquad
 \ell_{-m}^{\,f}f_{\mathrm{odd}}
       =\frac12\ell_{-2m}^{\,f_{\mathrm{odd}}}f_{\mathrm{odd}} .
\]
It follows, with $s=n+m$, that
\begin{align}
 R^b_{nm}(f)
 &=\frac13R^S_{2n,2m}(f_{\mathrm{odd}})-\frac83R^S_{nm}(f)
                         +E_{nm}S_{2s}(f_{\mathrm{odd}}),\label{eq:odd-reduction}
\shortintertext{where}
 E_{nm}
 &=\frac{2nm(n+2m)(2n+m)}
                    {s(4s^2-1)(s^2-1)}.\notag
\end{align}
Nehari's bound
$|{\operatorname{S}}f_{\mathrm{odd}}(z)|(1-|z|^2)^2\leq6$
\cite{Nehari} and Cauchy's estimate at $r^2=1-1/s$ give
\[
 |S_{2s}(f_{\mathrm{odd}})|\leq6es^2,\qquad
 E_{nm}\leq\frac{4nm}{s^3},\qquad
 \frac{\alpha_n}{\alpha_{2n}}\leq\sqrt{10}.
\]
Consequently the last normalized matrix in \eqref{eq:odd-reduction}
satisfies
\begin{align*}
 |\alpha_n\alpha_mE_{nm}S_{2(n+m)}(f_{\mathrm{odd}})|
       &\leq\frac{96e}{\sqrt{nm}(n+m)}.
\shortintertext{Its row sums obey}
 \sum_{m\geq2}\frac{96e}{\sqrt{nm}(n+m)}
       &\leq96e\int_0^\infty\frac{dx}{\sqrt{nx}(n+x)}
        =\frac{96\pi e}{n}\leq48\pi e.
\end{align*}
The Schur test and \eqref{eq:schwarzian-operator-bound} yield
\begin{equation}\label{eq:full-stress-block}
 \|(\alpha_n\alpha_m(a_\kappa R^S_{nm}+R^b_{nm}))\|
 \leq300\left(|a_\kappa-8/3|+10/3\right)+48\pi e.
\end{equation}
Equations \eqref{eq:mixed-block} and \eqref{eq:full-stress-block}
prove the horizontal operator bound.

The horizontal--low block equals
\[
 {H^\rho}(e_i,z_a)
 =D_i\sigma_a+\frac12\sum_j\Omega_{ij}^a\rho_j,\qquad
 \Omega_{ij}^a=([e_i,e_j],z_a)_{\widehat H}.
\]
The rotation coefficients are
$\Omega_{(n,c),(n,s)}^0=-2/(n^2-1)$; the other nonzero blocks
join adjacent frequencies and have absolute size
$(2n+1)/\sqrt{d_nd_{n+1}}$.
The Bers estimate in \eqref{eq:stress} gives
$|\rho_{n,\tau}|\leq C_\kappa\sqrt n$.
Each row and column of $\Omega^a$ has only a fixed number of
nonzero entries, all of size $O(n^{-2})$ at frequency $n$. Hence
\begin{equation}\label{eq:score-low-contraction}
 \sup_h\sum_{j,a}\left|\sum_i\Omega_{ij}^a\rho_i(h)\right|^2
 \leq C_\kappa\sum_{n\geq2}n^{-4}n<\infty.
\end{equation}
Since $t_1=2a_1$, the uniform bound
\eqref{eq:negative-response-energy} and the response identification
in Lemma~\ref{lem:coefficient-domain} give
\[
 \sup_h\sum_{i\in J_H}|D_it_1(h)|^2<\infty.
\]
The low-low block involves only finitely many coefficient polynomials
and the scores $\rho_{2,c},\rho_{2,s}$; each is uniformly bounded.
Together with the preceding two estimates, this proves
\eqref{eq:uniform-hessian}.

For the capacity score $2\widehat Dc$, retain the $R^b_{nm}$ block in
\eqref{eq:odd-reduction}, the tridiagonal mixed block, and the preceding
low-mode estimates. The same Schur and Hilbert--Schmidt bounds give
\begin{equation}\label{eq:capacity-hessian-bound}
 \|\operatorname{Hess}c\|_{\widehat H\to\widehat H}
 \leq C,
 \qquad
 \operatorname{Hess}c(v,w)=D_vD_wc-D_{\nabla_vw}c.
\end{equation}

Every matrix entry is measurable in $h$. Finite-vector approximation
and the uniform bound therefore imply that
$h\mapsto{H^\rho}_hv$ is measurable for each $v\in\widehat H$.
\end{proof}

\subsection{Proof of Lemma~\ref{lem:ricci-bound}}\label{app:ricci}
\begin{proof}
Its Fourier eigenvalues are
\begin{align}
 \Ric_0&=\frac{\pi^2}6-\frac78,&
 \Ric_1&=\frac38-\frac{\pi^2}6,\label{eq:ricci-low}\\
 \Ric_2&=-\frac{65}{16},&
 \Ric_n&=-\frac{13}6
       -\frac{10n^2-4}{n^2(n^2-1)(n^2-4)},\quad n\geq3.
       \label{eq:ricci-high}
\end{align}
In particular $\|\Ric\|\leq65/16$.

Write $c(p,q)$ for the coefficient in \eqref{eq:full-connection}.
The horizontal connection coefficient is
\[
 c_H(p,q)=
 \begin{cases}
 0,& |q|\leq1\text{ or }|p+q|\leq1,\\
 p-q,& |p|\leq1,\ |q|\geq2,\ |p+q|\geq2,\\
 c(p,q),& |p|,|q|,|p+q|\geq2.
 \end{cases}
\]
Let $R_H$ be defined by \eqref{eq:curvature-coefficients} with $c_H$.
For $n,k\geq2$ one has $R_H(k,n,-n)=0$. With
$F_n(x)=(2x-n)(x+n)/(x^3-x)$, substitution gives
\[
 \frac{R_H(-k,n,-n)}{d_n}
 =\begin{cases}
 F_n(k+n)-(k+n)(2n-k)/d_n,&2\leq k\leq n+1,\\
 F_n(k+n)-F_n(k),&k\geq n+2.
 \end{cases}
\]
The tail is $O_n(k^{-2})$. Summing and telescoping gives
\[
 \frac1{d_n}\sum_{|k|\geq2}R_H(k,n,-n)
 =-\frac1{d_n}\sum_{k=2}^{n+1}(k+n)(2n-k)=-\frac{13}6.
\]
For the full connection put
$J_{0,n}=\{0,\pm1,n-1,n,n+1,-n-1,-n,-n+1\}$.
Outside $J_{0,n}$ the difference of the two curvature summands is
$nk/\omega_k$. The set $J_{0,n}$ is symmetric, and
\[
 \sum_{k\notin J_{0,n}}\left|\frac{nk}{\omega_k}\right|<\infty,
 \qquad
 \frac{nk}{\omega_k}+\frac{-nk}{\omega_{-k}}=0.
\]
Put
\[
 \Delta_n(k)=R(k,n,-n)
       -\mathbf1_{\{|k|\geq2\}}R_H(k,n,-n).
\]
For example, the $0,\pm n$ contributions combine as
\begin{align*}
 \Delta_n(0)+\Delta_n(n)+\Delta_n(-n)
 &=\frac{n^2}{d_n}+0-\frac{3n^2}{d_n}\\
 &=-\frac{2n^2}{d_n}.
\end{align*}
Substituting \eqref{eq:full-connection} into
\eqref{eq:curvature-coefficients} for the remaining indices and
collecting over the common denominator gives
\begin{align*}
 8\sum_{k\in J_{0,2}}\Delta_2(k)+91&=0,\\
 n(n^2-4)\sum_{k\in J_{0,n}}\Delta_n(k)+10n^2-4&=0,
                                             \qquad n\geq3.
\end{align*}
Since $d_n=n(n^2-1)$, division by $d_n$ yields
\begin{align*}
 \frac1{d_n}\sum_{k\in J_{0,n}}\Delta_n(k)
 &=\begin{cases}
 -91/48,&n=2,\\
 -(10n^2-4)/[n^2(n^2-1)(n^2-4)],&n\geq3.
 \end{cases}
\end{align*}
This proves \eqref{eq:ricci-high}.

For $n=0$, $R(k,0,0)=k^2/\omega_k^2$, so
\[
 \Ric_0=\frac12+2\sum_{k\geq2}\frac1{(k^2-1)^2}
                         =\frac{\pi^2}6-\frac78.
\]
For $n=1$, put $s_k=(2k+1)^2/(d_kd_{k+1})$.
For $k\geq3$,
\[
 R(k,1,-1)=\frac3{k^2-1}+s_{k-1},\qquad
 R(-k,1,-1)=-\frac3{k^2-1}+s_k .
\]
After division by $\omega_1=2$, the three parts are
\begin{align*}
 \frac1{\omega_1}\sum_{|k|\leq1}R(k,1,-1)&=0,\\
 \frac1{\omega_1}\sum_{|k|=2}R(k,1,-1)&=-\frac{407}{288},\\
 \shortintertext{and}
 \frac1{\omega_1}\sum_{|k|\geq3}R(k,1,-1)
    &=\sum_{k\geq2}s_k-\frac{25}{288}.
\end{align*}
The decomposition
\[
 s_k=\frac34\left(\frac1{k-1}-\frac1{k+2}\right)
 -\frac54\left(\frac1k-\frac1{k+1}\right)
 -\frac12\left(\frac1{k^2}+\frac1{(k+1)^2}\right)
\]
therefore gives
\[
 \Ric_1=-\frac32+\sum_{k\geq2}s_k=\frac38-\frac{\pi^2}6.
\]
These series are absolutely convergent.
\end{proof}

\subsection{Proof of Lemma~\ref{lem:two-commutator-bounds}}\label{app:commutators}
\begin{proof}
In the positive horizontal block the first squared commutator norm
is diagonal in the input frequency, with jump kernel
\begin{equation}\label{eq:connection-jump-kernel}
 {k_\nabla}(n,m)=
 \frac{(2\max(n,m)-\min(n,m))^2d_{\min(n,m)}}
      {d_{|n-m|}d_{\max(n,m)}},
 \qquad n,m\geq2,\quad |n-m|\geq2.
\end{equation}
There are no entries between positive and negative horizontal
covector frequencies. For $2\leq k\leq n/2$,
\begin{align}
 {k_\nabla}(n,n\pm k)&\leq Cn^2/k^3,\label{eq:jump-near}
\shortintertext{and the paired difference satisfies}
 |{k_\nabla}(n,n+k)-{k_\nabla}(n,n-k)|
                         &\leq Cn/k^2.\notag
\end{align}
For the far parts,
\[
 {k_\nabla}(n,m)\leq
 \begin{cases}
 Cm^3/n^4,&m<n/2,\\
 Cn^3/m^4,&m>3n/2.
 \end{cases}
\]
Let $s(x)=\log\log(e+x)$.
Near $n$, $|s(n+k)-s(n)|\leq Ck/[n\log(e+n)]$;
in the far parts use
$|s(m)-s(n)|\leq|\log(m/n)|$ and dyadic summation.
These estimates give
\[
 \sup_n\sum_m{k_\nabla}(n,m)|s(m)-s(n)|^2<\infty,
\]
which proves the first bound in \eqref{eq:two-commutator-bounds}
for the horizontal block.
The low derivative directions add adjacent jumps with coefficients
$O(n^2)$; low covector couplings have squared size $O(n^{-4})$.
Their corresponding moments are bounded as well.

On $[n/2,3n/2]$, the slow cutoff satisfies
\[
 |w_R'|\leq\frac C{Rn\log(e+n)},\qquad
 |w_R''|\leq\frac C{Rn^2\log(e+n)}.
\]
For the score term use the actual spatial field
\[
 b(h,\theta)=B(h,\theta)+\sigma_c(h)\cos\theta
                              +\sigma_s(h)\sin\theta.
\]
It has uniformly bounded Zygmund norm by
\eqref{eq:drift-modulus} and \eqref{eq:low-scores}.
Writing $b_k$ for its Fourier coefficients, the positive horizontal
block of $\sum_A\widehat\rho_A\Xi_A$, up to its imaginary
sign, has entries
\[
 b_{n-m}\,a(n,m),\qquad
 a(n,m)=(2\max(n,m)-\min(n,m))
                         \sqrt{d_{\min(n,m)}/d_{\max(n,m)}} .
\]
For $|n-m|=1$, completing the metric adds only a weighted shift
of size $O(n^{-2})$ times a bounded first coefficient.
In the commutator the main entries acquire
$w_R(m)-w_R(n)$.

Let $(\Delta_q)_{q\geq0}$ be smooth dyadic Fourier projections, and let
$b_+$ and $b_-$ be the positive and negative Fourier parts of $b$.
Their dyadic kernels have uniformly bounded $L^1$ norms. Consequently
\begin{align*}
 \|\Delta_qb_\pm\|_\infty&\leq C_\kappa2^{-q}.
\shortintertext{Summing the differentiated blocks up to frequency $N$ gives}
 \left\|\partial_\theta
       \sum_{2^q\leq N}\Delta_qb_\pm\right\|_\infty
       &\leq C_\kappa\log(e+N).
\end{align*}
For $n,m\asymp N$ and $n>m$, smoothly extend
$(2n-m)\sqrt{d_m/d_n}$ to the comparable rectangle.
Let $\sigma_{R,N}$ be its product with
$(w_R(n)-w_R(m))/(n-m)$, rescaled to $n/N,m/N$ and compactly
supported in a fixed rectangle. Its derivatives satisfy
\begin{align*}
 \|\sigma_{R,N}\|_{C^4}
 &\leq\frac C{R\log(e+N)}.
\shortintertext{Four integrations by parts therefore give}
 \sum_{k,\ell\in\Z}|\widehat\sigma_{R,N}(k,\ell)|
 &\leq\frac C{R\log(e+N)}.
\end{align*}
Each Fourier term multiplies a matrix on the left and right by
diagonal unitaries. Hence, for every bounded matrix $A$,
\[
 \|(\sigma_{R,N}(n/N,m/N)A_{nm})\|
 \leq\left(\sum_{k,\ell}|\widehat\sigma_{R,N}(k,\ell)|\right)\|A\|.
\]
The matrix with entries $(n-m)b_{n-m}$ restricted to this rectangle
and to positive differences has norm at most
$C_\kappa\log(e+N)$, by the preceding low-pass derivative bound.
Thus the commutator block has norm at most $C_\kappa/R$.
The negative-difference block is treated by the adjoint expression.
Comparable dyadic blocks have finite overlap.

For separated dyadic blocks $N=2^q$, $M=2^\ell$ with $\max(n,m)\asymp N$,
$\min(n,m)\asymp M\leq N/4$, the same Fourier-multiplier estimate gives
\[
 \|(b_{n-m}a(n,m)(w_R(m)-w_R(n)))_{n\asymp N,m\asymp M}\|
 \leq\frac{C_\kappa}R
             (1+|\log(N/M)|)(M/N)^{3/2}.
\]
Indeed $|n-m|\asymp N$,
$\|\sum_{|a-q|\leq3}\Delta_ab\|_\infty\leq C_\kappa/N$, and the rescaled
symbol derivatives obey the displayed bound. Summing the dyadic blocks uses
\[
 \sum_{q\in\Z}(1+|q|)2^{-3|q|/2}<\infty.
\]
The low rows and columns are bounded by
$C_\kappa R^{-1}n^{-3/2}(1+s(n))$; their squared sum is finite.

For $0<r<1$, let $\widehat\rho_A^{(r)}$ be the Abel-regularized full
score obtained by replacing every score coefficient $b_k$ with
$r^{|k|}b_k$, and define on finite vectors
\[
 S_R^{(r)}v:=-\sum_A\widehat\rho_A^{(r)}[\Xi_A,Q_R]v,
\]
including the low rows and columns. The Zygmund estimates above are
uniform in $r$, hence
\begin{equation}\label{eq:score-commutator-uniform}
 \sup_{0<r<1}\|S_R^{(r)}\|_{\widehat H\to\widehat H}
                                      \leq \frac{C_\kappa}{R}.
\end{equation}
For every finite Fourier vector $v$, coefficientwise convergence and
the $O(n^{-3/2})$ column tail give
\[
 \lim_{r,s\uparrow1}\|(S_R^{(r)}-S_R^{(s)})v\|_{\widehat H}=0.
\]
Thus, on finite vectors,
\[
 S_Rv=\lim_{r\uparrow1}S_R^{(r)}v
      =-\sum_A\widehat\rho_A[\Xi_A,Q_R]v.
\]
Equation \eqref{eq:score-commutator-uniform} extends $S_R$ to
$\widehat H$ with $\|S_R\|\leq C_\kappa/R$.

For $N\geq2$, let $\mathcal I_N$ index $z_0,z_c,z_s$ and both real
directions of every horizontal frequency $2\leq k\leq N$, and define
\[
 D_{R,N}=\sum_{A\in\mathcal I_N}
                   [\Xi_A,[\Xi_A,Q_R]].
\]
Rotation equivariance makes $D_{R,N}$ diagonal. Since
$\nabla_{l_p}l_n$ has frequency $p+n$, its horizontal main block is
\[
 D_{R,N}e_n^+
 =2\sum_{\substack{m\geq2\\2\leq|m-n|\leq N}}
       {k_\nabla}(n,m)(w_R(m)-w_R(n))e_n^+
       +\text{the low corrections from }\mathcal I_N.
\]
The negative block is its conjugate. In the near range put
$K_\pm={k_\nabla}(n,n\pm k)$. Both jumps occur whenever
$2\leq k\leq\min(N,\lfloor n/2\rfloor)$, and Taylor's formula gives
\begin{align*}
 &\left|K_+(w_R(n+k)-w_R(n))
             +K_-(w_R(n-k)-w_R(n))\right|\\
 &\qquad\leq k|w_R'(n)|\,|K_+-K_-|
       +\frac{k^2}{2}(K_++K_-)
             \sup_{x\in[n/2,3n/2]}|w_R''(x)|\\
 &\qquad\leq\frac{C}{R\log(e+n)}\frac1k.
 \shortintertext{Summing over the retained pairs yields}
 \sum_{2\leq k\leq\min(N,n/2)}
       \frac{C}{R\log(e+n)}\frac1k&\leq\frac CR.
\end{align*}
In the far ranges, the absolute estimates use
\[
 |w_R(m)-w_R(n)|
       \leq\frac CR\left|\log\frac mn\right|,
\]
so every retained subsum is bounded by $C/R$.
The $p=\pm1$ derivative directions are both retained. Writing
$k_\pm^{(1)}(n)$ for their squared jump weights and $k_{\mathrm{low}}(n,m)$
for the total squared weight from a high input $n$ to a low
covector frequency $m$, the full connection gives
\begin{align*}
 k_\pm^{(1)}(n)&\leq Cn^2,
 &|k_+^{(1)}(n)-k_-^{(1)}(n)|&\leq Cn,\\
 \shortintertext{whereas}
 \sum_{|m|\leq1}k_{\mathrm{low}}(n,m)&\leq Cn^{-4}.
\end{align*}
The adjacent jumps satisfy the paired Taylor bound; the low covector
couplings and, by skew-adjointness, the low rows are absolutely
summable. These estimates are uniform in $N$, and each fixed input
column has an absolutely summable far tail. Hence
\begin{align*}
 \sup_N\|D_{R,N}\|&\leq\frac CR,\\
 \lim_{M,N\to\infty}\|(D_{R,N}-D_{R,M})v\|_{\widehat H}&=0
       \quad(v\text{ of finite support}).\\
 \shortintertext{By density,}
 D_R&=\mathrm{s}\!\operatorname{-lim}_{N\to\infty}D_{R,N},
       \qquad\|D_R\|\leq C/R.
\end{align*}

For a fixed normalized Fourier vector
$u_n=l_n/\sqrt{\omega_n}$ and large $k$, the connection and Zygmund estimates give
\[
 \|\Xi_{e_{k,\tau}}u_n\|\leq C_n k^{-2},
 \qquad
 |\rho_{k,\tau}|\leq C_\kappa\sqrt k,
 \qquad \tau\in\{c,s\}.
\]
Thus
\[
 \sum_{k\geq2,\tau}
       \|\rho_{k,\tau}\Xi_{e_{k,\tau}}u_n\|
 \leq C_{n,\kappa}\sum_{k\geq2}k^{-3/2}<\infty.
\]
On finite vectors set
\[
 A_N(h)=\sum_{A\in\mathcal I_N}\widehat\rho_A(h)\Xi_A,
 \qquad F_{R,N}(h)=A_N(h)Q_R.
\]
Finite rank of $Q_R$ gives $F_{R,N}\to F_R$ in operator norm.
Skew-symmetry of $A_N$ on finite vectors gives the bounded extensions
\begin{align*}
 -[A_N,Q_R]&=-F_{R,N}-F_{R,N}^*,\\
 \shortintertext{and, by the finite-column Abel limit,}
 -[A_N,Q_R]&\longrightarrow-F_R-F_R^*=S_R
                       \quad\text{in operator norm}.
\end{align*}
Consequently,
\begin{align*}
 {T_R^{\nabla}}&=D_R+S_R\\
  &=\mathrm{s}\!\operatorname{-lim}_{N\to\infty}
        \bigl(D_{R,N}-[A_N,Q_R]\bigr),\\
 \shortintertext{and}
 \mathop{\mathrm{ess\,sup}}_h\|{T_R^{\nabla}}(h)\|
 &\leq\|D_R\|+\mathop{\mathrm{ess\,sup}}_h\|S_R(h)\|
 \leq\frac{C_\kappa}{R}.
\end{align*}
\end{proof}

\subsection{Proof of Lemma~\ref{lem:feedback-chart}}\label{app:feedback-chart}
\begin{proof}
Fix $\Psi_*\in{U_s}$ and put $h_*=h_0\circ\Psi_*^{-1}$.
Sobolev embedding gives
\[
 0<\inf\Psi_*'\leq\sup\Psi_*'<\infty,
 \qquad (\Psi_*^{-1})'=\frac1{\Psi_*'\circ\Psi_*^{-1}}.
\]
Thus $\Psi_*^{-1}$ is a $C^1$ bi-Lipschitz quasisymmetric map,
and \cite[Lemma~3.1]{Response} makes $h_*$ removable.
Choose $(f_*,g_*)$ by dividing its $f'(0)=1$ pair by the complex
leading exterior coefficient. Equation~\eqref{eq:JL.2} gives
\[
 f_*(0)=0,\qquad
 g_*(w)=w+c_*+O(w^{-1}),\qquad f_*(\D)\subset4\overline\D,
\]
where $c_*$ is the exterior constant coefficient in this normalization.

For $\Psi$ near $\Psi_*$ write
$\eta=\Psi\circ\Psi_*^{-1}=\id+\varphi$. Then
\[
 \varphi=(\Psi-\Psi_*)\circ\Psi_*^{-1},\qquad
 \|\varphi\|_{H^s}\leq C_{\Psi_*}\|\Psi-\Psi_*\|_{H^s}.
\]
Choose $0<R_0<R_1<1$, with $r^2<R_0$ for Poisson weights.
Let $\chi_0$ be a smooth radial cutoff, zero on $[0,R_1]$
and one near $1$, and set
\[
 \Lambda_\varphi(ae^{i\theta})
       =a\exp(i\{\theta+\chi_0(a)\varphi(\theta)\}).
\]
Direct differentiation gives
\begin{align*}
 \partial\Lambda_\varphi
  &=e^{i\chi_0\varphi}
         \left(1+\frac{\chi_0\varphi'}2
                       +\frac{ia\chi_0'\varphi}2\right),\\
 \bar\partial\Lambda_\varphi
  &=\frac12e^{2i\theta+i\chi_0\varphi}
                        (-\chi_0\varphi'+ia\chi_0'\varphi).\\
 \shortintertext{Consequently,}
 \left\|\frac{\bar\partial\Lambda_\varphi}
                   {\partial\Lambda_\varphi}\right\|_\infty
  &\leq C\|\varphi\|_{C^1}\leq C\|\varphi\|_{H^s}.
\end{align*}
After shrinking the chart this norm is at most $q<1$,
$\Lambda_\varphi$ is a disk diffeomorphism, and
$\Lambda_\varphi=\id$ on $R_1\D$.
Transport the dilatation to the fixed physical domain:
\begin{equation}\label{eq:transported-beltrami}
 \lambda_\varphi=
 \left(\frac{f_*'}{\overline{f_*'}}
                  \frac{\bar\partial\Lambda_\varphi}
                       {\partial\Lambda_\varphi}\right)\circ f_*^{-1}
       \quad\text{on }f_*(\D),\qquad
 \lambda_\varphi=0\quad\text{elsewhere}.
\end{equation}
The derivative formulas and the fixed transport give, for
$\xi\in H^s_{\mathrm{per}}$,
\[
 \supp\lambda_\varphi\subset4\overline\D,\qquad
 \|d_\varphi\lambda_\varphi[\xi]\|_{L^\infty\cap L^p}
       \leq C_p\|\xi\|_{H^s},\qquad 1\leq p<\infty.
\]
This derivative is continuous in $\varphi$.

Let $\operatorname{Be}$ and $\operatorname{Ca}$ be respectively the Beurling
transform and the normalized Cauchy transform on $\C$.
Shrink the chart so that $q\|\operatorname{Be}\|_{L^p\to L^p}<1$
for some $p>2$. The principal solution of the Beltrami
equation is
\begin{equation}\label{eq:beltrami-neumann}
 v_\lambda=(\operatorname{Id}_{L^p}-\lambda\operatorname{Be})^{-1}\lambda,\qquad
 W_\lambda(z)=z+\operatorname{Ca} v_\lambda(z).
\end{equation}
The principal-solution construction of \cite{AIM} applies because
\[
 \|\lambda_\varphi\operatorname{Be}\|_{L^p\to L^p}
       \leq q\|\operatorname{Be}\|_{L^p\to L^p}<1,
 \qquad \|\operatorname{Ca} v\|_\infty\leq C_p\|v\|_p
             \quad(\supp v\subset4\overline\D).
\]
Put $R_\lambda=(\operatorname{Id}_{L^p}-\lambda\operatorname{Be})^{-1}$.
The derivative is
\[
 dR_\lambda[\dot\lambda]
       =R_\lambda\dot\lambda\operatorname{Be} R_\lambda.
\]
Thus $\varphi\mapsto W_{\lambda_\varphi}-\id$ is continuously
differentiable in the uniform norm.
The new welded pair is
\[
 f_\varphi=W_{\lambda_\varphi}\circ f_*\circ\Lambda_\varphi^{-1},
 \qquad g_\varphi=W_{\lambda_\varphi}\circ g_* .
\]
The complex-dilatation chain rule and \eqref{eq:transported-beltrami}
give
\[
 \bar\partial f_\varphi=\bar\partial g_\varphi=0,
 \qquad
 g_\varphi^{-1}\circ f_\varphi=h_*\circ\eta^{-1}=h_0\circ\Psi^{-1}
                                      \quad\text{on }S^1.
\]
On the smaller disk, the source deformation is the identity, so
\[
 f_\varphi=W_{\lambda_\varphi}\circ f_*
                                      \quad\text{on }R_1\D.
\]
For each integer $m\geq0$, Cauchy's estimates give
\[
 \|d_\varphi\partial_z^m f_\varphi[\xi]
                       \|_{C(R_0\overline\D)}
 \leq C_m\|\xi\|_{H^s}.
\]
The normalization
\[
 f_\varphi^\circ(z)
       =\frac{f_\varphi(z)-f_\varphi(0)}{f_\varphi'(0)}
\]
preserves this differentiability since $\inf_\varphi|f_\varphi'(0)|>0$
locally. Subtracting the two principal terms in \eqref{eq:stress}
then gives a continuously differentiable regular remainder
on each smaller closed disk.

For Fej\'er weights only finitely many such coefficients occur.
For Poisson weights the formula is
\[
 B_{\gamma^r}^{\mathrm{geom}}(h,\theta)
                    =2\Ima A_{f_h^\circ}(r^2e^{i\theta}).
\]
For the same parameter direction, Cauchy's estimates give
\begin{align*}
 \left|d_\varphi t_n(f_\varphi^\circ)[\xi]\right|
 &\leq C R_0^{2-n}\|\xi\|_{H^s}.\\
 \shortintertext{Summing the regularized coefficients yields}
 \left\|d_\varphi B_{\gamma^r}^{\mathrm{geom}}[\xi]\right\|_{C^{s+1}}
 &\leq C\left(\sum_{n\geq2}\frac{n^{s+1}}{d_n}
                          \left(\frac{r^2}{R_0}\right)^n\right)
                    \|\xi\|_{H^s}.
\end{align*}
The series converges because $r^2<R_0$.
Sobolev product and composition estimates now give, on each
smaller chart ball,
\[
 \|b_\gamma(\Psi)-b_\gamma(\widetilde\Psi)\|_{H^s}
                  \leq C\|\Psi-\widetilde\Psi\|_{H^s}.
\]
For the noise, the same estimates followed by square summation
give
\begin{align*}
 \|\sigma_\gamma(\Psi)-\sigma_\gamma(\widetilde\Psi)\|_{\HS}
 &\leq C_\gamma\|\Psi-\widetilde\Psi\|_{H^s},\\
 \shortintertext{where, for either fixed regularization,}
 C_\gamma^2&=C\sum_{n\geq2}\gamma_n^2\alpha_n^2n^{2(s+2)}<\infty.
\end{align*}
\end{proof}

\section{Gaussian operator and stress estimates}\label{app:gaussian-estimates}

\subsection{Gaussian exponential verification}\label{app:gaussian-exponential}

Retain the notation of Section~\ref{sec:gaussian-lift}. For finite
field coordinates $\mathbf y$, Wick's rule shows that
\[
 \frac{K(\Vzero_n e^{\mathbf u\cdot\mathbf y})
       -K[(\bar\Sigma_n-K\bar\Sigma_n)e^{\mathbf u\cdot\mathbf y}]}
       {M(\mathbf u)}
\]
is a polynomial of degree at most two in $\mathbf u$. The linear and
holomorphic quadratic coefficients are those in
\eqref{eq:KD.7}--\eqref{eq:KD.8}; the antiholomorphic coefficients are
$-\overline{\ell_{-n}m_r}$ and
$-\overline{\ell_{-n}C_{rs}}$ by \eqref{eq:KD.11}--\eqref{eq:KD.12}.
The mixed holomorphic--antiholomorphic coefficients vanish, since
\[
 K[(\Vzero_nY_r)\bar y_s+y_r\Vzero_n\bar Y_s]
 -K[(\bar\Sigma_n-K\bar\Sigma_n)y_r\bar y_s]=0,
 \qquad y_r=Y_r-m_r.
\]
Consequently
\begin{align*}
 &K(\Vzero_n e^{\mathbf u\cdot\mathbf y})
       -K[(\bar\Sigma_n-K\bar\Sigma_n)e^{\mathbf u\cdot\mathbf y}]\\
 &\qquad=M(\mathbf u)
       \left(\mathbf u\cdot\ell_n\mathbf m
                   +\tfrac12\mathbf u^\top(\ell_n\mathbf C)\mathbf u\right)
       =\ell_nM(\mathbf u).
\end{align*}

\subsection{Finite-flow estimates}\label{app:finite-flow-estimates}

Retain the notation of Section~\ref{sec:joint-finite-flow}, with
$\pi(dh,dY,dZ)=\nu_\beta(dh)K_h^Q(dY,dZ)$ before reconstruction.
Its field marginal is \eqref{eq:KD.26}.

Define standard real Gaussian coordinates by
\[
 \widehat\xi_{2n-1}=\sqrt{2n}\operatorname{Re}Y_n,\qquad
 \widehat\xi_{2n}=\sqrt{2n}\operatorname{Im}Y_n\quad(n\geq1).
\]

Choose an even integer $d$ large enough for the $Y$-coordinates used by
$\Sigma_v$ and $P$. In these standardized real coordinates,
$\Sigma_v$ is a finite quadratic polynomial
\[
 \Sigma_v(x)=x^\top A_vx+\ell_v^\top x+c_v,
 \qquad x\in\R^d,\quad A_v=A_v^\top.
\]
The first $d$ standardized coordinates of $T_{\psi_s}Y$ have the form
\[
 C_s\widehat\xi+d_s,
 \qquad C_s:\ell^2(\R)\longrightarrow\R^d.
\]
Each row of $C_s$ is rapidly decreasing in the Fourier index. Indeed,
after changing variables in a field pairing, its test function is
smooth; its Fourier coefficients decay faster than any power. This
decay is uniform for $s$ in a compact interval, also after differentiating
in $s$. The deterministic vector $d_s$ is obtained from
$QP_0\log\psi_s'$. Thus $C_s$ has uniformly bounded Hilbert--Schmidt norm
on compact time intervals. Define
\begin{align*}
 M_t&=\int_0^t C_s^\top A_vC_s\,ds,\\
 \shortintertext{and}
 b_t&=\int_0^t C_s^\top(2A_vd_s+\ell_v)\,ds,\\
 \shortintertext{and}
 a_t&=\int_0^t
      (d_s^\top A_vd_s+\ell_v^\top d_s+c_v)\,ds.
\end{align*}
These are Bochner integrals in the trace-class operators, $\ell^2(\R)$,
and $\R$, respectively. For example,
\[
 \|C_s^\top A_vC_s\|_1\leq\|A_v\|\,\|C_s\|_{\mathrm{HS}}^2.
\]
Consequently, $M_t$ is real symmetric and trace class, and in
$L^2(\mathbb G)$,
\begin{align*}
 J_t
 &=\lim_{N\to\infty}\left(
     \sum_{k,l\leq N}({M_t})_{kl}\widehat\xi_k\widehat\xi_l
       +\sum_{k\leq N}({b_t})_k\widehat\xi_k\right)+{a_t}.\\
 \shortintertext{For $|t|\leq1$, their norms satisfy}
 \|{M_t}\|_1+\|{b_t}\|_2+|{a_t}|
 &\leq C_v|t|.
\end{align*}
For $2|p|\|{M_t}\|<1$, finite Gaussian integration and its
trace-class limit give
\begin{align*}
 \mathbb E_{\mathbb G}e^{pJ_t}
 ={}&e^{p{a_t}}\det(\operatorname{Id}-2p{M_t})^{-1/2}\\
 &\times\exp\left(\frac{p^2}{2}
   \left\langle(\operatorname{Id}-2p{M_t})^{-1}
                      {b_t},{b_t}\right\rangle\right).
\end{align*}
For small $|t|$, the formula for both signs and H\"older's inequality
prove \eqref{eq:finite-flow-exponential}.

To prove \eqref{eq:finite-flow-graph}, let $\Pi_M$ retain the first $M$ Fourier frequencies of
$\widehat\xi$, and let $Y^{[M]}$ be the corresponding truncated field.
Put
\begin{align*}
 J_{t,M}(Y)&=\int_0^t\Sigma_v(T_{\psi_s}Y^{[M]})\,ds,\\
 H_{t,M}(Y,Z)&=P(T_{\psi_t}Y^{[M]},Z)e^{-J_{t,M}(Y)}.
\end{align*}
The quadratic coefficients of $J_{t,M}$ are
$\Pi_MM_t\Pi_M$, $\Pi_Mb_t$, and $a_t$. Their norm bounds, and hence the
preceding small-time exponential-moment bounds, are uniform in $M$.

For a rapidly decreasing real coefficient sequence $a$, write
$\lambda_a=\sum_j a_j\widehat\xi_j$. Since $v$ has finite Fourier
support, its field derivation has the form
\[
 \Vzero_v\lambda_a=\lambda_{L_va}+\eta_v(a),
\]
where $L_v$ is a finite-band matrix with entries of at most linear growth
in the Fourier index, and $\eta_v$ is a finite-coordinate linear
functional. If $K_v$ bounds the Fourier support of $v$, then
\[
 \|[L_v,\Pi_M]a\|_{\ell^2}
 \leq C_v\left\|
  \mathbf1_{\{n>M-K_v\}}(1+n)a\right\|_{\ell^2}
 \longrightarrow0.
\]
Here the weight $n$ denotes the Fourier frequency of each real coordinate.
Applying this estimate to the smooth test sequences defining the
coordinates of $T_{\psi_s}Y$, uniformly on compact time intervals,
proves convergence of the projected linear functionals and their field
derivatives in every finite Gaussian $L^p$ norm. Products and integration
in $s$ give the same convergence for $J_{t,M}$ and its field derivative.
H\"older's inequality and the uniform exponential bounds therefore give,
for sufficiently small $|t|$, convergence in $L^2(\mathbb G\otimes\mathbb G)$:
\begin{align*}
 H_{t,M}&\longrightarrow H_t:=P(T_{\psi_t}Y,Z)e^{-J_t(Y)},\\
 \shortintertext{and}
 \Vzero_vH_{t,M}&\longrightarrow\Vzero_vH_t.
\end{align*}

For fixed $M$, multiply $H_{t,M}$ by a smooth cutoff equal to one on
the radius-$L$ ball in its finite coordinate space. As $L\to\infty$,
the cutoff error and its field derivative tend to zero in $L^2$ by
the same exponential bounds. The resulting smooth compactly supported
function has finite Hermite sums converging with its derivatives and
with polynomial coordinate weights in Gaussian $L^2$: multiplication
and differentiation act as raising and lowering operators on the
Hermite coefficients. The coefficients decay faster than any power,
since multiplication by the square root of the Gaussian density gives
a Schwartz function. Thus a diagonal choice of these Hermite sums
gives finite field polynomials $P_{t,k}$ such that
\[
 \|P_{t,k}-H_t\|_2
 +\|\Vzero_vP_{t,k}-\Vzero_vH_t\|_2\longrightarrow0.
\]
The cylinder $F_0(h\circ\psi_t)$ has bounded value and $D_v$ derivative.
Define
\[
 G_{t,k}(h,Y,Z)=F_0(h\circ\psi_t)P_{t,k}(Y,Z)
\]
These tests lie in the algebra of \eqref{eq:KD.28}. Their $L^2$ graph
convergence, together with
\[
 \|\Sigma_v(G_{t,k}-G_t)\|_{L^1(\pi)}
 \leq\|\Sigma_v\|_{L^2(\mathbb G)}\|G_{t,k}-G_t\|_{L^2(\pi)},
\]
proves \eqref{eq:finite-flow-graph}.
Finally, the same bounds applied to the smooth time derivatives of
$C_s,d_s$ control $\partial_tG_t$ locally in $L^1(\pi)$.

\subsection{Proof of Proposition~\ref{prop:gaussian-core}}\label{app:gaussian-core}

\begin{proof}
Put $\xi_n=\sqrt nY_n$. Use annihilators $a_n=\partial_{\xi_n}$ and
$b_n=\partial_{\bar\xi_n}$, and their Gaussian adjoints
$a_n^*=\xi_n-b_n$, $b_n^*=\bar\xi_n-a_n$.  Direct substitution in \eqref{eq:KD.2} gives the
positive noise component, up to its fixed unitary phase, as
\begin{equation}\label{eq:GC.1}
d_n^{-1/2}(J_n-P_n+q_na_n),\qquad
 q_n=Qn\sqrt n,\qquad
 P_n=\sum_{k+\ell=n}\sqrt{k\ell}\,a_ka_\ell.
\end{equation}
The negative component is the conjugate expression.  Explicitly,
\[
 J_n=\sum_{k\geq1}\sqrt{k(k+n)}
       (a_k^*a_{k+n}-b_{k+n}^*b_k)
       -\sum_{k+\ell=n}\sqrt{k\ell}\,b_k^*a_\ell.
\]
The operator $J_n$ preserves Gaussian degree. Regard each particle
as having a signed frequency $p\in\mathbb Z\setminus\{0\}$. Its
one-particle matrix satisfies
\[
 p\longmapsto p-n,\qquad
 |(J_n)_{p-n,p}|=\sqrt{|p(p-n)|},\qquad p-n\ne0.
\]
The signed matrix enters $A_0$; its absolute coefficients enter the off-diagonal
bounds below. The negative noise component shifts in the opposite
direction.

Let ${N_{\mathrm{part}}}$ be Gaussian particle number, acting as multiplication
by $r$ on the $r$th homogeneous Wiener chaos. Thus $\Vzero=J+T$, where
$J$ preserves ${N_{\mathrm{part}}}$, while $T$ lowers it by one or two.
On degree $r$ the aggregate lower part satisfies
\begin{equation}\label{eq:GC.2}
\|Tf\|_{\ell^2\otimes L^2}
                   \le C(Q)(r+1)\|f\|.
\end{equation}
The pair and charge terms satisfy
\begin{align*}
 \sum_{k+\ell=n}\frac{k\ell}{d_n}&=\frac16,\qquad
 \|(d_n^{-1/2}P_nf)_{n\geq2}\|_{\ell^2\otimes L^2}
          \leq\sqrt{\frac{r(r-1)}6}\|f\|.
\shortintertext{For the charge term,}
 \sup_{n\geq2}\frac{q_n^2}{d_n}&\leq\frac{4Q^2}3,
 \qquad
 \|(d_n^{-1/2}q_na_nf)_{n\geq2}\|_{\ell^2\otimes L^2}
          \leq\frac{2|Q|}{\sqrt3}\sqrt r\|f\|.
\end{align*}
Their sum proves \eqref{eq:GC.2}.

To verify the operator domain, fix a normalized Wick monomial $f$.
Take $n$ larger than every occupied frequency of $f$. In a
normal-ordered term, write $p,q$ for occupied signed frequencies.
The elementary contributions from $J_n^*J_n$, $q_na_n^*J_n$ and
$P_n^*J_n$, after division by $d_n$, have coefficients bounded by
\begin{align*}
 |c^{JJ}_{n;p,q}|
 &\leq C_f\frac{\sqrt{|pq(p-n)(q+n)|}}{d_n}
                    \leq C_f n^{-2},\\
 |c^{qJ}_{n;p}|
 &\leq C_f\frac{|q_n|\sqrt{|p(p-n)|}}{d_n}
                    \leq C_{f,Q} n^{-1},\\
 \shortintertext{and, for $k,\ell\geq1$ with $k+\ell=n$,}
 |c^{PJ}_{n;p;k,\ell}|
 &\leq C_f\frac{\sqrt{k\ell\,|p(p-n)|}}{d_n}
                    \leq C_f n^{-3/2}.
\end{align*}
Occupation factors are absorbed in $C_f$. There are only a fixed
number of occupied input choices, whereas the last line has $n-1$
ordered splittings; conjugate components satisfy the same bounds.
Their squared tail sum is therefore bounded by
\[
 C_{f,Q}\sum_{n\geq2}\bigl(n^{-4}+n^{-2}+n\,n^{-3}\bigr)<\infty.
\]
The diagonal contractions are absolutely summable, since
\[
 \frac{\|J_nf\|^2}{d_n}\leq C_f\frac n{d_n}
                                  \leq C_f n^{-2}.
\]
The remaining terms have finite frequency support, and each output
in the displayed tails has only a bounded number of input choices.
Let $A_{\rm alg}$
denote the formal matrix on $\mathrm{Pol}_0$. Finite-noise integration
by parts, followed by these square-summable limits, gives
\[
 \langle\Vzero f,\Vzero g\rangle=\langle A_{\rm alg}f,g\rangle,
 \qquad f,g\in\mathrm{Pol}_0.
\]
For $g\in\Dom(\Vzero)$, choose $g_m\in\mathrm{Pol}_0$ with
$g_m\to g$ and $\Vzero g_m\to\Vzero g$ in $L^2$. Thus
\[
 \langle\Vzero f,\Vzero g\rangle
 =\lim_m\langle A_{\rm alg}f,g_m\rangle
 =\langle A_{\rm alg}f,g\rangle.
\]
Therefore
\[
 \mathrm{Pol}_0\subset\Dom(A_0),\qquad
 A_0f=A_{\rm alg}f\quad(f\in\mathrm{Pol}_0).
\]

Fix $M<\infty$, and let $\Pi_M$ project onto degrees at most $M$.
On $\Pi_M\mathrm{Pol}_0$ consider
\[
                 A_M^{\rm alg}=\Pi_M A_{\rm alg}\Pi_M.
\]
It is symmetric and nonnegative, with form $\|\Vzero f\|^2$; $\Vzero f$ also
has degree at most $M$.  We prove that $A_M^{\rm alg}$ is essentially
self-adjoint. Constants in this step may depend on $M$.

Index normalized Wick monomials by finite occupations
$x=(m_p)_{p\in\mathbb Z\setminus\{0\}}$ with
$r(x)=\sum_pm_p\le M$.  Set
\begin{equation}\label{eq:GC.3}
S(x)=\sum_p|p|m_p,\qquad
 c_x^0=\prod_p\frac{|p|^{-m_p/2}}{\sqrt{m_p!}},\qquad
 \upsilon(s)=\sqrt{\log(s+e)},\qquad c_x=c_x^0e^{\upsilon(S(x))}.
\end{equation}
Set $\omega_x=c_x^2$. Multiplication by $c_x$ identifies
$\ell^2(\omega)$ with counting $\ell^2$. For every finite $s$,
\[
 \#\{x:S(x)\leq s\}<\infty.
\]
Moreover,
\[
 S(x)=0\iff x=0,\qquad \Vzero1=A_01=0.
\]
The estimates below concern $S(x)\geq1$.

Write $a_{xy}$ for the matrix of $A_M^{\rm alg}$ in the normalized
Wick basis.  The key estimate is
\begin{equation}\label{eq:GC.4}
\sup_x\sum_{y\ne x}|a_{xy}|\frac{c_y}{c_x}
         [\upsilon(S(y))-\upsilon(S(x))]^2\le C_{M,Q}.
\end{equation}
Occupation factors and the
number of choices of the finitely many input particles are bounded
by constants depending only on $M$.  Before inserting
$e^{\upsilon(S(y))-\upsilon(S(x))}$, the weight $c_y^0/c_x^0$ cancels the
square roots of all output frequencies in a normal-ordered monomial.
The remaining frequency factors are the input frequencies.

For clarity, the cancellations giving the rates in
\eqref{eq:GC.5}--\eqref{eq:GC.8} are recorded below.  Here $N=|n|$,
$d_N=N^3-N$, and the displayed ratios omit only Wick occupation
factors bounded in terms of $M$:
\[
\begin{array}{c|c|c|c}
\text{block}&\text{raw coefficient}&c_y^0/c_x^0&
       \text{weighted rate}\\[2pt]\hline
J^*J&
 \dfrac{\sqrt{|p(p-n)q(q+n)|}}{d_N}&
 \sqrt{\dfrac{|pq|}{|p-n||q+n|}}&
 \dfrac{|pq|}{d_N}\\[9pt]
J^*P&
 \dfrac{\sqrt{k\ell\,|p(p+n)|}}{d_N}&
 \sqrt{\dfrac{k\ell\,|p|}{|p+n|}}&
 \dfrac{|p|k\ell}{d_N}\\[9pt]
P^*J&
 \dfrac{\sqrt{k\ell\,|p(p-n)|}}{d_N}&
 \sqrt{\dfrac{|p|}{k\ell\,|p-n|}}&
 \dfrac{|p|}{d_N}\\[9pt]
J^*q&
 \dfrac{|q_n|\sqrt{|p(p-n)|}}{d_N}&
 \sqrt{\dfrac{N|p|}{|p-n|}}&
 \dfrac{|q_n|\sqrt N\,|p|}{d_N}
       \le C_Q\dfrac{|p|}{N}\\[9pt]
q^*J&
 \dfrac{|q_n|\sqrt{|p(p-n)|}}{d_N}&
 \sqrt{\dfrac{|p|}{N|p-n|}}&
 \dfrac{|q_n|\,|p|}{\sqrt N\,d_N}
       \le C_Q\dfrac{|p|}{N^2}
\end{array}
\]
In the pair blocks, $k,\ell\ge1$, $k+\ell=N$; hence $J^*P$
has one input triple per transition, whereas $P^*J$ has $N-1$
ordered splittings for each fixed input frequency.  Since
$|q_n|=|Q|N^{3/2}$ and $d_N\asymp N^3$, the last two bounds follow
directly.  The negative-frequency blocks are their conjugates.

Normal ordering $J_n^*J_n$ produces a diagonal one-particle term and
a two-particle transition
\begin{equation}\label{eq:GC.5}
(p,q)\longmapsto(p-n,q+n),
 \qquad\text{absolute weighted rate}\quad
                 \le C_M\frac{|pq|}{d_{|n|}}.
\end{equation}
Only nonzero output frequencies are retained.  The change in $S$
has absolute value at most $2|n|$.  Summing $|pq|$ over the input
choices is at most $C_M S^2$.

The pair-decreasing block $J^*P$ removes signed frequencies $k,\ell,p$
with $k,\ell$ of the same sign and $n=k+\ell$, and inserts $p+n$.
Put $N=|n|=|k|+|\ell|$.  Its base weighted rate is at most
\begin{equation}\label{eq:GC.6}
C_M\frac{|p k\ell|}{d_N}.
\end{equation}
It never increases $S$, and its decrease $\delta$ is at most $2N$.
Consequently its total downward second moment is bounded by
\begin{equation}\label{eq:GC.7}
\sum_{y:\,S(y)<S(x)}|a^{J^*P}_{xy}|\frac{c_y^0}{c_x^0}
                      [S(x)-S(y)]^2
             \le C_M S(x)^2.
\end{equation}
Indeed each input triple contributes at most
$C|p|\,|k\ell|/N\le C|p|\min(|k|,|\ell|)\le CS^2$, and
there are only $C_M$ triples.

For its pair-increasing adjoint $P^*J$, fix an input frequency $p$.
There are $N-1$ ordered splittings of $n$ into same-sign nonzero
frequencies.  Each has base rate at most $C_M|p|/d_N$.
After the splittings and input particles are summed, the rate at
noise frequency $N$ is at most
\begin{equation}\label{eq:GC.8}
C_M\frac{S}{N^2}.
\end{equation}
The increase of $S$ is at most $2N$.

The degree-decreasing charge block $J^*q$ removes $p,n$ and inserts
$p+n$. Its base rate is at most $C_{M,Q}|p|/N$, and its decrease is at
most $2\min(|p|,N)$. Thus its downward second moment satisfies
\eqref{eq:GC.7}. The adjoint $q^*J$ has base rate at most
$C_{M,Q}|p|/N^2$, which gives \eqref{eq:GC.8}.

The remaining blocks $P^*P$, $P^*q$, $q^*P$ and $q^*q$ preserve $S$
and contribute zero to \eqref{eq:GC.4}. The normal-ordering contraction
in $J^*J$ is diagonal: each input has a unique output under a fixed
signed shift.

For $s'=s-\delta\geq0$,
\begin{equation}\label{eq:GC.9}
e^{\upsilon(s')-\upsilon(s)}[\upsilon(s')-\upsilon(s)]^2
                           \le C\delta^2/s^2,\quad s\ge1.
\end{equation}
For $\delta\le s/2$, use the derivative of $\upsilon$; for $\delta>s/2$,
use $\sup_{t\ge0}e^{-t}t^2<\infty$. Together with \eqref{eq:GC.7}, this
bounds all downward charge and pair blocks.

For an upward jump $s'=s+\delta$ with $0\leq\delta\leq2N$,
\[
 e^{\upsilon(s')-\upsilon(s)}[\upsilon(s')-\upsilon(s)]^2
 \leq
 \begin{cases}
 CN^2/[s^2\log(s+e)],&N\leq s,\\
 C_\varepsilon(N/s)^\varepsilon,&N>s,
 \end{cases}
 \qquad s\geq1.
\]
Here $\varepsilon>0$ is arbitrary. Taking $s=S(x)$ and
$\varepsilon<1$ in \eqref{eq:GC.8} gives
\[
 \frac{Cs}{s^2\log(s+e)}\sum_{2\leq N\leq s}1
 +Cs\sum_{N>s}\frac{(N/s)^\varepsilon}{N^2}\leq C.
\]

For \eqref{eq:GC.5}, divide the sum into three ranges. When
$2\leq N\leq s$ and both endpoints are at least $s/2$, the cost is
$CN^2/[s^2\log(s+e)]$.

If an endpoint falls below $s/2$, then $N\geq s/4$, and the
exponential square difference is bounded by a constant. The range
$N>s$ is upward, with cost at most $C_\varepsilon(N/s)^\varepsilon$.
For $\varepsilon<2$, summing these three contributions gives
\[
 \frac C{\log(s+e)}\sum_{2\leq N\leq s}\frac{N^2}{d_N}
 +Cs^2\sum_{N\geq\max(2,s/4)}N^{-3}
 +Cs^2\sum_{N>s}N^{-3}(N/s)^\varepsilon\leq C.
\]
Together with \eqref{eq:GC.7} and \eqref{eq:GC.9}, this proves
\eqref{eq:GC.4}; the finitely many small values of $s$ are absorbed
in the constant.

Choose $\chi_R(x)=\chi(\upsilon(S(x))/R)$ with $\chi$ smooth, one near
zero, and compactly supported.  It has finite support, and \eqref{eq:GC.4}
implies
\begin{equation}\label{eq:GC.10}
\sum_y|a_{xy}|\frac{c_y}{c_x}
                    |\chi_R(y)-\chi_R(x)|^2
                        \le C_{M,Q}/R^2.
\end{equation}
Let $u\in\ker((A_M^{\rm alg})^*+\lambda)$, $\lambda>0$, and
write $u_x=c_x f_x$. The weak equation can be tested against
$c\chi_R^2f$, a finite-support vector. Symmetry gives
\[
 \Rea\langle u,A_M^{\rm alg}(\chi_R^2u)\rangle
 =\langle A_M^{\rm alg}(\chi_Ru),\chi_Ru\rangle
   +\frac12\Rea\sum_{x,y}a_{xy}\bar u_xu_y
                         (\chi_R(x)-\chi_R(y))^2.
\]
Consequently,
\begin{equation}\label{eq:GC.11}
\lambda\|\chi_R f\|_{\ell^2(\omega)}^2
       +\langle A_M^{\rm alg}(c\chi_R f),c\chi_R f\rangle
                \le C_{M,Q}R^{-2}\|f\|_{\ell^2(\omega)}^2.
\end{equation}
The symmetry $|\omega_xa_{xy}c_y/c_x|=|\omega_ya_{yx}c_x/c_y|$
and $2|f_xf_y|\leq|f_x|^2+|f_y|^2$ bound its remainder by
\[
 \sum_x\omega_x|f_x|^2\sum_y|a_{xy}|\frac{c_y}{c_x}
                       |\chi_R(y)-\chi_R(x)|^2
 \leq C_{M,Q}R^{-2}\|f\|_{\ell^2(\omega)}^2.
\]
This is absolutely summable by \eqref{eq:GC.10}. Every
matrix row defines the weak adjoint equation because a single Wick
basis vector has an $\ell^2$ image, as checked above.

Since $A_M^{\rm alg}\geq0$ and $\chi_R\to1$ pointwise,
\[
 \lambda\|u\|^2
 =\lambda\lim_{R\to\infty}\|\chi_Rf\|_{\ell^2(\omega)}^2
 \leq\lim_{R\to\infty}C_{M,Q}R^{-2}\|u\|^2=0.
\]
Thus $\ker((A_M^{\rm alg})^*+\lambda)=\{0\}$ and
$A_M^{\rm alg}$ is essentially self-adjoint. Denote its closure by $A_M$.

To pass from finite sectors to the full Gaussian space, let
$u\in\ker(A_{\rm alg}^*+\lambda)$.
The hard particle projection preserves $\Dom(\Vzero)$, since
\[
 \Vzero\Pi_M f=\Pi_M\Vzero f+[T,\Pi_M]f,
 \qquad \|[T,\Pi_M]\|\leq C_{M,Q}.
\]
Here the projection on the noise side acts on its scalar-chaos
factor, and the commutator has finite particle support by
\eqref{eq:GC.2}. Projecting polynomial form approximations by
$\Pi_M$ now shows that
the form domain of $A_M$ is precisely
$\Pi_M L^2\cap\Dom(\Vzero)$, with form $\|\Vzero f\|^2$.

For $f\in\Pi_M\mathrm{Pol}_0$, the part of $A_0 f$ above degree $M$ is
\begin{equation}\label{eq:GC.12}
B_Mf=(1-\Pi_M)T^*\Vzero f,
                  \qquad \|B_Mf\|\le C_{M,Q}\|\Vzero f\|.
\end{equation}
Indeed $\Vzero f$ has degree at most $M$, and the number-preserving
adjoint $J^*$ cannot leave that sector.  The adjoint row bound in
\eqref{eq:GC.12} follows from \eqref{eq:GC.2} on degrees at most $M+2$.

The deficiency equation thus says
\begin{equation}\label{eq:GC.13}
\langle\Pi_Mu,(\lambda+A_M)f\rangle
                   =-\langle(1-\Pi_M)u,B_Mf\rangle.
\end{equation}
By \eqref{eq:GC.12},
\[
 |\langle(1-\Pi_M)u,B_Mf\rangle|
 \leq C_{M,Q}\|u\|\|A_M^{1/2}f\|.
\]
Form-resolvent duality supplies $w_M\in\Dom(A_M^{1/2})$
such that
\[
 \lambda\langle w_M,f\rangle
 +\langle A_M^{1/2}w_M,A_M^{1/2}f\rangle
 =-\langle(1-\Pi_M)u,B_Mf\rangle.
\]
For $f\in\Pi_M\mathrm{Pol}_0$, subtracting \eqref{eq:GC.13} yields
\[
 \langle\Pi_Mu-w_M,(\lambda+A_M)f\rangle=0.
\]
Essential self-adjointness gives
\[
 \overline{(\lambda+A_M)\Pi_M\mathrm{Pol}_0}=\Pi_ML^2.
\]
Hence $\Pi_Mu=w_M$ and
\begin{equation}\label{eq:GC.14}
\Pi_Mu\in\Dom(\Vzero)
                       \qquad\text{for every }M.
\end{equation}

Let $\pi_r=\Pi_r-\Pi_{r-1}$ be the projection onto chaos degree $r$,
with $\Pi_{-1}=0$, and use the same notation on the scalar-chaos
factor of the noise space. Define the local output components by
\[
 \pi_r(\Vzero u):=\pi_r\Vzero\Pi_{r+2}u,\qquad r\geq0.
\]
Only degrees $r,r+1,r+2$ can contribute to this component. Thus
\[
 \pi_r\Vzero\Pi_Mu=\pi_r\Vzero\Pi_{r+2}u,
                           \qquad M\geq r+2.
\]
For every finite-chaos form test $f$, the weak equation
has finite pairings and takes the form
\begin{equation}\label{eq:GC.15}
\lambda\langle u,f\rangle
                           +\langle \Vzero u,\Vzero f\rangle=0.
\end{equation}
For a polynomial test of degree at most $M$, $A_0 f$ has degree at
most $M+2$, and the identity uses only the finite-sector form domain:
\[
 \langle u,A_0 f\rangle
   =\langle\Pi_{M+2}u,A_0 f\rangle
   =\langle \Vzero\Pi_{M+2}u,\Vzero f\rangle.
\]
For a finite-chaos form test $f$ of degree at most $M$, finite-sector
form density gives $f_m\in\Pi_M\mathrm{Pol}_0$ with
\[
 f_m\longrightarrow f,\qquad \Vzero f_m\longrightarrow\Vzero f
 \quad\text{in }L^2.
\]
The preceding equality with $\Pi_{M+2}u$ therefore passes to $f$ and
gives \eqref{eq:GC.15}.

Take a smooth nonincreasing $\zeta:[0,\infty)\to[0,1]$, equal to
one on $[0,1]$ and zero on $[2,\infty)$, and put
\[
       \Theta_R=\zeta(R^{-1}\log({N_{\mathrm{part}}}+1))^2.
\]
The same spectral multiplier acts on the scalar-chaos part of the
noise space. Write $\theta_r=\zeta(R^{-1}\log(r+1))^2$ for its value
on degree $r$. From \eqref{eq:GC.2},
\begin{equation}\label{eq:GC.16}
C_R^{\mathrm G}=[\Vzero,\Theta_R],\qquad
                    \|C_R^{\mathrm G}\|\le C(Q)/R.
\end{equation}
The second commutator also factors as
\begin{equation}\label{eq:GC.17}
[[\Vzero,\Theta_R],\Theta_R]=\Theta_R H_R^{\mathrm G},
                    \|H_R^{\mathrm G}\|\le C(Q)/R^2.
\end{equation}
A lowering block has input degree
$r$ and output degree $s=r-d$, $d=1,2$, and its double-commutator
coefficient is $(\theta_r-\theta_s)^2$ times that block.  Since
$\zeta$ is nonincreasing, $\theta_s\ge\theta_r$, and
\begin{equation}\label{eq:GC.18}
\frac{(\theta_r-\theta_s)^2}{\theta_s}
       \le4|\zeta(R^{-1}\log(r+1))-
                 \zeta(R^{-1}\log(s+1))|^2
       \le\frac{C_d}{R^2(r+1)^2}.
\end{equation}
The quotient is taken as zero if $\theta_s=0$.  Multiplying by
the $O_Q(r+1)$ block norm from \eqref{eq:GC.2} proves \eqref{eq:GC.17}; the at most
two overlapping output degrees only change the constant.

By \eqref{eq:GC.14}, $\Theta_R^2u$ is an admissible test in \eqref{eq:GC.15}.
The exact first-order operator identity is
\begin{equation}\label{eq:GC.19}
\operatorname{Re}\langle \Vzero u,\Vzero\Theta_R^2u\rangle
   =\|\Vzero\Theta_Ru\|^2-\|C_R^{\mathrm G}u\|^2
       +\operatorname{Re}\langle \Vzero u,\Theta_RH_R^{\mathrm G}u\rangle.
\end{equation}
Since $\Theta_R\Vzero u=\Vzero\Theta_Ru-C_R^{\mathrm G}u$,
\eqref{eq:GC.16}--\eqref{eq:GC.19} and
Young's inequality give
\begin{equation}\label{eq:GC.20}
\lambda\|\Theta_Ru\|^2
              +\tfrac12\|\Vzero\Theta_Ru\|^2
                    \le C(Q) R^{-2}\|u\|^2,\qquad R\ge1.
\end{equation}
All occurrences of $\Vzero u$ in this calculation have finite output
particle support, so \eqref{eq:GC.14} applies. Since $\Theta_Ru\to u$
in $L^2$, letting $R\to\infty$ gives $u=0$. Thus $A_{\rm alg}$ is
essentially self-adjoint and its closure is $A_0$.
\end{proof}

\subsection{Proof of Lemma~\ref{lem:smooth-polynomial-approximation}}\label{app:smooth-polynomial-approximation}

\begin{proof}
By Hahn--Banach, it suffices to consider a continuous functional
annihilating all field polynomials. With $p'$ conjugate to $p$, write
\[
 \ell(f)=\sum_{|\eta|\leq k}\int h_\eta\partial^\eta f\,d\mathbb G_d,
 \qquad h_\eta\in L^{p'}(\mathbb G_d),\qquad
 \ell(x^\gamma)=0\quad(\gamma\in\mathbb N_0^d).
\]
The function
\[
 \widehat\ell(z)=\sum_{|\eta|\leq k}(iz)^\eta
             \int h_\eta(x)e^{iz\cdot x}\,d\mathbb G_d(x)
\]
is entire by Gaussian H\"older estimates, and
\[
 \partial_z^\gamma\widehat\ell(0)=i^{|\gamma|}\ell(x^\gamma)=0,
 \qquad \widehat\ell=0.
\]
Define the tempered distribution
\[
 T=\sum_{|\eta|\leq k}(-1)^{|\eta|}
                    \partial^\eta(h_\eta\,\mathbb G_d).
\]
For $\phi\in C_c^\infty(\mathbb R^d)$,
\begin{align*}
 \langle T,\phi\rangle
 &=\sum_{|\eta|\leq k}\int h_\eta\partial^\eta\phi\,d\mathbb G_d
   =\ell(\phi),\\
 \shortintertext{and, with the Fourier convention $e^{iz\cdot x}$,}
 \widehat T(z)&=\widehat\ell(z)=0,\qquad z\in\mathbb R^d.
\end{align*}
Fourier uniqueness gives $T=0$, hence $\ell=0$ on $C_c^\infty$.
Cutoff and convolution density in the Gaussian Sobolev norm then
give $\ell=0$ on the full space.

For the mark write $F(x,q)=\sum_{m\in\mathbb Z}F_m(x)e^{imq}$.
If $b>l+1$, integration by parts in $q$ gives
\begin{align*}
 \|F_m\|_{W^{k,p}(\mathbb G_d)}
 &\leq |m|^{-b}\|F\|_{k,b,p}\quad(m\ne0).\\
 \shortintertext{Summing the Fourier tail, we obtain}
 \left\|F-\sum_{|m|\leq N}F_me^{imq}\right\|_{k,l,p}
 &\leq C\|F\|_{k,b,p}\sum_{|m|>N}|m|^{l-b}\longrightarrow0.
\end{align*}
Approximating the finitely many $F_m$ by field polynomials completes
the proof.
\end{proof}

\subsection{Proof of the Gaussian stress modulus}\label{app:stress-modulus}

At each fixed angle, $T_N=\partial_\theta \Bstress_N$ is a linear plus a
centred quadratic polynomial in the standard complex Gaussians.
The linear variance is at most
\[
 C Q^2\sum_{n=2}^N\frac{n^5}{d_n^2}\le C(Q)\log(N+2).
\]
Up to realification, phases, and a universal factor, the quadratic
coefficient matrix is
\begin{equation}\label{eq:GO.3}
a^{(N)}_{kl}=\frac{\sqrt{kl}(k+l)}{d_{k+l}}
                        \mathbf1_{2\le k+l\le N}.
\end{equation}
Since $n/d_n\le C/n^2$ for $n\ge2$, Schur's test with weights
$p_k=k^{-1/2}$ gives
\[
 \sum_l|a^{(N)}_{kl}|p_l
 \le C\sqrt k\sum_l(k+l)^{-2}\le Cp_k.
\]
Thus $\sup_N\|a^{(N)}\|_{\rm op}<\infty$.  Also
\begin{equation}\label{eq:GO.4}
\|a^{(N)}\|_{\rm HS}^2
 \le C\sum_{n=2}^N n^{-4}\sum_{k=1}^{n-1}k(n-k)
 \le C\log(N+2).
\end{equation}
The Gaussian determinant identity, \eqref{eq:GO.4}, and
Cauchy--Schwarz for the linear term give
\[
 \log\mathbb E e^{\lambda T_N(\theta)}
 \leq C(Q)\lambda^2\log(N+2),\qquad |\lambda|\leq c(Q).
\]
Exponential Chebyshev bounds with
$|\lambda|=\min(c(Q)/2,u/[2C(Q)\log(N+2)])$ yield
\begin{equation}\label{eq:GO.5}
\mathbb G(|T_N(\theta)|>u)
 \le C(Q)\exp\left[-c(Q)\min\left\{
              \frac{u^2}{\log(N+2)},u\right\}\right].
\end{equation}
For a trigonometric polynomial of degree $N$, the Bernstein inequality
and a deterministic grid of at most $CN$ points imply
$\|T_N\|_\infty\le2\max_{\rm grid}|T_N|$.
Define
\[
 \Lambda_Q^0(Y):=\sup_{m\ge1}
     \frac{\|\partial_\theta \Bstress_{2^m}\|_\infty}{m+1}
\]
and set $\Bstress_1=0$ when forming the first dyadic block.
For sufficiently large $A$ depending only on $Q$,
\eqref{eq:GO.5} and the grid bound imply
\[
 \mathbb G(\Lambda_Q^0>A)
 \leq C(Q)\sum_{m\geq1}2^m e^{-c(Q)A(m+1)}
 \leq C(Q)e^{-c(Q)A}.
\]

Let $D_m^\Sigma=\Bstress_{2^m}-\Bstress_{2^{m-1}}$. Its Fourier support lies in a fixed
annulus of radius $2^m$.  A smooth Fourier multiplier equal to
$1/(in)$ on this annulus has convolution-kernel $L^1$ norm at most
$C2^{-m}$.  Hence
\begin{equation}\label{eq:GO.6}
\|D_m^\Sigma\|_\infty
 \le C2^{-m}\|(D_m^\Sigma)'\|_\infty
 \le C\Lambda_Q^0(Y)(m+1)2^{-m}.
\end{equation}
\begin{align*}
 \sum_{m\geq1}\|D_m^\Sigma\|_\infty
 &\leq C\Lambda_Q^0(Y)\sum_{m\geq1}(m+1)2^{-m}<\infty.\\
 \shortintertext{Hence}
 \Bstress_{2^m}&\longrightarrow\Bstress
                    \quad\text{in }C(S^1).
\end{align*}
Given $\delta=\operatorname{dist}_{S^1}(x,y)\leq1$, choose
$2^{-m}\asymp\delta$. Then
\begin{align*}
 |\Bstress(Y,x)-\Bstress(Y,y)|
 &\leq\delta\|\Bstress_{2^m}'\|_\infty
                   +2\sum_{j>m}\|D_j^\Sigma\|_\infty\\
 &\leq C\Lambda_Q^0\{\delta(m+1)+(m+1)2^{-m}\}\\
 &\leq C\Lambda_Q^0\delta\log(e/\delta).
\end{align*}
Taking $\Lambda_Q$ to be a fixed multiple of $1+\Lambda_Q^0$ proves
\eqref{eq:GO.2}.

\subsection{Proof of Lemma~\ref{lem:conditional-stress-bounds}}
\label{app:stress-covariance}
\begin{proof}
Let $Y(z)=\sum_{k\geq1}Y_kz^k=m(z)+\widetilde Y(z)$ be a Gaussian
analytic field with complex coefficients and mean $m$.
In real coordinates, assume the covariance of
$(\sqrt{k}\widetilde Y_k)_{k\geq1}$ is bounded by $C_0\operatorname{Id}$.
Assume also
\begin{equation}\label{eq:B.1}
\sup_{z\in\mathbb D}(1-|z|^2)|m'(z)|\le M.
\end{equation}
Use the stress coefficients $\Sigma_n$ of \eqref{eq:KD.3} and the
convention $\operatorname{Var}(U)=\mathbb E|U-\mathbb EU|^2$ for
complex random variables.
There is a universal numerical constant $C$ such that, for every finite
complex sequence $a$,
\begin{equation}\label{eq:B.2}
\operatorname{Var}\!\left(\sum_{n\ge2}a_n\alpha_n\Sigma_n\right)
 \le C\big(C_0^2+C_0M^2+C_0Q^2\big)\sum_n|a_n|^2.
\end{equation}
Real sine and cosine coordinates satisfy the same bound up to a
numerical factor. Centering and zero-mode terms also satisfy
\[
 \operatorname{Var}\!\left(\sum_{n\geq2}a_n\alpha_nn\widetilde Y_n\right)
 \leq CC_0\sum_{n\geq2}|a_n|^2\alpha_n^2n
 =CC_0\sum_{n\geq2}\frac{2|a_n|^2}{n^2-1}
 \leq CC_0\|a\|_{\ell^2}^2.
\]

In coordinates $\xi_k=\sqrt{k}\widetilde Y_k$, the quadratic-chaos
matrix $M_a$ has entries
$a_{k+\ell}\alpha_{k+\ell}\sqrt{k\ell}$. Hence
\begin{equation}\label{eq:stress-hilbert-schmidt}
 \|M_a\|_{\mathrm{HS}}^2
 =\sum_{k,\ell\ge1}|a_{k+\ell}|^2\alpha_{k+\ell}^2 k\ell
 =\frac13\|a\|_{\ell^2}^2.
\end{equation}
The identity uses $\sum_{k=1}^{n-1}k(n-k)=d_n/6$.
Wick's formula then gives
\[
 \operatorname{Var}\!\left(
   \sum_{k,\ell\geq1}a_{k+\ell}\alpha_{k+\ell}
           \sqrt{k\ell}\,\xi_k\xi_\ell\right)
       \leq CC_0^2\|a\|_{\ell^2}^2.
\]
The charge term satisfies
\[
 \operatorname{Var}\!\left(\sum_n a_n\alpha_nQn^2\widetilde Y_n\right)
 \leq C C_0Q^2\sum_n|a_n|^2\alpha_n^2n^3
       \le \tfrac83 C C_0Q^2\|a\|^2.
\]

For the linear term from the mean, the required operator is
\[
 (T_{\mathrm{mean}}a)_k=\sqrt{k}\sum_{\ell\ge1}
                   \ell m_\ell\alpha_{k+\ell}a_{k+\ell}.
\]
Set
\[
 A_a(z)=\sum_{n\ge2}\frac{a_n}{\alpha_n}z^{n-2},\qquad
 F_x'(z)=\sum_{k\ge1}\sqrt{k}\,x_kz^{k-1}.
\]
Parseval in $A^2_2$ gives
\[
 \|A_a\|_{A^2_2}=\|a\|_{\ell^2},
 \qquad \|F_x'\|_{L^2(dA/\pi)}=\|x\|_{\ell^2}.
\]
Indeed, $\int_0^1u^{n-2}(1-u)^2du=2/d_n=\alpha_n^2$.
For the convention $\langle a,x\rangle=\sum_k a_k\overline{x_k}$,
\begin{align*}
 \langle T_{\mathrm{mean}}a,x\rangle
 &=\overline{\int_{\mathbb D}A_{\overline a}(z)
      \overline{m'(z)F_{\overline x}'(z)}
       (1-|z|^2)^2\,\frac{dA(z)}\pi}.
\shortintertext{The Bloch bound and Cauchy--Schwarz yield}
 |\langle T_{\mathrm{mean}}a,x\rangle|
       &\le\|A_{\overline a}\|_{A^2_2}
                              \|m'F_{\overline x}'\|_{A^2_2}
       \le M\|a\|\|x\|.
\end{align*}
Consequently,
\[
 \operatorname{Var}\!\left(2\sum_{k,\ell\geq1}
       a_{k+\ell}\alpha_{k+\ell}k\ell m_\ell\widetilde Y_k\right)
       \leq CC_0M^2\|a\|_{\ell^2}^2.
\]
Orthogonality of the first and second chaoses yields \eqref{eq:B.2}.

Let $E$ be a real separable Hilbert space equipped with an isonormal
Gaussian process. Index the following families by $i\in\mathbb N$.
Let $M_i$ be self-adjoint Hilbert--Schmidt operators on $E$, and
let $\ell_i\in E$. Suppose that, for every finite real sequence $a$,
\begin{align}
 \left\|\sum_i a_i M_i\right\|_{\mathrm{HS}(E)}
       &\leq C_M\|a\|_{\ell^2},\qquad
 \left\|\sum_i a_i\ell_i\right\|_E\leq C_\ell\|a\|_{\ell^2},\notag
\shortintertext{and}
 \left\|\sum_i M_iM_i^*\right\|_{E\to E}&\leq C_{\mathrm{row}}^2.
       \label{eq:W.1}
\end{align}
The positive operator sum is the increasing strong limit of its finite
partial sums. Identify self-adjoint Hilbert--Schmidt operators with
symmetric tensors in $E\otimes E$. Let $I_r$ denote the associated
Wiener integral, normalized by
\[
 \mathbb E[I_r(f)I_s(g)]=\mathbf1_{\{r=s\}}r!\langle f,g\rangle,
 \qquad \eta_i=I_2(M_i)+I_1(\ell_i).
\]
For every finitely
supported Hilbert--Schmidt matrix $a=(a_{ij})$,
\begin{equation}\label{eq:W.2}
\left\|\sum_{i,j}a_{ij}
       \big(\eta_i\eta_j-\mathbb E[\eta_i\eta_j]\big)
       \right\|_{L^2}
       \le C(C_M,C_{\mathrm{row}},C_\ell)\|a\|_{\mathrm{HS}(\ell^2)}.
\end{equation}

Let $\operatorname{Sym}$ denote tensor symmetrization. The product
formula gives
\begin{align*}
 F_4&=\operatorname{Sym}\sum_{i,j}a_{ij}M_i\otimes M_j,\\
 F_3&=\operatorname{Sym}\sum_{i,j}a_{ij}
                (M_i\otimes\ell_j+\ell_i\otimes M_j),\\
 F_2&=\operatorname{Sym}\sum_{i,j}a_{ij}
                         (4M_iM_j+\ell_i\otimes\ell_j),\\
 F_1&=2\sum_{i,j}a_{ij}(M_i\ell_j+M_j\ell_i).
\shortintertext{The centered product therefore expands as}
 \sum_{i,j}a_{ij}(\eta_i\eta_j-\mathbb E\eta_i\eta_j)
    &=\sum_{r=1}^4I_r(F_r).
\end{align*}
For the contraction terms, set $N_i=\sum_j a_{ij}M_j$ and
$R:E\otimes\ell^2\to E$, $R((x_i)_i)=\sum_iM_ix_i$. Then
\begin{equation}\label{eq:W.3}
\left\|\sum_{i,j}a_{ij}M_iM_j\right\|_{\mathrm{HS}}
  =\left\|\sum_i M_iN_i\right\|_{\mathrm{HS}}
  \le \|R\|\left(\sum_i\|N_i\|_{\mathrm{HS}}^2\right)^{1/2}
  \le C_{\mathrm{row}} C_M\|a\|_{\mathrm{HS}}.
\end{equation}
The same calculation with $M_j$ replaced by $\ell_j$, and then with
$a$ replaced by its transpose, gives
\begin{align*}
 \|F_4\|&\leq C_M^2\|a\|_{\mathrm{HS}},\qquad
 \|F_3\|\leq2C_MC_\ell\|a\|_{\mathrm{HS}},\\
 \|F_2\|&\leq(4C_{\mathrm{row}}C_M+C_\ell^2)\|a\|_{\mathrm{HS}},\qquad
 \|F_1\|\leq4C_{\mathrm{row}}C_\ell\|a\|_{\mathrm{HS}}.
\shortintertext{Orthogonality of the chaoses now gives}
 \left\|\sum_{r=1}^4I_r(F_r)\right\|_2^2
 &=\sum_{r=1}^4r!\|F_r\|^2
 \leq C(C_M,C_{\mathrm{row}},C_\ell)^2\|a\|_{\mathrm{HS}}^2.
\end{align*}

By duality,
\[
 \sum_{i,j}\left|\mathbb E[(\eta_i\eta_j-
                    \mathbb E\eta_i\eta_j)g]\right|^2
 \leq C(C_M,C_{\mathrm{row}},C_\ell)^2\|g\|_2^2.
\]

For the complex holomorphic stress, before splitting real and
imaginary coordinates, take
\[
 (M_n)_{k\ell}=\alpha_n\sqrt{k\ell}\,
                     \mathbf1_{\{k+\ell=n\}},
 \qquad n\ge2.
\]
The Hilbert--Schmidt estimate \eqref{eq:stress-hilbert-schmidt}
gives the first bound in \eqref{eq:W.1} with $C_M=1/\sqrt3$.
Moreover $\sum_n M_nM_n^*$ is diagonal, and its $k$th entry is
\begin{equation}\label{eq:W.4}
\sum_{n>k}\frac{2k(n-k)}{n^3-n}
           \le \frac{8}{3}k\sum_{n>k}\frac1{n^2}
           \le \frac{8}{3}.
\end{equation}
The same estimates, up to universal realification factors, apply to
the real sine and cosine stress matrices.  The linear stress map is
bounded whenever its conditional Gaussian mean has uniformly bounded
Bloch seminorm, by the weighted Bergman calculation established in
\eqref{eq:B.2}.

If the conditional real covariance is a bounded positive operator
$S\le C_0\operatorname{Id}$, replace $M_i$ by $S^{1/2}M_iS^{1/2}$ and
$\ell_i$ by $S^{1/2}\ell_i$.  The Bessel bound persists, and
\[
 \sum_i(S^{1/2}M_iS^{1/2})(S^{1/2}M_iS^{1/2})^*
       \le C_0 S^{1/2}\Big(\sum_i M_iM_i^*\Big)S^{1/2}.
\]
Thus \eqref{eq:W.2} holds uniformly over the conditional Gaussian kernels
with bounded covariance and Bloch mean.

For the same Gaussian field and the connection of
Section~\ref{sec:geometry}, set
\[
 H^{\rm field}_{ij}=\Vzero_i\Sigma_j-
                                  \Sigma_{\nabla_{e_i}e_j}.
\]
The affine Witt action is the one in \eqref{eq:KD.2}.

The following estimate holds on the horizontal indices:
\begin{equation}\label{eq:GH.1}
\left\|\sum_{i,j}a_{ij}
       (H^{\rm field}_{ij}-\mathbb EH^{\rm field}_{ij})\right\|_{L^2}
       \le C(C_0,M,Q)\|a\|_{\mathrm{HS}},
\end{equation}
for finite matrices $a$, uniformly over the conditional states above.

For the same-sign block at zero mean, put $N=n+m$ and
$\beta_{mn}=d_m(m+2n)/d_N$. In coordinates $\xi_k=\sqrt k \widetilde Y_k$, the
quadratic coefficient of this block is, up to one fixed normalization,
\begin{equation}\label{eq:GH.2}
\alpha_n\alpha_m\sqrt{k\ell}\,
  [\beta_{mn}-(k-n)_+-(\ell-n)_+],\qquad k+\ell=n+m.
\end{equation}
The charge-linear coefficient is
\begin{equation}\label{eq:GH.3}
Q\alpha_n\alpha_m mn\sqrt N\,
                \frac{mn+2}{N^2-1}\,\xi_N.
\end{equation}
For \eqref{eq:GH.3}, assume $n\le m$; its coefficient is at most
$C|Q|\sqrt n/N$, and
\[
 \sum_{2\leq n\leq N/2}\frac{Q^2n}{N^2}\leq CQ^2.
\]

For \eqref{eq:GH.2}, the connection coefficient is
$c(n,m)=-d_m(m+2n)/d_N=-\beta_{mn}$.
By symmetry exchange $n,m$ so that $n\le m$. Since
$k+\ell=n+m$, we have
$(k-m)_+=(n-\ell)_+$ and $(\ell-m)_+=(n-k)_+$.
Thus, with $b=d_n(2m+n)/d_N$, the bracket becomes
\[
          b-(n-k)_+-(n-\ell)_+.
\]
We estimate the $b$ term and the two positive-part terms
separately. The bulk coefficient has bound
$C n^{3/2}\sqrt{k(N-k)}/N^{7/2}$, with
\begin{align*}
 \sup_{n\leq N/2}\sum_{k=1}^{N-1}
       \frac{n^{3/2}\sqrt{k(N-k)}}{N^{7/2}}&\leq C,
\shortintertext{and}
 \sup_{1\leq k<N}\sum_{n\leq N/2}
       \frac{n^{3/2}\sqrt{k(N-k)}}{N^{7/2}}&\leq C.
\end{align*}
For an edge $k<n$, the coefficient has bound
\[
       C\frac{(n-k)\sqrt k}{n^{3/2}N},
\]
and
\[
 \sum_{k<n}\frac{(n-k)\sqrt k}{n^{3/2}N}\leq\frac{Cn}N,
 \qquad
 \sum_{k<n\leq N/2}\frac{(n-k)\sqrt k}{n^{3/2}N}
                                      \leq C\sqrt{\frac kN}.
\]
The other edge follows by $k\leftrightarrow N-k$.
Schur's test bounds the map on every fixed total-frequency
level $N$ uniformly. Distinct levels give orthogonal
second-chaos components, so Wick's isometry and square
summation yield \eqref{eq:GH.1} on the same-sign block.

At zero mean, the mixed-sign block, apart from deterministic and
charge-linear terms, equals
\begin{equation}\label{eq:GH.4}
2T_Y^*T_Y,\qquad
 (T_Ya)_k=\sqrt k\sum_{\ell\ge1}\ell Y_\ell
                                  \alpha_{k+\ell}a_{k+\ell}.
\end{equation}
Its quadratic matrix map splits along a fixed difference $d\ge0$.
The associated nonnegative Volterra kernel is
\[
 H_d(k,n)=(n-k)\sqrt{k(k+d)}\,\alpha_n\alpha_{n+d},
                       \qquad n>k.
\]
The two bounds
\[
 H_d(k,n)\leq C\sqrt{k}\,n^{-3/2},\qquad
 H_d(k,n)\leq\frac{C(n-k)}{n(n+d)}
\]
give
\begin{align*}
 \sup_{k\geq1}\sum_{n>k}H_d(k,n)
    &\leq C\sup_k\sqrt k\sum_{n>k}n^{-3/2}\leq C,
\shortintertext{and}
 \sup_{n\geq2}\sum_{k<n}H_d(k,n)
    &\leq C\sup_n\frac n{n+d}\leq C.
\end{align*}
Schur's bound is uniform in $d$. The centered
second-chaos components at different differences $d$
are orthogonal; square summation therefore gives the
corresponding bound for the whole mixed-sign block.
The charge residual in the convention of \eqref{eq:KD.2} is a constant multiple of
$\alpha_n\alpha_{n+d}n^2\sqrt d\,\xi_d$. It also has a Bessel bound,
because
\[
 \sum_{n\ge2}\alpha_n^2\alpha_{n+d}^2n^4d
       \le C d\sum_{n\ge2}\frac{n}{(n+d)^3}\le C.
\]
At difference one, the additional low-connection charge coefficient is
a fixed multiple of $Q\alpha_n\alpha_{n+1}(2n+1)$. Its squares satisfy
\[
 Q^2\sum_{n\geq2}\alpha_n^2\alpha_{n+1}^2(2n+1)^2
       =4Q^2\sum_{n\geq2}\frac{(2n+1)^2}{d_nd_{n+1}}<\infty.
\]

For the Bloch mean, the weighted Bergman multiplier calculation gives
$\|T_{\mathrm{mean}}\|\le M$.
For a deterministic first-chaos coefficient vector $x$, put
$Y_k=x_k/\sqrt k$. Then
\begin{equation}\label{eq:GH.5}
\|T_x\|_{\mathrm{HS}}^2
 =\sum_{k,p\ge1}kp\alpha_{k+p}^2|x_p|^2\le C\|x\|^2.
\end{equation}
Consequently,
\[
 \|T_{\mathrm{mean}}^*T_x+T_x^*T_{\mathrm{mean}}\|_{\mathrm{HS}}
 \leq2\|T_{\mathrm{mean}}\|\,\|T_x\|_{\mathrm{HS}}
 \leq CM\|x\|.
\]

For the same-sign block, its bulk term factors through the bounded map
$x\mapsto T_{\mathrm{mean}}^*x$ and the level multiplier
\[
 \sum_{n+m=N,\ n\le m}
               |\alpha_n\alpha_m b/\alpha_N|^2\le C.
\]
For the two edges, take $n\leq m$, put $r=n-k$, and define
\[
 (E_x)_{r,m}=\mathbf1_{\{r\leq m\}}\alpha_m\sqrt{r(m+r)}\,x_{m+r},
 \qquad
 (E_{\mathrm{mean}})_{r,m}
       =\alpha_m\sqrt r\,(m+r)m_{m+r}\quad(r\leq m).
\]
The edge matrices are the upper-triangular parts of
$T_{\mathrm{mean}}^\top E_x$ and $T_x^\top E_{\mathrm{mean}}$,
respectively. Changing variables $\ell=m+r$
gives $\|E_x\|_{\mathrm{HS}}\leq C\|x\|$.
To bound $E_{\mathrm{mean}}$, choose
$\psi\in C_c^\infty((1/2,8))$ equal to one on $[1,6]$.
For $N\ge2$ set $r_N=1-N^{-1}$ and define
\[
 m^{[N]}(e^{i\theta})
        =\sum_{k\ge1}\psi(k/N)m_ke^{ik\theta}.
\]
Two summations by parts in the multiplier coefficients give
\begin{align*}
 \left|\sum_{k\ge1}
      \frac{\psi(k/N)}{k r_N^{k-1}}e^{ik\theta}\right|
   &\leq C\min\{1,(N|\sin(\theta/2)|)^{-2}\}.
\shortintertext{Integrating over the circle gives}
 \left\|\sum_{k\ge1}
      \frac{\psi(k/N)}{k r_N^{k-1}}e^{ik\theta}
          \right\|_{L^1(d\theta/(2\pi))}&\le\frac CN.
\end{align*}
Convolving this kernel with
$e^{i\theta}m'(r_Ne^{i\theta})$ proves
\[
 \|m^{[N]}\|_\infty
       \le\frac CN\|m'(r_Ne^{i\cdot})\|_\infty
       \le\frac{CM}{N(1-r_N^2)}\le CM.
\]
Extend $E_{\mathrm{mean}}$ by the same formula on dyadic rectangles
$r\asymp R$, $m\asymp M_0$ with $R\leq2M_0$, and set it to zero
outside these rectangles. Their frequencies satisfy $M_0\le r+m\le6M_0$,
so the corresponding Hankel coefficient block is a compression
of the Hankel operator with symbol $m^{[M_0]}$. Therefore
\[
 \|(m_{r+m})_{r\asymp R,\,m\asymp M_0}\|_{\ell^2\to\ell^2}
        \le\|m^{[M_0]}\|_\infty\le CM.
\]
Splitting the weight
as $\sqrt r\,m\alpha_m+r^{3/2}\alpha_m$ gives rectangle norm at most
\[
           CM[(R/M_0)^{1/2}+(R/M_0)^{3/2}].
\]
For dyadic $R,M_0$, the two Schur sums are bounded by
\begin{align*}
 \sup_{M_0}\sum_{R\leq2M_0}
       [(R/M_0)^{1/2}+(R/M_0)^{3/2}]&\leq C,
\shortintertext{and}
 \sup_R\sum_{M_0\geq R/2}
       [(R/M_0)^{1/2}+(R/M_0)^{3/2}]&\leq C.
\end{align*}
Thus $\|E_{\mathrm{mean}}\|\leq CM$. Since $T_x(r,n)=0$ for $r\geq n$,
the extra entries $r>m$ vanish when $n\leq m$. Projection onto this
triangle is contractive in the Hilbert--Schmidt norm, so
\[
 \|\mathbf1_{\{n\leq m\}}T_{\mathrm{mean}}^{\top}E_x\|_{\mathrm{HS}}
 +\|\mathbf1_{\{n\leq m\}}T_x^{\top}E_{\mathrm{mean}}\|_{\mathrm{HS}}
       \leq CM\|x\|.
\]
The conditional covariance and
Bloch-mean bounds from Section~\ref{sec:gaussian-lift} are uniform in
$h$. For $S\le C_0\operatorname{Id}$, the quadratic kernels are
compressed by $S^{1/2}\otimes S^{1/2}$, with norm at most $C_0$,
while the charge and mean-induced linear kernels are compressed
by $S^{1/2}$, with norm at most $\sqrt{C_0}$.
After realification, the preceding estimates give
\[
 \eqref{eq:B.2}\ \Longrightarrow\ \eqref{eq:conditional-score-bound},
 \qquad
 \eqref{eq:W.2}+\eqref{eq:GH.1}
      \ \Longrightarrow\ \eqref{eq:conditional-matrix-bound}.
\]
\end{proof}

\section{Approximation in generator and gradient norms}\label{sec:mixed-core}

The following abstract lemma converts $L^1$ range density into
simultaneous approximation of a function, its gradient, and its generator.
The applications in Section~\ref{sec:conditional-stress} verify its
hypotheses for the marked Gaussian form and the tilted welding form.

Let $m$ be a probability measure. Let
$D_m:\Dom(D_m)\subset L^2(m)\to L^2(m;H)$
be a closed, densely defined gradient whose associated form is
conservative, strongly local and Markovian. Set $A_m=D_m^*D_m$.

Let $C_{\mathrm{test}}\subset\Dom(A_m)$ be a linear space invariant
under complex conjugation. Assume $D_mg\in L^4(m;H)$ for its real
elements and
\[
 \overline{(\lambda+A_m)C_{\mathrm{test}}}^{\,L^1(m)}=L^1(m),
 \qquad \lambda>0.
\]
Equip $\Dom(D_m)$ with the form norm
\[
 \|v\|_{\Dom({D_m})}^2=\|v\|_2^2+\|{D_m}v\|_2^2.
\]
For real $g\in C_{\mathrm{test}}$ and bounded smooth $\chi$ with
bounded derivatives of every order, the operator chain rule gives
\begin{align*}
 \chi(g)&\in\Dom(A_m),\\
 \shortintertext{and}
 {A_m}\chi(g)&=\chi'(g){A_m}g
                         -\chi''(g)\|{D_m}g\|^2.
\end{align*}
For every real $g\in C_{\mathrm{test}}$ and every such $\chi$ with
$\chi(0)=0$ and $0\leq\chi'\leq1$, suppose there are real
$p_k\in C_{\mathrm{test}}$ with
\begin{equation}\label{eq:mixed-contraction-approximation}
 \|p_k-\chi(g)\|_2
 +\|{D_m}p_k-{D_m}\chi(g)\|_2
 +\|{A_m}p_k-{A_m}\chi(g)\|_1
                         \longrightarrow0.
\end{equation}
Let $\ell$ be a continuous real linear functional on the real form
domain $\Dom({D_m})$.

\begin{lemma}\label{lem:mixed-core}
For every bounded real $u\in\Dom(A_m)$ there are real
$p_n\in C_{\mathrm{test}}$ such that
\begin{align}
 p_n&\longrightarrow u&&\text{in }L^2(m),\notag\\
 D_mp_n&\longrightarrow D_mu&&\text{in }L^2(m;H),\notag
\shortintertext{and}
 A_mp_n&\longrightarrow A_mu&&\text{in }L^1(m).
       \label{eq:mixed-convergence}
\end{align}
If a real $u\in L^\infty(m)$ satisfies, for some $\lambda>0$,
\begin{equation}\label{eq:mixed-distribution}
 \langle u,(\lambda+{A_m})p\rangle=\ell(p),
                         \qquad p\in C_{\mathrm{test}}\text{ real},
\end{equation}
then $u\in\Dom({D_m})$.
\end{lemma}
\begin{proof}
For the operator chain rule, take bounded real $v\in\Dom(D_m)$.
The gradient chain rule and
integration by parts give
\begin{align*}
 \langle{D_m}\chi(g),{D_m}v\rangle
 &=\langle{A_m}g,\chi'(g)v\rangle
       -\int\chi''(g)\|{D_m}g\|^2v\,dm\\
 &=\left\langle\chi'(g){A_m}g
               -\chi''(g)\|{D_m}g\|^2,v\right\rangle.
\end{align*}
The right-hand coefficient is in $L^2(m)$, since
\[
 \|\chi'(g){A_m}g-\chi''(g)\|{D_m}g\|^2\|_2
 \le\|\chi'\|_\infty\|{A_m}g\|_2
          +\|\chi''\|_\infty\|{D_m}g\|_4^2.
\]
Density of bounded form tests proves the stated identity and
$\chi(g)\in\Dom(A_m)$.

Now fix $\lambda>0$. Range density and conjugation give real
$g_n\in C_{\mathrm{test}}$ such that
\[
 h_n=(\lambda+{A_m})(g_n-u)\longrightarrow0
                                      \quad\text{in }L^1(m).
\]
The Markov resolvent then gives
\[
 \|g_n-u\|_1\leq\lambda^{-1}\|h_n\|_1,\qquad
 \|{A_m}g_n-{A_m}u\|_1\leq2\|h_n\|_1.
\]
Put $v_n=g_n-u$.
For $C>0$, take a bounded increasing smooth function $\tau_C$ with
$\tau_C(0)=0$, $\tau_C'=1$ on $[-C,C]$, and
$\|\tau_C\|_\infty\leq2C$. The form chain rule gives
\begin{align}
 \int_{\{|v_n|\leq C\}}\|D_mv_n\|^2\,dm
 &\leq\lambda\int v_n\tau_C(v_n)\,dm
       +\int\tau_C'(v_n)\|D_mv_n\|^2\,dm\notag\\
 &=\langle h_n,\tau_C(v_n)\rangle
 \leq2C\|h_n\|_1.\label{eq:local-truncation-energy}
\end{align}
Choose $M>2\|u\|_\infty$ and a bounded smooth $\chi$ equal to the
identity on $[-M,M]$, with $0\leq\chi'\leq1$ and $\chi'=0$ outside
$[-2M,2M]$. Set $f_n=\chi(g_n)$ and
$E_n=\{|v_n|\geq M/2\}$. Then
\begin{align*}
 m(E_n)&\leq\frac2M\|v_n\|_1\longrightarrow0,
\shortintertext{and bounded truncation gives}
 \|f_n-u\|_2^2&\leq3M\|v_n\|_1\longrightarrow0.
\end{align*}
The gradient chain rule and \eqref{eq:local-truncation-energy} give
\begin{align*}
 {D_m}f_n-{D_m}u
   &=\chi'(g_n){D_m}v_n+(\chi'(g_n)-1){D_m}u.
\shortintertext{Consequently,}
 \|{D_m}f_n-{D_m}u\|_2^2
   &\leq2\int_{\{|v_n|\leq3M\}}\|{D_m}v_n\|^2\,dm
        +2\int_{E_n}\|{D_m}u\|^2\,dm\\
   &\leq12M\|h_n\|_1+2\int_{E_n}\|{D_m}u\|^2\,dm
      \longrightarrow0.
\end{align*}
For the generator,
\[
 {A_m}f_n=\chi'(g_n){A_m}g_n
                         -\chi''(g_n)\|{D_m}g_n\|^2.
\]
Since $\{\chi''(g_n)\ne0\}\subset E_n\cap\{|v_n|\leq3M\}$,
\begin{align*}
 \|{A_m}f_n-{A_m}u\|_1
 &\leq\|{A_m}g_n-{A_m}u\|_1+\int_{E_n}|{A_m}u|\,dm\\
 &\quad+2\|\chi''\|_\infty
       \left(6M\|h_n\|_1+\int_{E_n}\|{D_m}u\|^2\,dm\right)
       \longrightarrow0.
\end{align*}
By \eqref{eq:mixed-contraction-approximation}, choose real
$p_n\in C_{\mathrm{test}}$ with
\[
 \|p_n-f_n\|_2+\|{D_m}p_n-{D_m}f_n\|_2
                  +\|{A_m}p_n-{A_m}f_n\|_1\leq n^{-1}.
\]
The triangle inequality proves \eqref{eq:mixed-convergence}.

To recover the form domain from the distributional identity, take
real $h\in L^\infty(m)$ and set
$f=(\lambda+{A_m})^{-1}h$. Then
$f\in\Dom({A_m})$ and
$\|f\|_\infty\leq\lambda^{-1}\|h\|_\infty$.
Apply \eqref{eq:mixed-convergence} to $f$, obtaining real $p_n$.
Equation \eqref{eq:mixed-distribution} gives
\begin{align*}
 |\langle u,h\rangle-\ell(f)|
 &\leq\|u\|_\infty
       \left(\lambda\|p_n-f\|_1+\|{A_m}p_n-{A_m}f\|_1\right)
       +|\ell(p_n-f)|\longrightarrow0.
\end{align*}
Coercivity gives a unique $z\in\Dom({D_m})$ with
\[
 \lambda\langle z,v\rangle+\langle{D_m}z,{D_m}v\rangle
                          =\ell(v),\qquad v\in\Dom({D_m})\text{ real}.
\]
For the preceding $f$ and $h$,
\[
 \langle z,h\rangle
 =\lambda\langle z,f\rangle+\langle{D_m}z,{D_m}f\rangle
 =\ell(f)=\langle u,h\rangle.
\]
Density of $L^\infty(m)$ in $L^2(m)$ gives $u=z\in\Dom({D_m})$.
\end{proof}

\end{document}